\documentclass [11pt,reqno]{amsart}
\usepackage {amsmath,amssymb,verbatim,geometry,color,esint}
\usepackage[all]{xy}
\usepackage{mathrsfs}
\usepackage[backref,pagebackref,pdftex,hyperindex]{hyperref}
\usepackage[pdftex]{graphicx}
\usepackage{tikz}
\usetikzlibrary{math}

\usepackage{accents}
\newcommand{\dbtilde}[1]{\accentset{\approx}{#1}}

\newcommand{\A}{\mathbb{A}}
\newcommand{\B}{\mathbb{B}}
 \newcommand{\C}{\mathbb{C}}

\newcommand{\Gm}{\mathbb{G}_\mathrm{m}}

\newcommand{\N}{\mathbb{N}}
\renewcommand{\P}{\mathbb{P}}
 \newcommand{\Q}{\mathbb{Q}}
 \newcommand{\R}{\mathbb{R}}
 \newcommand{\Z}{\mathbb{Z}}

\newcommand{\fa}{\mathfrak{a}}
\newcommand{\fb}{\mathfrak{b}}

\newcommand{\fm}{\mathfrak{m}}

\newcommand{\fp}{\mathfrak{p}}

\newcommand{\hD}{\widehat{\Delta}}

\newcommand{\hf}{\widehat{\f}}

\renewcommand{\hom}{\mathrm{hom}}

\newcommand{\tX}{\widetilde{X}}
\newcommand{\tcF}{\widetilde{\cF}}

\newcommand{\cA}{\mathcal{A}}

\newcommand{\cC}{\mathcal{C}}

\newcommand{\cE}{\mathcal{E}}
\newcommand{\cF}{\mathcal{F}}

\newcommand{\cH}{\mathcal{H}}
\newcommand{\cI}{\mathcal{I}}

\newcommand{\cL}{\mathcal{L}}
\newcommand{\cM}{\mathcal{M}}

\newcommand{\cO}{\mathcal{O}}
\newcommand{\cP}{\mathcal{P}}

\newcommand{\cR}{\mathcal{R}}

\newcommand{\cU}{\mathcal{U}}

\newcommand{\cX}{\mathcal{X}}
\newcommand{\cY}{\mathcal{Y}}

\newcommand{\tcL}{\widetilde{\mathcal{L}}}
\newcommand{\tcX}{\widetilde{\mathcal{X}}}

\newcommand{\Xan}{X^{\mathrm{an}}}

\newcommand{\Xdiv}{X^{\mathrm{div}}}
\newcommand{\Xlin}{X^{\mathrm{lin}}}

\renewcommand{\a}{\alpha}
\renewcommand{\b}{\beta}
\renewcommand{\d}{\delta}
\newcommand{\e}{\varepsilon}
\newcommand{\f}{\varphi}
\newcommand{\unipar}{\varpi}
\newcommand{\g}{\gamma}
\newcommand{\la}{\lambda}
\newcommand{\om}{\omega}

\newcommand{\Ga}{\Gamma}
\newcommand{\La}{\Lambda}
\newcommand{\p}{\psi}

\newcommand{\cf}{{\rm cf.\ }} 
\newcommand{\eg}{{\rm e.g.\ }} 
\newcommand{\ie}{{\rm i.e.\ }}

\newcommand{\inter}{\cdot\ldots\cdot}
\newcommand{\winter}{\wedge\dots\wedge}

\newcommand{\hto}{\hookrightarrow}

\newcommand{\dom}{\mathrm{dom}}

\newcommand{\gf}{\mathrm{gf}}

\DeclareMathOperator{\dd}{{d}}

\DeclareMathOperator{\depth}{depth}

\newcommand{\an}{\mathrm{an}}

\DeclareMathOperator{\en}{E}

\DeclareMathOperator{\ii}{I}
\DeclareMathOperator{\jj}{J}
\DeclareMathOperator{\env}{P}

\DeclareMathOperator{\qq}{Q}

\DeclareMathOperator{\Cz}{C^0}

\newcommand{\redu}{\mathrm{red}}

\DeclareMathOperator{\Aff}{Aff}

\DeclareMathOperator{\Capa}{Cap}
\DeclareMathOperator{\charac}{char}

\DeclareMathOperator{\Hom}{Hom}

\DeclareMathOperator{\MA}{MA}

\DeclareMathOperator{\one}{\mathbf{1}}

\DeclareMathOperator{\Spec}{Spec}
\DeclareMathOperator{\Amp}{Amp}
\DeclareMathOperator{\supp}{supp}
\DeclareMathOperator{\vol}{vol}

\DeclareMathOperator{\Pic}{Pic}

\DeclareMathOperator{\ord}{ord}
\DeclareMathOperator{\Nef}{Nef}

\DeclareMathOperator{\Psef}{Psef}

\DeclareMathOperator{\PSH}{PSH}
\DeclareMathOperator{\CPSH}{CPSH}

\DeclareMathOperator{\Zun}{Z^1}
\DeclareMathOperator{\Zunb}{Z^1_{\mathrm{b}}}

\DeclareMathOperator{\Numb}{N^1_{\mathrm{b}}}
\DeclareMathOperator{\Carb}{Car_{\mathrm{b}}}
\DeclareMathOperator{\Carbplus}{Car^+_{\mathrm{b}}}
\DeclareMathOperator{\Nefb}{Nef_{\mathrm{b}}}

\DeclareMathOperator{\tii}{\overline{I}}

\DeclareMathOperator{\PL}{PL}

\DeclareMathOperator{\Hnot}{H^0}

\DeclareMathOperator{\te}{T}
\DeclareMathOperator{\teL}{T_L}

\DeclareMathOperator{\Num}{N^1}

\DeclareMathOperator{\Bs}{Bs}

\DeclareMathOperator{\Car}{Car}

\DeclareMathOperator{\VCar}{VCar}

\DeclareMathOperator{\enh}{\mathrm{E}^{\uparrow}}
\DeclareMathOperator{\enb}{\mathrm{E}^{\downarrow}}

\DeclareMathOperator{\ddc}{dd^c}

\newcommand{\supstar}{\mathrm{sup}^\star }

\newcommand{\envstar}{\mathrm{P}^\star }

\newcommand{\ssup}{\mathrm{sup}}

\renewcommand{\div}{\mathrm{div}}

\newcommand{\triv}{\mathrm{triv}}
\newcommand{\val}{\mathrm{val}}

\newcommand{\lin}{\mathrm{lin}}

\newcommand{\NA}{\mathrm{NA}}

\newcommand{\tor}{\mathrm{tor}}

\newcommand{\D}{\Delta}

\newcommand{\cro}[1]{[\![#1]\!]}
\newcommand{\lau}[1]{(\!(#1)\!)}
\newcommand{\simto}{\overset\sim\to}

\numberwithin{equation}{section}       

\newtheorem{prop} {Proposition} [section]
\newtheorem{thm}[prop] {Theorem} 
\newtheorem{defi}[prop] {Definition}
\newtheorem{lem}[prop] {Lemma}
\newtheorem{cor}[prop]{Corollary}
\newtheorem{prop-def}[prop]{Proposition-Definition}

\newtheorem*{thmA}{Theorem A} 
 
\newtheorem*{thmB}{Theorem B} 
\newtheorem*{thmC}{Theorem C}

\newtheorem{exam}[prop]{Example}
\newtheorem{rmk}[prop]{Remark}
\newtheorem{qst}[prop]{Question}
\newtheorem{conj}[prop]{Conjecture}
\theoremstyle{remark}

\title{Global pluripotential theory over a trivially valued field}
\date{\today}

\author{S{\'e}bastien Boucksom \and Mattias Jonsson}

\address{CNRS--CMLS\\
  \'Ecole Polytechnique\\
  F-91128 Palaiseau Cedex\\
  France}
\email{sebastien.boucksom@polytechnique.edu}

\address{Dept of Mathematics\\
  University of Michigan\\
  Ann Arbor, MI 48109-1043\\
  USA}
\email{mattiasj@umich.edu}

\begin{document}

\begin{abstract}
We develop global pluripotential theory in the setting of Berkovich geometry over a trivially valued field. Specifically, we define and study functions and measures of finite energy and the non-Archimedean Monge--Amp\`ere operator on any (possibly reducible) projective variety. We also investigate the topology of the space of valuations of linear growth, and the behavior of psh functions thereon. 

\medskip

\emph{Nous d\'eveloppons une th\'eorie du pluripotentiel global dans le contexte de la g\'eom\'etrie de Berkovich sur un corps trivialement valu\'e. Plus pr\'ecis\'ement, nous d\'efinissons et \'etudions des fonctions et mesures d'\'energie finie et un op\'erateur de Monge--Amp\`ere non-archim\'edien sur toute vari\'ete projective (\'eventuellement r\'eductible). Nous explorons \'egalement la topologie de l'espace des valuations \`a croissance lin\'eaire, et le comportement des fonctions psh sur celui-ci.}
\end{abstract}

\dedicatory{Dedicated to Ahmed Zeriahi for all his contributions to complex analysis and geometry}

\maketitle

\setcounter{tocdepth}{1}
\tableofcontents
 \section*{Introduction}
 %
%
%
%
The main purpose of the present paper is to lay the foundations of pluripotential theory in the setting of Berkovich geometry over a trivially valued field, paralleling as much as possible the known theory in the complex analytic case. 

Pluripotential theory is a crucial tool in complex analysis and geometry, and A.~Zeriahi has been a central protagonist in the story of its development. For an excellent introduction, see the book by V.~Guedj and A.~Zeriahi~\cite{GZBook}. This theory lies at the heart of the variational approach to complex Monge--Amp\`ere equations~\cite{BBGZ},
which proved particularly fruitful in relation to K\"ahler--Einstein metrics, and ultimately led to a new proof of the Yau--Tian--Donaldson (YTD) conjecture for Fano manifolds~\cite{YTD2}, later generalized to (log terminal) Fano varieties~\cite{LTW19,Li19}. 

There are several reasons for developing a corresponding theory over a trivially valued field, that is, any field equipped with the trivial (non-Archimedean) valuation. First, the Berkovich approach to non-Archimedean geometry is the one closest to the complex analytic intuition, and it is natural to see to what extent the complex analytic theory admits a counterpart in the setting of Berkovich spaces.
Second, Berkovich spaces over a trivially or discretely valued field can often be used to study degenerations in complex geometry~\cite{pshsing,valmul,eigenval,dyncomp,konsoib,BHJ2,Shi19}, and such degenerations  are central to the variational approach to the YTD conjecture~\cite{YTD2,Li20}, through the notion of geodesic rays.

The thrust of the YTD conjecture is to relate the existence of a solution to a non-linear PDE to an algebro-geometric condition known as \emph{K-stability}, which has recently also come to play a key role in the study of moduli spaces, especially for Fano varieties.~\cite{ABHX,BHLX,BLX,LXZ21,XZ20}. The relation of K-stability with spaces of valuations and non-Archimedean geometry over trivially valued fields, originally pointed out in~\cite{BHJ1}, is one major motivation to endeavor the present study, and will be further exploited in the companion papers~\cite{nakstab1,nakstab2}. 

K-stability of a polarized projective variety $(X,L)$ over an algebraically closed field $k$ is a condition phrased in terms of $\Gm$-equivariant degenerations $(\cX,\cL)\to\A^1$ of $(X,L)$ known as \emph{test configurations}~\cite{Don02}. Our basic proposal, which goes back to~\cite{BHJ1}, consists in interpreting $(\cX,\cL)$ in terms of a \emph{piecewise linear function} $\f_\cL$ on the Berkovich analytification $X^\an$ with respect to the trivial absolute value on the ground field $k$. To readers familiar with non-Archimedean geometry, this will sound very natural indeed: the base change of a test configuration $(\cX,\cL)$ by $\Spec k\cro{\unipar}\to\A^1=\Spec k[\unipar]$ provides a model for the base change $(X_K,L_K)$ to the non-Archimedean field $K=k\lau{\unipar}$, and hence a model/PL metric on the Berkovich analytification $L_K^\an\to X_K^\an$. This metric can further be canonically identified with a PL function on $X_K^\an$, thanks to the reference metric on $L_K^\an$ induced by the trivial model, and $\f_\cL$ is simply the restriction of this function to $X^\an\hookrightarrow X_K^\an$, the embedding being realized by Gauss extension.
 
While this point of view was basically the one adopted in a previous version of this article~\cite{trivval}, we have tried here to take a more elementary and self-contained approach, avoiding for the most part any explicit reference to general Berkovich geometry---which accounts in part for the length of the present article. 
%
%
%
%
\subsection*{Functions and measures of finite energy}
Let $k$ be an arbitrary algebraically closed field, and $X$ an irreducible\footnote{While the main body of the article deals with possibly reducible varieties, we assume here for simplicity that $X$ is irreducible.} projective variety over $k$, with function field $k(X)$.
In order to describe more precisely our main results, recall first that the Berkovich analytification $X^\an$ of $X$ with respect to the trivial valuation on $k$ is a compactification of the space $X^\val$ of valuations $v\colon k(X)^\times\to\R$, trivial on $k$, endowed with the topology of pointwise convergence. Points of $X^\an$ can be understood as \emph{semivaluations} $v$ on $X$, \ie valuations $v\in Y^\val$ for some subvariety $Y\subset X$. The set $X^\div\subset X^\val$ of \emph{divisorial valuations}, which are attached to prime divisors over $X$, is dense in $X^\an$. 

For every $v\in X^\an$ and every section $s\in\Hnot(X,mL)$, $m\in\N$, we can make sense of $v(s)\in[0,+\infty]$ by trivializing $L$ at the center of $v$. Setting $|s|(v):=\exp(-v(s))$ defines a continuous function $|s|\colon X^\an\to[0,1]$, and the topology of $X^\an$ is in fact defined by the set of such functions. Building on these, we introduce the following classes of functions on $X^\an$:

\begin{itemize}

\item \emph{Fubini--Study functions} for $L$ are continuous functions $\f\colon X^\an\to\R$ of the form 
$$
\f=m^{-1}\max_i\{\log|s_i|+\la_i\},
$$
where $(s_i)$ is a finite set in $\Hnot(X,mL)$ without common zeroes and $\la_i\in\Q$. This defines a subset $\cH=\cH(L)$ of $\Cz(X):=\Cz(X^\an,\R)$, and we show that the map $(\cX,\cL)\mapsto\f_\cL$ alluded to above sets up a 1--1 correspondence between $\cH(L)$ and the set of ample, integrally closed test configurations $(\cX,\cL)$ for $(X,L)$ (see~\S\ref{sec:intclosed} for the definition when $X$ is not normal). 

\smallskip

\item \emph{piecewise linear (PL) functions} on $X^\an$ are elements of the $\Q$-linear subspace $\PL(X)$ of $\Cz(X)$ spanned by $\cH=\cH(L)$; this subspace turns out to be independent of $L$, and is dense in $\Cz(X)$. In the present setting, PL functions play the role of smooth functions from the complex analytic case. It is proved in Appendix A that PL functions are induced by usual piecewise linear functions on dual complexes of snc test configurations, when $X$ is smooth and $\charac k=0$. 

\smallskip

\item \emph{$L$-psh functions} are usc functions $\f\colon X^\an\to\R\cup\{-\infty\}$, not identically $-\infty$, that can be obtained as limits of decreasing nets in $\cH(L)$. The set $\PSH=\PSH(L)$ of $L$-psh functions is stable under addition of a real constant, finite maxima, decreasing limits, and is the smallest such class of functions that contains all functions of the form $m^{-1}\log|s|$, $s\in\Hnot(X,mL)\setminus\{0\}$. It follows from Dini's lemma that $\Cz(X)\cap\PSH$ is the closure of $\cH$ in the topology of uniform convergence. 
We prove that the restriction of an $L$-psh function $\f$ to $X^\div$ is finite-valued, and determines $\f$. We endow the set $\PSH$ with the topology of pointwise convergence on $X^\div$. 

\end{itemize}
In the main body of the text, we actually work with $\om$-psh functions, where $\om\in\Num(X)$ is any (real) ample (and sometimes even arbitrary) numerical class on $X$. However, here we will stick to the above setting, for simplicity.

As in~\cite{BHJ1}, the \emph{Monge--Amp\`ere energy} $\en\colon\cH\to\R$ is defined on $\f\in\cH$ as the normalized intersection number (or height)
$$
\en(\f):=\frac{(\bar\cL^{n+1})}{(n+1)V}, 
$$
where $n=\dim X$, $V=(L^n)$, and $(\bar\cX,\bar\cL)\to\P^1$ is the canonical $\Gm$-equivariant compactification of the ample, integrally closed test configuration $(\cX,\cL)\to\A^1$ that represents $\f$. The above normalization ensures that $\en(\f+c)=\en(\f)+c$ for $c\in\Q$. 

Adapting to our setting the original approach of~\cite{CL06,Gub10}, we also attach to $\f$ its \emph{Monge--Amp\`ere measure} $\MA(\f)$, a Radon probability measure on $X^\an$ with finite support in $X^\div$, defined using intersection numbers computed on $\cX_0$. We then have
$$
\frac{d}{dt}\bigg|_{t=0}\en((1-t)\f+t\p)=\int_{X^\an}(\p-\f)\MA(\f)
$$
for all $\p\in\cH$, and this actually characterizes the measure $\MA(\f)$, since $\cH$ spans the dense subset $\PL(X)\subset\Cz(X)$.

This formula shows that $\en$ is increasing on $\cH$, and we canonically extend it by monotonicity to a usc functional $\en\colon\PSH\to\R\cup\{-\infty\}$, continuous along decreasing nets. We denote by
$$
\cE^1:=\left\{\f\in\PSH\mid\en(\f)>-\infty\right\}
$$
the set of $L$-psh functions of finite energy. Thus $\f\in\cE^1$ iff $\f\colon X^\an\to\R\cup\{-\infty\}$ is a decreasing limit of functions $\f_i\in\cH$ with $\en(\f_i)$ bounded.
Note that the complex analytic analogue of $\cE^1$ has been well studied, see~\cite{GZ2,BBGZ,BBEGZ,Dar15,DDL18}.

The \emph{weak topology} of $\cE^1$ is the subset topology from $\PSH$ (\ie the topology of pointwise convergence on $X^\div$), and the \emph{strong topology} is the coarsest refinement of the weak topology that makes $\en\colon\cE^1\to\R$ continuous. Decreasing nets in $\cE^1$ are strongly convergent, and $\cH$ is thus strongly dense in $\cE^1$. 

As in~\cite{BBGZ}, we dually introduce an energy functional $\en^\vee\colon\cM\to[0,+\infty]$ on the space $\cM$ of Radon probability measures $\mu$ on $X^\an$, by setting
$$
\en^\vee(\mu):=\sup_{\f\in\cH}\left(\en(\f)-\int\f\,\mu\right)=\sup_{\f\in\cE^1}\left(\en(\f)-\int\f\,\mu\right),
$$
where the right-hand equality follows from a simple approximation argument. Thus $\en^\vee$ is convex, and lsc with respect to the usual weak topology of $\cM$. The space
$$
\cM^1:=\left\{\mu\in\cM\mid\en^\vee(\mu)<\infty\right\}
$$
of measures of finite energy has a weak topology inherited from $\cM$, and a \emph{strong topology}, defined as the coarsest refinement of the weak topology for which $\en^\vee$ becomes continuous. In contrast to $\cE^1$, we prove that both $\cM^1$ and its strong topology are independent of $L$. 

For any two $\f,\f'\in\cH$, we further have $\en(\f)-\en(\f')\ge\int(\f-\f')\MA(\f)$, a reflection of the concavity of $\en$. This precisely means that $\f$ computes the supremum defining the energy of $\MA(\f)$, which thus lies in $\cM^1$. 

The main contribution of the present article can be summarized as follows. 

\begin{thmA} The Monge--Amp\`ere operator $\MA\colon\cH\to\cM^1$ admits a unique continuous extension $\MA\colon\cE^1\to\cM^1$, where both sides are equipped with the strong topology. 

It further induces a topological embedding with dense image $\cE^1/\R\hookrightarrow\cM^1$, which is onto if $X$ is smooth,  and either $\charac k=0$ or $\dim X\le 2$. 
\end{thmA}
Theorem~A can be viewed as a trivially valued analogue of the main result of~\cite{BBGZ}, which itself is a version of the celebrated result of Yau~\cite{Yau} and a later version by Ko{\l}odziej~\cite{kolodziej1}. In the non-Archimedean setting, earlier results include~\cite{Liu11,nama,YZ17}.

\medskip
Theorem A actually brings together several main steps, of various flavors, that we now proceed to describe. 

\subsubsection*{Monotone extension of the Monge--Amp\`ere operator to $\cE^1$}

In a first step, we prove that the Monge--Amp\`ere operator $\MA\colon\cH\to\cM^1$ admits a unique extension $\MA\colon\cE^1\to\cM^1$ that is continuous along decreasing nets (with respect to the weak topology of $\cM^1$). 

Since any function in $\cE^1$ is, by definition, the limit of a decreasing net in $\cH$, uniqueness is clear. Our proof of existence is rather different from previously used approaches~\cite{GZ2,BEGZ,nama}, and proceeds via a direct monotonicity argument. It is phrased in terms of a multilinear \emph{energy pairing}
$$
(L_0,\f_0)\inter(L_n,\f_n)\in\Q
$$
defined on tuples $(L_i,\f_i)$ with $L_i\in\Pic(X)_\Q$ and $\f_i\in\PL(X)$ (and which basically amounts to the Deligne pairing, see Remark~\ref{rmk:enDeligne} and compare~\cite[Definition 2.1]{SD} in the K\"ahler case), but we reformulate it here in a language that will perhaps be more familiar to some readers. The key point is to observe that, for any $\p\in\cH$, the functional $F_\p\colon\cH\to\R$ defined by
$$
F_\p(\f):=(n+1)\en(\f)+\int(\p-\f)\MA(\f)
$$
(which is equal to $V^{-1}(L,\p)\cdot(L,\f)^n$ in terms of the energy pairing) is increasing, and further satisfies 
$$
F_\p(\f)-F_0(\f)=\int\p\MA(\f),\quad F_\p(\f)\ge (n+1)\en(\f)+\inf\p-\sup\f.
$$
Like the energy $\en$, $F_\p$ can thus be monotonically extended to a finite-valued functional $F_\p\colon\cE^1\to\R$, continuous along decreasing nets. Since $\cH$ spans the dense subpace $\PL(X)\subset\Cz(X)$, it is then easy to infer the existence, for each $\f\in\cE^1$, of a unique measure $\MA(\f)\in\cM$ such that $\int\p\MA(\f)=F_\p(\f)-F_0(\f)$ for all $\p\in\cH$. 
 
 This provides an extension $\MA\colon\cE^1\to\cM$ that is continuous along decreasing nets. As above, the concavity of $\en$ further implies, for each $\f\in\cE^1$, that $\en^\vee(\MA(\f))=\en(\f)-\int\f\MA(\f)$, and hence $\MA(\f)\in\cM^1$.

\subsubsection*{Maximizing sequences}
Pick $\mu\in\cM^1$, and consider a sequence $(\f_i)$ in $\cE^1$ that computes $\en^\vee(\mu)=\sup_{\f\in\cE^1}\left(\en(\f)-\int\f\,\mu\right)$; we call this a \emph{maximizing sequence for $\mu$}. In the complex analytic case, it follows \emph{a posteriori} from the variational approach developed in~\cite{BBGZ} that $\MA(\f_i)\to\mu$ weakly in $\cM$. Here we show this directly, by relying on a uniform differentiability property of a natural monotone extension $\en^\downarrow$ of $\en$ to arbitrary usc functions. More precisely, we infer from a key estimate in~\cite{BGM} that we have, for all $\f\in\cE^1$, $\p\in\PL(X)$ and $\e>0$, 
$$
\en^\downarrow(\f+\e\p)=\en(\f)+\e\int\p\MA(\f)+O(\e^2),
$$
where the estimate is uniform with respect to $\f$. We can then apply the variational argument of~\cite{BBGZ} to infer that $\int\p\MA(\f_i)\to\int\p\,\mu$, which yields weak convergence $\MA(\f_i)\to\mu$. 

\subsubsection*{Quasimetrics}
Following~\cite{BBGZ,BBEGZ}, the next step is to show that Aubin's functional
$$
\ii(\f,\f'):=\int(\f-\f')(\MA(\f')-\MA(\f))
$$ 
satisfies a quasi-triangle inequality, and descends to a quasi-metric on $\cE^1/\R$ that defines the strong topology. As in~\cite{BBEGZ}, the quasi-triangle inequality is based on fairly sophisticated estimates obtained by repeated applications of a Cauchy--Schwarz inequality for the energy pairing. But a new feature here is that we also dually introduce 
$$
\ii^\vee(\mu,\mu'):=\inf_{\f\in\cE^1}\left(\jj_\mu(\f)+\jj_{\mu'}(\f)\right)
$$
on $\cM^1$, where $\jj_\mu(\f):=\en^\vee(\mu)-\en(\f)+\int\f\,\mu\ge 0$ tends to zero precisely along maximizing sequences for $\mu$. We prove that the Monge--Amp\`ere operator is bi-Lipschitz with respect to $\ii$ and $\ii^\vee$, \ie $\ii^\vee(\MA(\f),\MA(\f'))\approx\ii(\f,\f')$, and also that $\ii^\vee(\MA(\f_i),\MA(\f'_i))\to\ii^\vee(\mu,\mu')$ if $(\f_i)$, $(\f'_i)$ are maximizing sequences for $\mu$, $\mu'$. This allows us to show that $\ii^\vee$ is a quasi-metric on $\cM^1$ that defines the strong topology, and hence that $\MA\colon\cE^1/\R\to\cM^1$ is a topological embedding with dense image.

  In~\cite{nakstab1,nakstab2} we will show that $\cE^1$ and $\cM^1$ can be equipped with natural metrics (rather than quasimetrics) that define the strong topologies, in such a way that the induced (pseudo)metric on $\cE^1/\R$ is a metric, and $\MA\colon\cE^1/\R\to\cM^1$ is an isometry. The space $\cM^1$ and subspaces thereof will play a key role for the approach to K-stability in~\cite{nakstab2}.

\subsubsection*{  The envelope property}  
In the complex analytic case, it is a basic fact that the (usc) envelope $\supstar_i\f_i$ of any family $(\f_i)$ of psh functions on a (smooth) complex manifold remains psh. This remains true for psh functions on any complex space that is locally irreducible in the analytic topology, but fails in general without this assumption. In our setting, we say that $L$   has the \emph{envelope property} if the previous property holds.   We prove that it implies that $X$ is \emph{unibranch}, which means that the normalization $\nu\colon X^\nu\to X$ is a homeomorphism (in the Zariski topology), and is equivalent to $X$ being locally analytically irreducible in the complex case. Conversely, we conjecture that   the envelope property   holds as soon as $X$ is unibranch (\eg normal), and prove that it holds when $X$ is smooth and either $\charac k=0$ or $\dim X\le 2$, by adapting arguments from~\cite{siminag,GJKM17}. 

We further show that   the envelope property   is equivalent to the analogues of several fundamental known facts in the complex analytic case:
\begin{itemize}
\item weak compactness of $\PSH/\R$;
\item completeness of $\cE^1/\R$ with respect to the quasi-metric $\ii$; 
\item surjectivity of $\MA\colon\cE^1\to\cM^1$.
\end{itemize}
 
It is also equivalent to the property that
for any $\f\in\Cz(X)$, the \emph{psh envelope}
\begin{equation*}
\env(\f):=\sup\{\p\in\PSH, \p\le\f\}
\end{equation*}
is continuous (or, equivalently, $L$-psh).
 
%
%
%
%

\subsection*{Valuations of linear growth}
The results mentioned so far have been quite similar to the complex analytic picture; we now focus on some aspects that are specific to the non-Archimedean setting. We say that a subset $E\subset X^\an$ is \emph{pluripolar} if $E\subset\{\f=-\infty\}$ for some $\f\in\PSH$, this condition being independent of $L$. As opposed to the complex analytic case, a point $v\in X^\an$ can be nonpluripolar. This is for instance the case for a divisorial valuation $v\in X^\div$, and we show more generally that a point $v\in X^\an$ is nonpluripolar iff $v$ is a \emph{valuation of linear growth} in the sense of~\cite{BKMS}, \ie $v(s)\le Cm$ for all $s\in\Hnot(X,mL)\setminus\{0\}$, for a uniform constant $C>0$.

We turn the set $X^\lin$ of valuations of linear growth into a metric space---later shown to be complete---by setting 
$$
\dd_\infty(v,w):=\sup\left\{m^{-1}|v(s)-w(s)|\mid m\ge 1,\,s\in\Hnot(X,mL)\setminus\{0\}\right\}.
$$
We refer to the metric space topology of $X^\lin$ as the \emph{strong topology}, the \emph{weak topology} being the one inherited from $X^\an$. This interacts nicely with the space $\cM^1$, as follows: 

\begin{thmB}   A point $v\in X^\an$ lies in $X^\lin$ iff the Dirac mass $\d_v$ lies in $\cM^1$. Furthermore, the map $v\mapsto\d_v$ defines a closed embedding $X^\lin\hookrightarrow\cM^1$, with respect to both the weak and the strong topologies (on both sides),   onto the set of extremal points of the convex set $\cM^1$. 
\end{thmB}
Thus a net $(v_i)$ in $X^\lin$ converges strongly to $v$ iff $v_i\to v$ weakly and $\en^\vee(\d_{v_i})\to\en^\vee(\d_v)$. As we shall see in the companion paper~\cite{nakstab1}, the energy $\en^\vee(\d_v)$ coincides with the \emph{expected vanishing order} $\mathrm{S}(v)$~\cite{FO16,BlJ}, an invariant that appears in the definition of the stability threshold (or $\d$-invariant) that has come to play a key role in recent works on K-stability of Fano varieties, such as~\cite{BX,BLZ} to name just a few.

Our final main result analyzes the behavior of $L$-psh functions on $X^\lin$. 

\begin{thmC} Pick any $L$-psh function $\f$. Then: 
\begin{itemize}
\item[(i)] the restriction of $\f$ to $X^\lin$ is strongly continuous, and even $1$-Lipschitz with respect to $\dd_\infty$; 
\item[(ii)] the restriction of $\f$ to any $\dd_\infty$-bounded subset of $X^\lin$ is weakly continuous; 
\item[(iii)] if $\f\in\cE^1$, then $\f$ has sublinear growth with respect to $\dd_\infty$, \ie $|\f|\le A \dd_\infty^{1-\e}+B$ on $X^\lin$ for some $A,B>0$ and $0<\e<1$; 
\item[(iv)] if $\f\in\cE^1$ and $\MA(\f)$ is supported in a $\dd_\infty$-bounded subset of $X^\lin$, then $\f\in\Cz(X)$. 
\end{itemize}
\end{thmC}
  In~(iii), $\dd_\infty$ denotes the distance to any fixed point $v_0\in X^\lin$. The assumption in (iv) holds, for instance, if $\MA(\f)$ has finite support in $X^\div$. More generally, as   we show in Appendix~A, dual complexes of snc test configurations (when $X$ is smooth and $\charac k=0$) provide strongly compact subsets of $X^\lin$, and (iv) therefore applies to Monge--Amp\`ere measures with support in such a dual complex, as in~\cite{nama}.
%
%
%
%
\subsection*{Relation to other works and outlook}
This article is in part a continuation of our joint work with T.~Hisamoto~\cite{BHJ1}, which, along the work of K.~Fujita's~\cite{FujitaVolume,FujitaValcrit}  and C.~Li~\cite{Li17,Li18} , first emphasized the role of valuations in the study of K-stability. Besides the clear influence from previous works in the complex setting, especially~\cite{GZ2,BBGZ,BBEGZ}, this article owes a great debt to our joint work with C.~Favre~\cite{siminag,nama}. Inspired in part by the local analysis of~\cite{valtree,valmul,hiro}, and by the unpublished work of M.~Kontsevich and Y.~Tschinkel~\cite{KoTs}, it paved the way to non-Archimedean pluripotential theory, for smooth varieties over a discretely valued field of residue characteristic $0$. These developments built upon the notion of a semipositive continuous metric on a line bundle, as developed by S.-W. Zhang~\cite{Zha95}, Bloch--Gillet--Soul\'e~\cite{BGS}, Gubler~\cite{gublerlocal, gublerequi}, Chambert--Loir~\cite{CL06}, and others.

In the last few years, a number of works on non-Archimedean pluripotential theory have appeared, including~\cite{BGJKM20,GJKM17,GM,BE,BGM}, as well as the first version of this paper~\cite{trivval}. 
Among these,~\cite{BE,BGM} work over arbitrary non-Archimedean fields, including the trivially valued case, and thus have some amount of overlap with the present work. As mentioned above, the main result of~\cite{BGM} in fact plays a key role in our analysis of maximizing sequences. We conversely expect all results in the present article to extend (once properly formulated) to the case of an arbitrary non-Archimedean ground field. This is work in progress, but some results in this generality have in fact already appeared in~\cite{trivval}, and R.~Reboulet has initiated an interesting study of the metric geometry of $\cE^1$ in this generality~\cite{Reb20,Reb21}, perhaps more closely related to the companion paper~\cite{nakstab1} of the present work. 

As opposed to our approach to pluripotential theory, which adopts a global definition of psh functions justified by analogous results in the complex case, it is important to point out that a theory of a local nature is also emerging, thanks to the foundational work of A.~Chambert-Loir and A.~Ducros~\cite{CLD}; see also~\cite{GK17}. The one-dimensional situation was studied in detail in the thesis of A.~Thuillier~\cite{thuillierthesis} (see also~\cite{BRBook}).
In our theory, functions in $\PSH$ are defined as decreasing limits of nicer functions, namely those in $\cH$. In the complex setting, the global notions rely on the local ones, and what is a definition in the trivially valued case becomes an important theorem in the complex situation, see~\cite{Dem92,BK}.

In a different direction, the approach followed in this paper is likely to carry over to the general case of $(1,1)$-classes on compact K\"ahler manifolds, where the basic formalism of test configurations and K-stability was introduced in~\cite{DR,SD}. In that case, Berkovich analytification is of course not available, but a good replacement for it can be constructed as a limit of dual complexes of snc test configurations, as in~\cite[Section 4]{konsoib}. 
%
%
%
%
\subsection*{Structure of the paper}
This article is organized as follows. 

\begin{itemize}
\item Besides recalling some basic facts on Berkovich analytification, test configurations and valuations of linear growth, Section 1 extends the relation between divisorial valuations and test configurations that was drawn in~\cite{BHJ1} to possibly non-normal varieties, emphasizing the role of   what we call \emph{integrally closed}   test configurations. 

\item Section 2 introduces Fubini--Study and PL functions, and describes their relation to test configurations. `Almost trivial' test configurations are also revisited from this perspective (fixing, in particular, a minor issue in~\cite{BHJ1}). 

\item In Section 3 we define and study plurisubharmonicity for PL functions, and undertake a thorough investigation of the energy pairing, for which various estimates are derived from a basic Cauchy--Schwarz inequality. 

\item   Sections 4 and 5 are devoted to general psh functions. In the former section, we  establishing some basic properties and introduce the notion of pluripolar sets. In the latter, we make a detailed study of envelopes of psh functions, and the associated notion of negligible sets. We prove that divisorial points are negligible.

\item   In Section 6 we study psh functions that are homogeneous with respect to the scaling action of $\R_{>0}$. The relation with (nef) $b$-divisors is discussed, and a general decomposition theorem for psh functions is established.   

\item Section   7  extends the energy pairing to psh functions by monotonicity, and studies mixed Monge--Amp\`ere measures and various energy functionals on functions of finite energy. 

\item    In Section 8 we extend the Monge--Amp\`ere energy functional to more general functions, and prove a crucial uniform differentiability result, using~\cite{BGM}.  

\item  Section 9 is devoted to the Monge--Amp\`ere energy of a Radon probability measure, and a preliminary study of the space $\cM^1$ of  measures of finite energy. In particular, we introduce the important notion of a maximizing net for a measure. 
  
\item In Section   10,   we introduce the strong topology on $\cM^1$, and prove that it is defined by the quasi-metric $\ii^\vee$, with respect to which it is complete. We further show that $\cM^1$ and its strong topology are independent of the choice of polarization, when $X$ is irreducible. 

\item Section   11   is devoted to the space $X^\lin$ of valuations of linear growth. It establishes Theorem B above, as well as (most of) Theorem C. 

\item In Section   12   we turn to the strong topology and quasi-metric for $\cE^1$. We complete the proof of Theorem A, and also investigate the continuity of solutions to Monge--Amp\`ere equations, completing the proof of Theorem C. 

\item Assuming the   envelope property,   Section  13   endeavors a detailed study of the Bedford--Taylor capacity, and proves that negligible sets are pluripolar. 

\item Finally, Appendix A adapts to the trivially valued setting the well-known description of the Berkovich analytification as a limit of dual complexes, while Appendix B provides a condensed description of various objects considered in the paper in the case of toric varieties. 

\end{itemize}
%
%
\subsection*{Acknowledgement}
  We thank R.~Berman, B.~Berndtsson, A.~Chambert-Loir, A.~Ducros, C.~Favre, M.~Florence, 
  W.~Gubler, T.~Hisamoto, M.~Maculan, M.~Musta\c{t}\u{a} and J.~Poineau for fruitful discussions.
  The first author was partially supported by the ANR grant GRACK\@. 
  The second author was partially supported by NSF grants DMS-1600011 and DMS-1900025,  
  the Knut and Alice Wallenberg foundation, 
  and the United States---Israel Binational Science Foundation.
  Part of this work was carried out at IMS, Singapore.
   
  Finally we thank the referees for a careful reading of the manuscript, with many corrections and useful suggestions.
   
%
%
\subsection*{Notation and conventions}
\begin{itemize}
\item  We use the standard abbreviations \emph{usc} for `upper semicontinuous', \emph{lsc} for `lower semicontinuous', \emph{wlog} for `without loss of generality'', and \emph{iff} for `if and only if'. We also use \emph{PL} for `piecewise linear',  \emph{snc} for `simple normal crossing', and \emph{psh} for `plurisubharmonic'. 
\item Following the Bourbaki convention, all compact and locally compact topological spaces are required to be Hausdorff. 
\item A \emph{net} in a set $X$ is a family $(x_i)_{i\in I}$ of elements of $X$ indexed by a directed set, \ie a partially preordered set in which any two elements are dominated by a third one.   A \emph{subnet} of $(x_i)_{i\in I}$ is a net of the form $(x_{\f(j)})_{j\in J}$ where $\f\colon J\to I$ is increasing and final, \ie for each $i\in I$ there exists $j\in J$ with $\f(j)\ge i$. A (Hausdorff) topological space $X$ is compact iff every net in $X$ admits a convergent subnet.   
\item If $X$ is a Hausdorff topological space, and $\f\colon X\to\R\cup\{\pm\infty\}$ is any function, then the \emph{usc regularization} $\f^\star$ of $\f$ is the smallest usc function with $\f^\star\ge\f$. Concretely, $\f^\star(x)=\limsup_{y\to x}\f(y)$. The \emph{lsc regularization} is defined by $\f_\star=-(-\f)^\star$. 
   
\item We work over an algebraically closed field $k$, for the most part of arbitrary characteristic. In this paper, a \emph{variety} (over $k$) is a separated $k$-scheme of finite type that is reduced, but not necessarily irreducible, nor even equidimensional. 
\item  For any variety $X$, we denote by $\Pic(X)$ the Picard group of isomorphism classes of line bundles. Elements of the associated $\Q$-vector space $\Pic(X)_\Q:=\Pic(X)\otimes_\Z\Q$ can be viewed as isomorphism classes of $\Q$-line bundles.  
\item If $X$ is a projective variety, we denote by $\Num(X)$ the finite dimensional $\R$-vector space of numerical classes of $\R$-Cartier divisors on $X$.   It comes with a surjective linear map $\Pic(X)_\R\to\Num(X)$ induced by  $L\mapsto c_1(L)$ . 
\item Ample classes form a nonempty open convex cone $\Amp(X)\subset\Num(X)$, whose closure is the closed convex cone $\Nef(X)\subset\Num(X)$ of nef classes. For $\theta,\theta'\in\Num(X)$, we write $\theta\ge\theta'$ if $\theta-\theta'$ is nef.   We generally denote by $\theta$ an element of $\Num(X)$, and by $\om$ an element of $\Amp(X)$. 

\item The cone $\Psef(X)\subset\Num(X)$ of  \emph{pseudoeffective classes} is defined as the closed convex cone generated by classes of effective $\R$-Cartier divisors; its interior $\mathrm{Big}(X)$ is the cone of \emph{big classes}. 
\item A section of a line bundle on $X$ is \emph{regular} if it does not vanish identically along any irreducible component of $X$. Its zero scheme is then a Cartier divisor on $X$.  
 
\item An \emph{ideal} $\fa$ on $X$ is a coherent ideal sheaf $\fa\subset\cO_X$, and similarly for fractional ideals.
\item We denote by $\unipar$ the coordinate of the affine line $\A^1=\Spec k[\unipar]$ and of the torus $\Gm=\Spec k[\unipar^{\pm 1}]$. 
\item For $x,y\in\R_+$, $x\lesssim y$ means $x\le C_n  y$ for a constant $C_n>0$ only depending on $n$, and $x\approx y$ if $x\lesssim y$ and $y\lesssim x$. Here $n$ will be the dimension of a fixed variety $X$ over $k$.
%
%

\item A \emph{quasi-metric} on a set $Z$ means a function $d\colon Z\times Z\to\R_+$ that is symmetric, separates points, and satisfies the \emph{quasi-triangle inequality}
$$
d(x,y)\le C\left(d(x,z)+d(z,y)\right)
$$
for some constant $C>0$. This is equivalent to requiring the quasi-ultrametric inequality 
$d(x,y))\le C\max\{d(x,z),d(z,y)\}$ for some other constant $C>0$, and $d^\a$ is then also a quasi-metric for any $\a\in\R_{>0}$. A quasi-metric space $(Z,d)$ comes with a Hausdorff topology, and even a uniform structure. In particular, Cauchy sequences and completeness make sense for $(Z,d)$. Such uniform structures have a countable basis of entourages, and are thus metrizable, by general theory. A subset $E\subset Z$ is \emph{bounded} if $d$ is bounded on $E\times E$. 
\end{itemize}
%
%
%
%
%
%
\section{Berkovich analytification and test configurations}\label{sec:berk}
In what follows, $X$ denotes a (possibly reducible) projective variety over $k$ (see the conventions above). The main purpose of this section is to review the relation between the Berkovich analytification of $X$ and test configurations, following the approach of~\cite{BHJ1}. 
%
%
%
%
\subsection{The Berkovich analytification}\label{sec:Berk}
 Here we note some facts about the Berkovich analytification of $X$ with respect to the trivial absolute value on $k$.  
\subsubsection{  The space of valuations }
Assume first that $X$ is irreducible, with function field $k(X)$. In this paper, a \emph{valuation} on $X$ means a real-valued valuation $v:k(X)^\times\to\R$, trivial on $k$. The \emph{trivial valuation} $v_\triv=v_{X,\triv}\in X^\val$ is defined by $v_\triv(f)=0$ for all $f\in k(X)^\times$. 

We denote by $X^\val$ the space of valuations on $X$, endowed with the topology of pointwise convergence on $k(X)^\times$.  As soon as $\dim X\ge 1$, $X^\val$ is not locally compact (see \S\ref{sec:curve} below).   

By the valuative criterion of properness, each valuation $v$ on $X$ admits a \emph{center} $c(v)=c_X(v)\in X$, characterized as the unique (scheme) point $\xi\in X$ such that $v\ge 0$ on the local ring $\cO_{X,\xi}$ and $v>0$ on its maximal ideal. Note that $c(v)$ is the generic point of $X$ iff $v=v_\triv$. 

In this paper, a \emph{divisorial valuation} $v$ on $X$ is a valuation of the form $v=t\ord_E$, where $t\in\Q_{>0}$ and $E$ is a prime divisor on a normal, projective birational model $X'\to X$. The center of $v$ on $X$ is then the generic point of the image of $E$ in $X$. It will be convenient to also count the trivial valuation $v_\triv=\lim_{t\to 0} t\ord_E$ as a divisorial valuation.  We write $X^\div$ for the set of divisorial valuations on $X$.  
\subsubsection{  The Berkovich analytification }
Returning to the case of a possibly reducible variety $X$, we denote by $X^\an$ the Berkovich analytification of $X$ with respect to the trivial absolute value on $k$, as in~\cite{BerkBook}. For our purposes, it will be sufficient to view $X^\an$ as a compact\footnote{Recall that all compact spaces are required to be Hausdorff in this paper.} topological space, whose points can be interpreted as \emph{semivaluations} on $X$, \ie valuations $v$ on some irreducible subvariety $Y\subset X$, called the \emph{support} of $v$ and denoted $s(v)$. As a set, we thus have $X^\an=\coprod_Y Y^\val$ with  $Y$  ranging over all irreducible subvarieties of $X$, and the topology of $X^\an$ is the coarsest one such that for each (Zariski) open $U\subset X$ we have: 

\begin{itemize}
\item the set $U^\an\subset X^\an$ of semivaluations whose support meets $U$ is open; 
\item for each $f\in\cO(U)$, the function $|f|\colon U^\an\to\R_+$ defined by $|f|(v):=\exp(-v(f))$ is continuous. 
\end{itemize}
Sets of the form $U^\an$ are open for the \emph{Zariski topology} of $X^\an$.

The Berkovich analytification is functorial. Any morphism $h\colon Y\to X$ of varieties induces a continuous map $h^\an\colon Y^\an\to X^\an$. For simplicity, we will write $h$ instead of $h^\an$. The analytification functor satisfies various GAGA properties, see~\cite[\S3.4]{BerkBook}. For example, if $(X^\b)_\b$ are the connected components of $X$, then $X^\an=\coprod_\b( X^\b)^\an$, and each $(X^\b)^\an$ is connected.
\subsubsection{  Valuations and divisorial valuations }
Mapping $v\in X^\an$ to the generic point $s(v)$ of its support defines a continuous map 
$$
s\colon X^\an\to X.
$$
On the other hand, mapping $v$ to its center $c(v)\in X$ defines a map 
$$
c\colon X^\an\to X,
$$
which is this time \emph{anticontinuous}, in the sense that $c^{-1}(U)\subset U^\an$ is closed (hence compact) for any open subset $U\subset X$ (in the language of Berkovich geometry, $c^{-1}(U)$ is actually a $k$-analytic domain of $X^\an$). 

We say that $v\in X^\an$ is a \emph{valuation} on $X$ (as opposed to a semivaluation) if its support is an irreducible component of $X$,   \ie $s(v)$ is a generic point of $X$. Denoting by $(X_\a)$ the irreducible components of $X$, the set of valuations can be written as
$$
X^\val=\coprod_\a X_\a^\val.
$$
 
We define the set of \emph{divisorial valuations} on $X$ as 
$$
X^\div:=\coprod_\a X_\a^\div. 
$$
It is a dense subset of $X^\an$ (see Theorem~\ref{thm:divdense}).

\begin{rmk}\label{rmk:countable} Assume $\dim X\ge 1$. When $k$ is countable, the set $X^\div$ is countable as well. The compact space $X^\an$ is thus separable, and hence metrizable. This fails when $k$ is uncountable, as $X^\an$ is not even first countable in that case. However, it nevertheless follows from~\cite{Poi} that any point lying in the closure of a subset $E\subset X^\an$ is the limit of a sequence in $E$; in particular, any closed subset of $X^\an$ is sequentially compact. 
\end{rmk}

Given any ideal $\fb\subset\cO_X$ and $v\in X^\an$ with center $c(v)\in X$, one sets
\begin{equation}\label{equ:videal}
v(\fb):=\min\{v(f)\mid f\in\fb_{c(v)}\}\in [0,+\infty], 
\end{equation}
the minimum being achieved among any given set of generators of $\fb_{c(v)}$. Denoting by $Z\subset X$ the zero locus of $\fb$, we have 
$$
v(\fb)>0\Longleftrightarrow c(v)\in Z,\quad v(\fb)=+\infty\Longleftrightarrow s(v)\subset Z.
$$
For any two ideals $\fb,\fb'$, we have
\begin{equation}\label{equ:videal2}
v(\fb\cdot\fb')=v(\fb)+v(\fb'),\quad v(\fb+\fb')=\min\{v(\fb),v(\fb')\}. 
\end{equation} 
The map $v\mapsto v(\fb)$ is continuous on $X^\an$, and such functions generate the topology of $X^\an$. In fact, denoting by $\cI$ the set of ideals of $X$, it is easy to check that $X^\an$ can be identified with the `tropical spectrum' of $\cI$, \ie the space of all functions $\chi\colon \cI\to[0,+\infty]$ that satisfy~\eqref{equ:videal2}, endowed with the topology of pointwise convergence (compare~\cite[\S 1.2]{jonmus}). 
\subsubsection{  Partial order and scaling action }
The space $X^\an$ is endowed with a natural partial order relation, most easily described from the tropical spectrum perspective by
$$
v\ge v'\Longleftrightarrow v(\fb)\ge v'(\fb)\text{ for all ideals }\fb\subset\cO_X.
$$
 As in~\cite[Lemma 4.4]{jonmus}, one checks that $v\ge v'$ iff $c(v)\in\overline{\{c(v')\}}$ and $v(f)\ge v'(f)$ for all $f\in\cO_{X,c(v)}$. 

There is also a natural continuous, order preserving \emph{scaling action}
$$
\R_{>0}\times X^\an\to X^\an\quad (t,v)\mapsto tv,
$$
which induces, in turn, an action on functions $\f\colon X^\an\to\R\cup\{\pm\infty\}$ by setting for $t\in\R_{>0}$ and $v\in X^\an$ 
\begin{equation}\label{equ:scale}
(t\cdot\f)(v):=t\,\f(t^{-1} v).
\end{equation}
  Note that $t\cdot\f=\f$ for all $t$ iff $\f$ is \emph{homogeneous}, which means, in this paper, $\f(tv)=t\f(v)$ for all $t>0$ and $v\in X^\an$ . This action of course preserves the set 
$$
\Cz(X):=\Cz(X^\an,\R)
$$ 
of continuous functions $\f\colon X^\an\to\R$. The reason for adding a factor $t$ in~\eqref{equ:scale} will become clear with Proposition~\ref{prop:FS}~(iii) and Theorem~\ref{thm:PSH}~(ii) (see also Lemma~\ref{lem:basechange} for a geometric interpretation in terms of base change).

\subsubsection{Trivial semivaluations}\label{sec:triv} For each closed point $p\in X(k)$, $\{p\}^\val$ consists of a single (trivial) semivaluation $v_{p,\triv}$, and $X(k)$ is thus naturally realized as a subset of $X^\an$. More generally, to each irreducible subvariety $Y\subset X$ is associated the \emph{trivial semivaluation} $v_{Y,\triv}\in X^\an$ with support $Y$. The set $$
X^\triv\subset X^\an
$$
of trivial semivaluations can be identified with the set of functions $\chi\colon\cI\to\{0,+\infty\}$ satisfying~\eqref{equ:videal2} and equals the set of fixed points of the scaling action of $\R_{>0}$ on $X^\an$. One easily checks that $X^\triv$ is the closure of $X(k)\subset X^\an$. 

Any scheme point $\xi\in X$ is the generic point of a subvariety $Y\subset X$, and taking $\xi$ to $v_{Y,\triv}$ defines a bijection $X\simto X^\triv$, which
is a section of both maps $s\colon X^\an\to X$, $c\colon X^\an\to X$. The 
inverse $X^\triv\to X$ is continuous but not a homeomorphism.

\begin{exam}\label{exam:vYZ} Pick $v\in X^\an$, denote by $Y$ its support and by $Z$ the closure of its center $c(v)$. Then $Z\subset Y$, and $v_{Y,\triv}\le v\le v_{Z,\triv}$. Further, 
 $$
 \lim_{t\to 0_+} tv=v_{Y,\triv},\quad \lim_{t\to+\infty} tv=v_{Z,\triv}.
 $$
 \end{exam}

The set $X^{\triv}\cap X^\val$   of trivial valuations is in 1--1 correspondence with the irreducible components $X_\a$ of $X$;   its elements will   be denoted by 
$$
v_{\triv,\a}:=v_{X_\a,\triv}\in X^\div. 
$$
When $X$ is irreducible, there is only one such valuation, denoted by $v_\triv$.

The following condition will arise many a time in this paper. 
\begin{defi}\label{defi:genfin} A function $\f\colon X^\an\to\R\cup\{-\infty\}$ is \emph{generically finite} if $\f|_{X_\a^\an}\not\equiv-\infty\text{ for all }\a$.
\end{defi}

\begin{lem}\label{lem:maxprin} Let $\f\colon X^\an\to\R\cup\{-\infty\}$ be decreasing. Then:
\begin{itemize}
\item[(i)] for any irreducible subvariety $Y\subset X$, $\f|_{Y^\an}$ satisfies the `maximum principle'
\begin{equation}\label{equ:maxprin}
\sup_{Y^\an}\f=\f(v_{Y,\triv}); 
\end{equation} 
\item[(ii)] $\f$ is generically finite iff it is finite at $v_{\triv,\a}$ for each $\a$; 
\item[(iii)] if $\f$ is further usc and $Z\subset X$ denotes the closure of the center of $v\in X^\an$, then 
$$
\f(tv)\searrow\f(v_{Z,\triv})=\sup_{Z^\an}\f
$$
as $t\to+\infty$. 
\end{itemize}
\end{lem}
\begin{proof} For any $v\in Y^\an$, we have $v\ge v_{Y,\triv}$, see Example~\ref{exam:vYZ}. Thus $\f(v)\le\f(v_{Y,\triv})$, which yields (i), and hence (ii). To see (iii), note that for $t\in\R_{>0}$ we have
$$
tv\le v_{Z,\triv}\Longrightarrow\f(tv)\ge\f(v_{Z,\triv})=\sup_{Z^\an}\f,
$$
since $\f$ is decreasing, and $\limsup_{t\to+\infty}\f(tv)\le\f(v_{Z,\triv})$, since $\f$ is usc.
\end{proof}

\subsubsection{  The Berkovich analytification of a curve }\label{sec:curve} 
  Here we describe $X^\an$ in the case when $X$ is a curve, \ie of pure dimension one. We may and will assume $X$ is connected.

If $X$ is smooth,  then $X^\an$ is a star-shaped $\R$-tree rooted at $v_\triv$, see Figure~\ref{fig:curve}. Further, the Berkovich topology coincides with the weak tree topology as defined in~\cite[\S3.1.4]{valtree} or~\cite[\S2.1.1]{dynberko}, for example.

More specifically, each closed point $p\in X(k)$ determines a compactified ray $\iota_p:[0,+\infty]\hookrightarrow X^\an$, with $\iota_p(0)=v_\triv$, $\iota_p(t)=t\ord_p$ for $t\in(0,+\infty)$, and $\iota_p(+\infty)=v_{p,\triv}$. 
  The half-open rays $\iota_p((0,+\infty])=c^{-1}(p)$ form an open partition of $X^\an\setminus\{v_\triv\}$, and we have 
\begin{equation*}
X^\an=X^\val\sqcup X(k),\quad X^\triv=\overline{X(k)}=X(k)\sqcup\{v_\triv\},\quad X^\div=\{v_{\triv}\}\sqcup\coprod_{p\in X(k)}\iota_p(\Q_{>0}).
\end{equation*}

A neighborhood basis of $v_\triv$ is formed by complements of the union of finitely many segments $\iota_{p_j}([t_j,+\infty]$, where $t_j\in\R_{>0}$. In particular, every neighborhood of $v_\triv$ contains at least one (in fact, infinitely many) compactified rays;    this prevents $X^\val$ from being locally closed in $X^\an$, and $X^\val$ is therefore not locally compact.

In the general case, let $(X_\a)$ be the irreducible components of $X$.
Then the normalization morphism $X^\nu\to X$ induces a surjective map
\begin{equation*}
  X^{\nu,\,\an}=\coprod_\a X_\a^{\nu,\an}\to X^\an
\end{equation*}
 
that identifies the endpoints $v_{p_i,\triv}$ of all rays corresponding to a point $p_i\in\nu^{-1}(p)$,  $p\in X(k)$. 
 
\begin{figure}[ht]
  \centering
  \begin{tikzpicture}

    \tikzmath{\scpar = 0.5;}
    \tikzmath{\dang = 60;}
    \tikzmath{\x1 = 0; \y1 =0;}
    \tikzmath{\scale=2; \ang=0; }
    \tikzmath{\dotsize=0.8; }
    \tikzmath{\nobranches=10; }

    \node at (0,0) [circle,fill,inner sep=\dotsize]{};
    \foreach \i in {1,...,\nobranches} {%
      \tikzmath{\ang=360*\i/\nobranches;}
      \tikzmath{\x1 = \scale*cos(\ang); \y1 = \scale*sin(\ang);}
      \draw (0,0)--(\x1,\y1);
      \node at (\x1,\y1) [circle,fill,inner sep=\dotsize]{};
    }
    \node at (1.3,0) [circle,fill,inner sep=\dotsize]{};
    \draw (1.3,0.3) node {$\ord_p$};
   \draw (2.6,0) node {$v_{p,\triv}$};

  \end{tikzpicture}
  \hspace*{1cm}
  \begin{tikzpicture}

    \tikzmath{\scpar = 0.5;}
    \tikzmath{\dang = 60;}
    \tikzmath{\x1 = 0; \y1 =0;}
    \tikzmath{\scale=2; \ang=0; }
    \tikzmath{\dotsize=0.8; }
    \tikzmath{\nobranches=10; }
    \tikzmath{\nobranchesm=8; }

    \node at (0,0) [circle,fill,inner sep=\dotsize]{};
    \foreach \i in {2,...,\nobranchesm} {%
      \tikzmath{\ang=360*\i/\nobranches;}
      \tikzmath{\x1 = \scale*cos(\ang); \y1 = \scale*sin(\ang);}
      \draw (0,0)--(\x1,\y1);
      \node at (\x1,\y1) [circle,fill,inner sep=\dotsize]{};
    }
    \draw (0,0) .. controls (1.0,0.8) .. (2,0);
    \draw (0,0) .. controls (1.0,-0.8) .. (2,0);
    \node at (1.0,0.6) [circle,fill,inner sep=\dotsize]{};
    \node at (1.0,-0.6) [circle,fill,inner sep=\dotsize]{};
    \node at (2.0,0) [circle,fill,inner sep=\dotsize]{};
    \draw (2.6,0) node {$v_{p,\triv}$};
    \draw (1.2,0.9) node {$\ord_{p_1}$};
    \draw (1.2,-0.9) node {$\ord_{p_2}$};

  \end{tikzpicture}
  \caption{The Berkovich analytification of a smooth curve (left) and a nodal curve (right), see \S\ref{sec:curve}.}\label{fig:curve}
\end{figure}

 This simple piecewise linear picture admits a far-reaching generalization: as we shall see in Appendix~\ref{sec:dual}, if $X$ is smooth of dimension $n$ and $\charac k=0$, then $X^\an$ can be written as the projective limit of the family of simplicial complexes (of dimension at most $n$) attached to all snc test configurations. 
%
%
\subsection{Test configurations}\label{sec:tc}
A \emph{test configuration} $\cX$ for $X$ consists of: 
\begin{itemize}
\item[(i)] a flat, projective morphism of schemes $\pi\colon\cX\to\A^1$; 
\item[(ii)] a $\Gm$-action on $\cX$ lifting the canonical action on $\A^1$;  
\item[(iii)] an isomorphism $\cX_1\simeq X$.
\end{itemize}
By~\cite[Proposition 2.6]{BHJ1}, the scheme $\cX$ is reduced, and hence a variety. The central fiber $\cX_0$ is a principal Cartier divisor, defined by the regular function $\pi^\star\unipar$, with $\unipar$ denoting the coordinate on $\A^1$. The open set $\cX\setminus\cX_0$ is Zariski dense, and (iii) amounts to the data of a $\Gm$-equivariant isomorphism 
\begin{equation}\label{equ:equivisom}
\cX\setminus\cX_0\simeq X\times\Gm
\end{equation}
over $\Gm\subset\A^1$. As a result, $\pi\colon\cX\to\A^1$ admits a canonical $\Gm$-equivariant compactification $\pi\colon\bar\cX\to\P^1$, obtained by simply extending~\eqref{equ:equivisom} to $\bar\cX\setminus\cX_0\simeq X\times(\P^1\setminus\{0\})$ over $\P^1\setminus\{0\}$. Thus $\bar\cX$ is a projective variety, of dimension $\dim X+1$. 

\smallskip

For each subvariety $Y\subset X$, the closure $\cY\subset\cX$ of the image of $Y\times\Gm$ under~\eqref{equ:equivisom} is a test configuration for $Y$. This applies to the irreducible components $X_\a$ of $X$, and induces the irreducible decomposition $\cX=\bigcup_\a\cX_\a$. In particular, $\cX$ is irreducible iff $X$ is. 

\smallskip

If $L$ is a $\Q$-line bundle on $X$, a \emph{test configuration} $(\cX,\cL)$ for $(X,L)$ consists of a test configuration $\cX$ for $X$, a $\Gm$-linearized $\Q$-line bundle $\cL$ on $\cX$, and an identification $(\cX,\cL)_1\simeq (X,L)$ compatible with $\cX_1\simeq X$. We also say that $\cL$ is a test configuration for $L$, determined on $\cX$. We then have a canonical $\Gm$-equivariant isomorphism
\begin{equation}\label{equ:tciso}
(\cX\setminus\cX_0,\cL)\simeq (X,L)\times\Gm, 
\end{equation}
and a canonical extension $\bar\cL$ of $\cL$ to a $\Gm$-linearized $\Q$-line bundle on $\bar\cX$. 

When $\cL$ (and hence $L$) are honest line bundles,~\eqref{equ:tciso} induces an isomorphism of $k[\unipar^{\pm 1}]$-modules
$$
\Hnot(\cX,\cL)_{k[\unipar^{\pm 1}]}\simeq \Hnot(X,L)_{k[\unipar^{\pm 1}]},
$$
which allows to view $\Hnot(\cX,\cL)$ as a $k[\unipar]$-submodule of $\Hnot(X,L)_{k[\unipar^{\pm 1}]}$. 

\begin{exam}\label{exam:vertcar} Test configurations $(\cX,\cL)$ for $(X,\cO_X)$ are in 1--1 correspondence with \emph{vertical $\Q$-Cartier divisors} $D$ on $\cX$, by which we mean $\Gm$-invariant $\Q$-Cartier divisors on $\cX$ with support in $\cX_0$. 
\end{exam}

\begin{exam}\label{exam:trivtc} The \emph{trivial test configuration} $\cX_\triv$ for $X$ is the product 
$X\times\A^1$, with the trivial $\Gm$-action on $X$. If $L$ is a $\Q$-line bundle on $X$, the trivial test configuration $(\cX_\triv,\cL_\triv)$ for $(X,L)$ is defined by $\cL_\triv:=p_1^\star L$, with $p_1\colon \cX_\triv\to X$ the first projection. 
\end{exam}

\begin{exam}\label{exam:defcone} Consider a closed subscheme $Z\subset X$, with ideal $\fb\subset\cO_X$. The blowup $\mu\colon \cX\to\cX_\triv$ of $Z\times\{0\}\subset X\times\{0\}=\cX_\triv$ is a test configuration for $X$, known as the \emph{deformation to the normal cone of $Z$}. The central fiber splits into a sum of two effective Cartier divisors $\cX_0=\tX+P$, where $\tX$ is the strict transform of $X\times\{0\}$, which is isomorphic to the blowup of $Z$ in $X$, and $P$ is the exceptional divisor of $\mu$, which can be identified to the projective compactification $P(C_{Z/X}\oplus 1)$ of the normal cone $C_{Z/X}=\Spec_Z\left(\bigoplus_{m\in\N}\fb^m/\fb^{m+1}\right)$ (see~\cite[\S 5.1]{FultonInter}). 
\end{exam}

Test configurations for $X$ form a category, a morphism $\mu\colon\cX\to\cX'$ being a $\Gm$-equivariant morphism over $\A^1$, compatible with the isomorphisms $\cX'_1\simeq X\simeq \cX_1$. There is at most one morphism $\cX\to\cX'$ between any two given test configurations, and we say that $\cX$ \emph{dominates} $\cX'$ when it exists. Two test configurations that dominate each other are canonically isomorphic, and can thus safely be identified. Any two test configurations can be dominated by a third, and the set of (isomorphism classes of) test configurations for $X$ thus forms a directed poset. 

 \begin{lem}\label{lem:relamp} Let $\cX$ be a test configuration that dominates $\cX_\triv$ via a morphism $\mu\colon\cX\to\cX_\triv$. Then $\cX$ admits a vertical $\Q$-Cartier divisor $D$ that is $\mu$-ample. 
\end{lem}
\begin{proof} Since the structure morphism $\pi\colon\cX\to\A^1$ is $\Gm$-equivariant and projective (by definition of a test configuration), we can pick a $\Gm$-linearized, ample line bundle $\cL$ on $\cX$. Let $L$ be its restriction to $X\simeq\cX_1$, so that $(\cX,\cL)$ is a test configuration for $(X,L)$. If we denote by $L_\cX$ the pullback of $L$ by $\cX\to\cX_\triv\to X$, then $D:=\cL-\mu^\star L_\cX$ is a vertical $\Q$-Cartier divisor on $\cX$, and it is $\mu$-ample.
\end{proof}

%
%
%
%
\subsection{Gauss extension}\label{sec:Gauss}
Each test configuration $\cX$ for $X$ comes with a topological embedding 
$$
\sigma\colon X^\an\hookrightarrow(\cX\setminus\cX_0)^\an\subset\bar\cX^\an,
$$
called \emph{Gauss extension}, with image the set of $k^\times$-invariant semivaluations $w\in\cX^\an$ such that $w(\unipar)=1$ (and hence centered on $\cX_0$), and defined as follows. For each irreducible subvariety $Y\subset X$, the associated test configuration $\cY\subset\cX$ provides a canonical embedding of function fields
$$
k(Y)\subset k(Y)(\unipar)=k(\cY_\triv)\simeq k(\cY). 
$$
The Gauss extension of a valuation $v\colon k(Y)^\times\to\R$ is then defined as the unique valuation $\sigma(v)\colon k(\cY)^\times\to\R$ such that 
$$
\sigma(v)\left(\sum_{d\in\N} f_d\unipar^d\right)=\min_d(v(f_d)+d)
$$
on $k(Y)[\unipar]\hookrightarrow k(\cY)$. 

By~\cite[Lemma 4.2]{BHJ1}, $\sigma$ defines a 1--1 correspondence between $X^\val$ (resp.~$X^\div$) and the set of $k^\times$-invariant valuations (resp.~divisorial valuations) $w$ on $\cX$ such that $w(\unipar)=1$. 

\begin{rmk}\label{rmk:Gausslaurent} Gauss extension is in fact independent of the choice of $\cX$, in the following sense: the base change $X_K$ of $X$ to the non-Archimedean field $K:=k\lau{\unipar}$ admits a Berkovich analytification $X_K^\an$, whose points can again be interpreted as semivaluations $w$ on $X_K$, compatible with the given valuation on $K$, \ie trivial on $k$ and such that $w(\unipar)=1$. Given any test configuration $\cX$ for $X$, we thus have a canonical identification of $X_K^\an$ with the image of $\sigma\colon X^\an\to\cX^\an$, and the corresponding map $\sigma\colon X^\an\to X_K^\an$ is a continuous section of the natural projection $\pi\colon X_K^\an\to X^\an$. 

  For any non-Archimedean field extension $F/k$, \cite[Corollaires 3.7 \& 3.14]{Poi} more generally yields a canonical continuous section $\sigma\colon X^\an\to X_F^\an$ of the projection $\pi\colon X_F^\an\to X^\an$, that takes $v\in X^\an$ to the unique point $\sigma(v)$ in the Shilov boundary of $\pi^{-1}(v)=\cM(\cH(v)\widehat{\otimes} F)$.   
\end{rmk}

%
%
%

\subsection{Integrally closed test configurations and divisorial valuations}\label{sec:intclosed}

\begin{defi} We say that a test configuration $\cX$ for $X$ is \emph{integrally closed} if the scheme $\cX$ is integrally closed in the generic fiber $\cX_{k(\unipar)}$ of $\pi\colon\cX\to\A^1$. 
\end{defi}
In other words, $\cX$ is integrally closed iff it can be covered by affine open subsets $\cU=\Spec\cA$ such that the $k[\unipar]$-algebra $\cA$ is integrally closed in $\cA_{k(\unipar)}$. Note that $\cX$ is normal iff $X$ is normal and $\cX$ is integrally closed. 

\begin{exam} The trivial test configuration $\cX_\triv=X\times\A^1$ is integrally closed. 
\end{exam}

The integral closure of any test configuration $\cX$ in $\cX_{k(\unipar)}$ defines a finite, $\Gm$-equivariant morphism $\tcX\to\cX$ (because the scheme $\cX$, being of finite type over a field, is excellent), which induces an isomorphism on the generic fibers over $\A^1$. Thus $\tcX$ is an integrally closed test configuration for $X$, which we simply call the \emph{integral closure} of $\cX$. It is characterized as the unique integrally closed test configuration that dominates $\cX$ via a finite morphism. 

  When $X$ is normal, the integral closure $\tcX$ of any test configuration $\cX$ coincides with its normalization $\cX^\nu$ (see also Remark~\ref{rmk:normvsint}).

\begin{lem}[Zariski's main theorem]\label{lem:Zarmain} If $\cX$ is integrally closed and $\mu\colon\cX'\to\cX$ is a morphism of test configurations,  then $\mu_\star\cO_{\cX'}=\cO_\cX$. 
\end{lem}
\begin{proof} By coherence of direct images, $\mu_\star\cO_{\cX'}$ is a finite $\cO_\cX$-module. Sections of $\mu_\star\cO_{\cX'}$ on an open $\cU\subset\cX$ are thus rational functions on $\cU$ that are regular on the generic fiber and integral over $\cO_\cU$, and hence belong to $\cO_\cU$.
\end{proof}

Integral closedness admits the following characterization, a `vertical version' of the usual Serre criterion for normality. 

\begin{thm}\label{thm:serrecrit} A test configuration $\cX$ is integrally closed iff it is: 
\begin{itemize}
\item[(i)] \emph{vertically $R_1$}, in the sense that $\cX$ is regular at each generic points of $\cX_0$;
\item[(ii)] \emph{vertically $S_2$}, in the sense that $\depth\cO_{\cX,\xi}\ge\min\{2,\dim\cO_{\cX,\xi}\}$ for all $\xi\in\cX_0$.
\end{itemize}
\end{thm}

\begin{rmk} Condition (i) was called \emph{partially normal} in~\cite[Definition 3.7]{Oda13}. 
\end{rmk}

\begin{lem}\label{lem:S2} For any test configuration $\cX$ for $X$, we have: 

\begin{itemize}
\item[(i)]  $\cX$ is vertically $S_2\Longleftrightarrow\cX_0$ is $S_1$, \ie without embedded points; 
\item[(ii)] $\cX$ is $S_2\Longleftrightarrow X$ is $S_2$ and $\cX$ is vertically $S_2$. 
\end{itemize}
\end{lem}
\begin{proof} Since $\cX_0$ is a Cartier divisor, each $\xi\in\cX_0$ satisfies 
$$
\depth\cO_{\cX,\xi}=\depth\cO_{\cX_0,\xi}+1,\quad\dim\cO_{\cX,\xi}=\dim\cO_{\cX_0,\xi}+1. 
$$
This yields (i), while (ii) is a direct consequence of the isomorphism 
$$
\cX\setminus\cX_0\simeq X\times\Gm.
$$ 
\end{proof}

\begin{cor}\label{cor:redvertnorm} For any test configuration $\cX$, the following are equivalent:
\begin{itemize}
\item[(i)] $\cX_0$ is reduced; 
\item[(ii)] $\cX$ is integrally closed, and $\cX_0$ is generically reduced.
\end{itemize}
\end{cor}
\begin{proof} The scheme $\cX_0$ is reduced iff it is generically reduced and $S_1$. By Lemma~\ref{lem:S2}, $\cX_0$ is thus reduced iff it is generically reduced and $\cX$ is vertically $S_2$. Finally, $\cX_0$ generically reduced implies that $\cX$ is regular at each generic point of $\cX_0$, since the latter is a Cartier divisor. The result is now a consequence of Theorem~\ref{thm:serrecrit}. 
\end{proof}

Before entering the proof of Theorem~\ref{thm:serrecrit}, we introduce some terminology that will be used throughout this paper.
\begin{defi} Let $\cX$ be a test configuration for $X$.
\begin{itemize} 
\item[(i)] A \emph{vertical fractional ideal} $\fa$ on $\cX$ is a coherent fractional ideal sheaf that is $\Gm$-invariant, and trivial outside $\cX_0$. 
\item[(ii)] The \emph{polar scheme} $P\subset\cX$ of a vertical fractional ideal $\fa$ is the closed subscheme of $\cX$ defined by the \emph{ideal of poles} $\fp:=\{f\in\cO_\cX\mid f \fa\subset\cO_\cX\}$.
\item[(iii)] A \emph{vertical Cartier divisor} $D$ on $\cX$ is a $\Gm$-invariant Cartier divisor with support in $\cX_0$.
\end{itemize}
\end{defi}
Note that $P$ in (ii) is supported in $\cX_0$, $\fa$ being trivial outside $\cX_0$. 

Vertical Cartier divisors on $\cX$ are in 1--1 correspondence with locally principal vertical fractional ideals of $\cX$, via $D\mapsto\cO_\cX(D)$. If $D$ is further effective, then it coincides with the polar scheme of $\cO_\cX(D)$. 

\begin{lem}\label{lem:vertnormfrac} A test configuration $\cX$ is integrally closed iff every vertical fractional ideal $\fa$ on $\cX$ that is integral over $\cO_\cX$ satisfies $\fa\subset\cO_\cX$.
\end{lem}
\begin{proof} The `only if' part is obvious. Conversely, consider the integral closure $\mu\colon\tcX\to\cX$. Then $\fa:=\mu_\star\cO_{\tcX}$ is a vertical fractional ideal that is integral over $\cO_\cX$, and the `if part' follows. 
\end{proof}

The next result is the key step in the proof of Theorem~\ref{thm:serrecrit}. 

\begin{lem}\label{lem:vertnorm} Let $\cX$ be an integrally closed test configuration for $X$, and $\fa$ be a vertical fractional ideal on $\cX$. Then:
\begin{itemize}
\item[(i)] every associated point of the scheme of poles of $\fa$ is a generic point of $\cX_0$;
\item[(ii)] $\cX$ is vertically $R_1$. 
\end{itemize}
\end{lem}
\begin{proof} We follow the usual proof of Serre's criterion for normality. Let $\fp=\{f\in\cO_\cX\mid f \fa\in\cO_\cX\}$ be the ideal of poles, and pick an associated point $\xi\in\cX_0$ of the polar scheme $P\subset\cX$. By definition of an associated point, there exists $f\in\cO_{\cX,\xi}$ such that $f\notin\fp_\xi$ but $f\fm_\xi\subset\fp_\xi$, with $\fm_\xi$ the maximal ideal of $\cO_{\cX,\xi}$. Then $f\fm_\xi\cdot\fa_\xi$ is an ideal of $\cO_{\cX,\xi}$, and hence $f\fm_\xi\cdot\fa_\xi\subset\fm_\xi$ or $f\fm_\xi\cdot\fa_\xi=\cO_{\cX,\xi}$. In the former case, the usual determinant trick implies that $f\fa_\xi$ is integral over $\cO_{\cX,\xi}$. Since $\cX$ is integrally closed, we infer $f\fa_\xi\subset\cO_{\cX,\xi}$, \ie $f\in\fp_\xi$, a contradiction. We thus necessarily have $f\fm_\xi\cdot\fa_\xi=\cO_{\cX,\xi}$, which proves that $\fm_\xi$ is invertible, and hence that $\cO_{\cX,\xi}$ is a DVR\@. This proves (i), as well as (ii), taking $\fa=\cO_\cX(\cX_0)$.  
\end{proof}

\begin{proof}[Proof of Theorem~\ref{thm:serrecrit}] Assume that $\cX$ is integrally closed. By Lemma~\ref{lem:vertnorm}, $\cX$ is vertically $R_1$, and the polar scheme $\cX_0$ of $\cO_\cX(\cX_0)$ has no embedded points. By Lemma~\ref{lem:S2}, $\cX$ is thus vertically $S_2$, which proves (i)$\Longrightarrow$(ii). 

Assume, conversely, that $\cX$ is vertically $R_1$ and $S_2$, and pick a vertical fractional ideal $\fa$ on $\cX$ that is integral over $\cO_\cX$. According to Lemma~\ref{lem:vertnormfrac}, we need to show that $\fa\subset\cO_{\cX}$, which means that its ideal of poles $\fp:=\left\{f\in\cO_\cX\mid f\fa\subset\cO_\cX\right\}$ is trivial. Arguing by contradiction, suppose that the polar scheme $P$ is non-empty, and pick an associated point $\xi\in\cX_0$ of $P$. If $\dim\cO_{\cX,\xi}=1$, then $\xi$ is a generic point of $\cX_0$. Since $\cX$ is vertically $R_1$, $\cO_{\cX,\xi}$ is regular, and hence integrally closed. Thus $\fa_\xi\subset\cO_{\cX,\xi}$, which contradicts $\xi\in P$. We thus have $\dim\cO_{\cX,\xi}\ge 2$, and hence $\depth\cO_{\cX,\xi}\ge 2$, since $\cX$ is vertically $S_2$. 

Since $\xi$ is an associated point of $P$, we can find, as above, $f\in\cO_{\cX,\xi}$ such that $f\notin\fp_\xi$ but $f\fm_\xi\subset\fb_\xi$, \ie $f\fm_\xi\cdot\fa_\xi\subset\cO_{\cX,\xi}$. Since $f\notin\fp_\xi$, there exists $g\in\fa_\xi$ such that $h:=fg\notin\cO_{\cX,\xi}$. Write $h=a/b$ with $a,b\in\cO_{\cX,\xi}$ and $b$ a non-zerodivisor. Then $a\notin(b)$, but $a\fm_\xi\subset(b)$. Thus $\xi$ is an associated point of the Cartier divisor $D=(b=0)$, which contradicts $\depth\cO_{D,\xi}\ge 1$.
\end{proof}

Generalizing~\cite[\S4.2]{BHJ1}, we associate to every irreducible component $E$ of the central fiber $\cX_0$ of an integrally closed test configuration $\cX$ a divisorial valuation $v_E\in X^\div$, as follows. 

By Theorem~\ref{thm:serrecrit}, the local ring of $\cX$ at the generic point of $E$ is a DVR (compare~\cite[Lemme 2.1]{CLT}), and hence defines a divisorial valuation $\ord_E:k(\cX)^\times\to\Z$. As in~\cite[Definition 4.4]{BHJ1}, we set $b_E:=\ord_E(\unipar)=\ord_E(\cX_0)$, and define a valuation $v_E$ on $X$ as the restriction of $w_E:=b_E^{-1}\ord_E$ to $$
k(X)\subset k(X)(\unipar)=k(\cX_\triv)\simeq k(\cX).
$$
Since $w_E$ is $k^\times$-invariant and $w_E(\unipar)=1$, we have $\sigma(v_E)=w_E$. 

\begin{lem}\label{lem:tcdiv} A valuation $v$ on $X$ is divisorial iff $v=v_E$ for an irreducible component $E$ of some integrally closed test configuration $\cX$ for $X$. 
\end{lem}
This follows from~\cite[Theorem 4.6]{BHJ1}. While the latter assumes $X$ normal, its proof is a rather simple consequence of a theorem of Zariski~\cite[Lemma 2.45]{KM}, which does not depend on this assumption.

\begin{exam}\label{exam:defconeval} Assume $X$ is smooth, and let $Z\subset X$ be a smooth irreducible subvariety. The vanishing order at the generic point of $Z$ is then a divisorial valuation $\ord_Z\in X^\div$; indeed, denoting by $\pi\colon\tX\to X$ the blowup of $Z$ and $E$ its exceptional divisor, we have $\ord_Z=\ord_E$. Alternatively, $\ord_Z=v_P$ with $P$ the exceptional divisor of the blowup $\cX\to\cX_\triv$ of $Z\times\{0\}$, \ie the deformation to the normal cone of $Z$ in $X$ (see Example~\ref{exam:defcone}). 
\end{exam}

 \begin{rmk}\label{rmk:normvsint} The normalization $\cX^\nu$ of a test configuration $\cX$ for $X$ is a test configuration for the normalization $X^\nu$ of $X$. The normalization morphism $\cX^\nu\to\cX$ factors through the integral closure $\tcX$, and the induced morphism $\cX^\nu\to\tcX$ is an isomorphism over each generic point of $\tcX_0$, $\tcX$ being regular at such points (compare~\cite[Lemma 3.9]{Oda13}). As a consequence, the irreducible components of $\tcX_0$ and those of $\cX^\nu_0$ induce the same set of valuations in $X^\div=(X^\nu)^\div$. 
\end{rmk}

 \begin{lem}\label{lem:CLT} A vertical fractional ideal $\fa$ on an integrally closed test configuration $\cX$ satisfies $\fa\subset\cO_\cX$ iff $\ord_E(\fa)\ge 0$ for each irreducible component $E$ of $\cX_0$. 
\end{lem}

\begin{proof} If $\fa\subset\cO_\cX$, then trivially $w(\fa)\ge 0$ for all $w\in\cX^\an$. Conversely assume $\ord_E(\fa)\ge 0$ for all $E$. We need to show that the polar scheme $P$ of $\fa$ is empty. Suppose this is not the case, and pick an associated point $\xi$ of $P$. By Lemma~\ref{lem:vertnorm}~(ii), $\xi$ is the generic point of some component $E$ of $\cX_0$, and the assumption $\ord_E(\fa)\ge 0$ thus yields $\fa_\xi\subset\cO_{\cX,\xi}$, \ie $\xi\notin P$, a contradiction. 
\end{proof}

%
%
%
\subsection{Valuations of linear growth}\label{sec:lin}
Let $L$ be a line bundle on $X$. A semivaluation $v\in X^\an$ can be naturally evaluated on any section $s\in \Hnot(X,L)$, by defining $v(s)$ as the value of $v$ on the germ in $\cO_X$ corresponding to $s$ in any local trivialization of $L$ at the center of $v$. Thus $v(s)\in[0,+\infty]$, $v(s)=\infty$ iff $s$ vanishes along the support of $v$, and $v(s)>0$ iff $s$ vanishes at the center of $v$. Setting $|s|(v):=e^{-v(s)}$ defines a continuous function $|s|\colon X^\an\to[0,1]$.  

\begin{rmk}\label{rmk:trivmetric} This construction reflects the existence of the \emph{trivial metric} $|\cdot|$ of $L^\an$, a continuous metric characterized by $|\tau|\equiv 1$ on the compact set $c^{-1}(U)$ for any trivializing section $\tau\in \Hnot(U,L)$ on an open subset $U\subset X$. 
\end{rmk}

In what follows, we fix an ample line bundle $L$ on $X$. 

\begin{lem}\label{lem:Zar} A semivaluation $v\in X^\an$ is a valuation iff $v(s)<+\infty$ for all regular sections $s\in \Hnot(X,mL)$, $m\in\N$. 
\end{lem}
\begin{proof} Suppose $v\notin X^\val$, so that its support $Y\subset X$ is not an irreducible component of $X$.  Since $L$ is ample, we can find for $m\gg 1$ a nonzero section $s\in\Hnot(X,mL)$ that vanishes along $Y$, but not along any irreducible component of $X$. Then $s$ is regular, and $v(s)=\infty$. The converse direction is clear. 
\end{proof}

Following~\cite{BKMS,BlJ}, we introduce: 

\begin{defi}\label{defi:TT} The \emph{maximal vanishing order} of $L$ at $v\in X^\an$ is defined as 
$$
\teL(v)=\sup\left\{m^{-1}v(s)\mid m\ge 1,\, s\in \Hnot(X,mL)\text{ regular}\right\}\in[0,+\infty]. 
$$
We say that $v$ has \emph{linear growth} if $\teL(v)<+\infty$, and denote by $X^\lin$ the corresponding subset of $X^\an$. 
\end{defi}
The linear growth condition is easily seen to be independent of the choice of an ample line bundle $L$. By Lemma~\ref{lem:Zar}, we have $X^\lin\subset X^\val$, and the inclusion is strict in general~\cite[Example 2.19]{BKMS}. We will later interpret $X^\lin$ as the set of nonpluripolar points of $X^\an$, cf.~ Proposition ~\ref{prop:linnpp} below. 

\begin{lem}\label{lem:Tcomp} Pick an irreducible component $X_\a$ of $X$, and set $L_\a:=L|_{X_\a}$. For each $v\in X_\a^\an\subset X^\an$ we then have 
$\mathrm{T}_{L_\a}(v)=\teL(v)$. In particular, 
$$
X^\lin=\coprod_\a X_\a^\lin.
$$
\end{lem}
\begin{proof} Since the restriction to $X_\a$ of a regular section on $X$ is regular, we trivially have $\mathrm{T}_{L_\a}(v)\ge \teL(v)$. Conversely pick a nonzero section $s\in\Hnot(X_\a,mL_\a)$. For $r\gg 1$, $s^r$ extends to a regular section $\tilde s\in\Hnot(X,rmL)$. Then $r v(s)=v(s^r)=v(\tilde s)\le \teL(v) mr$, and hence $m^{-1} v(s)\le \teL(v)$, which proves that $\mathrm{T}_{L_\a}(v)\le \teL(v)$. 
\end{proof}

\begin{prop}\label{prop:divlin} Every divisorial valuation has linear growth, \ie $X^\div\subset X^\lin$. 
\end{prop}

\begin{proof} Pick $v\in X^\div$, and choose a projective birational morphism $\mu\colon X'\to X$ with $X'$ normal and a prime divisor $E\subset X'$ such that $v=t\ord_E$, $t\in\Q_{>0}$. Let $H$ be an ample line bundle on $X'$. Pick a regular section $s\in\Hnot(X,mL)$ and set $a:=\ord_E(s)$. Then $\div(\mu^\star s)-aE$ is an effective Weil divisor on $X'$, and hence 
$$
a (E\cdot H^{n-1})\le(\div(\mu^\star s)\cdot H^{n-1})=m(\mu^\star L\cdot H^{n-1}).
$$
Thus $a\le Cm$ for a uniform constant $C>0$, which proves that $\ord_E$, and hence also $v$, has linear growth. 
\end{proof}
More generally, any $v\in X^\an$ such that $v\le v'$ for some $v'\in X^\div$ is a valuation of linear growth. Conversely,~\cite[Theorem~2.16]{BKMS} implies:  

\begin{exam}\label{exam:BKMS} If $v\in X^\val$ is centered at a closed point of $X$, then $v\in X^\lin$ iff $v\le v'$ for some $v'\in X^\div$. 
\end{exam} 
%
%
%
%
 \section{Piecewise linear and Fubini--Study functions}
As before, $X$ denotes a projective variety. We introduce the classes of piecewise linear and Fubini--Study functions on the Berkovich space $X^\an$, and interpret them in terms of test configurations, along the lines of~\cite{BHJ1}. 
%
%
\subsection{ Flag \ ideals and piecewise linear functions}\label{sec:ideals}
Recall from~\S\ref{sec:Berk} that to any ideal $\fb\subset\cO_X$ is associated a continuous function $X^\an\to[0,+\infty]$ given by $v\mapsto v(\fb)$. For reasons that will become clearer later, we define
$$
\log|\fb|:X^\an\to[-\infty,0]
$$ 
by setting $\log|\fb|(v):=-v(\fb)$. The function $\log|\fb|$ is homogeneous with respect to the scaling action of $\R_{>0}$, and~\eqref{equ:videal2} yields
\begin{equation}\label{equ:psideal}
\log|\fb\cdot\fb'|=\log|\fb|+\log|\fb'|,\quad \log|\fb+\fb'|=\max\{\log|\fb|,\log|\fb'|\}
\end{equation}
for all ideals $\fb,\fb'\subset\cO_X$. 

Following~\cite{Oda13,BHJ1}, we define a \emph{flag ideal} $\fa$ as a vertical fractional ideal on $\cX_\triv=X\times\A^1$, \ie a $\Gm$-invariant, coherent fractional ideal sheaf that is trivial on $X\times\Gm$, according to our conventions. We then have a weight decomposition
\begin{equation}\label{equ:flagweight}
\fa=\sum_{\la\in\Z}\fa_\la\unipar^{-\la}
\end{equation}
for a decreasing sequence of ideals $\fa_\la\subset\cO_X$ such that $\fa_\la=\cO_X$ for $\la\ll 0$ and $\fa_\la=0$ for $\la\gg 0$. For any $v\in X^\an$ with Gauss extension $\sigma(v)$, we have
$$
\sigma(v)(\fa)=\min_\la\{v(\fa_\la)-\la\}. 
$$
We define a continuous function $\f_\fa\colon X^\an\to\R$ by setting $\f_\fa(v)=-\sigma(v)(\fa)$, \ie
\begin{equation}\label{equ:gaussflag}
\f_\fa=\max_\la\{\log|\fa_\la|+\la\}. 
\end{equation}
For any two flag ideals $\fa,\fa'$, we have 
\begin{equation}\label{equ:flagsum}
\f_{\fa\cdot\fa'}=\f_\fa+\f_{\fa'},\quad \f_{\fa+\fa'}=\max\{\f_\fa,\f_{\fa'}\}. 
\end{equation}

\begin{defi}\label{defi:PL} We define the space of \emph{piecewise linear (PL) functions on $X^\an$} as the $\Q$-linear subspace 
$$
\PL(X)\subset\Cz(X)
$$
generated by all functions $\f_\fa$ attached to flag ideals $\fa$.  
\end{defi}
By~\eqref{equ:flagsum}, the subset
\begin{equation}\label{equ:PLplus}
\PL^+(X):=\left\{m^{-1}\f_\fa\mid m\in\Z_{>0},\,\fa\ \text{flag ideal}\right\}\subset\PL(X)
\end{equation}
 
is stable under sums and finite maxima, and contains all $\Q$-valued constant functions. It is also stable under multiplication by $\Q_+$: if $c\in\Q_+$, then  $c\f_\fa=m^{-1}\f_{\fa^{mc}}$ for any $m\in\Z_{>0}$ is such that $mc\in\Z$. Further, it is stable  under the scaling action of $\Q_{>0}$: if $t\in\Q_{>0}$ and $\fa=\sum_\la\fa_\la\varpi^{-\la}$, then $t\cdot\f_\fa=m^{-1}\f_{\fa'}$, where
$m\in\Z_{>0}$ is such that $mt\in\Z$, and $\fa'=\sum_\la\fa_\la^m\varpi^{-mt\la}$.

It follows that the $\Q$-vector space $\PL(X)$ is stable under finite maxima and minima, and under the scaling action of $\Q_{>0}$, and it contains all constant $\Q$-valued functions. Further, any function in $\PL(X)$ can be written as a difference of functions in $\PL^+(X)$, and is $\Q$-valued on $X^\div$.

We refer to Theorem~\ref{thm:VCarPL} and Appendix~\ref{sec:dual} for an interpretation of PL functions in terms of test configurations and PL functions on simplicial complexes, respectively. In our setting, PL functions play the role of smooth functions in the complex analytic case, as illustrated by the next result. 

\begin{thm}\label{thm:PLdense} The space $\PL(X)$ is dense in $\Cz(X)$ for the topology of uniform convergence. 
\end{thm}
By the `lattice version' of the Stone--Weierstrass theorem, this is a direct consequence of the following result. 

\begin{lem}\label{lem:FSdense} The $\Q$-linear subspace $\PL(X)\subset\Cz(X)$ is stable under max (and hence min), contains the constants in $\Q$, and separates the points of $X^\an$. 
\end{lem}

\begin{proof} The first two properties are clear. As mentioned in~\S\ref{sec:Berk}, the topology of $X^\an$ is generated by the functions $\log|\fb|$ attached to ideals $\fb\subset\cO_X$. For any two $v,v'\in X^\an$ we can thus find an ideal $\fb$ such that $\log|\fb|(v)\ne\log|\fb|(v')$ (since $X^\an$ is Hausdorff), and it follows that $\f:=\max\{\log|\fb|,-m\}\in\PL(X)$ separates $v,v'$ for $m\gg 1$. 
\end{proof}

\begin{rmk} We will occasionally consider the $\R$-vector space $\PL(X)_\R\subset\Cz(X)$ generated by $\PL(X)$. As opposed to the latter, $\PL(X)_\R$ is not closed under max. 
\end{rmk}

\begin{exam}\label{exam:defconePL}
  Given a closed subscheme $Z\subset X$ with ideal $\fb\subset\cO_X$, the function $$\f_Z:=\max\{\log|\fb|,-1\}$$ lies in $\PL^+(X)$, since 
  $\f_Z=\f_\fa$ with $\fa=\fb+(\unipar)$. Note that
\begin{equation}\label{equ:defconescale}
t\cdot\f_Z=\max\{\log|\fb|,-t\}
\end{equation}
for all $t\in\Q_{>0}$.
\end{exam} 

\begin{exam}\label{exam:curvePL} Assume $X$ is a smooth irreducible curve, and recall the description of $X^\an$ in \S\ref{sec:curve}. As a special case of Example~\ref{exam:defconePL}, each $p\in X(k)$ defines a function $\f_p\in\PL^+(X)$ that satisfies for all $q\in X(k)$ and $t\in[0,+\infty]$ 
$$
\iota_q^\star\f_p(t)=\left\{
\begin{array}{ll}
\max\{-t,-1\} &\  \text{if } q=p; \\
0 & \text{ otherwise}. 
\end{array}
\right. 
$$
Using~\eqref{equ:defconescale} and the fact that every nontrivial ideal of $X$ is of the form $\fb=\cO_X(-\sum_i a_i p_i)$ for a finite set $(p_i)$ in $X(k)$ and $a_i\in\Z_{>0}$, one checks that the functions $t\cdot\f_p$ with $t\in\Q_{>0}$ and $p\in X(k)$ span $\PL(X)$. A function $\f\in\Cz(X)$ is PL iff $\f$ is constant on all but finitely many rays of $X^\an$ and $\Q$-PL on these rays, and $\f\in\PL^+(X)$ iff $\f$ is further convex (and hence decreasing, being bounded) on each ray. 
\end{exam}
%
%
\subsection{PL functions and test configurations}
Let $\cX$ be a test configuration for $X$, and recall that Gauss extension of (semi)valuations yields an embedding $\sigma_\cX\colon X^\an\hookrightarrow\cX^\an$ onto the set of $k^\times$-invariant semivaluations $w\in\cX^\an$ such that $w(\cX_0)=w(\unipar)=1$. 

Recall also that a \emph{vertical $\Q$-Cartier divisor} on $\cX$ means a $\Gm$-invariant $\Q$-Cartier divisor with support in $\cX_0$ (see Example~\ref{exam:vertcar}). Such divisors form a finite dimensional $\Q$-vector space, denoted by 
$$
\VCar(\cX)_\Q.
$$ 

To each $D\in\VCar(\cX)_\Q$ we associate a continuous function $\f_D\in\Cz(X)$ by setting
$$
\f_D(v):=\sigma_\cX(v)(D)
$$
for $v\in X^\an$, where the right-hand side is defined as $m^{-1}\sigma_\cX(v)(\cO_\cX(-mD))$ for any choice of $m\in\Z_{>0}$ such that $mD$ is a Cartier divisor (and hence $D\ge 0\Rightarrow \f_D\ge 0$). 

The map $D\mapsto\f_D$ is $\Q$-linear, and invariant under pull-back: if $\mu\colon\cX'\to\cX$ is a morphism of test configurations then $\sigma_{\cX'}(v)(\mu^\star D)=\sigma_\cX(v)(D)$ for all $v\in X^\an$, and hence $\f_{\mu^\star D}=\f_D$. It thus gives rise to a $\Q$-linear map 
\begin{equation}\label{equ:VCarC0}
\varinjlim_\cX\VCar(\cX)_\Q\to\Cz(X),
\end{equation}
where the direct limit ranges over the directed poset of (isomorphism classes of) test configurations (or merely integrally closed ones, since they form a cofinal subset).

\begin{thm}\label{thm:VCarPL} The map $D\mapsto\f_D$ in~\eqref{equ:VCarC0} induces a $\Q$-linear isomorphism 
$$
\varinjlim_\cX\VCar(\cX)_\Q\simeq\PL(X). 
$$
\end{thm}

\begin{lem}\label{lem:PLvanishing}
  Let $\cX$ be an integrally closed test configuration for $X$, and $D\in\VCar(\cX)_\Q$. Then $D$ is effective iff $\f_D(v_E)\ge 0$ for each irreducible component $E$ of $\cX_0$. 
\end{lem}
\begin{proof} After passing to a multiple, we may assume that $D$ is Cartier. Since $\sigma_\cX(v_E)=b_E^{-1}\ord_E$, Lemma~\ref{lem:CLT} shows that $-\f_D(v_E)=b_E^{-1}\ord_E(\cO_\cX(-D))$ is nonnegative for all $E$ iff $\cO_\cX(-D)\subset\cO_\cX$, which is also equivalent to $D$ being effective. 
\end{proof}

\begin{lem}\label{lem:domL} For any $\f\in\Cz(X)$, the following properties are equivalent:
\begin{itemize}
\item[(i)] $\f\in\PL^+(X)$; 
\item[(ii)] there exists a test configuration $\cX$ that dominates $\cX_\triv$ via $\mu\colon\cX\to\cX_\triv$, and a $\mu$-semiample vertical $\Q$-divisor $D\in\VCar(\cX)_\Q$ such that $\f=\f_D$.
\end{itemize}
\end{lem}
\begin{proof} Assume (i), and write $\f=m^{-1}\f_{\fa}$ for a flag ideal $\fa$ on $\cX_\triv$ and $m\in\Z_{>0}$. Denote by $\mu\colon \cX\to\cX_\triv$ the blowup of $\fa$, so that $\fa\cdot\cO_\cX=\cO_\cX(E)$ with $E\in\VCar(\cX)$. Then $\f=\f_D$ with $D:=m^{-1}E\in\VCar(\cX)_\Q$, which is $\mu$-semiample. This proves (i)$\Rightarrow$(ii). Conversely, assume (ii),  and pick $m$ sufficiently divisible such that $\cO_\cX(mD)$ is $\mu$-globally generated. Then $\cO_\cX(mD)=\fa\cdot\cO_\cX$ with $\fa:=\mu_\star\cO_\cX(mD)$, and hence $\f=m^{-1}\f_\fa$, which proves (ii)$\Rightarrow$(i). 
\end{proof}

\begin{proof}[Proof of Theorem~\ref{thm:VCarPL}] Lemma~\ref{lem:PLvanishing} implies that~\eqref{equ:VCarC0} is injective. Any $D\in\VCar(\cX)_\Q$ with $\mu\colon\cX\to\cX_\triv$ can be written as a difference of $\mu$-(semi)ample divisors; on the other hand, any $\f\in\PL(X)$ is a difference of functions in $\PL^+(X)$, and Lemma~\ref{lem:domL} thus shows that the image of~\eqref{equ:VCarC0} is precisely $\PL(X)$. 
\end{proof}

Next we prove a result that will be used in~\S\ref{sec:topopsh}.
\begin{defi}\label{defi:Rees}
  Given a flag ideal $\fa$, we define the set $\Sigma_\fa\subset X^\div$ of \emph{Rees valuations} of $\fa$ as the finite set of divisorial valuations associated to the irreducible components of $\cX_0$, where $\cX\to\cX_\triv$ is the integral closure of the blowup of $\fa$.
\end{defi}
As the blowup of any ideal is canonically isomorphic to the blowup of any power of that ideal, we have $\Sigma_{\fa^r}=\Sigma_\fa$ for any $r\ge1$.
\begin{rmk}\label{rmk:Rees}
  Let $\cX\to\cX_\triv$ be the integral closure of the blowup along a flag ideal $\fa\subset\cO_{\cX_\triv}$, and $(E_i)$ the irreducible components of $\cX_0$. The Rees valuations of $\fa$ are then the valuations $v_{E_i}\in X^\div$ on $X$. Note that the terminology is a bit abusive, as (assuming $X$ is normal) the divisorial valuations $\ord_{E_i}$ on $\cX$ can also be seen as the Rees valuations of the ideal $\fa$. 
\end{rmk}
\begin{lem}\label{lem:divRees}
  Every divisorial valuation is a Rees valuation of some flag ideal. More generally, for any finite subset $\Sigma\subset X^\div$  there exists a flag ideal $\fa$ of $X$ such that $\Sigma\subset\Sigma_\fa$.
\end{lem}
\begin{proof}
  It suffices to consider the case $\Sigma=\{v\}$, $v\in X^\div$.
  By Lemma~\ref{lem:tcdiv}, there exists an integrally closed test configuration $\cX$ for $X$, and an irreducible component $E$ of $\cX_0$ such that $v_E=v$. Passing to a higher test configuration, we may assume that $\cX$ is the integral closure of the blowup of $\cX_\triv$ along a flag ideal $\fa$, and then $v\in\Sigma_\fa$.
\end{proof}  
\begin{lem}\label{lem:Shilov}
  For any flag ideal $\fa$ and any $\f\in\PL^+(X)$, we have
  \begin{equation*}
    \sup_{X^\an}(\f-\f_\fa)=\max_{\Sigma_\fa}(\f-\f_\fa).
\end{equation*}
\end{lem}
\begin{proof} Write $\f=m^{-1}\f_{\fa'}$, where $m\ge1$ and $\fa'$ is a flag ideal. Then $\f-\f_\fa=m^{-1}(\f_{\fa'}-\f_{\fa^m})$. As $\Sigma_{\fa^m}=\Sigma_\fa$, we may assume $m=1$. 

  Let $\{E_i\}_i$ be the irreducible components of $\cX_0$, so that $\Sigma_\fa=\{v_{E_i}\}_i$.  Write $\fa\cdot\cO_\cX=\cO_\cX(D)$ with $D\in\VCar(\cX)_\Q$. For each $i$ we have $\sigma(v_{E_i})=\ord_{E_i}(\cX_0)^{-1}\ord_{E_i}$, and hence 
$$
c:=\max_{\Sigma_\fa}(\f-\f_\fa)=\max_i\left\{\frac{\ord_{E_i}(D)-\ord_{E_i}(\fa')}{\ord_{E_i}(\cX_0)}\right\}.
$$
Pick $r\in\Z_{>0}$ such that $rc\in\Z$. For all $i$ we then have 
$$
\ord_{E_i}\left((\fa')^r(rc\cX_0-rD)\right)\ge 0, 
$$
and hence $(\fa')^r(rc\cX_0-rD)\subset\cO_\cX$, by Lemma~\ref{lem:CLT}. This yields, in turn,
$$
rc+r\sigma(v)(\fa')\ge r\sigma(v)(\fa)
$$
for all $v\in X^\an$, and we conclude, as desired, $\sup_{X^\an}(\f-\f_\fa)=c$. 
\end{proof}

%
%
\subsection{Density of divisorial valuations}\label{sec:divsense}
 
Using PL functions, we establish some topological properties of $X^\an$. Most importantly, we prove
 
\begin{thm}\label{thm:divdense}
  The set $X^\div$ is dense in $X^\an$.
\end{thm}
\begin{proof} By density of $\PL(X)$ in $\Cz(X)$ (see Theorem~\ref{thm:PLdense}), it suffices to prove that if $\f\in\PL(X)$ vanishes on $\Xdiv$, then $\f\equiv0$. By Theorem~\ref{thm:VCarPL} we have $\f=\f_D$ for some vertical $\Q$-Cartier divisor $D$ on an integrally closed test configuration $\cX$. The result now follows from Lemma~\ref{lem:PLvanishing}.
\end{proof}

Using this, we further show:
\begin{lem}\label{lem:nonvaldense} For any projective variety $X$, we have:
  \begin{itemize}
  \item[(i)]
    the set of semivaluations $v\in X^\an$ with one-dimensional support is dense in $X^\an$ iff every irreducible component of $X$ has dimension at least $1$;
  \item[(ii)]
    $X^\an\setminus X^\val$ is dense in $X^\an$ iff every irreducible component of $X$ has dimension at least $2$.
  \end{itemize}
\end{lem}
\begin{proof}
  We may assume that $X$ is irreducible. Indeed, a subset of $X^\an$ is dense iff its intersection with $X_\a^\an$ is dense in $X_\a^\an$ for every irreducible component $X_\a$ of $X$. 

  The case $\dim X=0$ is trivial, and if $\dim X=1$, then the set of valuations with one-dimensional support equals $X^\val$, which is dense, whereas $X^\an\setminus X^\val$ is not, see~\S\ref{sec:curve}.

  We may therefore assume $\dim X>1$. In this case, we claim that the set $\cC$ of semivaluations with one-dimensional support is dense in $X^\an$. This will prove~(i), and hence also~(ii), since any such semivaluation is contained in $X^\an\setminus X^\val$.

  Pick a prime divisor $E$ on a normal birational model $\pi\colon Y\to X$. By Theorem~\ref{thm:divdense}, it is enough to show that $\ord_E\in X^\div$ lies in the closure of $\cC$. By definition of the topology of $X^\an$, this amounts to the following: given an affine open subvariety $U\subset X$ that intersects $\pi(E)$ and $f_1,\dots,f_r\in\cO(U)$, we need to exhibit $v\in \cC$ such that $v(f_i)$ is arbitrarily close to $\ord_E(f_i)$ for $i=1,\dots,r$. We claim that we can actually find $v\in\cC$ such that $v(f_i)=\ord_E(f_i)$ for all $i$. To see this, denote by $Z$ the union of the irreducible components of $\sum_i\pi^\star\div(f_i)$ that are distinct from $E$. Since $E$ and $Y$ are smooth at the generic point of $E$, we can find an irreducible curve $C\subset Y$, not contained in the exceptional locus of $\pi$, and a closed point $p\in C$ that does not lie on $Z$, such that $E$ and $C$ intersect transversely at $p$. Setting $\ord_{(C,p)}(f):=\ord_p(f|_C)$ for $f\in\cO_{Y,p}$ defines a semivaluation $\ord_{(C,p)}\in Y^\an$, which lies outside $Y^\val$ since $Y$ has dimension at least $2$ at $p$, by assumption. The image of $\ord_{(C,p)}$ in $X^\an$ is a semivaluation with support $\pi(C)$, and hence $v\in\cC$. By construction, we further have $\ord_E(f_i)=\ord_E(\pi^\star f_i)=\ord_{(C,p)}(\pi^\star f_i)=v(f_i)$ for all $i$, and we are done.
\end{proof}
 
%
%
\subsection{Fubini--Study functions}\label{sec:FS}
 
We now introduce classes of functions defined by global sections of line bundles.

Recall that the \emph{base ideal} of a line bundle $L$ on $X$ is the ideal $\fb_L\subset\cO_X$ locally generated by the global sections $\Hnot(X,L)$. The corresponding closed subscheme of $X$ is called the \emph{base scheme} of $L$, while its \emph{base locus} $\Bs(L)$ is the underlying Zariski closed set, \ie
$$
\Bs(L)=\{x\in X\mid s(x)=0\ \text{for all }s\in\Hnot(X,L)\}.
$$
\begin{lem}\label{lem:effective} A line bundle $L$ admits a regular section $s\in\Hnot(X,L)$ iff $\Bs(L)$ is nowhere dense. 
\end{lem}
\begin{proof} The set of regular sections in $\Hnot(X,L)$ is the complement of the union of linear subspaces $V_\a:=\{s\in\Hnot(X,L)\mid s|_{X_\a}\equiv 0\}$. Thus $L$ admits a regular section iff each $V_\a$ is a strict subspace; this is also equivalent to saying that $\Bs(L)$ does not contain any component $X_\a$, \ie $\Bs(L)$ is nowhere dense. 
\end{proof}

Consider now a $\Q$-line bundle $L$. The \emph{asymptotic base locus} of $L$ is the Zariski closed subset $\B(L):=\Bs(mL)$ for $m$ sufficiently divisible. One says that $L$ is \emph{effective} (resp.~\emph{semiample}) if $\B(L)$ is nowhere dense (resp.~empty). Note that 
$$
L\ \text{effective }\Longrightarrow L\ \text{pseudoeffective, \ie }c_1(L)\in\Psef(X);
$$
$$
L\ \text{semiample }\Longrightarrow L\ \text{nef, \ie }c_1(L)\in\Nef(X).
$$
Indeed, the first implication follows from Lemma~\ref{lem:effective}, which yields a regular section $s\in\Hnot(X,mL)$ for $m$ sufficiently divisible, so that $m^{-1}\div(s)$ is an effective $\Q$-Cartier divisor in the numerical class of $L$. 

\medskip

Consider next an additive subgroup $\La\subset\R$ (the main cases being $\{0\}$, $\Q$, or $\R$), and a function $\f\colon X^\an\to\R\cup\{-\infty\}$ of the form 
\begin{equation}\label{equ:FS}
\f=m^{-1}\max_j\{\log|s_j|+\la_j\}
\end{equation}
with $m\in\Z_{>0}$ such that $mL$ is an honest line bundle, $(s_j)$ a finite set of sections of $\Hnot(X,mL)$, and $\la_j\in\La$. (Recall that~\eqref{equ:FS} means $\f(v)=m^{-1}\max_j\{-v(s_j)+\la_j\}$ for all $v\in X^\an$). Using Lemma~\ref{lem:maxprin}, the next result is straightforward: 

\begin{lem}\label{lem:FSgenfin} Every function $\f\colon X^\an\to\R\cup\{-\infty\}$ of the form~\eqref{equ:FS} is continuous, decreasing, and satisfies $\B(L)(k)\subset\{\f=-\infty\}$. Furthermore, the following are equivalent: 
\begin{itemize}
\item[(i)] $\f$ is generically finite, \ie $\f|_{X_\a^\an}\not\equiv-\infty$ for all $\a$ (see Definition~\ref{defi:genfin}); 
\item[(ii)] $\f$ is finite at $v_{\triv,\a}$ for all $\a$ (see~\S\ref{sec:triv}); 
\item[(iii)] $\f$ is finite valued on $X^\val$. 
\end{itemize}
\end{lem}

\begin{defi}\label{defi:FS} Given a subgroup $\La\subset\R$ and a $\Q$-line bundle $L$ on $X$, we say that a function $\f\colon X^\an\to\R\cup\{-\infty\}$ is 
\begin{itemize}
\item[(i)] a \emph{$\La$-rational, generically finite Fubini--Study function for $L$} if $\f$ is of the form~\eqref{equ:FS} and generically finite, \ie finite valued on $X^\val$ (see Lemma~\ref{lem:FSgenfin}); 
\item[(ii)] a \emph{$\La$-rational Fubini--Study function for $L$} if $\f$ is further finite valued on all of $X^\an$.
\end{itemize}
We denote by $\cH_\La(L)\subset\cH^\gf_\La(L)$ the spaces so defined. 
\end{defi}
Note that these sets only depend on the isomorphism class of $L$, \ie its image in $\Pic(X)_\Q$. Further,
$$
\cH_\La(L)=\cH^\gf_\La(L)\cap\Cz(X).
$$
For any $\f$ as in~\eqref{equ:FS}, we have 
\begin{equation}\label{equ:FSpower}
\f=(rm)^{-1}\max_j\{\log|s_j^r|+r\la_j\}
\end{equation}
for all $r\in\Z_{>0}$. Thus
\begin{equation}\label{equ:FSdivi}
\cH^\gf_\La(L)=\cH^\gf_{\Q\La}(L),\quad \cH_\La(L)=\cH_{\Q\La}(L),
\end{equation}
which means that the subgroup $\La\subset\R$ can always be assumed to be \emph{divisible} in the above definition. Finally, we trivially have
$$
\cH_0^\gf(L)\subset\cH^\gf_\La(L)\subset\cH^\gf_\R(L),\quad\cH_0(L)\subset\cH_\La(L)\subset\cH_\R(L).
$$
These sets can be empty; more precisely, it is straightforward to check that 
\begin{equation}\label{equ:cHvalnempty}
\cH^\gf_\R(L)\ne\emptyset\Longleftrightarrow\cH_0^\gf(L)\ne\emptyset\Longleftrightarrow L\ \text{effective}; 
\end{equation}
\begin{equation}\label{equ:cHnempty}
\cH_\R(L)\ne\emptyset\Longleftrightarrow\cH_0(L)=\{0\}\Longleftrightarrow L\ \text{semiample}.
\end{equation}
The next result summarizes further properties that are also readily checked (compare~\cite[Proposition 5.4]{BE}).

\begin{prop}\label{prop:FS} Pick any $L\in\Pic(X)_\Q$. Then: 
\begin{itemize}
\item[(i)] each $\f\in\cH^\gf_\La(L)$ is decreasing  on $X^\an$, and hence satisfies the maximum principle~\eqref{equ:maxprin} (see Lemma~\ref{lem:maxprin}); further, $\f\equiv-\infty$ on $\B(L)^\an$; 
\item[(ii)] $\cH^\gf_\La(L)$ and $\cH_\La(L)$ are both invariant under the scaling action~\eqref{equ:scale} restricted to the subgroup 
$\{t\in\R_{>0}\mid t\La\subset\La\}$ of $\R_{>0}$;
\item[(iii)] for all $L'\in\Pic(X)_\Q$ and $a\in\Q_{>0}$ we have 
$$
\cH^\gf_\La(aL)=a\cH^\gf_\La(L),\quad\cH^\gf_\La(L)+\cH^\gf_\La(L')\subset\cH^\gf_\La(L+L'),
$$
and similarly for $\cH_\La$;
\item[(iv)] for any morphism $f\colon Y\to X$ from a projective variety we have 
$$
f^\star\cH_\La(L)\subset\cH_\La(f^\star L),
$$ 
and the same holds for $\cH^\gf_\La$ if $f$ is surjective. 
\end{itemize}
\end{prop}
When $\f\in\cH^\gf_\La(L)$ is written as in~\eqref{equ:FS}, the scaling action in (ii) is simply given by
$$
t\cdot\f=m^{-1}\max_j\{\log|s_j|+t\la_j\}.
$$
For any $v\in X^\an$, we similarly have 
$$
\f(tv)=m^{-1}\max_j\{-t v(s_j)+\la_j\},
$$
which yields: 
\begin{lem}\label{lem:FSconvex} For any $\f\in\cH^\gf_\R(L)$ and $v\in X^\an$, $t\mapsto\f(tv)$ is convex and decreasing on $\R_{>0}$.
\end{lem}

Proposition~\ref{prop:FS}~(iv) admits the following partial converse, a key ingredient in the proof of Theorem~\ref{thm:descentpsh} below. 

\begin{lem}\label{lem:descentFS}  Let $\pi\colon Y\to X$ be the blowup of an ideal $\fb\subset\cO_X$, with exceptional divisor $E$. Denote by $s_E\in\Hnot(Y,\cO_Y(E))$ the canonical section, so that $\pi^\star\log|\fb|=\log|s_E|\in\cH^\gf_0(E)$. For any $L\in\Pic(X)_\Q$ we then have 
$$
\cH^\gf_\La(\pi^\star L-E)+\log|s_E|\subset \pi^\star\cH^\gf_\La(L). 
$$
\end{lem}
 
\begin{lem}\label{lem:Hart} In the notation of Lemma~\ref{lem:descentFS}, we have $\fb^m\subset\pi_\star\cO_Y(-mE)$ for all $m\in\N$, and equality holds for all $m$ large enough. 
\end{lem} 
\begin{proof} This follows from the fact that $Y$ is the relative Proj of the graded $\cO_X$-algebra $\cR=\bigoplus_{m\in\N}\fb^m$, which is generated in degree $1$, and that $\cO_Y(-E)=\cO_Y(1)$ (see for instance~\cite[Exercise II.5.9]{Har}). 
\end{proof}

\begin{proof}[Proof of Lemma~\ref{lem:descentFS}] Pick $\f\in\cH^\gf_\La(\pi^\star L-E)$, and write it as in~\eqref{equ:FS}, with $s_i\in\Hnot(Y,m(\pi^\star L-E))$. After replacing $m$ and the $s_i$ with $rm$ and $s_i^r$ for $r$ large enough as in~\eqref{equ:FSpower}, we may assume that $\pi_\star\cO_Y(-mE)=\fb^m$, by Lemma~\ref{lem:Hart}. For each $i$, $s_i s_E^m\in\Hnot(Y,m\pi^\star L)$ locally belongs to the ideal $\cO_Y(-mE)$, and hence $s_i s_E^m=\pi^\star \sigma_i$ with $\sigma_i\in\Hnot(X,\cO_Y(mL)\otimes\fb^m)$. This yields $\f+\log|s_E|=\pi^\star\p$ with $\p:=m^{-1}\max_i\{\log|\sigma_i|+\la_i\}$. Since $\pi$ is birational, $\p$ is finite valued at each $v_{\triv,\a}$ iff $\pi^\star \p=\f+\log|s_E|$ satisfies the analogous condition on $Y$, which is indeed the case since $s_E$ is a regular section. Thus $\p\in\cH^\gf_\La(L)$, and we are done. 
\end{proof}

%
%

For later use (see Theorem~\ref{thm:pshdiv}), we also establish the following `division' property: 

\begin{lem}\label{lem:FShomdiv} Assume that $X$ is normal. Pick $L\in\Pic(X)_\Q$, an effective $\Q$-Cartier divisor $E$ on $X$, and $\f\in\cH^\gf_\La(L)$. Then 
$$
\f\le\log|s_E|+O(1)\Longleftrightarrow\f-\log|s_E|\in\cH^\gf_\La(L-E).
$$
\end{lem}
We have somewhat abusively set $\log|s_E|:=m^{-1}\log|s_{mE}|\in\cH^\gf_0(E)$ for $m$ sufficiently divisible, where $s_{mE}\in\Hnot(X,mE)$ denotes the canonical section. 

\begin{proof} Write $\f$ as in~\eqref{equ:FS}. Replacing $m$ and $s_i$ with $rm$ and $s_i^r$ for $r$ large enough, we may assume that $mE$ is a Cartier divisor. For each $i$, we have $\log|s_i|+\la_i\le m\f\le\log|s_{mE}|+O(1)$, \ie $v(s_i)\ge v(mE)-C$ for all $v\in X^\an$ and a uniform constant $C$. Replacing $v$ with $tv$ and letting $t\to+\infty$, we infer $v(s_i)\ge v(mE)$ for all $v$. This holds in particular with $v=\ord_F$ for any irreducible component $F$ of $E$, which shows that $s_i$ locally belongs to the ideal $\cO_X(-mE)\subset\cO_X$, since $X$ is normal. For each $i$, we thus have $s_i=\sigma_i s_{mE}$ with $\sigma_i\in\Hnot(X,m(L-E))$, which yields $\f-\log|s_E|=m^{-1}\max_i\{\log|\sigma_i|+\la_i\}$. This function is further generically finite, since so are $\f$ and $\log|s_E|$, and we conclude, as desired, $\f-\log|s_E|\in\cH^\gf_\La(L-E)$. 
\end{proof}

We further observe that generically finite Fubini--Study functions are automatically constant on a substantial part of $X^\an$: 
\begin{lem}\label{lem:almostconstant} Pick $L\in\Pic(X)_\Q$ and $\f\in\cH^\gf_\R(L)$. Then there exists a non-empty Zariski open subset $U\subset X$ such that $\f\equiv\sup\f$ on $c^{-1}(U)$. 
\end{lem}
Recall that the center map $c:X^\an\to X$ is anticontinuous. Thus $c^{-1}(U)$ is closed in $X^\an$ (in fact, a $k$-analytic domain), but it has non-empty interior, as it contains the non-empty open subset $c^{-1}(\{p\})$ for any closed point $p\in U(k)$. 

\begin{proof} Write $\f=m^{-1}\max_i\{\log|s_i|+\la_i\}$ as in~\eqref{equ:FS}, with $s_i\ne 0$ for all $i$. For any $v\in X^\an$, $v(s_i)=-\log|s_i|(v)$ is nonzero iff $s_i$ vanishes at the center $c(v)$. If $c(v)\notin Z:=\bigcup_i (s_i=0)$, we thus have $\f(v)=\max_i\la_i=\sup\f$, which yields the result with $U:=X\setminus Z$. 
\end{proof}

The space
$$
\cH(L):=\cH_\Q(L)
$$ 
of \emph{(rational) Fubini--Study functions} plays a central role in this paper. As we shall see, when $L$ is ample, $\cH(L)$ is in 1--1 correspondence with the set of integrally closed, ample test configurations (see Corollary~\ref{cor:tcPL} below). On the other hand, $\cH_\R(L)$ is closely related to the notion of an `$\R$-test configuration' as considered for instance in~\cite{DS}; see~\cite{nakstab1} for details.

%
%
 
\subsection{Fubini--Study functions and PL functions}
We now study the relation between $\cH(L)$ and the space $\PL(X)$ introduced in~\S\ref{sec:ideals}.   

\begin{prop}\label{prop:FSideal} For any $L\in\Pic(X)_\Q$ and $\f\colon X^\an\to\R\cup\{-\infty\}$, we have $\f\in\cH(L)$ iff $\f=m^{-1}\f_\fa$ for a flag ideal $\fa$ and $m\in\Z_{>0}$ such that $mL$ is an honest line bundle and $m\cL_\triv\otimes\fa$ is globally generated on $\cX_\triv$. 
\end{prop}
Here $(\cX_\triv,\cL_\triv)$ is the trivial test configuration for $(\cX,\cL)$, see Example~\ref{exam:trivtc}.
In terms of the weight decomposition~\eqref{equ:flagweight}, note that $m\cL_\triv\otimes\fa$ is globally generated iff $mL\otimes\fa_\la$ is globally generated for all $\la\in\Z$. 

\begin{proof} Assume first $\f\in\cH(L)$ and write $\f=m^{-1}\max_i\{\log|s_i|+\la_i\}$ for a finite set $(s_i)$ in $\Hnot(X,mL)$ and $\la_i\in\Z$. The rational sections $(s_i\unipar^{-\la_i})$ of $m\cL_\triv$ on $\cX_\triv$ generate a flag ideal $\fa$ such that $m\cL_\triv\otimes\fa$ is globally generated on $\cX_\triv$, and $\f=m^{-1}\f_\fa$.

  Assume, conversely, that $\f=m^{-1}\fa$ with $\fa$ a flag ideal such that $m\cL_\triv\otimes\fa$ is globally generated. As above, write $\fa=\sum_{\la\in\Z}\fa_\la\unipar^{-\la}$. For each $\la\in\Z$, $\cO_X(mL)\otimes\fa_\la$ is globally generated by a finite set $(s_{\la,i})_i$ of sections in $\Hnot(X,mL)$, and~\eqref{equ:gaussflag} shows that
  
$$
\f=m^{-1}\max_{\la,i}\{\log|s_{\la,i}|+\la\},
$$
which proves that $\f\in\cH(L)$. 
\end{proof}

\begin{cor}\label{cor:FSideal} For any $L\in\Pic(X)_\Q$, we have $\Q_+\cH(L)\subset\PL^+(X)$, and equality holds if $L$ is ample. In particular, $\cH(L)$ spans the $\Q$-vector space $\PL(X)$ whenever $L$ is ample. 
\end{cor}
In view of Proposition~\ref{prop:FS}~(iii), this implies
\begin{cor}\label{cor:FSideal2} We have $\PL^+(X)=\bigcup_L\cH(L)$, where $L$ ranges over ample classes in $\Pic(X)$.
\end{cor}
Given a flag ideal $\fa$, we get an evaluation map $\PL(X)\to\Q^{\Sigma_\fa}$. We now show that this map is surjective. Using Lemma~\ref{lem:divRees}, this will imply that the evaluation map $\PL(X)\to\Q^\Sigma$ is surjective for any finite subset $\Sigma\subset X^\div$. For later purposes, we prove a more precise result.
\begin{lem}\label{lem:PLinterpol}
  For any flag ideal $\fa$, and $L\in\Pic(X)_\Q$ ample, the following property holds for $m$ sufficiently divisible. For any $c\in\Q^{\Sigma_\fa}$, there exists $r\ge1$ and $\rho\in\cH(L)$ such that $\p:=r(\f_\fa-m\rho)\in\PL(X)$ satisfies $\p(v)=c_v$ for all $v\in\Sigma_\fa$. 
\end{lem}
Before proving the lemma, we establish the following version of~\cite[Theorem 1.10]{BHJ1}.
\begin{lem}\label{lem:Rees} Let $\fa$ be a flag ideal. Denote by $\mu\colon\cX\to\cX_\triv$ the integral closure of the blowup along $\fa$, and by $E_1,\dots,E_N$ the irreducible components of $\cX_0$. Then
\begin{equation}\label{equ:Rees}
\bigcap_i\left\{f\in\cO_{\cX_\triv}\mid\ord_{E_i}(f)\ge\ord_{E_i}(\fa^m)\right\}\varsubsetneq\bigcap_{i>1}\left\{f\in\cO_{\cX_\triv}\mid\ord_{E_i}(f)\ge\ord_{E_i}(\fa^m)\right\}
\end{equation}
for all $m\in\N$ large enough.
\end{lem}
\begin{proof} The (antieffective) vertical Cartier divisor $D$ on $\cX$ such that $\cO_{\cX}(D)=\fa\cdot\cO_{\cX}$ is $\mu$-ample, and the left-hand side of~\eqref{equ:Rees} coincides with $\mu_\star\cO_{\cX}(mD)$. For $m\gg 1$, 
$$
\cO_{\cX}(mD)\varsubsetneq\cO_{\cX}(mD)\cdot\cO_{\cX}(E_1)
$$
are both $\mu$-globally generated, and taking $\mu_\star$ yields the result.
\end{proof}
\begin{proof}[Proof of Lemma~\ref{lem:PLinterpol}]
Let $\cX\to\cX_\triv$ be the integral closure of the blowup of $\cX_\triv$ along $\fa$, and $\{E_i\}_i$ the irreducible components of $\cX_0$, so that $\Sigma_\fa=\{v_i\}_i$, where $v_i=v_{E_i}$.
  By Lemma~\ref{lem:Rees}, after replacing $\fa$ with some power, we may assume that for any $j$, we have 
$$
\bigcap_i\left\{f\in\cO_{\cX_\triv}\mid\ord_{E_i}(f)\ge\ord_{E_i}(\fa)\right\}\varsubsetneq\fa'_j:=\bigcap_{i\ne j}\left\{f\in\cO_{\cX_\triv}\mid\ord_{E_i}(f)\ge\ord_{E_i}(\fa)\right\}.
$$
For $m$ sufficiently divisible, $m\cL_\triv\otimes\fa'_j$ is  globally generated for any $j$. After replacing $L$ by $mL$, we may assume $m=1$. Then, for each $j$,  there exists $s_j\in\Hnot\left(\cX_\triv,\cL_\triv\right)$ such that $\ord_{E_j}(s_j)<\ord_{E_j}(\fa)$ and $\ord_{E_i}(s_j)\ge\ord_{E_i}(\fa)$ for all $i\ne j$. Write $s_j=\sum_{\la\in\N}s_{j,\la}\varpi^\la$ with $s_{\la,j}\in\Hnot(X,L)$, and define functions $\rho'_j\colon X^\an\to\R\cup\{-\infty\}$ by  
$$
\rho'_j:=\max_\la\{\log|s_{j,\la}|-\la\}.
$$
Then $\rho'_j(v_i)=-b_i^{-1}\ord_{E_i}(s_j)$ and $\f_\fa(v_i)=-b_i^{-1}\ord_{E_i}(\fa)$; hence $\rho'_j(v_j)>\f_{\fa}(v_j)$ while $\rho'_j(v_i)\le\f_\fa(v_i)$ for $i\ne j$. If we pick $a=\min_i\f_\fa(v_i)$ and set $\rho_j:=\max\{\rho'_j,a\}$, then $\rho_j\in\cH(L)$ and 
$$
\max_{i\ne j}(\rho_j(v_i)-\f_\fa(v_i))\le0<\e_j:=\rho_j(v_j)-\f_\fa(v_j).
$$
For $r\ge1$, we now set
$$\rho:=\max_j\{\rho_j-\e_j-c_j/r\}\in\cH(L)\subset\PL^+(X)
\quad\text{and}\quad\p=r(\f_\fa-\rho)\in\PL(X).$$
If $r\gg1$, then $r\e_i+c_i\ge c_j$ for all $i,j$, which easily implies $\p(v_j)=c_j$ for all $j$.
\end{proof}

%
%
\subsection{Equivalence of test configurations}\label{sec:equivtc}
Let $L$ be a $\Q$-line bundle on $X$. As in~\cite[Definition 6.1]{BHJ1}, we say that that two test configurations $(\cX,\cL)$, $(\cX',\cL')$ for $(X,L)$ are \emph{equivalent} if $\cL$ and $\cL'$ agree after pulling back to a test configuration $\cX''$ dominating both $\cX$ and $\cX'$. 

For simplicity, we say that a test configuration $(\cX,\cL)$ is integrally closed (resp.~semiample, ample) if $\cX$ (resp.~$\cL$) is. 
Slightly generalizing~\cite[Lemma 6.3]{BHJ1} (which assumed $X$ normal), we have: 

\begin{prop}\label{prop:tcunique} If $L$ is an ample $\Q$-line bundle, then every semiample test configuration for $(X,L)$ is equivalent to a unique ample, integrally closed test configuration.
\end{prop}

\begin{proof} By~\cite[Proposition 2.17]{BHJ1}, every semiample test configuration is equivalent to an ample test configuration, which can further be assumed to be integrally closed after passing to the integral closure. This proves existence. 

To prove uniqueness, let $(\cX,\cL)$, $(\cX',\cL')$ be two ample, integrally closed test configurations for $(X,L)$ that are equivalent. After replacing $L$ with a multiple, we may assume that $\cL,\cL'$ are honest line bundles. By ampleness, it will then be enough to show that $\Hnot(\cX,m\cL)= \Hnot(\cX',m\cL')$ for all $m\in\N$, as $k[\unipar]$-submodules of $\Hnot(X,mL)_{k[\varpi^{\pm 1}]}$ (see~\S\ref{sec:tc}). Choose a test configuration $\cX''$ dominating $\cX$ and $\cX'$ via $\mu\colon\cX''\to\cX$, $\mu'\colon\cX''\to\cX'$, such that $\cL'':=\mu^\star \cL=\mu'^\star \cL'$.  Since $\cX$ and $\cX'$ are integrally closed, Lemma~\ref{lem:Zarmain} yields $\mu_\star\cO_{\cX''}=\cO_\cX$, $\mu'_\star\cO_{\cX''}=\cO_{\cX'}$, and the projection formula shows that
$$
\Hnot(\cX,m\cL)=\Hnot\left(\cX'',m\cL''\right)=\Hnot(\cX',m\cL').
$$
The proof is complete.
\end{proof}
%
%
\subsection{Fubini--Study functions and test configurations}\label{sec:FStest}
Let $L$ be a $\Q$-line bundle on $X$. For any test configuration $(\cX,\cL)$ for $(X,L)$, we can choose a test configuration $\cX'$ with two morphisms $\mu\colon\cX'\to\cX$, $\rho\colon\cX'\to\cX_\triv$.
Then 
$$
D:=\mu^\star\cL-\rho^\star\cL_\triv\in\VCar(\cX)_\Q,
$$
and mapping $\cL$ to $D$ yields a 1--1 correspondence between the set of equivalence classes of test configurations for $(X,L)$ and the $\Q$-vector space
 $\varinjlim_\cX\VCar(\cX)_\Q$,
which in turn is isomorphic to $\PL(X)$, by Theorem~\ref{thm:VCarPL}.
We write $\f_\cL$ for the element of $\PL(X)$ associated to $\cL$.
The restriction of $\f_\cL$ to $X^\div$ coincides with the function associated to the equivalence class of $(\cX,\cL)$ in~\cite[\S6]{BHJ1}.
 
\begin{thm}\label{thm:testPL} For any $\Q$-line bundle $L$ on $X$, the map $\cL\mapsto\f_\cL$ sets up a 1--1 correspondence between:
\begin{itemize}
\item[(i)] the set of equivalence classes of test configurations for $L$ and $\PL(X)$; 
\item[(ii)] the set of equivalence classes of semiample test configurations for $L$ and $\cH(L)$. 
\end{itemize}
\end{thm}
Note that the sets in (ii) are nonempty only if $L$ is semiample. Combined with Proposition~\ref{prop:tcunique}, we infer: 

\begin{cor}\label{cor:tcPL} If $L$ is an ample $\Q$-line bundle, then $\cL\mapsto\f_\cL$ defines a 1--1 correspondence between the set of ample, integrally closed test configurations for $L$ and $\cH(L)$. 
\end{cor}
The set $\cH(L)$ thus corresponds to $\cH^\NA(L)$ in the notation of~\cite[Definition 6.2]{BHJ1}. 

\begin{proof}[Proof of Theorem~\ref{thm:testPL}] By construction, the map $\cL\mapsto\f_\cL$ is the composition of $D\mapsto\f_D$ with the canonical bijection between equivalence classes of test configurations for $L$ and $\varinjlim_\cX\VCar(\cX)_\Q$. Thus (i) is a direct consequence of Theorem~\ref{thm:VCarPL}. 

To prove (ii), let $(\cX,\cL)$ be a semiample test configuration for $(X,L)$. After replacing $L$ with a multiple and pulling-back to a higher test configuration, we may assume without loss that $\cL$ is a globally generated line bundle, and that $\cX$ dominates the trivial test configuration via $\mu\colon\cX\to\cX_\triv$. Set $D:=\cL-\mu^\star \cL_\triv\in\VCar(\cX)$, so that $\f_\cL=\f_D$. Since $\cL$ is globally generated, $\cO_\cX(D)$ is $\mu$-globally generated, and the flag ideal $\fa:=\mu_\star\cO_\cX(D)$ thus satisfies $\fa\cdot\cO_{\cX}=\cO_\cX(D)$, and hence $\f_\fa=\f_D=\f_\cL$. Denote by $\fa'\subset\fa$ the flag ideal locally generated by 
$$
\Hnot\left(\cX_\triv,\cL_\triv\otimes\fa\right)\simeq\Hnot(\cX,\cL).
$$
Since $\cL$ is globally generated, we have $\fa'\cdot\cO_\cX=\cO_\cX(D)=\fa\cdot\cO_\cX$; hence $\f_{\fa'}=\f_\cL$. By construction, $L\otimes\fa'$ is globally generated, and Proposition~\ref{prop:FSideal} thus yields, as desired, $\f_\cL\in\cH(L)$. 

Conversely pick $\f\in\cH(L)$. After replacing $L$ with a multiple, Proposition~\ref{prop:FSideal} yields $\f=\f_\fa$ for a flag ideal $\fa$ on $\cX_\triv$ such that $\cL_\triv\otimes\fa$ is globally generated. Denoting by $\mu\colon\cX\to \cX_\triv$ the blowup along $\fa$, we have $\fa\cdot\cO_\cX=\cO_\cX(D)$ for a vertical Cartier divisor $D$ such that $\cL:=\mu^\star \cL_\triv+D$ is globally generated, and $\f_\cL=\f_D=\f_\fa=\f$.  
\end{proof}
The next result will be needed in~\S\ref{sec:pshPL}.
\begin{lem}\label{lem:amphigh} If $L$ is an ample $\Q$-line bundle on $X$ and $\cX$ any test configuration for $X$ that dominates $\cX_\triv$, then $L$ admits an ample test configuration $\cL$ determined on $\cX$. 
\end{lem}
\begin{proof} By Lemma~\ref{lem:relamp}, we can find $D\in\VCar(\cX)_\Q$ that is relatively ample over $\cX_\triv$. Denoting by $L_\cX$ the pullback of $L$ to $\cX$, it follows that $\cL:=L_\cX+\e D$ is ample for $\e\in\Q_{>0}$ small enough, and we are done. 
\end{proof}
%
%
\subsection{Geometric interpretation of the scaling action}\label{sec:geomscaling}
 
As noted in~\S\ref{sec:ideals}, the spaces $\PL(X)$ and $\PL^+(X)$ admit a natural 
scaling action by $\Q_{>0}$.
We now give a geometric interpretation of the scaling action of $\Z_{>0}$ in terms of base change.

Using flag ideals, this is easy: any function $\f\in\PL^+(X)$ can be written $\f=m^{-1}\f_\fa$ for some $m\in\Z_{>0}$ and some flag ideal $\fa$ on $\cX_\triv=X\times\A^1$. It then follows from the discussion after~\eqref{equ:PLplus} that if $d\in\Z_{>0}$, then $d\cdot\f=m^{-1}\f_{\fa^{(d)}}$, where $\fa^{(d)}$ is the pullback of $\fa$ under the map $\cX_\triv\to\cX_\triv$ given by $\unipar\mapsto\unipar^d$.

Following~\cite[\S6.3]{BHJ1}, we can give a smilar interpretation of the scaling action when the functions are associated to test configurations.
 
Let $L$ be a $\Q$-line bundle on $X$.
For any test configuration $(\cX,\cL)$ for $(X,L)$ and $d\in\Z_{>0}$, denote by $\cX_d\to\A^1$ the base change of $\cX\to\A^1$ by $\unipar\mapsto\unipar^d$, and let $\cL_d$ be the pull-back of $\cL$ to $\cX_d$. 

\begin{lem}\label{lem:basechange} For each $d\in\Z_{>0}$ we have $d\cdot\f_\cL=\f_{\cL_d}$.
\end{lem}
\begin{proof}
  After passing to a higher test configuration, we may assume, by linearity, that $L=\cO_X$ and $\cL=\cO_\cX(D)$ with $D\in\VCar(\cX)$. Denote by $\rho\colon\cX_d\to\cX$ the natural morphism.
  For any $v\in X^\an$, we then need to show that 
$$
\sigma_{\cX_d}(v)(\rho^\star D)=d\,\sigma_\cX(d^{-1} v)(D),
$$
where $\sigma_\cX$ and $\sigma_{\cX_d}$ denote the Gauss extensions. Let $Y\subset X$ be the support of $v$. Then $w:=\rho_\star\sigma_{\cX_d}(v)$ is a valuation on the induced test configuration $\cY\subset\cX$, and it will be enough to show that $w=d\,\sigma_\cX(d^{-1}v)$. To this end, note that $w$ is $k^\times$-invariant, and satisfies $w(\unipar)=\sigma_{\cX_d}(v)(\unipar^d)=d$. We thus have $d^{-1} w=\sigma_\cX(v')$ for a unique valuation $v'$ on $Y$. For all $f\in k(Y)$ we have $\rho^\star f=f$, and hence 
$$
v'(f)=d^{-1} w(f)=d^{-1} v(f),
$$
which proves, as desired, that $w=d\,\sigma_\cX(d^{-1} v)$. 
\end{proof}
\begin{cor}\label{cor:reducedfiber} Let $L$ be an ample $\Q$-line bundle, and pick $\f\in\cH(L)$. For all $d\in\Z_{>0}$ divisible enough, the unique ample, integrally closed representative of $d\cdot\f$ has a reduced central fiber. 
\end{cor}
\begin{proof} Let $(\cX,\cL)$ be the ample, integrally closed representative of $\f$. By Lemma~\ref{lem:basechange}, the ample, integrally closed representative of $d\cdot\f$ is the integral closure $(\tcX_d,\tcL_d)$ of $(\cX_d,\cL_d)$, and we thus need to show that the central fiber $\tcX_{d,0}$ is reduced for $d$ divisible enough. 

  By~\cite[Tag 09IJ]{Stacks}, the central fiber of the normalization $\cX^\nu_d$ is generically reduced for $d$ divisible enough. Since $\tcX_d$ is integrally closed, the normalization morphism $\cX^\nu_d\to\tcX_d$ is an isomorphism over the generic points of $\tcX_{d,0}$ (see Remark~\ref{rmk:normvsint}), which is thus also generically reduced. We conclude thanks to Corollary~\ref{cor:redvertnorm}. 
\end{proof}
%
%
\subsection{Almost trivial test configurations}\label{sec:almosttriv}
We end this section with an analysis of `almost trivial' test configurations. Following~\cite{Sto2,Oda4,BHJ1}, we introduce: 

\begin{defi} We say that a test configuration $\cX$ for $X$ is 
\begin{itemize}
\item[(i)] \emph{almost trivial} if the normalization $\cX^\nu$ of $\cX$ is trivial;
\item[(ii)] \emph{trivial in codimension $1$} if the canonical $\Gm$-equivariant birational map $\cX\dashrightarrow\cX_\triv$ is an isomorphism in codimension $1$.
\end{itemize}
\end{defi}
Note that (i) corresponds to~\cite[Definition 2.9]{BHJ1}, while (ii) corresponds to~\cite[Definition 1]{Sto2} and~\cite[Definition 3.3]{Oda4}. 

Recall that $\cX^\nu$ is a test configuration for the normalization $X^\nu$ of $X$. If $\cX$ is trivial in codimension $1$, then $\cX_0$ is generically reduced, and $\cX$ is thus regular at each generic point of $\cX_0$. 

As we shall see, (ii) implies (i), but the converse fails in general, even when $X$ is smooth, despite what was claimed in~\cite[Proposition 2.11]{BHJ1}. This was kindly pointed to us by Masafumi Hattori, together with the following simple example: 

\begin{exam} Let $X\subset\P^2$ be a smooth conic such that $[0:0:1]\notin X$, and consider the test configuration $\cX$ defined the $1$-parameter subgroup $\rho(t)[x_0:x_1:x_2]=[x_0:x_1:t x_2]$, which degenerates $X=\cX_1$ to a double line $\cX_0$. It comes with a morphism of test configurations $\cX_\triv\to\cX$, which is finite, and hence coincides with the normalization. Thus $\cX$ is almost trivial, but not trivial in codimension $1$, as $\cX_0$ is generically non-reduced. 
\end{exam}

The next result clarifies the situation, and also provides a simple interpretation of almost triviality in terms of PL functions. 

\begin{thm}\label{thm:almosttriv} Let $L$ be a $\Q$-line bundle on $X$ and $(\cX,\cL)$ a test configuration for $(X,L)$, with associated function $\f_\cL\in\PL(X)$. Then:
\begin{itemize}
\item[(i)] if $\cX$ is almost trivial, then $\f_\cL$ is locally constant. When $\cL$ is further ample, the converse holds as well; 
\item[(ii)] $\cX$ is almost trivial iff its integral closure $\tcX$ is trivial; 
\item[(iii)] $\cX$ is trivial in codimension $1$ iff $\cX$ is almost trivial and $\cX_0$ is generically reduced.
\end{itemize}
\end{thm}
When $X$ is equidimensional and $\cL$ is ample, (i) is also equivalent to $\jj(\f_\cL)=0$, cf.~Corollary~\ref{cor:MAinj} below (compare also~\cite[Theorem A~(ii)]{BHJ1}). 

\begin{proof} Assume that $\cX$ is almost trivial, \ie $\cX^\nu=\coprod_\a\cX^\nu_\a$ is trivial. Since $\cX^\nu_\a$ is integrally closed and its central fiber is irreducible, Lemma~\ref{lem:PLvanishing} implies that $\VCar(\cX^\nu_\a)_\Q$ is $1$-dimensional, and hence 
$$
\nu^\star(\cL-\cL_{\triv})|_{\cX_\a}=c_\a\cX^\nu_{\a,0}
$$
for some $c_\a\in\Q$. This implies that $\f_{\cL}\equiv c_\a$ on $X_\a^\div\simeq (X^\nu_\a)^\div$; hence $\f_\cL\equiv c_\a$ on $X_\a^\an$, by density of divisorial points, and we infer that $\f_\cL$ is constant on each connected component of $X^\an$, which proves the first part of (i). 

Assume conversely that $\cL$ is ample, and $\f_\cL$ is locally constant. Arguing on each connected component of $X^\an$, we may assume that $X$ is connected, and hence $\f_\cL\equiv c\in\Q$. After replacing $\cL$ with $\cL-c\cX_0$, we may assume $c=0$. Then $\f_{\nu^\star\cL}=\nu^\star\f_\cL\equiv 0$, and Corollary~\ref{cor:tcPL} implies that the ample, integrally closed test configuration $(\cX^\nu,\nu^\star\cL)$ is trivial, which concludes the proof of (i). 

We next turn to (ii). If $\tcX$ is trivial, then its normalization $\cX^\nu$ is trivial as well. To prove the converse, we may assume wlog that $\cL$ is ample (as any test configuration admits a $\Gm$-equivariant, ample line bundle). Assuming that $\cX^\nu$ is trivial, (i) shows that $\f_\cL$ is locally constant, and Corollary~\ref{cor:tcPL} yields as above that $\tcX$ is trivial as well, which proves (ii).  

Finally we prove (iii). If $\cX_0$ is generically reduced, then $\cX$ is regular at each generic point of $\cX_0$. Thus $\cX^\nu\to\cX$ is an isomorphism over these points, and $\cX$ is therefore trivial in codimension $1$ iff $\cX^\nu$ is. To prove (iii), we may thus assume wlog that $\cL$ is ample and $\cX$ is normal, irreducible, and trivial in codimension $1$, and we then need to show that $\cX$ is trivial. Since $\cX$ is trivial in codimension $1$, $\cX_0$ is irreducible, and $\ord_{\cX_0}=\ord_{\cX_{\triv,0}}=\sigma(v_\triv)$. After adding to $\cL$ a multiple of $\cX_0$, we assume $\f_\cL(v_\triv)=0$, and we then need to show that $\f_\cL\equiv 0$, by Corollary~\ref{cor:tcPL}. To see this, 
pick a test configuration $\cX'$ that dominates both $\cX$ and $\cX_\triv$, with morphisms $\mu\colon\cX'\to\cX$ and $\rho\colon\cX'\to\cX_\triv$, and set 
$$
D:=\mu^\star\cL-\rho^\star\cL_\triv\in\VCar(\cX')_\Q,
$$
so that $\f_\cL=\f_D$. Since $\ord_{\cX_0}(D)=\ord_{\cX_{\triv,0}}(D)=\f_D(v_\triv)=0$, $D$ is exceptional with respect to both $\mu$ and $\rho$. For $m$ divisible, we thus have 
$$
\fa_m:=\mu_\star\cO_{\cX'}(mD)\subset\cO_\cX(m\mu_\star D)\subset\cO_\cX,
$$
and similarly $\fa'_m:=\rho_\star\cO_{\cX'}(-mD)\subset\cO_{\cX_\triv}$. On the other hand, $D$ (resp.~$-D$) is semiample with respect to $\mu$ (resp.~$\rho$), and hence $\fa_m\cdot\cO_{\cX'}=\cO_{\cX'}(mD)$, $\fa'_m\cdot\cO_{\cX'}=\cO_{\cX'}(-mD)$ for $m$ large and divisible enough. It follows that $\cO_{\cX'}(\pm mD)\subset\cO_{\cX'}$; hence $D=0$, which concludes the proof.  
\end{proof}
%
%
%
%
 \section{Plurisubharmonic functions and energy pairing: the PL case}
In this section, $X$ is any projective variety, with irreducible components $(X_\a)$. We set $n:=\dim X$.   We introduce and study the class of PL functions that are $\theta$-psh for a numerical class $\theta\in\Num(X)$. We also introduce the energy pairing, defined, for the moment, on $(n+1)$-tuples of pairs $(\theta,\f)\in\Num(X)\times\PL(X)$, and study its finer properties when $\f$ is $\theta$-psh.

%
%
\subsection{Plurisubharmonic PL functions}\label{sec:pshPL}
In what follows, it will be convenient to allow real coefficients. We thus denote by 
$$
\PL_\R:=\PL(X)_\R\subset\Cz(X)
$$ 
the $\R$-vector space generated by $\PL:=\PL(X)$. Theorem~\ref{thm:VCarPL} induces an isomorphism
\begin{equation}\label{equ:PLinjlim}
\varinjlim_\cX\VCar(\cX)_\R\simeq\PL_\R.
\end{equation}
For any test configuration $\pi\colon\cX\to\A^1$, we denote by 
$$
\Num(\cX/\A^1)
$$ 
the space of $\pi$-numerical equivalence classes of all (not necessarily vertical) $\R$-Cartier divisors $D$ on $\cX$, \ie the quotient of $\Car(\cX)_\R$ by the subspace defined by $D\cdot C=0$ for all irreducible curves $C\subset\cX$ contained in some fiber of $\pi$. Since $\pi$ is projective, it follows from general theory that the $\R$-vector space $\Num(\cX/\A^1)$ is finite dimensional (this is also a consequence of Lemma~\ref{lem:relnef} below). We equip it with the corresponding $\R$-vector space topology. 

\smallskip

A class $\a\in\Num(\cX/\A^1)$ is \emph{nef} if $\a\cdot C\ge 0$ for each irreducible curve $C$ in a fiber of $\pi\colon\cX\to\A^1$. Nef classes form a closed convex cone 
$$
\Nef(\cX/\A^1)\subset\Num(\cX/\A^1).
$$
Given a morphism $\mu\colon\cX'\to\cX$ of test configurations and $\a\in\Num(\cX/\A^1)$, we have 
$$
\a\in\Nef(\cX/\A^1)\Longleftrightarrow\mu^\star\a\in\Nef(\cX'/\A^1).
$$
\begin{lem}\label{lem:relnef} For a class $\a\in\Num(\cX/\A^1)$, the following conditions are equivalent:
\begin{itemize}
\item[(i)] $\a\in\Nef(\cX/\A^1)$; 
\item[(ii)] $\a|_{\cX_0}\in\Nef(\cX_0)=\Nef(\cX_{0,\redu})$; 
\item[(iii)] $\a\cdot C\ge 0$ for all $\Gm$-invariant irreducible curves $C\subset\cX_0$. 
\end{itemize}
\end{lem}
In particular, the restriction map $\Num(\cX/\A^1)\to\Num(\cX_0)$ is injective. 

\begin{proof} Trivially, (i)$\Rightarrow$(ii)$\Rightarrow$(iii). Assume (ii), and pick a $\Gm$-linearized, ample line bundle $\cA$ on $\bar\cX$. For each $\e>0$, $(\a+\e\cA)|_{\cX_0}$ is ample, and $\a+\e\cA$ is thus relatively ample over a neighborhood of $0\in\A^1$, see~\cite[1.2.7]{Laz}. Thus $(\a+\e\cA)|_{\cX_t}$ is ample for some $t\in\Gm$, and hence for all $t\in\Gm$, thanks to the $\Gm$-equivariant isomorphism $\cX\setminus\cX_0\simeq X\times\Gm$ over $\Gm$ (note that $\Gm$, being connected, acts trivially on $\Num(\cX/\A^1)$). It follows that $\a+\e\cA$ is  nef for all $\e>0$, and hence $\a$ is  nef as well, proving (ii)$\Rightarrow$(i). Finally assume (iii). To prove (ii), we need to show that $\a\cdot C\ge 0$ for every (not necessarily $\Gm$-invariant) irreducible curve $C\subset\cX_0$, which we accomplish by way of a standard degeneration argument: denote by $C_t$ the image of $C$ under $t\in\Gm$, and note that $\a\cdot C=\a\cdot C_t$, by $\Gm$-invariance of $\a$. By properness of the components of the Chow scheme of $\cX_0$, $C_0:=\lim_{t\to 0} C_t$ exists as an effective $1$-cycle. It is $\Gm$-invariant, and hence a positive linear combination of $\Gm$-invariant irreducible curves. Thus 
$$
\a\cdot C=\lim_{t\to 0}\a\cdot C_t=\a\cdot C_0\ge 0,
$$
and we are done. 
\end{proof}

If a test configuration $\cX$ dominates $\cX_\triv=X\times\A^1$, then each class $\theta\in\Num(X)$ pulls back to a class $\theta_\cX\in\Num(\cX/\A^1)$ via the composition 
$$
\cX\to\cX_\triv=X\times\A^1\to X.
$$
In line with~\cite{Zha95,gublerlocal,CL06,siminag}, we introduce: 

\begin{defi}\label{defi:PLpsh}  We say that a function $\f\in\PL_\R$ is \emph{$\theta$-plurisubharmonic} (\emph{$\theta$-psh} for short) if it is determined by a divisor $D\in\VCar(\cX)_\R$ on a test configuration $\cX$ dominating $\cX_\triv$ such that $\theta_\cX+D\in\Nef(\cX/\A^1)$.
\end{defi}
We sometimes also say that the \emph{pair $(\theta,\f)\in\Num(X)\times\PL_\R$ is psh}. Here $\theta_{\cX}+D$ denotes, slightly abusively, the sum of $\theta_{\cX}$ and of the image of $D$ in $\Num(\cX/\A^1)$. 


This definition is independent of the choice of $\cX$. Indeed, for any morphism of test configurations $\mu\colon\cX'\to\cX$, we have 
$$
\theta_{\cX'}+\mu^\star D=\mu^\star\left(\theta_\cX+D\right),
$$
which is thus nef iff $\theta_\cX+D$ is  nef. Restricting to $\cX_1\simeq X$ further shows that a $\theta$-psh function in $\PL_\R$ can only exist when $\theta\in\Nef(X)$. 

We will use the notation
\begin{equation*}
  \PL_\R\cap\PSH(\theta)
  \quad\text{and}\quad
  \PL\cap\PSH(\theta)
\end{equation*}
for the sets of $\theta$-psh functions in $\PL_\R$ and $\PL$, respectively. At this point, the notation is purely formal, but it will be justified in~\S\ref{sec:psh}, by defining the set $\PSH(\theta)$.

When $\theta=c_1(L)\in\Num(X)$ is the numerical class of $L\in\Pic(X)_\Q$ on $X$, we simply speak of \emph{$L$-psh functions}. By Theorem~\ref{thm:testPL}, we have: 

\begin{exam}\label{exam:FSpsh} For any $\Q$-line bundle $L$ on $X$, $\cL\mapsto\f_\cL$ sets up a 1--1 correspondence between the set of semipositive non-Archimedean metrics in the sense of~\cite[Definition 6.4]{BHJ1}, \ie equivalence classes of  nef test configurations $\cL$ for $L$, and $\PL\cap\PSH(L)$. By Theorem~\ref{thm:testPL}~(ii), we thus have $\cH(L)\subset\PL\cap\PSH(L)$, the inclusion being strict in general. 
\end{exam}

\begin{rmk}\label{rmk:closedform}
Following~\cite{siminag,nama} we could define a closed $(1,1)$-form on $\Xan$ to be an element $\eta\in\varinjlim_\cX\Num(\cX/\A^1)$, and declare a function $\f\in\PL_\R$ to be $\eta$-psh if it is of the form $\f=\f_D$, with $D\in\VCar(\cX)_\R$ such that $\eta+D$ is  nef. In this paper we only consider the case when $\eta\in\Num(\cX_\triv/\A^1)$ is determined by the pullback of a class $\theta\in\Num(X)$. 
\end{rmk}

As a direct consequence of the fact that $\Nef(\cX/\A^1)\subset\Num(\cX/\A^1)$ is closed, we have: 
\begin{lem}\label{lem:limitpsh} Let $\theta_i\to\theta$ be a convergent sequence in $\Num(X)$. Pick $\f,\p\in\PL_\R$, and assume that $\f+c_i\p$ is $\theta_i$-psh for a sequence $c_i\to 0$. Then $\f$ is $\theta$-psh. 
\end{lem}

\begin{prop}\label{prop:pshPL} Pick $\theta,\theta'\in\Num(X)$, and assume that $\f,\f'\in\PL_\R$ are $\theta$-psh and $\theta'$-psh, respectively. Then:  
\begin{itemize}
\item[(i)] $\f+\f'$ is $(\theta+\theta')$-psh, and $t\f$ is $t\theta$-psh for all $t\in\R_{>0}$; 
\item[(ii)] for each $c\in\R$ and $t\in\Q_{>0}$, $\f+c$ and $t\cdot\f$ are $\theta$-psh; 
\item[(iii)] if $\theta=\theta'$ and $\f,\f'\in\PL$, then $\max\{\f,\f'\}\in\PL$ is $\theta$-psh; 
\item[(iv)] $\f$ is decreasing with respect to the partial order on $X^\an$, and hence satisfies~\eqref{equ:maxprin} (see Lemma~\ref{lem:maxprin}); 
\item[(v)] for any morphism $f\colon Y\to X$ from a projective variety, $f^\star\f$ is $f^\star\theta$-psh; 
\item[(vi)] if we further assume that $f\colon Y\to X$ is surjective, then
$$
\f\,\,\theta\text{-psh}\Longleftrightarrow\f|_{X_\a^\an}\,\,\theta|_{X_\a}\text{-psh for all }\a\Longleftrightarrow f^\star\f\,\,f^\star\theta\text{-psh}. 
$$
\end{itemize}
\end{prop}
In particular, (i) shows that $\PL_\R\cap\PSH(\theta)$ is a convex subset of $\PL_\R$.

\begin{lem}\label{lem:pshPL} For any $\theta\in\Num(X)$ and $\f\in\PL_\R\cap\PSH(\theta)$, there exists a sequence $\f_m\in\cH(L_m)$ with $L_m\in\Pic(X)_\Q$ such that $\f_m\to\f$ uniformly on $X^\an$ and $c_1(L_m)\to\theta$. 
\end{lem} 

\begin{proof} Pick a test configuration $\cX$ dominating $\cX_\triv$ and $D\in\VCar(\cX)_\R$ such that $\f=\f_D$. By openness of the (relatively) ample cone of $\Num(\cX/\A^1)$, we can find a sequence of $\Q$-line bundles $L_m$ on $X$ and $D_m\in\VCar(\cX)_\Q$ such that $c_1(L_m)\to\theta$ in $\Num(X)$, $D_m\to D$ in $\VCar(\cX)_\R$, and $c_1(L_m)_\cX+D_m\in\Num(\cX/\A^1)$ is ample for all $i$. In particular, $\cL_m:=(L_m)_\cX+D_m$ is a semiample test configuration for $L$, and hence $\f_m:=\f_{\cL_m}=\f_{D_m}\in\cH(L_m)$. Furthermore,  $D_m\to D$ implies $\f_m\to\f$ uniformly. 
\end{proof}

\begin{proof}[Proof of Proposition~\ref{prop:pshPL}] Pick a test configuration $\cX$ dominating $\cX_\triv$ and $D,D'\in\VCar(\cX)_\R$ that determine $\f,\f'$. By assumption, $\theta_\cX+D$ and $\theta'_\cX+D'$ are  nef. Thus 
$$
(\theta+\theta')_\cX+D+D' =\left(\theta_\cX+D \right)+\left(\theta'_\cX+D' \right)\quad\text{and}\quad(t\theta)_\cX+ tD =t\left(\theta_\cX+D \right)
$$
are nef as well, proving (i). 

For each $c\in\R$ we have $\f+c=\f_{D+c\cX_0}$. Since $\cX_0=\pi^\star(0)$ vanishes in $\Num(\cX/\A^1)$, $\theta_\cX+D+c\cX_0$ is  nef, and $\f+c$ is thus $\theta$-psh. To prove the remaining part of~(ii), it suffices to prove that if $d\in\Z_{>0}$, then $\f$ is $\theta$-psh iff $d\cdot\f$ is $\theta$-psh. In view of Lemma~\ref{lem:basechange}, this follows from the fact that $\theta_\cX+D\in\Num(\cX/\A^1)$ is nef iff its pullback under the base change of $\cX$ with respect to $\unipar\mapsto\unipar^d$ is nef.

Assume now $\theta=\theta'$ and $D,D'\in\VCar(\cX)_\Q$. By openness of $\Amp(X)$, we can find an $\Q$-line bundle $L$ on $X$ such that $c_1(L)-\theta$ is ample and arbitrarily small. Then $\f,\f'$ are $L$-psh, and it will be enough to show that $\max\{\f,\f'\}$ is $L$-psh as well, by Lemma~\ref{lem:limitpsh}. After perhaps passing to a higher test configuration, we may assume that $\cX$ supports an ample test configuration $\cL$ for $L$ (Lemma~\ref{lem:amphigh}). For each $\e\in\Q_{>0}$, the $\Q$-line bundle $L_\cX+D+\e\cL$ is  ample, and hence  semiample. Thus $\f+\e\f_\cL$ is Fubini--Study for $(1+\e)L$, and similarly for $\f'+\e\f_\cL$. By Proposition~\ref{prop:FS}, it follows that 
$$
\max\{\f+\e\f_\cL,\f'+\e\f_\cL\}=\max\{\f,\f'\}+\e\f_\cL
$$
is Fubini--Study for $(1+\e)L$, and hence $(1+\e)L$-psh. Using Lemma~\ref{lem:limitpsh} we conclude, as desired, that $\max\{\f,\f'\}$ is $L$-psh. This proves (iii). 

By Proposition~\ref{prop:FS}, any Fubini--Study function is decreasing, and (iv) thus follows from Lemma~\ref{lem:pshPL}. 

To prove (v), note that the $\Gm$-equivariant rational map $Y\times\Gm\dashrightarrow\cX$ induced by $f$ admits a $\Gm$-equivariant resolution of indeterminacies, which is thus a test configuration $\cY$ for $Y$ dominating $\cY_\triv$. Also denoting by $f\colon\cY\to\cX$ the corresponding morphism, we have $f^\star\f=\f_{f^\star D}$, and the result thus follows from the fact that $\theta_{\cY}+f^\star D=f^\star(\theta_\cX+D)$ is  nef, as the pullback of a  nef class. If $f\colon Y\to X$ is surjective, then $f\colon\cY\to\cX$ is surjective as well, and (vi) follows. 
\end{proof}

Given a test configuration $\pi\colon\cX\to\A^1$, we denote by $\Amp(\cX/\A^1)\subset\Nef(\cX/\A^1)$ the set of $\pi$-ample classes $\a\in\Num(\cX/\A^1)$, that is, classes whose restriction to $\cX_0$ is ample.

\begin{defi}\label{defi:Hdom} For any $\om\in\Amp(X)$, we denote by $\cH^\dom(\om)\subset\PL\cap\PSH(\om)$ the set of $\om$-psh PL functions of the form $\f_D$ with $D\in\VCar(\cX)_\Q$ for a test configuration $\cX$ dominating $\cX_\triv$ and $\om_\cX+D\in\Amp(\cX/\A^1)$.
\end{defi} 
After pulling back by the (finite) integral closure morphism $\tcX\to\cX$, one may always arrange that $\cX$ is integrally closed in this definition. 

\begin{exam} If $\om=c_1(L)$ for an ample $\Q$-line bundle, then $\cH(L)$ is in 1--1 correspondence with the set of ample, integrally closed test configurations $(\cX,\cL)$ for $(X,L)$ (see Corollary~\ref{cor:tcPL}), and 
$$
\cH^\dom(L):=\cH^\dom(c_1(L))\subset\cH(L)
$$ 
corresponds to the subset such that $\cX$ dominates $\cX_\triv$.
\end{exam}

As in Corollary~\ref{cor:FSideal}, we have
\begin{equation}\label{equ:HdomPL}
\Q_+\cH^\dom(\om)=\PL^+(X)
\end{equation}
for any $\om\in\Amp(X)$, so that $\cH^\dom(\om)$ spans the $\Q$-vector space $\PL(X)$. 

\begin{prop}\label{prop:PLdom} If $\om\in\Amp(X)$, then every $\f\in\PL_\R\cap\PSH(\om)$ is a uniform limit of functions in $\cH^\dom(\om)$.
\end{prop}

\begin{lem}\label{lem:amplemodel} For any $\om\in\Amp(X)$, the set of test configurations $\cX$ dominating $\cX_\triv$ such that there exists $H\in\VCar(\cX)_\Q$ with $\om_\cX+H\in\Amp(\cX/\A^1)$ is cofinal in the set of all test configurations. 
\end{lem}

\begin{proof} By Lemma~\ref{lem:relamp}, any test configuration for $X$ is dominated by a test configuration $\cX$ with a morphism $\mu\colon\cX\to\cX_\triv$ such that $\VCar(\cX)_\Q$ contains a $\mu$-ample divisor $A$. Since $\om_{\cX_\triv}\in\Num(\cX_{\triv}/\A^1)$ is ample, $\om_\cX+\e A= \mu^\star\om_{\cX_\triv}+\e A\in\Num(\cX/\A^1)$ is ample for $\e\in\Q_{>0}$ small enough, and setting $H:=\e A$ yields the result. 
\end{proof}

\begin{proof}[Proof of Proposition~\ref{prop:PLdom}] Pick any $\f\in\PL_\R\cap\PSH(\om)$, and write $\f=\f_D$ for some $D\in\VCar(\cX)_\R$, where $\cX$ is an integrally closed test configuration. Thus $\om_\cX+D \in\Num(\cX/\A^1)$ is  nef. By Lemma~\ref{lem:amplemodel} we may assume that there exists $H\in\VCar(\cX)_\R$ such that $\om_\cX+H$ is ample.   Pick a basis $(C_i)$ for $\VCar(\cX)_\Q$ and write $D=\sum_i c_iC_i$, $H=\sum_i c'_i C_i$ with $c_i,c'_i\in\R$. For $m\gg1$, the class  $\om_\cX+D +\frac1m H\in\Num(\cX/\A^1)$ is ample, and we can pick $\e_{m,i}\in(0,\frac1m)$ such that if $D_m:=\sum_i(c_i+\frac{c'_i}m+\e_{m,i})C_i$, then $D_m\in\VCar(\cX)_\Q$ and $\om_\cX+D_m$ is ample. It follows that $\f_{D_m}\in\cH^\dom(\om)$, and that $\f_{D_m}-\f=\tfrac 1m\f_H+\sum_i\e_{m,i}\f_{C_i}$ tends to $0$ uniformly on $X^\an$.  
\end{proof}
%
%
\subsection{The energy pairing}\label{sec:enpairing}
Recall that every test configuration $\pi\colon\cX\to\A^1$ admits a canonical compactification $\pi\colon\bar\cX\to\P^1$. If $\cX$ dominates $\cX_\triv$, then any  $\theta\in\Num(X)$ pulls back to a class $\theta_{\bar\cX}\in\Num(\bar\cX)$, whose image in $\Num(\cX/\A^1)$ coincides with $\theta_\cX$ considered above. 

Pick an $(n+1)$-tuple of pairs $(\theta_i,\f_i)\in\Num(X)\times\PL_\R$, $i=0,\dots,n$, and choose a test configuration $\cX$ dominating $\cX_\triv$ and divisors $D_i\in\VCar(\cX)_\R$ that determined the $\f_i$. Following~\cite[Definition 6.11]{BHJ1}, we introduce: 

\begin{defi}\label{defi:intpairing} The \emph{energy pairing} takes an $(n+1)$-tuple of pairs $(\theta_i,\f_i)\in\Num(X)\times\PL_\R$, $i=0,\dots,n$ to
\begin{equation}\label{equ:intpairing}
(\theta_0,\f_0)\inter(\theta_n,\f_n):=\left(\theta_{0,\bar\cX}+ D_0 \right)\inter\left(\theta_{n,\bar\cX}+ D_n \right)\in\R. 
\end{equation}
\end{defi}
The right-hand side is an intersection number against the fundamental class 
\begin{equation}\label{equ:fundclass}
[\bar\cX]=\sum_{\dim X_\a=n}[\bar\cX_\a]
\end{equation}
of the $(n+1)$-dimensional projective variety $\bar\cX$, with $ \theta_{i,\bar\cX}+D_i \in\Num(\bar\cX)$ now denoting (slightly abusively again) the sum of $\theta_{i,\bar\cX}$ and of the image of $D_i$ in $\Num(\bar\cX)$. By the projection formula, this definition is independent of the choice of $\cX$.

\begin{rmk}\label{rmk:enDeligne} Assume $\theta_i=c_1(L_i)$ with $L_i$ a line bundle on $X$, $i=0,\dots,n$. In that case, the energy pairing can be interpreted as a metric on the Deligne pairing $\langle L_0,\dots,L_n\rangle$, as follows. Pick $\f_i\in\PL_\R$, $i=0,\dots,n$. Using the trivial metric on $L_i$, $\f_i$ can be identified with a continuous dpsh metric $\phi_i$ on $L_i$ in the sense of~\cite[Definition 8.7]{BE}. By~\cite[Theorem 8.16]{BE}, one gets an induced metric $\langle\phi_0,\dots,\phi_n\rangle$ on the line $\langle L_0,\dots,L_n\rangle$, that can in turn be viewed as a real number, thanks to the trivial metric. Using~\cite[Theorem 8.18]{BE}, one can check that this number coincides with $(\theta_0,\f_0)\inter(\theta_n,\f_n)$. 
\end{rmk}

%
\begin{prop}\label{prop:intpairing} The energy pairing is a symmetric, multilinear form on $\Num(X)\times\PL_\R$. For all tuples $(\theta_i,\f_i)\in\Num(X)\times\PL_\R$, $i=0,\dots,n$, we further have: 
\begin{itemize}
\item[(i)] $(0,1)\cdot(\theta_1,\f_1)\inter(\theta_n,\f_n)=\left(\theta_1\inter\theta_n\right)_X$;
\item[(ii)]   for all $c_0,\dots,c_n\in\R$ we have 
\begin{align*} 
(\theta_0,\f_0+c_0)\inter(\theta_n,\f_n+c_n) & =(\theta_0,\f_0)\inter(\theta_n,\f_n)\\
& +\sum_{i=0}^n c_i(\theta_0\inter\widehat{\theta_i}\inter\theta_n)_X; 
\end{align*}  
\item[(iii)] $(\theta_0,0)\inter(\theta_n,0)=0$; 
\item[(iv)] $(\theta_0,t\cdot\f_0)\inter(\theta_n,t\cdot\f_n)=t(\theta_0,\f_0)\inter(\theta_n,\f_n)$ for all $t\in\Q_{>0}$; 
\item[(v)] denoting by $\nu\colon X^\nu\to X$ the normalization morphism, we have 
\begin{align*} 
(\theta_0,\f_0)\inter(\theta_n,\f_n) & =\sum_{\dim X_\a=n}(\theta_0,\f_0)|_{X_\a}\inter(\theta_n,\f_n)|_{X_\a}\\
& =(\nu^\star\theta_0,\nu^\star\f_0)\inter(\nu^\star\theta_n,\nu^\star\f_n); 
\end{align*}
\item[(vi)] pick an integrally closed test configuration $\cX$ dominating $\cX_\triv$, $D_1,\dots,D_n\in\VCar(\cX)_\R$, and set $\f_i:=\f_{D_i}$ for $i=1,\dots,n$. For all $\f\in\PL_\R$ we then have 
\begin{equation}\label{equ:intnumb}
(0,\f)\cdot(\theta_1,\f_1)\inter(\theta_n,\f_n)=\sum b_E\,\f(v_E)(\theta_{1,\cX}+ D_1 )|_E\inter (\theta_{n,\cX}+ D_n )|_E
\end{equation}
where the sum runs over the irreducible components $E$ of $\cX_0$ and $b_E:=\ord_E(\cX_0)=\ord_E(\unipar)$. 
\end{itemize}
\end{prop}
The right-hand side in (i) is an intersection number against the fundamental class $[X]=\sum_{\dim X_\a=n}[X_\a]$. We emphasize that $\cX$ in (vi) does not depend on $\f$. 

\begin{proof} The first assertion is an immediate consequence of the multilinearity and symmetry of the right-hand side of~\eqref{equ:intpairing}. In the notation of Definition~\ref{defi:intpairing},
\begin{multline*}
 (0,1)\cdot(\theta_1,\f_1)\inter(\theta_n,\f_n)=[\cX_0]\cdot (\theta_{1,\bar\cX}+ D_1 )\inter (\theta_{n,\bar\cX}+ D_n )\\
=[\cX_1]\cdot (\theta_{1,\bar\cX}+ D_1 )\inter (\theta_{n,\bar\cX}+ D_n )=(\theta_{1,\bar\cX})|_{\cX_1}\inter(\theta_{n,\bar\cX})|_{\cX_1}=(\theta_1\inter\theta_n)_X,
\end{multline*}
where the second equality holds by flatness of $\pi\colon\bar\cX\to\P^1$. This proves (i), which implies (ii), by multilinearity, while (iii) follows from the projection formula applied to $X\times\P^1\to X$. 

To prove (iv), we may assume that $t\in\Z_{>0}$, since we are dealing with a group action of $\Q_{>0}$. The result is now a simple consequence of the geometric description of $t\cdot\f$ in terms of base change in Lemma~\ref{lem:basechange} and the projection formula (compare~\cite[Lemma 6.13]{BHJ1}). Next, (v) follows~\eqref{equ:fundclass}, $\nu_\star[\bar\cX^\nu]=[\bar\cX]$, and the projection formula. Finally, for (vi), pick a morphism $\mu\colon\cX'\to\cX$ of test configurations and $D_0\in\VCar(\cX')_\R$ such that $\f=\f_{D_0}$. Then
\begin{align*}
(0,\f)\cdot(\theta_1,\f_1)\inter(\theta_n,\f_n)
  &= D_0 \cdot\left(\theta_{1,\bar\cX'}+\mu^\star D_1 \right)\inter\left(\theta_{n,\bar\cX'}+\mu^\star  D_n \right)\\
  &= D_0 \cdot\mu^\star\left(\left(\theta_{1,\bar\cX}+ D_1 \right)\inter\left(\theta_{n,\bar\cX}+ D_n \right)\right)\\
  &=\mu_\star [D_0] \cdot\left(\theta_{1,\bar\cX}+ D_1 \right)\inter\left(\theta_{n,\bar\cX}+ D_n \right)\\
  &=\sum_E\ord_E(D_0)[E]\cdot\left(\theta_{1,\bar\cX}+ D_1 \right)\inter\left(\theta_{n,\bar\cX}+ D_n \right)\\
&=\sum_Eb_E\,\f(v_E)[E]\cdot\left(\theta_{1,\bar\cX}+ D_1 \right)\inter\left(\theta_{n,\bar\cX}+ D_n \right),
\end{align*}
where we used the projection formula in the third equality.
\end{proof}

For psh pairs,~\eqref{equ:intnumb} implies the following crucial monotonicity property: 
\begin{lem}\label{lem:enmono} Consider psh pairs $(\theta_i,\f_i)\in\Nef(X)\times\PL_\R$, $i=1,\dots,n$, and pick also $\theta_0\in\Num(X)$ and $\f_0,\f'_0\in\PL_\R$. Then 
$$
\f_0\le\f'_0\Longrightarrow(\theta_0,\f_0)\cdot(\theta_1,\f_1)\inter(\theta_n,\f_n)\le(\theta_0,\f'_0)\cdot(\theta_1,\f_1)\inter(\theta_n,\f_n). 
$$
\end{lem}

Combined with Proposition~\ref{prop:intpairing}~(iii), this yields: 
\begin{cor}\label{cor:enmono} For all psh pairs $(\theta_i,\f_i)\in\Nef(X)\times\PL_\R$, $i=0,\dots,n$, we have
$$
\forall i\,\,\f_i\le 0\Longrightarrow(\theta_0,\f_0)\inter(\theta_n,\f_n)\le 0;
$$
$$
\forall i\,\,\f_i\ge 0\Longrightarrow(\theta_0,\f_0)\inter(\theta_n,\f_n)\ge 0.
$$
 \end{cor}

Anticipating on the more general construction of \S\ref{sec:mixedMA}, it is convenient to interpret~\eqref{equ:intnumb} by attaching to any tuple $(\theta_i,\f_i)\in\Num(X)\times\PL_\R$, $i=1,\dots,n$ its \emph{mixed Monge--Amp\`ere measure}
\begin{equation}\label{equ:mixedMAPL}
\bigwedge_{i=1}^n(\theta_i+\ddc\f_i):=\sum_E c_E\d_{v_E}
\end{equation}
with $c_E:=b_E(\theta_{1,\cX}+ D_1 )|_E\inter (\theta_{n,\cX}+ D_n )|_E$. It is thus a signed Radon measure on $X^\an$, with support a finite subset of $X^\div$, and characterized by
\begin{equation}\label{equ:mixedMA2PL}
\int_{X^\an} \f\,\bigwedge_{i=1}^n(\theta_i+\ddc\f_i)=(0,\f)\cdot (\theta_1,\f_1)\inter(\theta_n,\f_n)
\end{equation}
for all $\f\in\PL$ (recall that $\PL$ is dense in $\Cz(X)$, see Theorem~\ref{thm:PLdense}). By Proposition~\ref{prop:intpairing}~(i), 
\begin{equation}\label{equ:mixedMAtotal}
\int_{X^\an}\bigwedge_{i=1}^n(\theta_i+\ddc\f_i)=(\theta_1\inter\theta_n)_X.
\end{equation}
The symmetry of the intersection pairing yields the `integration by parts' formula
\begin{equation}\label{equ:intpartPL}
\int_{X^\an} \f_0\,\ddc \f_1\wedge\bigwedge_{i=2}^n(\theta_i+\ddc\f_i)=\int_{X^\an}\f_1\,\ddc\f_0\wedge\bigwedge_{i=2}^n(\theta_i+\ddc\f_i)
\end{equation}
for all $\f_i\in\PL_\R$, while Proposition~\ref{prop:intpairing}~(iv) yields 
\begin{equation}\label{equ:MAscale} 
\bigwedge_i(\theta_i+\ddc(t\cdot\f_i))=t_\star\bigwedge_i(\theta_i+\ddc\f_i)
\end{equation}
for all $t\in\Q_{>0}$. When $\f_i=0$ for some $i$, we drop the term $\ddc\f_i$ from the notation.
\begin{exam} We have
\begin{equation}\label{equ:mixedMA0}
\theta_1\winter\theta_n=\sum_{\dim X_\a=n}\left(\theta_1\inter\theta_n\right)_{X_\a}\d_{v_{\triv,\a}}. 
\end{equation}
\end{exam}

\begin{exam}\label{exam:defconeMA} Assume $Z\subset X$ are both smooth and irreducible, with associated valuation $\ord_Z\in X^\div$. Set $d:=\dim Z$, and consider the function $\f_Z$ of Example~\ref{exam:defconePL}. Denoting by $\mu\colon \cX\to\cX_\triv$ the blowup of $Z\times\{0\}$ and $P$ its exceptional divisor, we have $\cX_0=\tX+P$ (see Example~\ref{exam:defcone}), $\f_Z=-\f_P$, $v_{\tX}=v_\triv$ and $v_P=\ord_Z$ (see Example~\ref{exam:defconeval}). For any $\theta\in\Num(X)$, we have
$$
(\theta+\ddc\f_Z)^n=c\,\d_{\ord_Z}+\left((\theta^n)-c\right)\d_{v_{\triv}}
$$
where 
$$
c:=(\theta_\cX-P)|_P^n=\left(\mu^\star\theta|_Z+\cO_P(1)\right)^n=\sum_{j=0}^{d}{n\choose j}\left(\theta|_Z^j\cdot s_{d-j}(Z)\right),
$$
with $s_i(Z)$ the $i$-th Segre class of $Z$ (see~\cite[Chapter 4]{FultonInter}). 
\end{exam}

\begin{exam}\label{exam:curvePLMA} Assume that $X$ is a smooth irreducible curve. Then each $\f\in\PL_\R$ is a PL function on any ray $(t\ord_p)_{t\ge 0}$ with $p\in X(k)$, constant on all but finitely many of these. Moreover, for any $\theta\in\Num(X)$, we have
\begin{equation}\label{eqi:curvePLMA}
\theta+\ddc\f=(\deg\theta)\d_{v_\triv}+\D\f,
\end{equation}
where
$$
\D\f:=\sum_{p\in X(k)}\left[\frac{d}{dt}\bigg|_{0+}\f(t\ord_p)\d_{v_\triv}+\frac{d^2}{dt^2}\f(t\ord_p)\right]
$$
is (up to a sign) the \emph{tree Laplacian}. To see all this, note that the case $\f=0$ follows from~\ref{equ:mixedMA0}, so by linearity we may assume $\theta=0$.
We may further assume $\f=a\cdot\f_p$ with $a\in\Q_{>0}$ and $p\in X(k)$, since these functions span $\PL_\R$ (see Example~\ref{exam:curvePL}) and the operators $\ddc$ and $\D$ are $\R$-linear. By~\eqref{equ:MAscale}, we can also assume $a=1$. Then $\f(t\ord_p)=\max\{-t,-1\}$, and $\f(t\ord_q)=0$ for $q\ne0$. This implies $\D\f=\d_{\ord_p}-\d_{v_\triv}$, which coincides with $\ddc\f$, by Example~\ref{exam:defconeMA}. 
\end{exam}

When $(\theta_i,\f_i)\in\Num(X)\times\PL_\R$ is psh for $i=1,\dots,n$, the measure~\eqref{equ:mixedMAPL} is positive (which is equivalent to Lemma~\ref{lem:enmono}). Conversely: 

\begin{thm}\label{thm:goodman} Pick $\om\in\Amp(X)$ and $(\theta,\f)\in\Num(X)\times\PL_\R$. Then $\f$ is $\theta$-psh iff 
$$
(\theta+\ddc\f)\wedge(\om+\ddc\p)^{n-1}\ge 0
$$ 
for all $\p\in\cH^\dom(\om)$.
\end{thm}

\begin{exam}\label{exam:curvePLpsh} If $X$ is a smooth irreducible curve and $(\theta,\f)\in\Num(X)\times\PL_\R$, then $\f$ is $\theta$-psh iff $\theta+\ddc\f\ge 0$. In view of Example~\ref{exam:curvePLMA}, this amounts to: 
\begin{itemize} 
\item[(i)] for each $p\in X(k)$, $t\mapsto\f(t\ord_p)$ is convex---and hence decreasing, being bounded above; 
\item[(ii)] $\deg\theta+\sum_{p\in X(k)}\frac{d}{dt}\big|_{0+}\f(t\ord_p)\ge 0$. 
\end{itemize}
\end{exam}

As an important consequence of Theorem~\ref{thm:goodman}, we infer the following analogue of~\cite[Theorem 5.11]{siminag} and~\cite[Theorem 5.5]{GM}:

\begin{cor}\label{cor:goodman} For each $\theta\in\Num(X)$, $\PL_\R\cap\PSH(\theta)$ is closed in $\PL_\R$ with respect to the topology of pointwise convergence on $X^\div$. 
\end{cor}

\begin{proof}[Proof of Theorem~\ref{thm:goodman}] By Lemma~\ref{lem:amplemodel}, we can pick a test configuration $\cX$ dominating the trivial one   such that there exists $H\in\VCar(\cX)_\Q$ with $\om_\cX+H\in\Amp(\cX/\A^1)$   and $D\in\VCar(\cX)_\R$ with $\f=\f_{D}$. By Lemma~\ref{lem:relnef}, we need to show that $(\theta_{\cX}+D )\cdot C\ge 0$ for each $\Gm$-invariant irreducible curve $C\subset\cX_0$. Consider, as in the proof of~\cite[Proposition 8]{Good}, the blowup $\mu\colon\cX'\to\cX$ along $C$, with exceptional divisor $F$. Then $\cX'$ is a test configuration for $X$, and $-F$ is $\mu$-ample, so that  
$$
\a:=\om_{\cX'}+\mu^\star H-\e F\in\Amp(\cX'/\A^1)
$$ 
for $0<\e\ll 1$. Since $F$ is an effective Cartier divisor dominating $C$, there exists $a\in\Q_{>0}$ such that 
$$
\mu_\star\left(F\cdot\a^{n-1}\right)=a C
$$ 
in $\mathrm{N}_1(\cX/\A^1)$, and the projection formula thus yields
$$
\left(\theta_{\cX}+D\right)\cdot C=a^{-1} F\cdot \left(\theta_{\cX'}+\mu^\star D\right)\cdot \a^{n-1}. 
$$
Now
$$
F\cdot\left(\theta_{\cX'}+\mu^\star D \right)\cdot\a^{n-1}=\int\f_F(\theta+\ddc\f)\wedge(\om+\ddc\p)^{n-1}
$$
with $\p:=\f_{\mu^\star H-\e F}\in\cH^\dom(\om)$. By assumption, the right-hand integral is nonnegative. Thus $\left(\theta_{\cX}+D\right)\cdot C\ge 0$, and we are done.  
\end{proof}

In view of Theorem~\ref{thm:goodman}, Corollary~\ref{cor:goodman} is a direct consequence of the following continuity result. 

\begin{lem}\label{lem:contdivMA} Consider a tuple $(\theta_i,\f_i)\in\Num(X)\times\PL_\R$, $i=1,\dots,n$, and assume that $\f_1$ is the pointwise limit on $X^\div$ of a net $(\f_{1j})_j$ in $\PL_\R$. Then
$$
(\theta_1+\ddc\f_{1j})\wedge\bigwedge_{i=2}^n(\theta_i+\ddc\f_i)\to (\theta_1+\ddc\f_1)\wedge\bigwedge_{i=2}^n(\theta_i+\ddc\f_i)
$$
weakly as measures on $X^\an$. 
\end{lem}

\begin{proof} By density of $\PL$ in $\Cz(X)$ (\cf Theorem~\ref{thm:PLdense}), we need to show 
\begin{equation}\label{equ:contdivMA}
\int\f_0\,(\theta_1+\ddc\f_{1j})\wedge\bigwedge_{i=2}^n(\theta_i+\ddc\f_i)\to\int\f_0\,(\theta_1+\ddc\f_1)\wedge\bigwedge_{i=2}^n(\theta_i+\ddc\f_i)
\end{equation}
for all $\f_0\in\PL$. Now 
$$
\int\f_0\,(\theta_1+\ddc\f_{1j})\wedge\bigwedge_{i=2}^n(\theta_i+\ddc\f_i)
$$
$$
=\int\f_0\,\theta_1\wedge\bigwedge_{i=2}^n(\theta_i+\ddc\f_i)+\int\f_0\,\ddc\f_{1j}\wedge\bigwedge_{i=2}^n(\theta_i+\ddc\f_i),
$$
and~\eqref{equ:intpartPL} yields
$$
\int\f_0\,\ddc\f_{1j}\wedge\bigwedge_{i=2}^n(\theta_i+\ddc\f_i)=\int\f_{1j}\,\ddc\f_0\wedge\bigwedge_{i=2}^n(\theta_i+\ddc\f_i)
$$
$$
\to\int\f_1\,\ddc\f_0\wedge\bigwedge_{i=2}^n(\theta_i+\ddc\f_i)=\int\f_0\,\ddc\f_1\wedge\bigwedge_{i=2}^n(\theta_i+\ddc\f_i),
$$
since $\ddc\f_0\wedge\bigwedge_{i=2}^n(\theta_i+\ddc\f_i)$ is supported in a finite subset of $X^\div$. This proves~\eqref{equ:contdivMA}, and concludes the proof. 
\end{proof} 
%
%
\subsection{Convexity and Hodge-type estimates}\label{sec:estimates} 
The following consequence of the Hodge index theorem is one main building block of all the results to follow. 

\begin{lem}\label{lem:Hodge} For all psh pairs $(\theta_i,\f_i)\in\Num(X)\times\PL_\R$, $i=1,\dots,n-1$ and $\f\in\PL_\R$, we have 
$$
(0,\f)^2\cdot(\theta_1,\f_1)\inter(\theta_{n-1},\f_{n-1})\le 0. 
$$
\end{lem}
\begin{proof} When $(\theta_i,\f_i)\in\Num(X)_\Q\times\VCar(\cX)_\Q$ and $\f\in\PL$, this follows from~\cite[Lemma 6.14]{BHJ1}, itself a consequence of the Hodge index theorem. The general case easily follows by approximation, arguing as in the proof of Corollary~\ref{cor:goodman} and using Proposition~\ref{prop:PLdom}.
\end{proof}

\begin{thm}\label{thm:enconc} Pick $p\in\{0,\dots,n+1\}$. For $i=0,\dots,n-p$, let $(\theta_i,\f_i)\in\Num(X)\times\PL_\R$ be a psh pair, and write for brevity $\Ga:=(\theta_0,\f_0)\inter(\theta_{n-p},\f_{n-p})$. For any class $\theta\in\Num(X)$, the function 
$$
\f\mapsto(\theta,\f)^p\cdot\Ga
$$ 
is then concave on $\PL_\R\cap\PSH(\theta)$. 
\end{thm}
Here $\Ga$ is a purely notational device, and we are not trying to make sense of it as a cycle class, for example. 

\begin{lem}\label{lem:en} For all $\f,\p\in\PL_\R$ we have 
\begin{equation}\label{equ:en}
(\theta,\f)^p\cdot\Ga-(\theta,\p)^p\cdot\Ga
=\sum_{j=0}^{p-1}(0,\f-\p)\cdot(\theta,\f)^j\cdot(\theta,\p)^{p-1-j}\cdot\Ga
\end{equation}
and
\begin{equation}\label{equ:ender}
\frac{d}{dt}\bigg|_{t=0}(\theta,\f+t\p)^p\cdot\Ga=p(0,\p)\cdot(\theta,\f)^{p-1}\cdot\Ga. 
\end{equation}
\end{lem} 
\begin{proof} This follows from straightforward computations based on the multilinearity and symmetry of the energy pairing.
\end{proof}

\begin{lem}\label{lem:enincr} For all $\f,\p\in\PL_\R\cap\PSH(\theta)$, the sequence 
$$
a_j:= (0,\f-\p)\cdot(\theta,\f)^j\cdot(\theta,\p)^{p-1-j}\cdot\Ga
$$ 
is decreasing on $\{0,\dots,p-1\}$, and 
\begin{equation}\label{equ:enconc}
p a_{p-1}\le (\theta,\f)^p\cdot\Ga-(\theta,\p)^p\cdot\Ga\le p a_0. 
\end{equation}
\end{lem}
\begin{proof} For $j=0,\dots,p-1$, we have 
$$
a_{j+1}-a_j
=(0,\f-\p)^2\cdot(\theta,\f)^j\cdot(\theta,\p)^{p-j-2}\cdot\Ga\le 0,
$$
by Lemma~\ref{lem:Hodge}. In view of~\eqref{equ:en} this implies~\eqref{equ:enconc}, 
\end{proof}

\begin{proof}[Proof of Theorem~\ref{thm:enconc}] By~\eqref{equ:ender} and~\eqref{equ:enconc}, we have, for any two $\f,\p\in\PL_\R\cap\PSH(\theta)$,
$$
(\theta,\f)^p\cdot\Ga\le(\theta,\p)^p\cdot\Ga+\frac{d}{dt}\bigg|_{t=0}\left(\theta,(1-t)\p+t\f\right)^p\cdot\Ga. 
$$
This is equivalent to the concavity of $\f\mapsto(\theta,\f)^p\cdot\Ga$ on $\PL_\R\cap\PSH(\theta)$. 
\end{proof}

By Lemma~\ref{lem:Hodge}, we may introduce: 

\begin{defi} To each $(n-1)$-tuple of psh pairs $(\theta_i,\f_i)\in\Num(X)\times\PL_\R$, $i=1,\dots,n-1$, we associate a seminorm on $\PL_\R$ by setting
$$
\|\f\|_{(\theta_1,\f_1)\inter(\theta_{n-1},\f_{n-1})}:=\sqrt{-(0,\f)^2\cdot(\theta_1,\f_1)\inter(\theta_{n-1},\f_{n-1})}.
$$
\end{defi}

For all $\f,\p\in\PL_\R$, we then have the Cauchy--Schwarz inequality
\begin{equation}\label{equ:CS}
\left|(0,\f)\cdot(0,\p)\cdot(\theta_1,\f_1)\inter(\theta_{n-1},\f_{n-1})\right|\le\|\f\|_{(\theta_1,\f_1)\inter(\theta_{n-1},\f_{n-1})}\|\p\|_{(\theta_1,\f_1)\inter(\theta_{n-1},\f_{n-1})}. 
\end{equation}
For the remainder of this section we fix a nef class $\theta\in\Nef(X)$. 

\begin{defi}\label{defi:dtheta} For all $\f,\p\in\PL_\R\cap\PSH(\theta)$ we set 
\begin{equation}\label{equ:dtheta}
d_\theta(\f,\p):=\max_{0\le j\le n-1}\|\f-\p\|^2_{(\theta,\f)^j\cdot(\theta,\p)^{n-1-j}}. 
\end{equation}
\end{defi} 

When $\p=0$ we simply set 
$$
d_\theta(\f):=d_\theta(\f,0).
$$
We first note the following basic monotonicity property.

\begin{lem}\label{lem:dmono} Assume $\theta'\in\Num(X)$ satisfies $\theta'\ge\theta$. For all 
$$
\f,\p\in\PL_\R\cap\PSH(\theta)\subset\PL_\R\cap\PSH(\theta'),
$$
we then have 
$$
d_{\theta'}(\f,\p)\ge d_\theta(\f,\p).
$$
\end{lem}

\begin{proof} By assumption, $\theta'':=\theta'-\theta$ is nef. For $j=0,\dots,n-1$, we thus have 
$$
\|\f-\p\|^2_{(\theta',\f)^j\cdot(\theta',\p)^{n-1-j}}=-(0,\f-\p)^2\cdot\left((\theta,\f)+(\theta'',0)\right)^j\cdot\left((\theta,\p)+(\theta'',0)\right)^{n-1-j}. 
$$
Expanding out, Lemma~\ref{lem:Hodge} yields 
$$
(0,\f-\p)^2\cdot\left((\theta,\f)+(\theta'',0)\right)^j\cdot\left((\theta,\p)+(\theta'',0)\right)^{n-1-j}\le (0,\f-\p)^2\cdot(\theta,\f)^j\cdot(\theta,\p)^{n-1-j},
$$
and the result follows. 
\end{proof}

In analogy to~\cite[Theorem~1.8]{BBEGZ}, we shall prove that $d_\theta$ satisfies a quasi-triangle inequality. 

\begin{thm}\label{thm:triangle}
  For all $\f_1,\f_2,\f_3\in\PL_\R\cap\PSH(\theta)$ we have
  \begin{equation*}
    d_\theta(\f_1,\f_2)\lesssim d_\theta(\f_1,\f_3)+d_\theta(\f_3,\f_1).
  \end{equation*}
\end{thm}

Recall that we write $x\lesssim y$ if $x\le C_n  y$ for a constant $C_n>0$ only depending on $n$, and $x\approx y$ if $x\lesssim y$ and $y\lesssim x$. 

As we shall see later, if $\theta$ is ample, then $d_\theta(\f,\p)=0$ iff $\f-\p$ is constant, cf.~Corollary~\ref{cor:MAinj}. 

\begin{lem}\label{lem:dLconv} For all $\f,\f',\p\in\PL_\R\cap\PSH(\theta)$ we have 
\begin{equation}\label{equ:dtheta2}
d_\theta(\f,\p)\approx\|\f-\p\|^2_{\left(\theta,\frac{\f+\p}{2}\right)^{n-1}},
\end{equation}
\begin{equation}\label{equ:dJ}
d_\theta(\f,\p)\approx(\theta,\p)^{n+1}-(\theta,\f)^{n+1}+(n+1)(\f-\p)\cdot(\theta,\p)^n,
\end{equation}
and
\begin{equation}\label{equ:dconv}
d_\theta\left((1-t)\f+t\f',\p\right)\lesssim (1-t)d_\theta(\f,\p)+t\,d_\theta(\f',\p) 
\end{equation}
for all $t\in[0,1]$. 
\end{lem}

\begin{proof} Expanding out $\|\f-\p\|^2_{\left(\theta,\frac{\f+\p}{2}\right)^{n-1}}=-2^{1-n}(0,\f-\p)^2\cdot\left((\theta,\f)+(\theta,\p)\right)^{n-1}$ directly proves~\eqref{equ:dtheta2}. On the other hand, an elementary computation yields
\begin{multline}\label{equ:Jquad}
(\theta,\p)^{n+1}-(\theta,\f)^{n+1}+(n+1)(0,\f-\p)\cdot(\theta,\p)^n\\
=-\sum_{j=0}^{n-1}(j+1)(0,\f-\p)^2\cdot(\theta,\f)^j\cdot(\theta,\p)^{n-1-j}\\
=\sum_{j=0}^{n-1}(j+1)\|\f-\p\|^2_{(\theta,\f)^j\cdot(\theta,\p)^{n-1-j}},
\end{multline}
which implies~\eqref{equ:dJ}. By Theorem~\ref{thm:enconc}, the right-hand side of~\eqref{equ:dJ} is a convex function of $\f$, and~\eqref{equ:dconv} follows.  
\end{proof}

The key estimate is as follows. 

\begin{lem}\label{lem:CS} For all $\f,\f',\p\in\PL_\R\cap\PSH(\theta)$ we have
$$
\|\f-\f'\|^2_{(\theta,\p)^{n-1}}\lesssim\,d_\theta(\f,\f')^{\a_n}\max\{d_\theta(\f,\p),d_\theta(\f',\p)\}^{1-\a_n}
$$
with $\a_n:=2^{1-n}\in(0,1]$. 
\end{lem}
\begin{proof} Set $\tau:=\tfrac12(\f+\f')$, $f:=\f-\f'$,
$A:=d_\theta(\f,\f')$, $B:=\max\{d_\theta(\f,\p),d_\theta(\f',\p)\}$,
and, for $j=0,\dots,n-1$,
$$
b_j:=\|f\|^2_{(\theta,\p)^j\cdot(\theta,\tau)^{n-1-j}}=-(0,f)^2\cdot(\theta,\p)^j\cdot(\theta,\tau)^{n-1-j}.
$$
By~\eqref{equ:dtheta2}, we have $b_0\approx A$, and our goal is to show that
$$
b_{n-1}\lesssim A^{\frac 1{2^{n-1}}} B^{1-\frac 1{2^{n-1}}}.
$$ 
Since $f=(\f-\p)+(\p-\f')$, the triangle inequality for the seminorm $\|\cdot\|_{(\theta,\p)^{n-1}}$ yields 

  \begin{equation}\label{e213}
    b_{n-1}
    \le\left(\|\f-\p\|_{(\theta,\p)^{n-1}}+\|\f'-\p\|_{(\theta,\p)^{n-1}}\right)^2
      \le 4B.
  \end{equation}
If $A\ge B$, then $b_{n-1}\le 4 B\le 4 A^{\frac1{2^{n-1}}}B^{1-\frac1{2^{n-1}}}$, providing the desired estimate. We thus henceforth assume $A\le B$, and prove by induction on $j=0,\dots,n-1$ that 
  \begin{equation}\label{e215}
    b_j\lesssim  A^{\frac1{2^j}}B^{1-\frac1{2^j}}. 
  \end{equation}
For $j=0$ this holds since $b_0\approx A$. Now fix $0\le j\le n-2$, and note that
  \begin{align}\label{e216}
    b_{j+1}-b_j
    =& - (0,f)^2\cdot(\p-\tau)\cdot(\theta,\p)^j\cdot(\theta,\tau)^{n-2-j}\notag\\
    =&- (0,f)\cdot(0,\p-\tau)\cdot(\theta,\f)\cdot(\theta,\p)^j\cdot(\theta,\tau)^{n-2-j}\notag\\
     &+ (0,f)\cdot(0,\p-\tau)\cdot(\theta,\f')\cdot(\theta,\p)^j\cdot(\theta,\tau)^{n-2-j}.
    \end{align}
  Here we can use Cauchy--Schwarz to bound the the last two terms. For example,
$$
\left| (0,f)\cdot(0,\p-\tau)\cdot(\theta,\f)\cdot(\theta,\p)^j\cdot(\theta,\tau)^{n-2-j}\right|\le\|f\|_{(\theta,\f)\cdot(\theta,\p)^j\cdot(\theta,\tau)^{n-2-j}}\|\p-\tau\|_{(\theta,\f)\cdot(\theta,\p)^j\cdot(\theta,\tau)^{n-2-j}}.
$$
Using $(\theta,\f)\le 2(\theta,\tau)$ we can bound the first factor by $\sqrt{2b_j}$, and the second factor by
$$
\sqrt{2}\|\p-\tau\|_{(\theta,\p)^j\cdot(\theta,\tau)^{n-1-j}}\le\sqrt{2}\,d_\theta(\tau,\p).
 $$
 We have a similar bound for the last term in~\eqref{e216}. Adding the two bounds yields
 $b_{j+1}-b_j\le 4d_\theta(\tau,\p)\sqrt{b_j}$.
 By~\eqref{equ:dconv} $d_\theta(\tau,\p)\lesssim\sqrt{B}$, and we conclude 
\begin{equation}\label{e214}
    b_{j+1}-b_j\lesssim\sqrt{Bb_j}
  \end{equation} 
  for $0\le j\le n-2$. 
 Using the induction hypothesis~\eqref{e215}, we get
$$
b_{j+1}\lesssim b_j+\sqrt{Bb_j}\lesssim A^{\frac1{2^j}}B^{1-\frac1{2^j}}+A^{\frac1{2^{j+1}}}B^{1-\frac1{2^{j+1}}}
$$
$$
=A^{\frac1{2^{j+1}}}B^{1-\frac1{2^{j+1}}}\left(\left(\frac{A}{B}\right)^{\frac1{2^{j+1}}}+1\right).
$$
  By assumption, $A\le B$, so
  $b_{j+1}\lesssim A^{\frac1{2^{j+1}}}B^{1-\frac1{2^{j+1}}}$. 
  The proof is complete.
\end{proof}

\begin{proof}[Proof of Theorem~\ref{thm:triangle}]
Set $\tau:=(\f_1+\f_2)/2$. By~\eqref{equ:dtheta2}, we have 
  \begin{equation*}
    d_\theta(\f_1,\f_2)
    \lesssim\|\f_1-\f_2\|^2_{(\theta,\tau)^{n-1}}
    \lesssim\max_{i=1,2}\|\f_i-\f_3\|^2_{(\theta,\tau)^{n-1}}. 
  \end{equation*}
By  Lemma~\ref{lem:CS} we have, for $i=1,2$,
  \begin{equation*}
    \|\f_i-\f_3\|^2_{(\theta,\tau)^{n-1}}
    \lesssim d_\theta(\f_i,\f_3)^{\a_n}
    \max\{d_\theta(\f_i,\tau),d_\theta(\f_3,\tau)\}^{1-\a_n}.
  \end{equation*}

  By Lemma~\ref{lem:dLconv} $d_\theta(\f_i,\tau)\lesssim \,d_\theta(\f_1,\f_2)$ and $d_\theta(\f_3,\tau)\lesssim \max_{i=1,2}d_\theta(\f_i,\f_3)$. Altogether, this yields
  \begin{equation}\label{e401}
    d_\theta(\f_1,\f_2)
    \lesssim \max_{i=1,2}d_\theta(\f_i,\f_3)^{\a_n}
    \max\{d_\theta(\f_1,\f_2),\max_{i=1,2}d_\theta(\f_i,\f_3)\}^{1-\a_n}. 
  \end{equation}
 
When $d_\theta(\f_1,\f_2)\ge\max_{i=1,2}d_\theta(\f_i,\f_3)$,~\eqref{e401} yields
  $d_\theta(\f_1,\f_2)\lesssim \max_{i=1,2}d_\theta(\f_i,\f_3)$.
  But the same inequality trivially holds when 
  $d_\theta(\f_1,\f_2)\le\max_{i=1,2}d_\theta(\f_i,\f_3)$.
  \end{proof}

\begin{cor}\label{cor:CS} For all $\f_0,\dots,\f_n\in\PL_\R\cap\PSH(\theta)$ we have
$$
\|\f_0-\f_1\|^2_{(\theta,\f_2)\inter(\theta,\f_n)}\lesssim \,d_\theta(\f_0,\f_1)^{\a_n}\max_{0\le i\le n} d_\theta(\f_i)^{1-\a_n}
$$
with $\a_n=2^{1-n}$. 
\end{cor}

Recall that $d_\theta(\f)=d_\theta(\f,0)$. 

\begin{proof} Set $\tau:=\frac1{n-1}\sum_{i=2}^n\f_i\in\PL_\R\cap\PSH(\theta)$. Then
$$
\|\f_0-\f_1\|^2_{(\theta,\f_2)\inter(\theta,\f_n)}\lesssim \|\f_0-\f_1\|^2_{(\theta,\tau)^{n-1}}\lesssim \,d_\theta(\f_0,\f_1)^{\a_n}\max\{d_\theta(\f_0,\tau),d_\theta(\f_1,\tau)\}^{1-\a_n},
$$
by Lemma~\ref{lem:CS}. The result follows since 
\begin{align*}
  \max\{d_\theta(\f_0,\tau),d_\theta(\f_1,\tau)\}
  &\lesssim\max\left\{\max_i d_\theta(\f,\f_i),\max_i d_\theta(\f',\f_i)\right\}\\
  &\lesssim \max\left\{d_\theta(\f),d_\theta(\f'),\max_i d_\theta(\f_i)\right\}
\end{align*}
by quasi-convexity of $d_\theta$ (cf.~\eqref{equ:dconv}) and the quasi-triangle inequality.
\end{proof}
%
%
%
%
 \section{Plurisubharmonic functions}\label{sec:psh}
As before, $X$ is a projective variety, and $(X_\a)$ denotes its set of irreducible components.
 
In this section, we introduce the class $\PSH(\theta)=\PSH(X,\theta)$ of $\theta$-psh functions associated to an arbitrary numerical class $\theta\in\Num(X)$. We also study pluripolar sets, \ie loci where $\theta$-psh functions take the value $-\infty$.

%
%
\subsection{The class of $\theta$-psh functions}
Recall that:
\begin{itemize}
\item a function $\f\in X^\an\to\R\cup\{-\infty\}$ is called \emph{generically finite} if $\f|_{X_\a^\an}\not\equiv-\infty$ for all $\a$ (see Definition~\ref{defi:genfin});
\item for any $L\in\Pic(X)_\Q$, $\cH^\gf_\R(L)$ denotes the set of generically finite (continuous) functions $\f\colon X^\an\to\R\cup\{-\infty\}$ of the form 
$$
\f=m^{-1}\max_i\{\log|s_i|+\la_i\}
$$
where $(s_i)$ is a finite set of sections of $mL$ with $m$ sufficiently divisible, $\la_i\in\R$ (see Definition~\ref{defi:FS}). This set is non-empty iff $L$ is effective.  
\end{itemize}
Given any numerical class $\theta\in\Num(X)$, we are seeking to define a class $\PSH(\theta)$ of \emph{$\theta$-plurisubharmonic} functions (\emph{$\theta$-psh} functions for short) $\f\colon X^\an\to\R\cup\{-\infty\}$. By analogy with the complex analytic setting, these functions should be usc and generically finite, and the following properties should hold: 
\begin{itemize}
\item[(PSH1)] for any $L\in\Pic(X)_\Q$, we have $\cH^\gf_\R(L)\subset\PSH(L):=\PSH(c_1(L))$;
\item[(PSH2)] if a generically finite function $\f\colon X^\an\to\R\cup\{-\infty\}$ arises as the pointwise limit of a decreasing net $(\f_i)$ with $\f_i\in\PSH(\theta_i)$ and $\lim_i\theta_i=\theta$, then $\f\in\PSH(\theta)$. 
\end{itemize}
Note that (PSH2) implies the following two properties:
\begin{itemize}
\item[(PSH2a)] for any function $\f\colon X^\an\to\R\cup\{-\infty\}$, the set of $\theta\in\Num(X)$ such that $\f\in\PSH(\theta)$ is closed (possibly empty); 
\item[(PSH2b)] for any $\theta\in\Num(X)$, the set $\PSH(\theta)$ is closed under decreasing limits (provided the limit function is generically finite).
\end{itemize}
Taking (PSH1) and (PSH2) as minimal requirements directly leads to the following definition: 

\begin{defi}\label{defi:psh} For any $\theta\in\Num(X)$, we define a \emph{$\theta$-psh function} $\f\colon X^\an\to\R\cup\{-\infty\}$ as a generically finite, usc function that can be written as the pointwise limit of a decreasing net $\f_i\in\cH^\gf_\R(L_i)$ with $L_i\in\Pic(X)_\Q$ such that $\lim_i c_1(L_i)=\theta$ in $\Num(X)$. 
\end{defi}
We denote by 
$\PSH(\theta)$ the set of $\theta$-psh functions, and by 
$$
\CPSH(\theta):=\Cz(X)\cap\PSH(\theta)\subset\cE^\infty(\theta)\subset\PSH(\theta)
$$
the subsets of continuous and bounded $\theta$-psh functions, respectively. When needed, we also write $\PSH(X,\theta)=\PSH(\theta)$ etc. Using the notation in the definition, we have:
\begin{itemize}
\item since $\cH^\gf_\R(L_i)$ is non-empty, $L_i$ is necessarily effective, by~\eqref{equ:cHvalnempty}, and hence 
\begin{equation}\label{equ:pshpsef}
\PSH(\theta)\ne\emptyset\Longrightarrow\theta\in\Psef(X);
\end{equation}
\item if $\f$ is bounded, then $\f_i$ is finite valued, \ie $\f_i\in\cH(L_i)$; this implies that $L_i$ is semiample, by~\eqref{equ:cHnempty}, and hence 
\begin{equation}\label{equ:pshnef}
\cE^\infty(\theta)\ne\emptyset\Longrightarrow\theta\in\Nef(X);
\end{equation}
\item if $\f\in\CPSH(\theta)$ is further continuous, then $\f_i\to\f$ uniformly on $X^\an$, by Dini's lemma. 
\end{itemize}

\begin{rmk}\label{rmk:DPS} In the complex analytic case, psh functions and metrics are defined in a different way, which is local in nature. However, when $X$ is smooth, well-known arguments due to Demailly~\cite{Dem92} and relying on the Ohsawa--Takegoshi extension theorem show that any psh metric on a (necessarily pseudoeffective) line bundle $L$ can be approximated from above by metrics attached to sections of $mL+A$, where $A$ is a fixed ample line bundle and $m$ is large (compare for instance~\cite[Theorem~5.4]{Bou} or~\cite[Appendix]{GZ1}). As a result, any such metric satisfies the analogue of Definition~\ref{defi:psh}. 
\end{rmk} 

\begin{lem}\label{lem:consistent} For functions in $\PL(X)_\R$, Definition~\ref{defi:psh} is consistent with Definition~\ref{defi:PLpsh}. In other words, given $\theta\in\Num(X)$, a test configuration $\cX$ that dominates $\cX_\triv$ and $D\in\VCar(\cX)_\R$, we have
$$
\f_D\in\PSH(\theta)\Longleftrightarrow\theta_\cX+D\in\Nef(\cX/\A^1).
$$
\end{lem}
Indeed, this follows from Lemma~\ref{lem:limitpsh}, Lemma~\ref{lem:pshPL} and Corollary~\ref{cor:goodman}. For any $\theta\in\Num(X)$, we get, in particular,
\begin{equation}\label{equ:cstpsh}
\PSH(\theta)\ \text{contains all constant functions}\Longleftrightarrow\theta\in\Nef(X),
\end{equation}
which provides a converse to~\eqref{equ:pshnef}. 

\begin{rmk}\label{rmk:DPS2} If $L$ is a $\Q$-line bundle, then functions in $\PSH(L)=\PSH(c_1(L))$ can be interpreted as psh metrics on $L$, using the trivial metric. By~\eqref{equ:cstpsh}, any nef line bundle $L$ admits a bounded (continuous) psh metric, to wit the trivial metric. This is in stark contrast with the complex analytic case, where a nef and big line bundle may not admit any bounded psh metric (see~\cite[Example 5.4]{BEGZ}, based on~\cite[Example 1.7]{DPS}). 
\end{rmk}

\begin{thm}\label{thm:psh12} Properties (PSH1) and (PSH2) above are satisfied. 
\end{thm} 
As a consequence, Definition~\ref{defi:psh} is indeed the minimal one for which (PSH1) and (PSH2) are satisfied. 

\begin{lem}\label{lem:PSH} Let $K$ be a compact topological space, and assume we are given: 
\begin{itemize}
\item a set $S$ and, for each $s\in S$, a family $\cF_s$ of continuous functions $f\colon K\to\R\cup\{-\infty\}$, stable under addition of a constant;
\item a map $\rho\colon S\to T$ to a metrizable topological space $T$.
\end{itemize}
For each $t\in T$, denote by $\tcF_t$ the set of all (usc) functions $g\colon K\to\R\cup\{-\infty\}$ that can be written as the pointwise limit of a decreasing net $(f_i)$ with $f_i\in\cF_{s_i}$ for a net $(s_i)$ of $S$ such that $\lim_i \rho(s_i)=t$ in $T$. Then $\dbtilde{\cF}_t=\widetilde{\cF}_t$
for any $t\in T$. In other words, if a function $g\colon K\to\R\cup\{-\infty\}$ can be written as the pointwise limit of  a decreasing net $(g_i)_{i\in I}$ with $g_i\in\tcF_{t_i}$ for a net $(t_i)$ of $T$ such that $\lim_i t_i=t$, then $g\in\tcF_t$. 
\end{lem}
\begin{proof} Pick a metric $d$ on $T$ inducing the given topology. For each $i\in I$, $g_i$ is the pointwise limit of a decreasing net $(f_{i,j})_{j\in J_i}$ with $f_{i,j}\in\cF_{s_{i,j}}$ for a net $(s_{i,j})_{j\in J_i}$ in $S$ such that $t_{i,j}:=\rho(s_{i,j})$ satisfies $\lim_j t_{i,j}=t_i$. Denote by $A$ the set of triples $\a=(i,j,\e)$ with $i\in I$, $j\in J_i$ and $\e\in\R_{>0}$, and define a partial preorder on $A$ by setting 
\begin{equation}\label{equ:directed}
(i,j,\e)\ge(i',j',\e')\Longleftrightarrow\left\{
\begin{array}{ll}
i\ge i'\\
f_{i,j}+\e\le f_{i',j'}+\e'\ \text{on }K
\\
d(t_{i,j},t_i)+\e\le d(t_{i',j'},t_{i'})+\e',
\end{array}
\right. 
\end{equation}
We claim that $A$ is directed. To see this, pick $\a_1=(i_1,j_1,\e_1)$, $\a_2=(i_2,j_2,\e_2)$ in $A$. Since $I$ is directed, we can choose $i_3\in I$ with $i_3\ge i_1,i_2$. Pick also $\e_3>0$ with $\e_3<\min\{\e_1,\e_2\}$. Since $(g_i)_i$ is decreasing, we have $g_{i_3}\le\min\{g_{i_1},g_{i_2}\}\le\min\{f_{i_1,j_1},f_{i_2,j_2}\}$, and hence 
\begin{equation}\label{equ:i3i1}
g_{i_3}+\e_3<\min\{f_{i_1,j_1}+\e_1,f_{i_2,j_2}+\e_2\}
\end{equation}
pointwise on $K$. For each $j\in J_{i_3}$, define $h_j:K\to\R$ by 
$$
h_j:=\exp\left(f_{i_3,j}+\e_3\right)-\exp\left(\min\{f_{i_1,j_1}+\e_1,f_{i_2,j_2}+\e_2\}\right). 
$$
Then $(h_j)$ is a decreasing net of continuous functions such that $\lim_j h_j<0$ pointwise on $K$, by~\eqref{equ:i3i1}. The compact space $K$ can thus be written as the increasing union of the open sets $\{h_j<0\}$, and hence $K=\{h_{j_3}<0\}$ for some $j_3\in J_{i_3}$, \ie 
$$
f_{i_3,j_3}+\e_3<\min\{f_{i_1,j_1}+\e_1,f_{i_2,j_2}+\e_2\}\text{ on }K. 
$$
Since $\lim_j t_{i_3,j}=t_{i_3}$ and $\e_3<\min\{\e_1,\e_2\}$, we can further arrange, after increasing $j_3$ if needed, that 
$$
d(t_{i_3,j_3},t_{i_3})+\e_3<\min\{\e_1,\e_2\}\le\min\left\{d(t_{i_1,j_1},t_{i_1})+\e_1,d(t_{i_2,j_2},t_{i_2})+\e_2\right\}. 
$$
Thus $\a_3:=(i_3,j_3,\e_3)\in A$ satisfies $\a_3\ge \a_1,\a_2$, which proves, as desired, that $A$ is directed. 

Now, for each $\a=(i,j,\e)\in A$, set
$$
s_\a:=s_{i,j}\in S,\quad t_\a=\rho(s_\a)\in T,\quad\text{and}\quad g_\a:=f_{i,j}+\e\in\cF_{s_\a},
$$
where the last containment holds because $\cF_{s_\a}$ is assumed to be stable under addition of a constant. 
By construction, $(g_\a)_{\a\in A}$ is a decreasing net, and it will suffice to show $\lim_\a g_\a=g$ pointwise on $K$ and $\lim_\a t_\a=t$. Pick $\d>0$ and $x\in K$. Since $\lim_i g_i(x)=g(x)$ and $\lim_i t_i=t$, we can find $i_0\in I$ such that $g_{i_0}(x)<g(x)+\d$ and $d(t_i,t)<\d$ for all $i\ge i_0$. Using $\lim_j f_{i_0,j}(x)=g_{i_0}(x)$ and $\lim_j t_{i_0,j}=t_{i_0}$, we next pick $j_0\in J_{i_0}$ and $0<\e_0\ll 1$ such that $f_{i_0,j_0}(x)+\e_0<g(x)+\d$ and $d(t_{i_0,j_0},t_{i_0})+\e_0<\d$. Set $\a_0:=(i_0,j_0,\e_0)\in A$. For each $\a=(i,j,\e)\in A$ such that $\a\ge\a_0$, we then have, by~\eqref{equ:directed},
$$
g_\a(x)=f_{i,j}(x)+\e\le f_{i_0,j_0}(x)+\e_0<g(x)+\d,
$$
$$
d(t_\a,t_i)=d(t_{i,j},t_i)\le d(t_{i,j},t_i)+\e\le d(t_{i_0,j_0},t_{i_0})+\e_0<\d, 
$$
as well as $i\ge i_0$, which implies $d(t_i,t)<\d$ and hence $d(t_\a,t)<2\d$. This shows, as desired, $\lim_\a g_\a(x)=g(x)$ and $\lim_\a t_\a=t$, which concludes the proof. 
\end{proof}

\begin{proof}[Proof of Theorem~\ref{thm:psh12}] That (PSH1) holds is clear from Definition~\ref{defi:psh}. As for (PSH2), it corresponds to the special case of Lemma~\ref{lem:PSH} in which $K:=X^\an$, $\cF_L:=\cH^\gf_\R(L)$ for $L\in S:=\Pic(X)_\Q$, and $\rho\colon S\to T:=\Num(X)$ maps $L\in\Pic(X)_\Q$ to $c_1(L)$. Indeed, for any $\theta\in\Num(X)$, $\PSH(\theta)$ is precisely  the set of functions in $\tcF_\theta$ that are generically finite. 
\end{proof}

%
%
\subsection{Basic properties} 
The next few results provide some basic but important properties of $\theta$-psh functions. 

\begin{thm}\label{thm:PSH} For all $\theta,\theta'\in\Num(X)$ and $t\in\R_{>0}$, we have: 
\begin{itemize}
\item[(i)]
$\PSH(\theta)+\PSH(\theta')\subset\PSH(\theta+\theta')$ and $\PSH(t\theta)=t\PSH(\theta)$;
\item[(ii)] the set $\PSH(\theta)$ is convex, and stable under uniform limits, finite maxima, the additive action of $\R$, and the scaling action of $\R_{>0}$; 
\item[(iii)] for any morphism $f\colon Y\to X$ and $\f\in\PSH(X,\theta)$, either $f^\star\f\in\PSH(Y,f^\star\theta)$ or $f^\star\f\equiv-\infty$ on some irreducible component of $Y^\an$;
\item[(iv)]
  for any $\f\in\PSH(\theta)$ and any $\a$, we have $\f|_{X_\a^\an}\in\PSH(X_\a,\theta|_{X_\a})$.
\end{itemize}
\end{thm}
\begin{rmk}\label{lem:pshdisconnect}
  If $X$ is disconnected, with connected components $X^\b$, then it is easy to see that $\PSH(X,\theta)=\prod_\beta\PSH(X^\b,\theta|_{X^\b})$.
\end{rmk}
\begin{lem}\label{lem:twopsh} For any $\theta\in\Num(X)$ and $\f,\f'\in\PSH(\theta)$, we can find a net $(L_i)$ of $\Q$-line bundles and $\f_i,\f'_i\in\cH^\gf_\R(L_i)$ such that $\lim_i c_1(L_i)\to\theta$, $\f_i\searrow\f$, $\f'_i\searrow\f'$. 
\end{lem} 
\begin{proof} By definition, we can find decreasing nets $(\f_i)_{i\in I}$, $(\f'_j)_{j\in I'}$ and $\Q$-line bundles $M_i$, $M'_j$ such that $\f_i\in\cH^\gf_\R(M_i)$, $\f'_j\in\cH^\gf_\R(M'_j)$,  $\lim_i c_1(M_i)=\lim_j c_1(M'_j)=\theta$ and $\lim_i\f_i=\f$, $\lim_j\f'_j=\f'$. Replacing $I$, $I'$ with the directed set $I\times I'$  (equipped with the product preorder), we may assume $I=I'$. Fix an ample line bundle $A$. Since $\lim_i(c_1(M_i)-c_1(M'_i))=0$, we can find $\e_i\in\Q_{>0}$ such that $M_i-M'_i+\e_i A$ is ample and $\lim_i\e_i=0$. If $L_i:=M_i+\e_i A$, then $c_1(L_i)\to\theta$, and since $L_i-M_i$ and $L_i-M'_i$ are both ample, Proposition~\ref{prop:FS} yields $\f_i\in\cH^\gf_\R(M_i)\subset\cH^\gf_\R(L_i)$, $\f'_i\in\cH^\gf_\R(L_i)$, and we are done.  
\end{proof}

\begin{proof}[Proof of Theorem~\ref{thm:PSH}] We will repeatedly use Proposition~\ref{prop:FS}. Pick $\f\in\PSH(\theta)$, $\f'\in\PSH(\theta')$. As $\f$ and $\f'$ are finite on $v_{\triv,\a}$ for all $\a$, so is the function $\f+\f'$, which is therefore generically finite. Further, $\f$ and $\f'$ can be written as the limits of decreasing nets $(\f_i)_{i\in I}$ and $(\f'_j)_{j\in I'}$ with $\f_i\in\cH^\gf_\R(L_i)$, $\f'_j\in\cH^\gf_\R(L'_j)$ for $L_i,L'_j\in\Pic(X)_\Q$ such that $\lim_i c_1(L_i)=\theta$ and $\lim_j c_1(L'_j)=\theta'$, respectively. The net $(\f_i+\f'_j)_{(i,j)\in I\times I'}$ is decreasing, and converges pointwise to $\f+\f'$. Now  $\f_i+\f'_j\in\cH^\gf_\R(L_i+L'_j)$, and hence $\f+\f'\in\PSH(\theta+\theta')$, since $\lim_{(i,j)}(c_1(L_i)+c_1(L'_j))=\theta+\theta'$. This proves the first inclusion in (i).  

For any $c\in\R$ and $t\in\R_{>0}$, $\f_i+c$ and $t\cdot\f_i$ belong to $\cH^\gf_\R(L_i)$, and decrease, respectively, to $\f+c$ and $t\cdot\f$, which therefore belong to $\PSH(\theta)$. 

We next show that $t\f\in\PSH(t\theta)$. If $t$ is rational, this is clear, since each $tL_i$ is a $\Q$-line bundle, $t\f_i\in\cH^\gf_\R(tL_i)$ and $t\f_i\searrow t\f$. In the general case, we may assume, after replacing $\f$ with $\f-\sup\f$, that $\f\le 0$. Writing $t$ as the limit of an increasing sequence $t_m\in\Q_{>0}$, we then have $t_m\f\in\PSH(t_m\theta)$ and $t_m\f\searrow t\f$, $t_m\theta\to t\theta$, and hence $t\f\in\PSH(t\theta)$, by (PSH2). This concludes the proof of (i), which implies that $\PSH(\theta)$ is convex. 

To prove that $\PSH(\theta)$ is stable under maxima, pick $\f,\f'\in\PSH(\theta)$. By Lemma~\ref{lem:twopsh}, we can find a net $(L_i)$ of $\Q$-line bundles and $\f_i,\f'_i\in\cH^\gf_\R(L_i)$ such that $\lim_i c_1(L_i)\to\theta$, $\f_i\searrow\f$, $\f'_i\searrow\f'$. Then $\max\{\f_i,\f'_i\}\in\cH^\gf_\R(L_i)$ decreases to $\max\{\f,\f'\}$, which therefore belongs to $\PSH(\theta)$.

If a net in $\PSH(\theta)$ converges uniformly to a function $\p\colon X^\an\to\R\cup\{-\infty\}$, then $\p$ is generically finite, and is also the uniform limit of a sequence $(\p_m)$ in $\PSH(\theta)$. After passing to a subsequence, it is then easy to find constants $c_m\to 0$ such that $(\p_m+c_m)$ is decreasing, and we infer $\p\in\PSH(\theta)$, by (PSH2). This concludes the proof of (ii). 

Next we prove~(iii). Write $\f\in\PSH(\theta)$ as the decreasing limit of a net $\f_i\in\cH^\gf_\R(L_i)$ such that $\lim_i c_1(L_i)\to\theta$. Then each $\f_i\in\cH^\gf_\R(L_i)$ satisfies either $f^\star\f_i\in\cH^\gf_\R(f^\star L_i)$ or $f^\star\f_i\equiv-\infty$ on some irreducible component of $Y^\an$. If $f^\star\f$ is generically finite, then the latter case never occurs, since $f^\star\f\le f^\star\f_i$. Since $f^\star\f_i\searrow f^\star\f$ and $c_1(f^\star L_i)\to f^\star\theta$, (PSH2) yields $f^\star\f\in\PSH(f^\star\theta)$.

Finally,~(iv) follows from~(iii) since $\f$ is generically finite.
\end{proof}

\begin{cor}\label{cor:pshcvx} Pick $\theta\in\Num(X)$, $\f,\p\in\PSH(\theta)$ with $\f\le\p$, and let $\chi\colon(-\infty,0]\to\R$ be a continuous convex function such that $0\le\chi'\le1$. Then $\chi\circ(\f-\p)+\p\in\PSH(\theta)$.
\end{cor}

\begin{proof} Note that $\chi\circ(\f-\p)+\p$ is generically finite. Now, $\chi$ can be written as a decreasing limit of
  functions $\chi_m$ that are finite maxima of affine functions of the form $\ell(t)=a t+b$, where $a\in[0,1]$ and $b\in\R$. For each $m$, $\chi_m\circ(\f-\p)+\p$ is a finite maximum of functions of the form $\ell\circ(\f-\p)+\p=(1-a)\p+a\f+b$, and hence $\chi_m\circ(\f-\p)+\p\in\PSH(\theta)$, by Theorem~\ref{thm:PSH}. Now $\chi_m\circ(\f-\p)+\p\searrow\chi\circ(\f-\p)+\p$, and the result follows, by (PSH2). 
 \end{proof}

\begin{cor}\label{cor:exppsh} If $\theta\in\Nef(X)$ and $0\ge\f\in\PSH(\theta)$, then $\exp(\f)\in\cE^\infty(\theta)$. 
\end{cor}
Indeed, we can pick $\p=0$ (since $\theta$ is nef, see~\eqref{equ:pshnef}).
This corollary will allow us to reduce certain statements about general $\theta$-psh functions to the case of bounded ones.

\begin{prop}\label{prop:pshconvex} For any $\theta\in\Num(X)$ and $\f\in\PSH(\theta)$, we have: 
\begin{itemize}
\item[(i)] $\f$ is decreasing and usc on $X^\an$; 
\item[(ii)] for each irreducible subvariety $Y\subset X$, the `maximum principle' $$\sup_{Y^\an}\f=\f(v_{Y,\triv})$$ holds; in particular, 
\begin{equation}\label{equ:suppsh} 
\sup\f=\max_\a\f(v_{\triv,\a}).
\end{equation}
\item[(iii)] for any $v\in X^\an$, $t\mapsto\f(tv)$ is convex and decreasing on $[0,+\infty)$, and $\f(tv)\searrow\f(v_{Z,\triv})=\sup_{Z^\an}\f$ as $t\to+\infty$, with $Z$ the closure of the center of $v$. 
\end{itemize}
\end{prop}
\begin{proof} By Proposition~\ref{prop:FS}~(i), each $\f_i\in\cH^\gf_\R(L_i)$ is decreasing, and $t\mapsto\f_i(tv)$ is convex, by Lemma~\ref{lem:FSconvex}. Since $(\f_i)$ converges pointwise to $\f$, the latter inherits these properties. The rest follows from Lemma~\ref{lem:maxprin}. 
\end{proof}


\begin{exam}\label{exam:curvepsh} Let $X$ be a smooth irreducible curve, and recall the description of $X^\an$ given in \S\ref{sec:curve}. For any $\theta\in\Num(X)$, a usc function $\f\colon X^\an\to\R\cup\{-\infty\}$ is then $\theta$-psh iff 
\begin{itemize} 
\item[(i)] for each $p\in X(k)$, $\f_p(t):=\f(t\ord_p)$ is convex on $[0,+\infty)$ (and hence decreasing, being bounded above);
\item[(ii)] $\deg\theta+\sum_{p\in X(k)}\f_p'(0_+)\ge 0$ (which implies $\f_p'(0_+)=0$, and hence $t\mapsto\f(t\ord_p)$ constant, for all but countably many $p\in X(k)$).
\end{itemize}
Indeed, when $\f$ is PL, this follows from Example~\ref{exam:curvePLpsh}, and the general case is obtained by writing a function satisfying (i) and (ii) as a decreasing limit of PL functions with the same properties (compare~\cite[\S7]{valtree} and~\cite[\S2.5]{dynberko}). 
\end{exam}

This characterization implies $\PSH(0)=\R$ when $X$ is a smooth irreducible curve. As in the complex analytic case, this holds in fact on any connected variety, see Corollary~\ref{cor:pshcst} below. 

\begin{rmk}\label{rmk:qpsh} Any $\f\in\PL(X)_\R$ is $\eta$-psh with respect to some closed $(1,1)$-form $\eta$, in the sense of Remark~\ref{rmk:closedform}. However, $\f$ is in general not $\theta$-psh for any $\theta\in\Num(X)$, since such functions must be decreasing in the partial ordering on $\Xan$ (see Proposition~\ref{prop:pshconvex} ). As a result, several potential definitions of \emph{quasi-psh functions} coexist in our context. 
\end{rmk}

%
%
\subsection{The regularization theorem} 
Many properties of $\theta$-psh functions are obtained by regularization, that is, approximation by nicer $\theta$-psh functions. In our present approach, regularization of $\theta$-psh functions is built into their definition. One can, however, improve the properties of the approximants as follows: 

\begin{thm}\label{thm:pshample} Pick $\theta\in\Num(X)$ and $\f\in\PSH(\theta)$. Then:
\begin{itemize}
\item[(i)] $\f$ can be written as the limit of a decreasing net $\f_i\in\cH^\gf_\Q(L_i)$ with $L_i\in\Pic(X)_\Q$ such that $c_1(L_i)-\theta$ is ample and $\lim_i c_1(L_i)=\theta$; 
\item[(ii)] if $\theta\in\Nef(X)$ then (i) holds with $\f_i\in\cH(L_i)=\cH_\Q(L_i)$ and $L_i$ ample; 
\item[(iii)] if $\om:=\theta\in\Amp(X)$, then $\f$ can be written as the limit of a decreasing net in $\cH^\dom(\om)$ (see Definition~\ref{defi:Hdom}). 
\end{itemize}
\end{thm}
We shall later prove that these results actually hold with (countable) sequences instead of nets, see Corollary~\ref{cor:countreg}.

As already indicated by~(iii), the class of $\theta$-psh functions has particularly good properties when $\theta\in\Amp(X)$ is ample. In this case, we typically write $\om$ instead of $\theta$.
\begin{cor} For any $\om\in\Amp(X)$, the space $\CPSH(\om)=\Cz(X)\cap\PSH(\om)$ is the closure of $\cH^\dom(\om)\subset\Cz(X)$ in the topology of uniform convergence.
\end{cor}
\begin{proof} By Theorem~\ref{thm:PSH}, $\CPSH(\om)$ is closed in the topology of uniform convergence. By Theorem~\ref{thm:pshample} and Dini's lemma, any $\f\in\CPSH(\om)$ is the uniform limit of a net in $\cH^\dom(\om)$, and the result follows. 
\end{proof}
\begin{cor}\label{cor:pshonXlin}
  For any $\theta\in\Num(X)$ and any $\f\in\PSH(\theta)$ we have $\f>-\infty$ on $X^\lin$.
\end{cor}
\begin{proof}
  Pick an ample line bundle $L$ such that $c_1(L)-\theta$ is nef. Then $\f\in\PSH(L)$. Using Theorem~\ref{thm:PSH} and $X^\lin=\coprod_\a X_\a^\lin$, we may assume $X$ is irreducible. Now pick $v\in X^\lin$. Then $\teL(v)<\infty$, and we claim that $\f(v)\ge\f(v_\triv)-\teL(v)$. Indeed, $\teL(v)$ is the supremum of $m^{-1}v(s)$ over $m\in\Z_{>0}$ and $s\in \Hnot(X,mL)\setminus\{0\}$. Equivalently, $\p:=m^{-1}\log|s|$ satisfies $\p(v)\ge\p(v_\triv)-\teL(v)$. Adding constants and taking finite maxima shows that this also holds for all $\p\in \cH(L)$, and then for all $\f\in\PSH(L)$, by Theorem~\ref{thm:pshample}.
\end{proof}
As another consequence, we get: 
 
\begin{cor}\label{cor:charLpsh} Assume that $X$ is irreducible, and let $L$ be an ample $\Q$-line bundle. Then $\PSH(L)$ is the smallest class of functions $\f\colon X^\an\to\R\cup\{-\infty\}$, not identically $-\infty$, that contains all functions of the form $m^{-1}\log|s|$ with $m\in\Z_{>0}$ and $s\in \Hnot(X,mL)\setminus\{0\}$, and is stable under addition of a real constant, finite maxima, and decreasing limits. 
\end{cor}

\begin{proof}[Proof of Theorem~\ref{thm:pshample}] By definition, $\f$ is the limit of a decreasing net $\f_i\in\cH^\gf_\R(L_i)$ with $c_1(L_i)\to\theta$. Using~\eqref{equ:FS} it is straightforward to see that each $\f_i$ can further be written as a decreasing limit of functions in $\cH^\gf_\Q(L_i)$. By Lemma~\ref{lem:PSH}, we may thus assume that $\f_i\in\cH^\gf_\Q(L_i)$. We can also arrange that $c_1(L_i)-\theta\in\Amp(X)$. Indeed, given any ample line bundle $A$, we can find $\e_i\in\Q_{>0}$ such that $\lim_i\e_i=0$ and $L_i-\theta+\e_i A$ is ample. By Proposition~\ref{prop:FS}, we have $\f_i\in\cH^\gf_\Q(L_i)\subset\cH^\gf_\Q(L_i+\e_i A)$, and replacing $L_i$ with $L_i+\e_i A$ yields the claim. This proves (i). 

Assume next that $\theta$ is nef. In the notation of (i), it follows that each $L_i$ is ample. Further, each $\f_i$ is the decreasing limit of $\max\{\f_i,-m\}\in\cH(L_i)$, and Lemma~\ref{lem:PSH} yields (ii). 

Finally assume that $\om:=\theta$ is ample. Replacing $\f$ with $\f-\sup\f$, we may assume $\f\le 0$. By Proposition~\ref{prop:PLdom}, any function in $\PL\cap\PSH(\om)$ is a decreasing limit of functions in $\cH^\dom(\om)$. Relying again on Lemma~\ref{lem:PSH}, it will thus be enough to show that $\f$ is the decreasing limit of functions in $\PL\cap\PSH(\om)$. Use the notation of (ii), and pick $t\in\Q_{>1}$. For all $i$ large enough we then have $t\om-c_1(L_i)\in\Amp(X)$ for $i$ large enough, and hence $\cH(L_i)\subset\PL\cap\PSH(t\om)$. It follows that $t^{-1}\f$ is the decreasing limit of a net in $\PL\cap\PSH(\om)$. Since $\f\le 0$, we have $t^{-1}\f\searrow\f$ as $t\searrow 1$, and another application of Lemma~\ref{lem:PSH} shows that $\f$ is also the decreasing limit of functions in $\PL\cap\PSH(\om)$. This proves (iii). 
\end{proof}

\begin{rmk}\label{rmk:PSHsiminag} Assume that $\charac k=0$ and that $X$ is smooth and connected. Consider the field of Laurent series $K=k\lau{t}$, and denote by $X_K^\an$ the Berkovich analytification of the base change $X_K$; this comes with a natural projection $\pi\colon X_K^\an\to X^\an$. The base change $\om_K$ of any $\om\in\Amp(X)$ defines a closed $(1,1)$-form on $X_K^\an$ with ample de Rham class, in the sense of~\cite{siminag}. Using Theorem~\ref{thm:pshample}~(iii) and Theorem~\ref{thm:pshdual} in the Appendix, one can show that a function $\f\colon X^\an\to\R\cup\{-\infty\}$ is $\om$-psh in the present sense iff $\pi^\star\f$ is $\om_K$-psh on $X_K^\an$ in the sense of~\cite{siminag}. 
\end{rmk}
%
%
\subsection{The topology of $\PSH(\theta)$}\label{sec:topopsh}
Let us fix $\theta\in\Num(X)$ for the moment. By  Corollary~\ref{cor:pshonXlin} , any $\f\in\PSH(\theta)$ is finite-valued on $X^\div=\coprod_\a X_\a^\div$.  

\begin{defi} We endow $\PSH(\theta)$ with the topology of pointwise convergence on $X^\div$.
\end{defi}
As we shall later see in Theorem~\ref{thm:weaklin}, this topology coincides with the topology of pointwise convergence on $X^\lin=\coprod_\a X_\a^\lin$.  In view of~\eqref{equ:suppsh}, we have:

\begin{exam}\label{exam:supcont} The map  $\f\mapsto\sup\f=\max_\a\f(v_{\triv,\a})$ is continuous on $\PSH(\theta)$.
\end{exam}

\begin{thm}\label{thm:suppsh} For any $\f\in\PSH(\theta)$ and any usc function $\p\colon X^\an\to\R\cup\{-\infty\}$, we have  
$$
\f\le\p\text{ on }X^\div\Longleftrightarrow\f\le\p\text{ on }X^\an.
$$
In particular, $\f$ is the smallest usc extension of $\f|_{X^\div}$ to $X^\an$. 
\end{thm}

\begin{cor}\label{cor:Hausdorff} Every $\theta$-psh function is uniquely determined by its restriction to $X^\div$, and the topology of $\PSH(\theta)$ is Hausdorff. 
\end{cor}
 
Because of this, we will occasionally say that a function $\f\colon X^\div\to\R$ is $\theta$-psh if it admits a (necessarily unique) extension to a $\theta$-psh function on $X^\an$.

\begin{cor}\label{cor:pshcst} If $X$ is connected, then any $0$-psh function is constant, \ie $\PSH(0)=\R$. 
\end{cor}
\begin{proof} By GAGA, $X^\an$ is connected, and after replacing $X$ with an irreducible component, we may assume that $X$ is irreducible. By Corollary~\ref{cor:Hausdorff}, it is then enough to show that any $\f\in\PSH(0)$ with $\sup\f=\f(v_\triv)=0$ satisfies $\f(v)=0$ for any $v\in X^\div$. Consider the decreasing sequence $(m\f)_{m\in\N}$ and its pointwise limit $\p:=\lim_m (m\f)$. Since $\p(v_\triv)=0$, we have $\p\in\PSH(0)$ (see Theorem~\ref{thm:psh12}). By  Corollary ~\ref{cor:pshonXlin}, we thus have $\p(v)>-\infty$, \ie $\f(v)=0$.
\end{proof}

\begin{cor}\label{cor:decrweak} Pick a decreasing net $(\f_i)$ in $\PSH(\theta)$, and assume that it converges to $\f\in\PSH(\theta)$. Then $\f=\lim_i\f_i$ pointwise on $X^\an$.
\end{cor}
\begin{proof} The pointwise limit $\p$ of $(\f_i)$ on $X^\an$ coincides with $\f$ on $X^\div$. In particular, $\p(v_{\triv,\a})$ is finite for any $\a$, and hence $\p\in\PSH(\theta)$, by Theorem~\ref{thm:psh12}. By Corollary~\ref{cor:Hausdorff}, we infer $\f=\p$, which means that $\f_i$ converges to $\f$ pointwise on $X^\an$ 
\end{proof}

\begin{lem}\label{lem:supPL} For each $\p\in\PL(X)$, there exists a finite set $\Sigma\subset X^\div$ such that 
\begin{equation}\label{equ:supPL}
\sup_{X^\an}(\f-\p)=\max_\Sigma(\f-\p)
\end{equation}
for all $\f\in\PL^+(X)$, and the same result then holds for all $\f\in\PSH(\theta)$, $\theta\in\Num(X)$. Moreover, if $\p=r(\f_\fa-\rho)$, where $r\in\Q_{>0}$, $\fa$ is a flag ideal, and $\rho\in\PL^+(X)$, then we can pick $\Sigma=\Sigma_\fa$, the set of Rees valuations of $\fa$.
\end{lem}
\begin{proof}
  We can always write $\p=r(\f_\fa-\rho)$ as above, and by homogeneity we may assume $r=1$.
  The case $\f\in\PL^+(X)$ now follows from Lemma~\ref{lem:Shilov}, since $\f+\rho\in\PL^+(X)$. 

  Now assume $\f\in\PSH(\theta)$ for some $\theta\in\Num(X)$. After replacing $\theta$ with $\om\in\Amp(X)$ such that $\om-\theta$ is nef, we may assume $\theta=\om$ is ample. By Theorem~\ref{thm:pshample}, we can find a decreasing net $(\f_i)$ in $\cH^\dom(\om)\subset\PL^+(X)$ such that $\f_i\searrow\f$. For each $v\in X^\an$, we then have $\f_i(v)-\p(v)\le\max_{\Sigma_\fa}(\f_i-\p)$ for all $i$, and hence $\f(v)-\p(v)\le\max_{\Sigma_\fa}(\f-\p)$.
\end{proof}

\begin{proof}[Proof of Theorem~\ref{thm:suppsh}] If $\p\in\PL(X)$, the result follows from Lemma~\ref{lem:supPL}. Assume next $\p\in\Cz(X)$. By Theorem~\ref{thm:PLdense}, $\p$ is a uniform limit of functions in $\PL(X)$, and the result thus holds for $\p$ as well. In the general case, we can find a decreasing net $(\p_i)$ in $\Cz(X)$ with $\p_i\to\p$. Then $\f\le\p\le\p_i$ on $X^\div$ implies $\f\le\p_i$ on $X^\an$, and hence $\f\le\p$ on $X^\an$ as well. 
\end{proof}

\begin{rmk}\label{rmk:Hausdorff} Assume $\charac k=0$. By Corollary~\ref{cor:apprbelow}, any $v\in X^\an$ is the limit of a net $(v_i)$ in $X^\div$ such that $v_i\le v$ for all $i$. As any $\f\in\PSH(\theta)$ is decreasing and usc, it follows that $\f(v_i)\to\f(v)$, which provides a slightly more precise version of Corollary~\ref{cor:Hausdorff} in that case. 
\end{rmk}
  
We record another useful consequence of Lemma~\ref{lem:supPL}. 

\begin{prop}\label{prop:suppsh}
  For any $\p\in \Cz(X)$, $\f\mapsto\sup_{X^\an}(\f-\p)$ is continuous on $\PSH(\theta)$.
\end{prop}
\begin{proof} When $\p$ is PL, this is a direct consequence of Lemma~\ref{lem:supPL}. As above, the general case follows by density of $\PL(X)$ in $\Cz(X)$ with respect to uniform convergence.
\end{proof}

\begin{cor}\label{cor:pointusc} For any $v\in X^\an$, the evaluation map $\f\mapsto\f(v)$ is usc on $\PSH(\theta)$. 
\end{cor}
In other words, each convergent net $\f_i\to\f$ in $\PSH(\theta)$ satisfies $\f\ge\limsup_i\f_i$ pointwise on $X^\an$.

\begin{proof} Pick a convergent net $\f_i\to\f$ in $\PSH(\theta)$. Since $\f$ is usc and $X^\an$ is compact, we have $\f=\inf\{\p\in\Cz(X)\mid\p\ge\f\}$ pointwise. For each candidate $\p$ and $\e>0$, Proposition~\ref{prop:suppsh} implies that $\f_i\le\p+\e$ on $X^\an$ for all $i$ large enough. Thus $\limsup_i\f_i\le\p+\e$ pointwise on $X^\an$, and the result follows by taking the infimum over $\p$ and letting $\e\to 0$. 
\end{proof}

\begin{rmk}
  The evaluation map $\f\mapsto\f(v)\in\R\cup\{-\infty\}$ can fail to be continuous in general. For example, take $X=\P^1$ and $L=\cO(1)$. In homogeneous coordinates $[z_0:z_1]$, set $\f_m:=\max\{m^{-1}\log|z_1|,-1\}$ for $m\ge1$. Then $\f_m\to0$ in $\PSH(L)$, but for $p:=\{z_1=0\}\in X(k)$ we have $\f_m(v_{p,\triv})=-1$ for all $m$, see~\S\ref{sec:curve}.
\end{rmk}
 
We also note: 
\begin{lem}\label{lem:pullbackcont} Let $f\colon Y\to X$ be a morphism from a projective variety, such that any irreducible component of $Y$ is mapped onto a component of $X$. For any $\theta\in\Num(X)$, the map $f^\star\colon\PSH(\theta)\to\PSH(f^\star\theta)$ is then continuous. 
\end{lem} 
\begin{proof} The assumption guarantees that $f$ maps $Y^\div$ to $X^\div$, and the result follows. 
\end{proof} 

We conclude this section with an `almost birational invariance' property of $\theta$-psh functions. Recall that any projective birational morphism is (isomorphic to) the blowup of some (generically trivial) ideal, see~\cite[Theorem II.7.17]{Har}.

\begin{thm}\label{thm:descentpsh} For any birational morphism $\pi\colon Y\to X$, there exists $\om_X\in\Amp(X)$ and $\f_X\in\PSH(\om_X)$ such that
$$
\PSH(\pi^\star\theta)+\pi^\star\f_X\subset\pi^\star\PSH(\theta+\om_X). 
$$
for all $\theta\in\Num(X)$. If $\pi$ is realized as the blowup of a generically trivial ideal $\fb\subset\cO_X$, then one can take $\f_X=\log|\fb|$. 
\end{thm}
As we shall see in Lemma~\ref{lem:descentpsh2} below, we actually have $\PSH(\pi^\star\theta)=\pi^\star\PSH(\theta)$ when the \emph{envelope property} holds for $\theta$. However, the `error term' $\pi^\star\f_X$ is necessary in general (for instance if $\pi$ is the normalization of a reducible variety). 

\begin{cor}\label{cor:descentpsh} In the notation of Theorem~\ref{thm:descentpsh}, there exists $\om_Y\in\Amp(Y)$ such that $\PSH(\om_Y)+\pi^\star\f_X\subset\pi^\star\PSH(\om_X)$. 
\end{cor}

\begin{proof}[Proof of Theorem~\ref{thm:descentpsh}] Write $\pi$ as the blowup of an ideal $\fb\subset\cO_X$. Denote by $E$ the exceptional divisor, and choose an ample line bundle $H$ on $X$ such that $H\otimes\fb$ is globally generated. As $-E$ is $\pi$-ample, we can further arrange that $H':=\pi^\star H-E$ is ample. We claim that the result holds with $\f_X:=\log|\fb|$ and $\om_X:=c_1(H)$. Since $H\otimes\fb$ is globally generated, we can choose a finite subset $(s_i)$ of $\Hnot(X,H\otimes\fb)$ that locally generates $\fb$, and hence $\log|\fb|=\max_i\log|s_i|\in\cH^\gf_0(H)\subset\PSH(\om_X)$. 

Assume first that $\theta=c_1(L)$ with $L\in\Pic(X)_\Q$. Pick $\f'\in\PSH(\pi^\star L)$, and write $\f'$ as the limit of a decreasing net $\f'_i\in\cH^\gf_\R(L'_i)$ with $L'_i\in\Pic(X)_\Q$ and $c_1(L'_i)\to c_1(\pi^\star L)$. Since $H'$ is ample, $H'+(\pi^\star L-L'_i)=(\pi^\star(L+H)-E)-L'_i$ is ample for all $i$ large enough, and hence $\f'_i\in\cH^\gf_\R(\pi^\star (L+H)-E)$, by Proposition~\ref{prop:FS}. By Lemma~\ref{lem:descentFS}, there exists $\f_i\in\cH^\gf_\R(L+H)$ such that $\f'_i+\pi^\star\f_X=\pi^\star\f_i$. As $\pi$ induces an isomorphism $Y^\div\simto X^\div$, $\f_i$ is uniquely determined, and the net $(\f_i)$ is decreasing on $X^\an$, by Theorem~\ref{thm:suppsh}. Its pointwise limit $\f:=\lim_i\f_i$ satisfies $\f'+\pi^\star\f_X=\pi^\star\f$ on $X^\div$. In particular, $\f(v_{\triv,\a})$ is finite for all $\a$. This shows that $\f$ is generically finite, and hence $\f\in\PSH(L+H)$ (see Theorem~\ref{thm:psh12}), which proves the result in that case. 

In the general case, we can choose a sequence $L_m\in\Pic(X)_\Q$ such that $c_1(L_m)-\theta$ is nef and converges to $0$. By Theorem~\ref{thm:PSH}, we have $\PSH(\theta)\subset\PSH(L_m)$, and the previous step thus shows that any $\f'\in\PSH(\pi^\star\theta)\subset\PSH(\pi^\star L_m)$ satisfies $\f'+\pi^\star\f_X=\pi^\star\f_m$ with $\f_m\in\PSH(L_m+H)$. As noted above, $\f_m$ is uniquely determined by $\f$, and hence independent of $m$. Since $c_1(L_m)\to\theta$, Theorem~\ref{thm:psh12} yields $\f\in\PSH(\theta)$, which concludes the proof. 
\end{proof}

\begin{proof}[Proof of Corollary~\ref{cor:descentpsh}] Use the same notation as above, and set $\om_Y:=\tfrac 12(\pi^\star\om_X-E)\in\Amp(Y)$. Since $\pi^\star\f_X=\log|s_E|\subset\PSH(E)$, we have 
$\PSH(2\om_Y)+\pi^\star\f_X\subset\PSH(\pi^\star\om_X)$, and Theorem~\ref{thm:descentpsh} yields
$$
\PSH(2\om_Y)+2\pi^\star\f_X\subset\pi^\star\PSH(\pi^\star(2\om_X)),
$$
which is equivalent to the desired the result. 
\end{proof}
  
%
%
\subsection{Pluripolar sets}\label{sec:pluripolar}
Mimicking classical pluripotential theory, we introduce: 
 
\begin{defi} A subset $E\subset X^\an$ is \emph{pluripolar} if there exist   $\theta\in\Num(X)$   and $\f\in\PSH(\theta)$ such that $E\subset\{\f=-\infty\}$.
\end{defi}

\begin{lem}\label{lem:ppindep} Pick $\om\in\Amp(X)$. Then $E\subset X^\an$ is pluripolar iff $E\subset\{\f=-\infty\}$ for some $\f\in\PSH(\om)$. 
\end{lem}
\begin{proof} Assume $E\subset\{\p=-\infty\}$ for some $\p\in\PSH(\theta)$, $\theta\in\Num(X)$. For $t>0$ large enough we have $t\om\ge\theta$. Thus $\f:=t^{-1}\p$ is $\om$-psh (see Theorem~\ref{thm:PSH}), and satisfies $E\subset\{\f=-\infty\}$. 
\end{proof}

\begin{lem}\label{lem:Zarpp}
  For any Zariski closed subset $Z\subset X$, $Z^\an$ is pluripolar iff $Z$ is nowhere dense, \ie does not contain any irreducible component $X_\a$.
\end{lem}
\begin{proof}
  If $Z^\an\subset\{\f=-\infty\}$ with $\f\in\PSH(\theta)$, then $Z$ does not contain any component $X_\a$, since $\f$ is generically finite. Conversely, if $Z$ is nowhere dense then we can find a section $s$ of some ample line bundle $L$ that vanishes along $Z$, but not along any irreducible component of $X$. Then $\f:=\log|s|$ is $L$-psh, and $Z^\an\subset\{\f=-\infty\}$. 
\end{proof}
 
By the next result, any countable subset of $X(k)\subset X^\an$ is therefore also pluripolar, but some are Zariski dense.
\begin{lem}\label{lem:ppcount} Any countable union of pluripolar sets is pluripolar.
\end{lem}
   
\begin{proof}
  Let $(E_m)_{m\ge 1}$ be a sequence of pluripolar subsets, and fix $\om\in\Amp(X)$. For each $m$, pick $\f_m\in\PSH(\om)$ such that $E_m\subset\{\f_m=-\infty\}$, and set $a_{m,\a}:=\sup_{X_\a^\an}\f_m>-\infty$. We can pick $c_m>0$ such that $\sum_mc_m\le 1$ and $\sum_{m,\a}c_m|a_{m,\a}|<\infty$. Then $\f:=\sum c_m\f_m\in\PSH(\om)$ and $\bigcup_mE_m\subset\{\f=-\infty\}$.
\end{proof}

\begin{prop}\label{prop:linnpp} A point $v\in X^\an$ is nonpluripolar iff $v\in X^\lin$, \ie $v$ is a valuation of linear growth. 
\end{prop}
 
\begin{proof}
  Corollary~\ref{cor:pshonXlin} shows that no valuation of linear growth is pluripolar. Now suppose $v\in X^\an\setminus X^\lin$ and pick any ample line bundle $L$. For any integer $j\ge1$ there exist $m_j\ge1$ and a regular section $s_j\in\Hnot(X,m_jL)$ such that $m_j^{-1}v(s_j)\ge 2^j$. If we set $\f_j=m_j^{-1}\log|s_j|$ and $\p_j=\sum_{l=1}^j2^{-l}\f_l$, then $(\p_j)_j$ is a decreasing sequence in $\PSH(L)$ with $\p_j(v_{\triv,\a})=0$ and $\p_j(v)\le-j$. It follows that $\p=\lim\p_j=\sum_{l=1}^\infty2^{-l}\f_l$ satisfies $\p\in\PSH(L)$ and $\p(v)=-\infty$, so that $v$ is pluripolar.
\end{proof}
\begin{cor}\label{cor:ppint} Any pluripolar set has empty interior.
\end{cor}
\begin{proof}
  Indeed, a pluripolar set must be disjoint from the dense subset $X^\div\subset X^\lin$. \end{proof}
\begin{lem}\label{lem:ppbir}
  If $\pi\colon Y\to X$ is a birational morphism, then a subset $E\subset X^\an$ is pluripolar iff $\pi^{-1}(E)\subset Y^\an$ is pluripolar.
\end{lem}
\begin{proof}
  Write $\pi$ as the blowup of a generically trivial ideal $\fb\subset\cO_X$. Let $Z\subset X$ be the cosupport of $\fb$, and set $W:=\pi^{-1}(Z)$. Then $Z$ and $W$ are Zariski closed, nowhere dense subsets of $X$ and $Y$, respectively, and $\pi\colon Y\setminus W\to X\setminus Z$ is an isomorphism. Then $Z^\an\subset X^\an$ and $W^\an\subset X^\an$ are pluripolar, see Lemma~\ref{lem:Zarpp}. Set $\f_X:=\log|\fb|$ and pick $\om_X\in\Amp(X)$ as in Theorem~\ref{thm:descentpsh}.

  First suppose $E$ is pluripolar, and pick $\f\in\PSH(\om_X)$ such that $\f=-\infty$ on $E$. Then $\p:=\pi^\star\f\in\PSH(\pi^\star\om_X)$ and $\p=-\infty$ on the set $\pi^{-1}(E)$, which is therefore pluripolar.

  If instead $\pi^{-1}(E)$ is pluripolar, then pick $\om_Y\in\Amp(Y)$ as in Corollary~\ref{cor:descentpsh}, and $\p\in\PSH(\om_Y)$ such that $\p=-\infty$ on $\pi^{-1}(E)$. By Corollary~\ref{cor:descentpsh} we can find $\f\in\PSH(\om_X)$ such that $\p+\pi^\star\f_X=\pi^\star\f$. It follows that $\f=-\infty$ on $E\setminus Z^\an$, so $E\setminus Z^\an$ is pluripolar. As $Z^\an$ is also pluripolar, we conclude that $E$ is pluripolar, see Lemma~\ref{lem:ppcount}.
\end{proof}
Applying Lemma~\ref{lem:ppbir} to the canonical birational map $\coprod_\a X_\a\to X$, we get
\begin{cor}\label{cor:ppcomp}
  A set $E\subset X^\an$ is pluripolar iff $E\cap X_\a^\an$ is a pluripolar subset of $X_\a^\an$ for each $\a$.
\end{cor}
For later use, we note: 
\begin{exam}\label{exam:ppdense} Any trivial semivaluation $v\in X^\triv$ lies in the closure of some pluripolar subset $E\subset X^\an$, that can be chosen as a countable subset of $X(k)\subset X^\an$. When $X$ is a connected smooth curve, this follows from \S\ref{sec:curve}, by choosing any infinite sequence of distinct closed points in $X$. In the general case, $v$ lies in the closure of $X(k)$ (see~\S\ref{sec:triv}), and hence in the closure of a countable subset thereof (see Remark~\ref{rmk:countable}). 
\end{exam}
\begin{rmk}\label{rmk:ppnotdense}
If $\dim X=0$, then there are no nonempty pluripolar subsets of (the finite set) $X^\an$. If $\dim X=1$, then any pluripolar set must be contained in $X(k)\subset X^\an$, whose closure is the strict subset $X^\triv\subset X^\an$. In particular, pluripolar sets are never dense in that case (in stark contrast with the complex analytic case). 

When $X$ has no irreducible component of dimension $\le1$, the situation is more subtle: the set $X^\an\setminus X^\val$ is then dense (see Lemma~\ref{lem:nonvaldense}), and it is pluripolar if $k$ is countable, being the (countable) union of all strict irreducible subvarieties (see Lemma~\ref{lem:ppcount}). However, when $k$ is uncountable, pluripolar sets are again never dense: see Corollary~\ref{cor:ppnotdense} below. 
\end{rmk}

\subsection{The Alexander--Taylor capacity}\label{sec:AT}
In order to detect pluripolar sets, we introduce the following variant\footnote{\   More precisely, $\exp(-\!\te_\om)$ corresponds to the Alexander--Taylor capacity. }  of the classical Alexander--Taylor capacity~\cite{AT84}.

  \begin{defi} For any $\om\in\Amp(X)$ and any subset $E\subset X^\an$, we define $\te_\om(E)\in[0,+\infty]$ as follows:
    \begin{itemize}
    \item[(i)] if $X$ is irreducible, then 
        $\te_\om(E):=\sup_{\f\in\PSH(\om)}\left(\sup_{X^\an}\f-\sup_E\f\right)$.
      \item[(ii)]
        in general, $\te_\om(E)=\min_\a\te_{\om|_{X_\a}}(E\cap X_\a^\an)$.
      \end{itemize}
\end{defi}
 
Recall that $\sup_{X^\an}\f=\f(v_\triv)$ when $X$ is irreducible, see Proposition~\ref{prop:pshconvex}.

\begin{thm}\label{thm:ppT} A subset $E\subset X^\an$ is pluripolar iff $\te_\om(E)=\infty$ for some (equivalently, any) $\om\in\Amp(X)$.
\end{thm}
\begin{proof}
  Pick any ample class $\om\in\Num(X)$.
   By Corollary~\ref{cor:ppcomp} $E$ is pluripolar iff $E\cap X_\a^\an$ is pluripolar for all $\a$, and by definition, $T_\om(E)=\infty$ iff $T_{\om|_{X_\a}}(E\cap X_\a^\an)=\infty$ for all $\a$. We may therefore assume that $X$ is irreducible.   
  First assume that $E$ is pluripolar. By Lemma~\ref{lem:ppindep}, there exists $\f\in\PSH(\om)$ with $\f|_E\equiv-\infty$, and hence $\te_\om(E)=\infty$. Conversely, suppose $\te_\om(E)=\infty$. For each $m\in\N$ we can then find $\f_m\in\PSH(\om)$ such that $\sup\f_m=\f_m(v_\triv)=0$ and $\sup_E\f_m\le -2^m$. By convexity of $\PSH(\om)$, setting for each $m$ 
$$
\p_m:=\sum_{l=1}^m 2^{-l}\f_l=\sum_{i=1}^m 2^{-l}\f_l+2^{-m}\cdot 0
$$ 
defines a decreasing sequence in $\PSH(\om)$ such that $\p_m(v_\triv)=0$ and $\sup_E\p_m(v)\le -m$. By Theorem~\ref{thm:PSH}, 
$$
\p:=\lim_m\p_m=\sum_{l=1}^\infty 2^{-l}\f_l
$$ 
is $\om$-psh, and $\sup_E\p=-\infty$. 
\end{proof}

For line bundles, the chosen notation is compatible with Definition~\ref{defi:TT}: 

\begin{lem}\label{lem:TT} If $L\in\Pic(X)_\Q$ is ample and $\om=c_1(L)$, then $\te_\om(\{v\})=\teL(v)$ for any $v\in\Xan$.
\end{lem}
\begin{proof}
   
  By Lemma~\ref{lem:Tcomp} and the definition of the Alexander--Taylor capacity, we may assume that $X$ is irreducible.
  By definition, $\teL(v)$ is then the supremum of $m^{-1}v(s)$ over $m$ sufficiently divisible and $s\in \Hnot(X,mL)\setminus\{0\}$. Equivalently, it is the smallest constant such that all functions of the form $\f:=m^{-1}\log|s|$ satisfy $\f(v)\ge\f(v_\triv)-\teL(v)$. As $\f\in\PSH(L)$, we get $\te_L(v)\le\te_\om(\{v\})$. But adding constants and taking finite maxima shows that the inequality $\f(v)\ge\f(v_\triv)-\teL(v)$ also holds for all $\f\in \cH(L)$, and then for all $\f\in\PSH(L)$, by Theorem~\ref{thm:pshample}.
   
\end{proof}


%
%
%
%
 
\section{Envelopes and negligible sets}\label{sec:envnegl}
By our definitions, $\theta$-psh functions are well-behaved under decreasing limits. As in the complex analytic case, many important constructions involve increasing limits or envelopes of $\theta$-psh functions. Consider a bounded-above family $(\f_i)_i$ of $\theta$-psh functions. In general, the supremum $\f:=\ssup_i\f_i$ may fail to usc, and is therefore not $\theta$-psh. In this section, we study whether the usc regularization $\f^\star$ is $\theta$-psh. We conjecture that this is true when $X$ is unibranch, and we prove it when $X$ is smooth, $\theta$ is nef, and either $\dim X\le 2$ or $\charac k=0$, see Theorem~\ref{thm:contenvlisse}.
%
%
\subsection{Negligible sets}\label{sec:negl}
We start by studying the following notion, imported from classical pluripotential theory. 
For the time being $X$ is an arbitrary projective variety, with irreducible components $X_\a$.
 
\begin{defi}\label{defi:negl} A subset $E\subset X^\an$ is \emph{negligible} if there exists   $\theta\in\Num(X)$  and a bounded-above family $(\f_i)$ in $\PSH(\theta)$ such that $E\subset\{\ssup_i\f_i<\supstar_i\f_i\}$. 
\end{defi}
 
As in Lemma~\ref{lem:ppindep}, we may always choose $\theta$ as a fixed ample class.
 
Replacing $\f_i$ by $\exp(\f_i-\sup_j\sup_{X^\an}\f_j)-1$, see Corollary~\ref{cor:exppsh}, we may also assume $-1\le\f_i\le 0$.

Any subset of a negligible set is trivially negligible. We also have
\begin{lem}\label{lem:neglcount}
  A countable union of negligible sets is negligible.
\end{lem}
\begin{proof}
  Let $(E_m)_1^\infty$ be a sequence of negligible sets, and set $E:=\bigcup E_m$. Fix $\om\in\Amp(X)$. For each $m$ we can find a family $(\f_{m,i})_{i\in I_m}$ in $\PSH(\om)$ such that $-1\le\f_{m,i}\le0$ on $X^\an$ and $\f_m<\f_m^\star$ on $E_m$, where $\f_m=\sup_{i\in I_m}\f_{m,i}$. Now set $I=\prod_mI_m$ and, for $i=(i_m)_m\in I$, $\p_i=\sum_{m=1}^\infty2^{-m}\f_{m,i_m}$. Then $\p_i\in\PSH(\om)$ and $\p:=\sup_i\p_i=\sum_m2^{-m}\f_m$. This gives $\p<\p^\star$ on $E$, so $E$ is negligible.
\end{proof}

\begin{prop}\label{prop:plurineg} Every pluripolar subset is negligible.
\end{prop}
 
The converse implication is much more subtle. We will prove it when $\charac k=0$ or $\dim X\le 2$, see Corollary~\ref{cor:neglpp} below. 
\begin{proof} Suppose $E\subset X^\an$ is pluripolar, and pick $\p\in\PSH(\om)$ such that $\p\le 0$ and $\p=-\infty$ on $E$. We may assume $E=\{\p=-\infty\}$.
  For $m\ge 1$, set $\f_m:=m^{-1}\p\in\PSH(\om)$ and $\f:=\sup_m\f_m$. Then $\f\equiv-\infty$ on $E$ and $\f\equiv 0$ on $X^\an\setminus E$. Since $E$ has empty interior (see Corollary~\ref{cor:ppint}),
  $\f^\star\equiv 0$ on $X^\an$, and $E=\{\f<\f^\star\}$ is thus negligible.
\end{proof}
Just like pluripolar sets, the class of negligible sets is birationally invariant.
\begin{lem}\label{lem:neglbir}
  If $\pi\colon Y\to X$ is a birational morphism, then a subset $E\subset X^\an$ is negligible iff $\pi^{-1}(E)\subset Y^\an$ is negligible.
\end{lem}
\begin{proof}
  We follow the setup and notation of the proof of Lemma~\ref{lem:ppbir}.

  First suppose $E$ is negligible, and pick a family $(\f_i)_i$ in $\PSH(\om_X)$ uniformly bounded above such that $\ssup_i\f_i<\supstar_i\f_i$ on $E$. If we set $\p_i:=\pi^\star\f_i$, then $(\p_i)_i$ is a family in $\PSH(\pi^\star\om_X)$ that is uniformly bounded above. Moreover, $\ssup_i\p_i<\supstar_i\p_i$ on $\pi^{-1}(E)\setminus W^\an$, so $\pi^{-1}(E)\setminus W^\an$ is negligible. As $W^\an$ is pluripolar, and hence also negligible, it follows that $\pi^{-1}(E)$ is negligible, see Lemma~\ref{lem:neglcount}.

  If instead $\pi^{-1}(E)$ is negligible, then there exists a family $(\p_i)_i$ in $\PSH(\om_Y)$ uniformly  bounded above, such that $\ssup_i\p_i<\supstar_i\p_i$ on $\pi^{-1}(E)$. We can find $\f_i\in\PSH(\om_X)$ such that $\p_i+\pi^\star\f_X=\pi^\star\f_i$. Then $(\f_i)_i$ is uniformly bounded above, and $\ssup_i\f_i<\supstar_i\f_i$ on $E\setminus Z^\an$, so $E\setminus Z^\an$ is negligible, and we conclude as above.
\end{proof}
\begin{cor}\label{cor:neglcomp}
  A subset $E\subset X^\an$ is negligible iff $E\cap X_\a^\an$ is a negligible subset of $X_\a^\an$ for all $\a$.
\end{cor}
The following result will play a crucial role in what follows.

\begin{thm}\label{thm:divnegl}
Divisorial points are nonnegligible. 
\end{thm}
Thus non-empty open sets are also nonnegligible, as they contain divisorial points, by density of $X^\div$. The key ingredient is the following partial converse to Lemma~\ref{lem:supPL}. 

\begin{lem}\label{lem:divnegl} Assume that $X$ is irreducible. Then, for each $\theta\in\Num(X)$ and $v\in X^\div$, there exists $\p\in\PL(X)$ such that 
$$
\sup_{X^\an}(\f-\p)=(\f-\p)(v)
$$ 
for all $\f\in\PSH(\theta)$.
\end{lem}
 
\begin{proof}
  By Lemma~\ref{lem:divRees} we may find a flag ideal $\fa$ for which $v$ is a Rees valuation, that is, $v\in\Sigma_\fa$. Pick $L\in\Pic(X)_\Q$ ample such that  $c_1(L)\ge\theta$. Then $\PSH(\theta)\subset\PSH(L)$, so we may assume $\theta=c_1(L)$.
  As $X$ is irreducible, there exists $C>0$ such that $|\f(v)-\f(w)|\le C$ for all $w\in\Sigma_\fa$ and all $\f\in\PSH(L)$; for example, we can take $C=2\max_{w\in\Sigma_\fa}\teL(w)$, as follows from the proof of Corollary~\ref{cor:pshonXlin}. By Lemma~\ref{lem:PLinterpol} there exists $m,r\ge1$ and $\rho\in\cH(L)$ such that if we set $\p:=r(\f_\fa-m\rho)$, then $\p(v)=0$ and $\p(w)=2C$ for $w\in\Sigma_\fa$, $w\ne v$. Then $\p\in\PL(X)$, and for any $\f\in\PSH(L)$, the max of $\f-\p$ over $\Sigma_\fa$ is attained at $v$. On the other hand, it follows from Lemma~\ref{lem:supPL} that the supremum of $\f-\p$ on $X^\an$ is attained on $\Sigma_\fa$. The proof is complete.
\end{proof}
 
\begin{proof}[Proof of Theorem~\ref{thm:divnegl}]\phantom{}  In view of Corollary~\ref{cor:neglcomp} we may assume that $X$ is irreducible. Consider $\theta\in\Num(X)$   and a bounded-above family $(\f_i)$ in $\PSH(\theta)$, and set $\f:=\sup_i\f_i$. Pick $v\in X^\div$, and choose a PL function $\p$ as in Lemma~\ref{lem:divnegl}. For each $i$ we have $\f_i\le\f_i(v)-\p(v)+\p$. Thus $\f\le\f(v)-\p(v)+\p$, and hence $\f^\star\le\f(v)-\p(v)+\p$, by continuity of $\p$. This yields $\f^\star(v)=\f(v)$, which proves that $v$ is nonnegligible. 
\end{proof}
%
%
\subsection{The envelope property and unibranch varieties}\label{sec:envprop}
  In the rest of~\S\ref{sec:psh} we assume that $X$ is \textbf{irreducible} unless stated otherwise. (See Remark~\ref{rmk:envred}.)

\begin{defi} We say that a class $\theta\in\Num(X)$ has the \emph{envelope property} if, for any bounded-above family $(\f_i)$ in $\PSH(\theta)$, the usc upper envelope $\supstar_i\f_i$ is $\theta$-psh.
\end{defi}
In the complex analytic case, it is a basic fact that psh functions on domains in $\C^n$ satisfy the analogue of the envelope property. As a consequence, $\theta$-psh functions satisfy the envelope property for any projective complex manifold $X$ and any closed $(1,1)$-form $\theta$. More generally, the latter property holds when $X$ is unibranch (\ie locally analytically irreducible), \cf~\cite[Th\'eor\`eme 1.7]{DemSMF}, but fails in general, for the same reason as Theorem~\ref{thm:unibranch} below. 

With our global definition of $\theta$-psh functions, the envelope property turns out to be especially delicate. The purpose of what follows is to explore various implications, and establish it in an important special case. We first observe:

\begin{lem}\label{lem:envppty} Assume that $\theta\in\Num(X)$ can be written as the limit of classes $\theta_j\ge\theta$ that satisfy the envelope property. Then $\theta$ has the envelope property as well. 
\end{lem}
\begin{proof} Pick a bounded-above family $(\f_i)$ in $\PSH(\theta)$. Since $\theta\le\theta_j$, we have $\PSH(\theta)\subset\PSH(\theta_j)$.
  Thus $\supstar_i\f_i$ is $\theta_j$-psh for all $j$, and hence $\theta$-psh, by Theorem~\ref{thm:psh12}.
\end{proof}
\begin{cor}
  If the continuity property holds for all $\om\in\Amp(X)$, then it holds for all $\theta\in\Nef(X)$.
\end{cor}
We next prove that the envelope property is equivalent to the analogue of a basic compactness property in the complex analytic case. 

\begin{thm}\label{thm:envppty} For any $\theta\in\Num(X)$, the following properties are equivalent:
\begin{itemize}
\item[(i)] the envelope property holds for $\theta$; 
\item[(ii)] the space $\PSH_{\sup}(\theta):=\left\{\f\in\PSH(\theta)\mid\sup\f=0\right\}$ is compact;
\item[(iii)] every bounded-above, increasing net $(\f_j)$ in $\PSH(\theta)$ converges in $\PSH(\theta)$. 
\end{itemize}
\end{thm}
\begin{proof} Assume (i), and pick a net $(\f_i)$ in $\PSH_{\sup}(\theta)$. By Proposition~\ref{prop:linnpp}, $(\f_i(v))$ is a bounded net in $\R$ for each $v\in X^\div$. By Tychonoff's theorem, after passing to a subnet, we may thus assume that $(\f_i)$ converges pointwise on $X^\div$ to a function $\f\colon X^\div\to\R$, and it suffices to show that $\f$ extends to a function in $\PSH(\theta)$ (necessarily unique, by Corollary~\ref{cor:Hausdorff}). By (i), $\p_i:=\supstar_{j\ge i}\f_j$ is $\theta$-psh for each $i$. The net $(\p_i)$ is further decreasing, and its pointwise limit $\p:=\lim_i\p_i$ is thus either $\theta$-psh, or identically $-\infty$, by Theorem~\ref{thm:psh12}. Now Theorem~\ref{thm:divnegl} implies that $\p_i=\sup_{j\ge i}\f_j$ on $X^\div$. Thus $\p=\f$ on $X^\div$. In particular, $\p\not\equiv-\infty$, hence $\p\in\PSH(\theta)$. This proves (i)$\Rightarrow$(ii). 

Next, assume (ii), and consider a bounded-above, increasing net $(\f_j)$ in $\PSH(\theta)$. Then a subnet $\f_{j_i}-\sup\f_{j_i}$ converges to some $\f\in\PSH(\theta)$. Since $(\sup\f_j)$ is increasing and bounded-above, it converges to some $c\in\R$. Thus $\f_{j_i}\to\f+c$ on $X^\div$, and hence $\f_j\to\f+c$ on $X^\div$ since $(\f_j)$ is increasing. This proves (ii)$\Rightarrow$(iii).

Finally, assume (iii). Let $(\f_i)_{i\in I}$ be a bounded-above family of $\theta$-psh functions, and consider the increasing net $\p_J:=\max_{i\in J}\f_i$ parametrized by finite subsets $J\subset I$. Then $\p_J$ admits a limit $\p$ in $\PSH(\theta)$, and we claim that $\p=\supstar_i\f_i$. On $X^\div$ we have 
$$
\p=\lim_J\max_{i\in J}\f_i=\ssup_i\f_i=\supstar_i\f_i,
$$
by Theorem~\ref{thm:divnegl}. Since $\p$ is $\theta$-psh and $\supstar_i\f_i$ is usc, Theorem~\ref{thm:suppsh} yields $\p\le\supstar_i\f_i$ on $X^\an$. Similarly, Theorem~\ref{thm:suppsh} yields $\f_i\le\p$, and hence $\supstar_i\f_i\le\p$, since $\p$ is usc. This proves the claim, and hence (iii)$\Rightarrow$(i).
\end{proof}

Recall that  the variety $X$ is said to be  \emph{unibranch} if the following equivalent conditions hold (see~\cite[IV.7.6.3]{EGA} and~\cite[Corollary 32]{Kol}): 
\begin{itemize}
\item the normalization morphism $\nu\colon X^\nu\to X$ is a homeomorphism; 
\item the formal completion of $X$ at each of its points is irreducible.
\end{itemize}
In particular, any normal variety is unibranch. When $k=\C$, $X$ is unibranch iff the associated complex analytic space is locally irreducible in the analytic topology.   This is also true for the Berkovich analytification as considered in this paper: the variety $X$ is unibranch iff the $k$-analytic space $X^\an$ is locally irreducible in the analytic topology (see~\cite[Lemme 5.19]{Duc09}).

\begin{thm}\label{thm:unibranch} If the envelope property holds for some $\om\in\Amp(X)$, then $X$ is necessarily unibranch. 
\end{thm} 
The analogue of this result is known in the complex analytic setting, too. We have not been able to locate a precise reference, but the proof below of Theorem~\ref{thm:unibranch} can easily be adapted to that setting.
\begin{lem}\label{lem:descentpsh2} Let $\pi\colon Y\to X$ be a birational morphism. Pick $\theta\in\Num(X)$, and assume that $\theta$ can be written as the limit of a sequence of classes $\theta_m>\theta$ with the envelope property. Then $\pi^\star\colon\PSH(\theta)\simto\PSH(\pi^\star\theta)$ is a homeomorphism. 
\end{lem}
Recall that $\theta_m>\theta$ means that $\theta_m-\theta$ is ample. 

\begin{proof} After passing to a subsequence, we may assume that $\theta_m>\theta_{m+1}$ for all $m$, and hence $\PSH(\theta_{m+1})\subset\PSH(\theta_m)$. Since $\pi$ induces a bijection $Y^\div\simto X^\div$, $\pi^\star\colon\PSH(\theta)\to\PSH(\pi^\star\theta)$ is a topological embedding, and it suffices to show that it is onto. Pick $\p\in\PSH(\pi^\star\theta)$, and write $\pi$ as the blowup of an ideal $\fb\subset\cO_X$. By Theorem~\ref{thm:descentpsh}, there exists a sequence $\e_m\searrow 0$ such that $\p+\e_m\pi^\star\log|\fb|=\pi^\star\f_m$ with $\f_m\in\PSH(\theta_m)$. By Theorem~\ref{thm:suppsh}, $\f_m$ is uniquely determined, and $(\f_m)_{m\ge m_0}$ is an increasing sequence in $\PSH(\theta_{m_0})$ for any given $m_0$. By Theorem~\ref{thm:envppty}, $\f_m$ converges pointwise on $X^\div$ to $\f\in\PSH(\theta)$, which satisfies $\p=\pi^\star\f$ on $Y^\div$. The result follows. 
\end{proof}
 
\begin{proof}[Proof of Theorem~\ref{thm:unibranch}] Since the normalization morphism $\nu\colon X^\nu\to X$ is finite, $\nu^\star\om$ is ample, and hence $\cH^\dom(\nu^\star\om)$ spans $\PL(X^\nu)$, see~\eqref{equ:HdomPL}. By Lemma~\ref{lem:FSdense}, $\PL(X^\nu)$ separates the points of $(X^\nu)^\an$, and it is thus already the case of $\cH^\dom(\nu^\star\om)$. On the other hand, Lemma~\ref{lem:descentpsh2} implies that all functions in $\cH^\dom(\nu^\star\om)\subset\PSH(\nu^\star\om)$ descend to $X^\an$ (take $\om_m=(1+m^{-1})\om$). As a consequence, $\nu^\an\colon (X^\nu)^\an\to X^\an$ is injective. By the non-Archimedean GAGA principle, this implies that $\nu:X^\nu\to X$ is injective as well (see~\cite[\S3.4]{BerkBook}); this is enough to conclude that $\nu$ is a homeomorphism, and hence that $X$ is unibranch. 
\end{proof}

 As mentioned above,   in the complex analytic case the envelope property conversely holds as soon as $X$ is unibranch. It is thus natural to conjecture: 

\begin{conj}\label{conj:contenv} If $X$ is unibranch, then the envelope property holds  for all $\theta\in\Num(X)$.  
\end{conj}
  We will establish this conjecture in the important special case where $\theta$ is ample, $X$ is smooth, and $\charac k=0$ (see Theorem~\ref{thm:contenvlisse}). 
\begin{rmk}\label{rmk:envred}
  The envelope property makes sense also when $X$ is reducible. It is easy to see that $\theta\in\Num(X)$ has the envelope property iff $\theta|_Y\in\Num(Y)$ has the envelope property for every connected component $Y$ of $X$. Moreover, if $X$ is connected, then the proof of Theorem~\ref{thm:unibranch} shows that if the envelope property holds for some ample class in $\Num(X)$, then the normalization morphism $\nu\colon X^\nu\to X$ is a homeomorphism, which implies that $X$ is irreducible (and unibranch).
\end{rmk}
 
%
%
 
\subsection{Envelopes}\label{sec:pshenv}
For ample classes, the envelope property admits a useful reformulation. We continue to assume that $X$ is \textbf{irreducible} unless stated otherwise. Fix a class $\om\in\Amp(X)$.

\begin{defi}\label{defi:pshenv} The \emph{$\om$-psh envelope} of a function $\f\colon X^\an\to\R\cup\{\pm\infty\}$ is the function $\env_\om(\f)\colon X^\an\to\R\cup\{\pm\infty\}$ defined as the pointwise supremum
$$
\env_\om(\f):=\sup\left\{\p\in\PSH(\om)\mid\p\le\f\right\}.
$$
\end{defi}
Thus $\env_\om(\f)\equiv-\infty$ iff there is no $\p\in\PSH(\om)$ with $\p\le\f$. Despite the name, $\env_\om(\f)$ is not always $\om$-psh (and indeed not even usc in general).  
However, it is clear that 
\begin{itemize}
\item $\f\mapsto\env_\om(\f)$ is increasing;
\item $\env_\om(\f+c)=\env_\om(\f)+c$ for all $c\in\R$.
\end{itemize}
These properties formally imply the Lipschitz estimate
\begin{equation}\label{equ:envlip}
\sup\left|\env_\om(\f)-\env_\om(\f')\right|\le\sup|\f-\f'|. 
\end{equation}
when $\f,\f'$ are bounded.
\begin{lem}\label{lem:contenv} For any $\om\in\Amp(X)$, the following statements are equivalent:
\begin{itemize}
\item[(i)] $\om$ has the envelope property;
\item[(ii)] for any function $\f\colon X^\an\to\R\cup\{\pm\infty\}$, we have 
$$
\env_\om(\f)\equiv-\infty,\,\env_\om(\f)^\star\equiv+\infty,\text{ or }\env_\om(\f)^\star\in\PSH(\om).
$$
\item[(iii)] for any $\f\in\Cz(X)$, $\env_\om(\f)$ is continuous. 
\end{itemize}
\end{lem}
We refer to the property in~(iii) as \emph{continuity of envelopes}.
\begin{proof} First assume (i). Pick any $\f\colon X^\an\to\R\cup\{\pm\infty\}$, and suppose that the set $\cF:=\left\{\p\in\PSH(\theta)\mid\p\le\f\right\}$ is nonempty, so that $\env_\om(\f)\not\equiv-\infty$. If the functions in $\cF$ are uniformly bounded above, then $\env_\om(\f)^\star\in\PSH(\om)$, by (i). If not, then, by the definition of the Alexander--Taylor capacity we have 
$$
\env_\om(\f)(v)=\sup\left\{\p(v)\mid \p\in\cF\right\}\ge\sup\left\{\sup\p\mid\p\in\cF\right\}-\te_\om(v)=+\infty
$$ 
for all $v\in X^\div$, and hence $\env_\om(\f)^\star\equiv+\infty$, by density of $X^\div$. This proves (i)$\Rightarrow$(ii). 

Next we prove (ii)$\Rightarrow$(iii), so pick $\f\in\Cz(X)$. Then $\env_\om(\f)^\star\in\PSH(\om)$ is a competitor in the definition of $\env_\om(\f)$, and hence $\env_\om(\f)^\star\le\env_\om(\f)$. We conclude that $\env_\om(\f)^\star=\env_\om(\f)$ is usc.
We claim that $\env_\om(\f)$ is also lsc, and hence continuous. To prove the claim, it suffices to show that $\env_\om(\f)$ is a supremum of continuous functions, and for this it suffices to prove that for any $\p\in\PSH(\om)$ with $\p\le\f$ and any $\e>0$, there exists $\p'\in\Cz(X)$ with $\p\le\p'\le\f+\e$. Pick a decreasing net $\p_i\in\cH^\dom(\om)$ converging pointwise to $\p$. For each $v\in X^\an$, the set $\{\p_i<\f+\e\}$ is an open neighborhood of $v$ for $i$ large enough, by lower semicontinuity of $\f-\p_i$. By compactness of $X^\an$, it follows that $\p_i<\f+\e$ for all $i$ large enough, so we can take $\p'=\p_i$.

Finally, we prove (iii)$\Rightarrow$(i), following~\cite[Lemma 7.30]{BE}. Let $(\f_i)$ be a bounded-above family in $\PSH(\om)$, and set $\f:=\supstar_i\f_i$. Since $\f$ is usc and $X^\an$ is compact, we can find a decreasing net of continuous functions $(\p_j)$ such that $\p_j\to\f$. For each $i,j$, we have $\f_i\le\p_j$, and hence $\f_i\le \env_\om(\p_j)$, which in turn yields $\f\le \env_\om(\p_j)\le\p_j$. We have thus written $\f$ as the limit of the decreasing net of $\om$-psh functions $\env_\om(\p_j)$, which shows that $\f$ is $\om$-psh.
\end{proof}

\begin{cor}\label{cor:envusc} Assume that $\om\in\Amp(X)$ has the envelope property, and consider a usc function $\f\colon X^\an\to\R\cup\{-\infty\}$. Then: 
\begin{itemize}
\item[(i)] $\env_\om(\f)$ is $\om$-psh, or $\env_\om(\f)\equiv-\infty$; 
\item[(ii)] if $\f$ is the limit of a decreasing net $(\f_j)$ of bounded-above, usc functions, then $\env_\om(\f_j)\searrow\env_\om(\f)$.
\end{itemize}
\end{cor}
\begin{proof} By Lemma~\ref{lem:contenv}, either $\p:=\env_\om(\f)^\star$ is $\om$-psh, or $\env_\om(\f)\equiv-\infty$. Since $\env_\om(\f)\le\f$ and $\f$ is usc, we also have  $\p\le\f$. If $\p$ is $\om$-psh, then $\p\le\env_\om(\f)$, which proves (i). 

To see (ii), note that $\rho:=\lim_j\env_\om(\f_j)$ satisfies either $\rho\in\PSH(\om)$ or $\rho\equiv-\infty$, by Theorem~\ref{thm:PSH}. Furthermore, $\env_\om(\f_j)\le\f_j$ yields in the limit $\rho\le\f$, and hence $\rho\le\env_\om(\f)$ (by definition of $\env_\om(\f)$ if $\rho\in\PSH(\om)$, and trivially if $\rho\equiv-\infty$). Thus $\lim_j\env_\om(\f_j)=\rho=\env_\om(\f)$. On the other hand, $\env_\om(\f_j)\ge\env_\om(\f)$ implies $\rho\ge\env_\om(\f)$, which completes the proof of (ii). 
\end{proof}

For any function $\f\colon X^\an\to\R\cup\{\pm\infty\}$, we also introduce the pointwise envelope
$$
\qq_\om(\f):=\sup\left\{\p\in\CPSH(\om)\mid\p\le\f\right\};
$$
this is lsc and bounded below if $\f$ is bounded below, and $\equiv-\infty$ otherwise. Since each $\p\in\CPSH(\om)$ is the uniform limit of functions in $\cH^\dom(\om)$ (see Theorem~\ref{thm:PSH}), one easily checks that
$$
\qq_\om(\f)=\sup\{\p\in\cH^\dom(\om)\mid\p\le\f\}.
$$

\begin{lem}\label{lem:envlsc} Suppose $\f\colon X^\an\to\R\cup\{+\infty\}$ is bounded below, with lsc regularization $\f_\star\colon X^\an\to\R\cup\{+\infty\}$. Then:
\begin{itemize}
\item[(i)] $\qq_\om(\f)=\env_\om(\f_\star)$; 
\item[(ii)] $\qq_\om(\f)\in\Cz(X)\Longleftrightarrow\qq_\om(\f)\in\PSH(\om)$; 
\item[(iii)] if $\f$ is the pointwise limit of an increasing net $(\f_j)$ of bounded-below, lsc functions (and hence $\f$ is lsc), then 
$\env_\om(\f_j)\nearrow\env_\om(\f)$. 
\end{itemize}
\end{lem}

\begin{proof} A function $\p\in\CPSH(\om)$ satisfies $\p\le\f$ iff $\p\le\f_\star$. Thus $\qq_\om(\f)=\qq_\om(\f_\star)$, and we may therefore assume wlog that $\f$ is lsc. Trivially, $\env_\om(\f)\ge\qq_\om(\f)$. Pick $\p\in\PSH(\om)$ such that $\p\le\f$, and let $(\p_i)$ be a decreasing net in $\CPSH(\om)$ converging pointwise to $\p$. For each $\e>0$ and $v\in X^\an$, we can find $i$ such that $\{\p_i<\f+\e\}$ is an open neighborhood of $v$, by lower semicontinuity of $\f-\p_i$. By compactness of $X^\an$, it follows that $\p_i<\f+\e$ for all $i$ large enough. Thus $\p\le\p_i\le\qq_\om(\f+\e)=\qq_\om(\f)+\e$, which proves (i).

As noted above, $\qq_\om(\f)$ is lsc and bounded below. If $\qq_\om(\f)$ is $\om$-psh, then it is in particular usc, and hence $\qq_\om(\f)\in\Cz(X)$. Assume, conversely, that $\qq_\om(\f)\in \Cz(X)$. Since $\CPSH(\om)$ is stable under finite maxima, $\qq_\om(\f)$ is the pointwise limit of an increasing net $(\p_i)$ in $\CPSH(\om)$. By Dini's lemma, $\p_i\to\qq_\om(\f)$ uniformly on $X^\an$, and hence $\qq_\om(\f)\in\PSH(\om)$, by Theorem~\ref{thm:PSH}. This proves (ii). 

Consider finally a net $(\f_j)$ as in (iii). We trivially have $\lim_j\env_\om(\f_j)=\sup_j\env_\om(\f_j)\le\env_\om(\f)$. Pick $\e>0$ and $\p\in\CPSH(\om)$ such that $\p\le\f$. As above, for all $j$ large enough we have $\p\le\f_j+\e$, and hence $\p\le\env_\om(\f_j)+\e$. Thus $\env_\om(\f)\le\sup_j\env_\om(\f_j)$, and we are done. 
\end{proof}
%
%
 
\subsection{The envelope property on smooth varieties}\label{sec:contenvlisse}
The results earlier in this section would not be particularly useful unless we have examples of classes where the envelope property holds.
Arguing along the lines of~\cite[Theorem 8.5]{siminag}, we will prove: 

\begin{thm}\label{thm:contenvlisse} Assume that $X$ is smooth and connected, and that either $\charac k=0$ or $\dim X\le 2$. Then any $\theta\in\Nef(X)$ has the envelope property.  
\end{thm}
We do not know whether a class $\theta\in\Num(X)$ that is not nef has the envelope property.

\begin{cor}\label{cor:contenvcomp} Under the assumptions of Theorem~\ref{thm:contenvlisse}, the set
$$
\PSH_{\sup}(\theta)=\left\{\f\in\PSH(\theta)\mid\sup\f=0\right\}
$$
is compact. 
\end{cor}
In~\cite{siminag}, the above compactness property was established in the discretely valued case, by relying on much more involved arguments based on dual complexes and toroidal modifications.

\begin{lem}\label{lem:envline} Let $L$ be an ample line bundle on $X$, and $(\cX,\cL)$ a test configuration for $(X,L)$, with $\cL$ an honest line bundle. Denote by $\fa_m\subset\cO_\cX$ the base ideal of $m\cL$, which is a vertical ideal of $\cX$ for all $m\gg 1$. Then
$$
\env_L(\f_\cL)-\f_\cL=\sup_{m\ge 1} m^{-1}\f_{\fa_m}=\lim_{m\to\infty} m^{-1}\f_{\fa_m}
$$
pointwise on $X^\an$. Furthermore,  $\env_L(\f_\cL)$ is continuous iff $m^{-1}\f_{\fa_m}$ converges uniformly on $X^\an$. 
\end{lem}
 
As with flag ideals, see~\S\ref{sec:ideals}, the function $\f_{\fa_m}\in\PL(X)$ is defined by $\f_{\fa_m}(v)=-\sigma(v)(\fa_m)$, where $\sigma=\sigma_\cX$ denotes Gauss extension (see Remark~\ref{rmk:Gausslaurent}).

\begin{proof} Set $\f_m:=\f_\cL+m^{-1}\f_{\fa_m}$. Note first that $(\fa_m)_m$ is a graded sequence of ideals, \ie $\fa_m\cdot\fa_{m'}\subset\fa_{m+m'}$ for all $m,m'\in\N$. The sequence $m\mapsto\f_{\fa_m}$ is thus superadditive, and Fekete's lemma yields $\sup_{m\ge 1} m^{-1} \f_{\fa_m}=\lim_{m\to\infty}m^{-1} \f_{\fa_m}$. 

Next write $\cL=L_\cX+D$ with a vertical Cartier divisor $D$ on $\cX$, so that $\f_\cL=\f_D$. For each $m$, denote by $\mu_m\colon\cX_m\to\cX$ the integral closure of the blowup along $\fa_m$, and by $D_m$ the (antieffective) divisor on $\cX_m$ such that $\cO_{\cX_m}\cdot\fa_m=\cO_{\cX_m}(D_m)$; thus $\f_m=\f_D+m^{-1}\f_{D_m}$. By construction, $\cO_{\cX}(m\cL)\otimes\fa_m$ is  globally generated. This implies that
$$
\mu^\star(m\cL)+D_m=mL_{\cX_m}+\mu^\star(mD)+D_m
$$ 
is  nef, and hence $\f_m\in\PL\cap\PSH(L)$. Furthermore,  $\f_{D_m}\le 0$, so $\f_m\le\f_\cL$, and hence $\f_m\le\env_L(\f_\cL)$, which proves
$$
\sup_{m\ge 1} \f_m=\lim_{m\to\infty}\f_m\le \env_L(\f_\cL).
$$
Conversely, Lemma~\ref{lem:envlsc} implies 
$$
\env_L(\f_\cL)=\qq_L(\f_L)=\sup\left\{\p\in\cH(L)\mid\p\le\f_\cL\right\}.
$$
Pick $\p\in\cH(L)$ with $\p\le\f_\cL$. After replacing $\cX$ with a higher test configuration, we may and do assume, for the sake of notational simplicity, that $\p$ is determined by a vertical $\Q$-Cartier divisor $E$ on $\cX$. For $m$ large and divisible enough, Theorem~\ref{thm:testPL} shows that $mL_\cX+mE=m\cL+m(E-D)$ is  globally generated, and hence $\cO_\cX(m(E-D))\subset\fa_m$. This yields 
$$
m(\p-\f_\cL)=\f_{m(E-D)}\le\f_{\fa_m},
$$
and hence $\p\le\f_m$, which yields, as desired, the converse inequality
$$
\env_L(\f_\cL)=\sup\left\{\p\in\cH(L)\mid\p\le\f_\cL\right\}\le\sup_m \f_m.
$$
The final assertion is a simple consequence of Dini's lemma, using the superadditivity of $m\mapsto m\f_m$. 
\end{proof}

\begin{proof}[Proof of Theorem~\ref{thm:contenvlisse}] By Lemma~\ref{lem:envppty}, we may assume that $\theta=c_1(L)$ with $L\in\Pic(X)_\Q$ ample. Pick $\f\in\Cz(X)$. By Lemma~\ref{lem:contenv}, we need to show that $\env_L(\f)$ is continuous. Since $\PL(X)$ is dense in $\Cz(X)$ with respect to uniform convergence (see Theorem~\ref{thm:PLdense}), we may assume $\f\in\PL(X)$, by~\eqref{equ:envlip}. By Theorem~\ref{thm:testPL}, we have $\f=\f_\cL$ for a test configuration $(\cX,\cL)$ for $(X,L)$. After replacing $L$ with a multiple, we may further assume that $L$ and $\cL$ are honest line bundle.  

Using the notation of Lemma~\ref{lem:envline}, we need to show that $\f_m:=\f_\cL+m^{-1}\f_{\fa_m}$ converges uniformly to $\env_L(\f_\cL)$. Since we assume that $\charac k=0$ or $\dim X\le 2$ (and hence $\dim\cX\le 3$), we can rely on resolution of singularities and assume that $\cX$ is smooth and $\cX_0$ has simple normal crossings support. We may also assume that there exists an effective vertical $\Q$-divisor $E$ on $\cX$ such that $\cA:=\cL-E$ is  ample. 

Assume first that $\charac k=0$, and let $\fb_m$ be the multiplier ideal of the graded sequence $\fa_\bullet^m$. The inclusion $\fa_m\subset\fb_m$ is elementary, 
  and we have $\fb_{ml}\subset\fb_m^l$ for all $m,l$ by the subadditivity of 
  multiplier ideals. This implies that 
  \begin{equation*}
    m^{-1}\f_{\fb_m}\ge(ml)^{-1}\f_{\fb_{ml}}\ge(ml)^{-1}\f_{\fa_{ml}}\ge m^{-1}\f_{\fa_m}
  \end{equation*}
  for all sufficiently divisible $m$ and $l$. Letting $l\to\infty$ shows that
\begin{equation}\label{e508}
    \f_\cL+m^{-1}\f_{\fb_m}\ge \env_L(\f_\cL)\ge\f_m
  \end{equation}    
  for all sufficiently divisible $m$. 
    By the uniform global generation of multiplier ideals there exists $m_0\in\N$ such that 
   
$$\cO_\cX(m\cL+m_0\cA)\otimes\fb_m=\cO_\cX((m+m_0)\cL)\otimes\cO_\cX(-m_0E)\otimes\fb_m$$
 
  is  globally generated for all sufficiently divisible $m$.
  This implies that
  \begin{equation*}
    \f_{\fb_m}\le\f_{\fa_{m+m_0}}+C
  \end{equation*}
  for a constant $C$ independent of $m$. Combining this with~\eqref{e508} yields
   
  \begin{equation*}
    \f_m\le \env_L(\f_\cL)\le(1+\frac{m_0}{m})\f_{m+m_0}-\frac{m_0}{m}\f_\cL+\frac{C}{m},
  \end{equation*}
  for all sufficiently divisible $m$. This shows that $\env_L(\f_\cL)$ is a uniform limit of continuous functions, and hence continuous.

When $\charac k>0$, the very same argument applies with test ideals in place of multiplier ideals, see~\cite{GJKM17} for details. 
\end{proof}
\begin{rmk}
  In view of Remark~\ref{rmk:envred}, Theorem~\ref{thm:contenvlisse} is valid also when $X$ is smooth but possibly disconnected. A suitable version of Corollary~\ref{cor:contenvcomp} is also true.
\end{rmk}

%
%
%
%
  
\section{Homogeneous functions and $b$-divisors}\label{sec:homog}
In this section, we assume for simplicity that $X$ is \textbf{irreducible}. We study homogeneous $\theta$-psh functions and their relation to nef $b$-divisors. Inspired by the work of Ross and Witt Nystr\"om, we express an arbitrary $\theta$-psh function in terms of homogeneous ones, and use this to establish a version of Siu's decomposition theorem in our setting. 
%
%
\subsection{Homogeneous PL and homogeneous Fubini--Study functions}
Recall that any ideal $\fb\subset\cO_X$ determines a homogeneous, decreasing function $\log|\fb|\colon X^\an\to [-\infty,0]$.

\begin{defi}\label{defi:PLhom} The space of \emph{homogeneous PL functions on $X$} is defined as the $\Q$-vector space 
$$
\PL_\hom(X)\subset\Cz(X^\val)
$$
generated by the restriction to $X^\val$ of all functions $\log|\fb|$ attached to nonzero ideals $\fb\subset\cO_X$. 
\end{defi}

\begin{rmk} The terminology is slightly abusive, as a function $\f\in\PL_\hom(X)$ is not a PL function in the sense of Definition~\ref{defi:PL}, except if $\f=0$. Indeed, while functions in $\PL(X)\subset\Cz(X)$ are always bounded, $\f$ can only be bounded if $\f=0$, by homogeneity. 
\end{rmk}

By~\eqref{equ:psideal}, the set 
$$
\PL^+_\hom(X):=\left\{m^{-1}\log|\fb|\mid m\in\Z_{>0},\,0\ne\fb\subset\cO_X\right\}
$$
is a $\Q_+$-semivector subspace of $\PL_\hom(X)$ that is stable under finite maxima, and any function in $\PL_\hom(X)$ can be written as a difference of functions on $\PL^+_\hom(X)$.

\begin{rmk} For any $\p\in\PL^+_\hom(X)$ and $c\in\Q$, we have $\max\{\p,c\}\in\PL^+(X)$. 
\end{rmk}

 \begin{exam}\label{exam:homfrac} If $\fa$ is a nonzero fractional ideal on $X$, then setting $\log|\fa|(v):=-v(\fa)$ for $v\in X^\val$ defines a function $\log|\fa|\in\PL_\hom(X)$. Indeed, $\log|\fa|=\log|\fb'|-\log|\fb|$ with $\fb:=\{f\in\cO_X\mid f\fa\subset\cO_X\}$ and $\fb':=\fb\cdot\fa\subset\cO_X$. 
\end{exam}

\begin{rmk} As mentioned in~\S\ref{sec:Berk}, the space $X^\val$ (resp.~$X^\an$) can be reconstructed from $\PL^+_\hom(X)$ (resp.~$\PL_\hom(X)$) as its `tropical spectrum', \ie the set of all $\Q_+$-linear (resp.~$\Q$-linear) maps $\chi\colon\PL^+_\hom(X)\to\R\cup\{-\infty\}$ (resp.~$\chi\colon\PL_\hom(X)\to\R$) that commute with taking max (compare~\cite[\S 1.2]{jonmus}). In particular, $\PL_\hom(X)$ is a birational invariant of $X$, while $\PL^+_\hom(X)$ is not. 
\end{rmk}

\begin{defi} For any $L\in\Pic(X)_\Q$, we define the space of \emph{homogeneous Fubini--Study functions} for $L$ as 
$$
\cH_\hom(L):=\cH^\gf_0(L). 
$$
\end{defi} 
By Definition~\ref{defi:FS}, a homogeneous Fubini--Study function is thus a function $\f\colon X^\an\to [-\infty,0]$ of the form 
\begin{equation}\label{equ:FShom}
\f=m^{-1}\max_i\log|s_i|
\end{equation} 
for a finite set of nonzero sections $(s_i)$ of $mL$ with $m$ sufficiently divisible. Clearly, $\f=m^{-1}\log|\fb|$ for the ideal $\fb\subset\cO_X$ locally generated by the $(s_i)$, and hence 
$$
\cH_\hom(L)\subset\PL^+_\hom(X).
$$
Arguing as in Proposition~\ref{prop:FSideal}, we conversely have: 

\begin{lem}\label{lem:homFSideal} For any $L\in\Pic(X)_\Q$ and $\f\colon X^\an\to [-\infty,0]$ the following are equivalent:
\begin{itemize} 
\item[(i)] $\f\in\cH_\hom(L)$; 
\item[(ii)] $\f=m^{-1}\log|\fb|$ for a nonzero ideal $\fb\subset\cO_X$ and $m\in\Z_{>0}$ such that $mL$ is an honest line bundle and $mL\otimes\fb$ is globally generated on $X$. 
\end{itemize}
\end{lem}

\begin{cor}\label{cor:homFSideal} For any $L\in\Pic(X)_\Q$ ample, we have $\Q_+\cH_\hom(L)=\PL^+_\hom(X)$, and $\cH_\hom(L)$ spans the $\Q$-vector space $\PL_\hom(X)$. 
\end{cor}
\begin{cor}\label{cor:FSideal3} We have $\PL^+_\hom(X)=\bigcup_L\cH_\hom(L)$, where $L$ ranges over ample line bundles on $X$.
\end{cor}

\begin{exam}\label{exam:logsE} For any effective $\Q$-Cartier divisor $E$, set (as in Lemma~\ref{lem:FShomdiv}) $\log|s_E|:=m^{-1}\log|s_{mE}|$ for $m$ sufficiently divisible, where $s_{mE}\in\Hnot(X,\cO_X(mE))$ denotes the canonical section. Then $\log|s_E|\in\cH_\hom(E)$. 
\end{exam}

%
\subsection{Homogenization}\label{sec:homogpsh}
Recall that $\R_{>0}$ acts on functions $\f\colon X^\an\to\R\cup\{\pm\infty\}$ by 
$$
(t\cdot\f)(v)=t\f(t^{-1} v),
$$
so that $\f$ is homogeneous iff $t\cdot \f=\f$ for all $t$.  

\begin{defi} We define the \emph{homogenization} of a function $\f\colon X^\an\to\R\cup\{-\infty\}$ as the function $\hf\colon X^\an\to\R\cup\{-\infty\}$ such that
$$
\hf(v):=\inf_{t>0}(t\cdot\f)(v)
$$
for $v\in X^\an$. 
\end{defi}
Obviously, $\hf$ is homogeneous, $\hf\le\f$, and $\hf$ is the largest function with these properties.

\begin{exam}\label{exam:homogtriv} For any irreducible subvariety $Y\subset X$, $v_{Y,\triv}\in X^\triv$ is a fixed point under the action of $\R_{>0}$, and hence 
$$
\hf(v_{Y,\triv})=\left\{
\begin{array}{lll}
0 & \text{ if } & \f(v_{Y,\triv})\ge 0,\\
-\infty & \text{ if } & \f(v_{Y,\triv})<0.
\end{array}
\right. 
$$
\end{exam}

\begin{lem}\label{lem:homogdec} For any decreasing net $(\f_i)$ of functions $\f_i\colon X^\an\to\R\cup\{\pm\infty\}$ with pointwise limit $\f=\lim_i\f_i$, we have $\hf=\lim_i\hf_i$. 
\end{lem}
\begin{proof} Since $\f\le\f_i$, we have $\hf\le\hf_i\le\f_i$. Thus $\p:=\lim_i\hf_i$ is homogeneous and satisfies $\hf\le\p\le\f$, and hence $\p=\hf$, by the maximality property of $\hf$. 
\end{proof}
Using Lemma~\ref{lem:maxprin} and Example~\ref{exam:homogtriv}, the next result is straightforward to check. 

\begin{lem}\label{lem:genpos} If $\f\colon X^\an\to\R\cup\{-\infty\}$ is decreasing, then $\hf$ is decreasing as well. Further, $\hf\not\equiv-\infty$ iff $\sup\f\ge 0$. 
\end{lem}

For any $\theta\in\Num(X)$, we denote by 
$$
\PSH_\hom(\theta)\subset\PSH_{\sup}(\theta)
$$
the set of homogeneous $\theta$-psh functions. 

%

\begin{thm}\label{thm:homogFS} For any $\f\in\cH^\gf_\Q(L)$ with $L\in\Pic(X)_\Q$, we have $\hf\in\cH_\hom(L)$ if $\sup\f\ge 0$, and $\hf\equiv-\infty$ otherwise. 
\end{thm}

\begin{cor}\label{cor:PL+hom} For any $\f\in\PL^+(X)$, we have $\hf\in\PL^+_\hom(X)$ if $\sup\f\ge 0$, and $\hf\equiv-\infty$ otherwise. 
\end{cor}

\begin{cor}\label{cor:pshhom} Pick $\theta\in\Num(X)$ and $\f\in\PSH(\theta)$ such that $\sup\f\ge 0$. Then:
\begin{itemize}
\item[(i)] $\hf\in\PSH_\hom(\theta)$; 
\item[(ii)] we can find a decreasing net $(\p_i)$ such that $\p_i\in\cH_\hom(L_i)$ with $L_i\in\Pic(X)_\Q$ and $\lim_i c_1(L_i)=\theta$ and $\p_i\searrow\hf$;
\item[(iii)] when $\theta=c_1(L)$ with $L\in\Pic(X)_\Q$ ample, (ii) holds with $L_i=L$ for all $i$. 
\end{itemize}
\end{cor}

\begin{rmk}\label{rmk:homocount} By Theorem~\ref{thm:moncount} below, Corollary~\ref{cor:pshhom} is actually valid with (countable) sequences instead of nets.
\end{rmk}

\begin{cor}\label{cor:pphom} For any $\om\in\Amp(X)$ and $E\subset X^\an$ pluripolar, there exists $\p\in\PSH_\hom(\om)$ such that $E\subset\{\p=-\infty\}$.
\end{cor}

%


The proof of Theorem~\ref{thm:homogFS} relies on the following elementary result. 

\begin{lem}\label{lem:legendre} Pick $\la_1,\dots,\la_r\in\R$ such that $\max_i\la_i\ge 0$, and set for $x\in\R^r$ 
$$
g(x):=\inf_{t>0}\max_i\{x_i+\la_i t\}.
$$
\begin{itemize}
\item[(i)] If $\max_i\la_i=0$ (resp.~$\min_i\la_i\ge 0$) then $g(x)=\max_{\la_i=0} x_i$ (resp.~$g(x)=\max_i x_i$). 
\item[(ii)] In general, $g(x)=\max_{w\in W}\langle w,x\rangle$ for a finite subset $W$ of the $(r-1)$-simplex
$$
\sigma=\left\{w\in\R_+^r\mid\sum_i w_i=1\right\}. 
$$
\item[(iii)] If the $\la_i$ are rational, then we can chose $W\subset\sigma\cap\Q_+^r$.  
\end{itemize}
\end{lem}

\begin{proof} Note first that $\max_i\la_i\ge 0$ implies that $g(x)\ge\min_i x_i$ is finite for all $x\in\R^r$. The proof of (i) is straightforward. To see (ii), note that the epigraph of $g$ is the projection to $\R^r$ of the epigraph of the convex, homogeneous PL function $f(x,t):=\max_i\{x_i+\la_i t\}$ on $\R^r\times\R_+$. This implies that $g$ is a convex, homogeneous PL function as well, and hence 
$$
g(x)=\sup_{w\in W}\langle y,x\rangle.
$$
with $W$ the (finite) set of vertices of the Newton polyhedron
$$
P:=\left\{w\in\R^r\mid \langle w,x\rangle\le g(x)\right\}
$$
Finally, $g$ is increasing in each variable, and 
$$
g(x_1+c,\dots,x_r+c)=g(x)+c
$$
for all $c\in\R$. This implies that $P\subset\sigma$, which proves (ii). If the $\la_i$ are rational, then $f$ is $\Q$-PL\@. Thus $g$ is $\Q$-PL as well, and $P$ is then a rational polyhedron, whose set $W$ of vertices is thus rational. This proves (iii). 
\end{proof}

\begin{proof}[Proof of Theorem~\ref{thm:homogFS}]
   By Lemma~\ref{lem:genpos} we may assume $\sup\f\ge0$.   Write
$$
\f=m^{-1}\max_i\{\log|s_i|+\la_i\}
$$ 
with $m\in\Z_{>0}$, $s_1,\dots,s_r\in\Hnot(X,mL)$ and $\la_i\in\Q$, and note that $\max_i\la_i\ge\sup\f\ge 0$. Then
$$
\hf=m^{-1}\inf_{t>0}\max_i\left\{\log|s_i|+t\la_i\right\},
$$
and Lemma~\ref{lem:legendre} thus yields a finite subset $W\subset\sigma\cap\Q_+^r$ such that 
$$
\hf=m^{-1}\max_{w\in W}\left\{\sum_i w_i\log|s_i|\right\}.
$$
Pick $b\in\Z_{>0}$ such that $b w\in\N^r$ for all $w\in W$. Then 
$$
\sum_i w_i\log|s_i|=b^{-1}\log|s_w|
$$
with $s_w:=\prod_i s_i^{b w_i}\in\Hnot(X,mbL)$, and hence $\hf=(mb)^{-1}\max_{w\in W}\log|s_w|$. Thus $\hf\in\cH_\hom(L)$, and we are done.
\end{proof}

\begin{proof}[Proof of Corollary~\ref{cor:PL+hom}] By Corollaries~\ref{cor:FSideal} and~\ref{cor:homFSideal}, for any ample line bundle $L$ we have $\PL^+(X)=\Q_+\cH(L)$ and $\PL^+_\hom(X)=\Q_+\cH_\hom(L)$. We conclude using Theorem~\ref{thm:homogFS}. 
\end{proof}

\begin{proof}[Proof of Corollary~\ref{cor:pshhom}] By Lemma~\ref{lem:genpos}, (i) follows from (ii). By Theorem~\ref{thm:pshample}~(i), we can write $\f$ as the limit of a decreasing net $\f_i\in\cH^\gf_\Q(L_i)$ with $L_i\in\Pic(X)_\Q$ and $c_1(L_i)\to\theta$. For each $i$, we have $\sup\f_i\ge\sup\f\ge 0$, and hence $\hf_i\in\cH_\hom(L_i)$, by Theorem~\ref{thm:homogFS}. By Lemma~\ref{lem:homogdec}, we have $\hf_i\searrow\hf$, which proves (ii). Finally, if $\theta=c_1(L)$ with $L\in\Pic(X)_\Q$ ample, Theorem~\ref{thm:pshample}~(iii) shows we can take $L_i=L$ in the above argument, and (iii) follows. 
\end{proof}

\begin{proof}[Proof of Corollary~\ref{cor:pphom}] By Lemma~\ref{lem:ppindep}, we can find $\f\in\PSH(\om)$ such that $E\subset\{\f=-\infty\}$. After adding a constant to $\f$, we may assume $\sup\f=0$. By Corollary~\ref{cor:pshhom}, we then have $\hf\in\PSH_\hom(\om)$, and $E\subset\{\hf=-\infty\}$ since $\hf\le\f$. 
\end{proof}

%


Relying on Corollary~\ref{cor:pshhom}, we now establish the following version of Siu's decomposition theorem (\cf~\cite[III.8.16]{DemBook}) for homogeneous psh functions---see Theorem~\ref{thm:Siu} below for a statement in the general case. 

\begin{thm}\label{thm:pshdiv} Assume that $X$ is normal. Pick $\theta\in\Num(X)$ and $E$ an effective $\Q$-Cartier divisor. For any $\p\in\PSH_\hom(\theta)$, we then have 
$$
\p\le\log|s_E|\Longleftrightarrow\p-\log|s_E|\in\PSH_\hom(\theta-E).
$$
\end{thm}
See Example~\ref{exam:logsE} for the notation. 

\begin{proof} By Corollary~\ref{cor:pshhom}, we can write $\p$ as the limit of a decreasing net $\p_i\in\cH_\hom(L_i)$ with $L_i\in\Pic(X)_\Q$ and $c_1(L_i)\to\theta$. Since $K:=\{\log|s_E|=-1\}$ is compact in $X^\an$ and $\sup_K\p\le -1$, we can find $t_i\in\Q_{>0}$ such that $\lim_i t_i=1$ and $\sup_K\p_i\le-t_i$, by Dini's lemma. By homogeneity, this implies $\p_i\le t_i\log|s_E|$ on $\R_{>0}K=\{\log|s_E|>-\infty\}$, and hence on $X^\an$, since $\{\log|s_E|>-\infty\}\supset X^\div$ is dense and $\p_i$ and $\log|s_E|$ are continuous (or by Theorem~\ref{thm:suppsh}). 

By Lemma~\ref{lem:FShomdiv}, $\p'_i:=\p_i-t_i\log|s_E|=\p_i-\log|s_{t_i E}|$ lies in $\cH_\hom(L_i-t_i E)$. Since $(\p_i)$ is a decreasing net, $(\p'_i)$ is decreasing on $X^\div\subset\{\log|s_E|>-\infty\}$, and hence on $X^\an$, by Theorem~\ref{thm:suppsh}. Further, $\p_i(v_{\triv,\a})=0$ for all $\a$ and all $i$. Since $c_1(L_i)-t_i E\to \theta-E$, Theorem~\ref{thm:psh12} shows that $\p':=\lim_i\p'_i$ is $(\theta-E)$-psh, with $\p=\p'+\log|s_E|$ on $X^\div$, and hence on $X^\an$, by Corollary~\ref{cor:Hausdorff}. 
\end{proof}

\begin{cor}\label{cor:homcompact} If all classes in $\Amp(X)$ have the envelope property, then $\PSH_\hom(\theta)$ is compact for any $\theta\in\Num(X)$.
\end{cor}
Recall that the assumption holds if $X$ is smooth and $\charac k=0$, see Corollary~\ref{cor:contenvcomp}.

\begin{proof} Pick an effective Cartier divisor $E$ such that $\theta':=\theta+E\in\Num(X)$ is ample. Since $\theta'$ has the envelope property, $\PSH_{\sup}(\theta')$ is compact (see Theorem~\ref{thm:envppty}). Now pick a net $(\p_i)$ in $\PSH_\hom(\theta)$. For each $i$, we have $\p'_i:=\p_i+\log|s_E|\in\PSH_\hom(\theta')$. After passing to a subnet, we may thus assume $\p'_i\to\p'\in\PSH_\hom(\theta')$. Since $\p'_i\le\log|s_E|$, we get in the limit $\p'\le\log|s_E|$ on $X^\div$, and hence on $X^\an$, by Theorem~\ref{thm:suppsh}. By Theorem~\ref{thm:pshdiv}, we get $\p:=\p'-\log|s_E|\in\PSH_\hom(\theta)$. Further, $\p_i\to\p$ on $X^\div$, hence in $\PSH_\hom(\theta)$, which proves that the latter space is compact.
\end{proof}

%
%
\subsection{The homogeneous decomposition of a psh function}\label{sec:homogdecomp}
For any function $\f\colon X^\an\to\R\cup\{-\infty\}$ and $\la\in\R$, we set 
\begin{equation}\label{equ:hfla}
\hf^\la:=\widehat{\f-\la}=\inf_{t>0}\left\{t\cdot\f-t\la\right\}. 
\end{equation}
Thus $\hf^\la$ is largest homogeneous function such that $\hf^\la+\la\le\f$.   By Corollary~\ref{cor:pshhom}, we have, for any $\f\in\PSH(\theta)$ and $\la\in\R$,
$$
\left\{
\begin{array}{l}
\la\le\inf\f\Longrightarrow\hf^\la=0,\\
\la\le\sup\f\Longrightarrow\hf^\la\in\PSH_\hom(\theta),\\
\la>\sup\f\Longrightarrow\hf^\la\equiv-\infty. 
\end{array}
\right. 
$$

\begin{lem} For any $\theta\in\Num(X)$ and $\f\in\PSH(\theta)$, $(\hf^\la)_{\la<\sup\f}$ is a concave and decreasing family of functions in $\PSH_\hom(\theta)$, in the sense that $\lambda\mapsto\hf^\la(v)$ is concave and decreasing for all $v\in X^\an$. Moreover, the map
$$
(-\infty,\sup\f]\ni\la\mapsto\hf^\la\in\PSH_\hom(\theta)
$$ 
is continuous.
\end{lem}

\begin{proof} The extremal characterization of $\hf^\la$ shows that it is a decreasing function of $\la$, and concavity follows directly from~\eqref{equ:hfla}. For any $v\in X^\div$, we need to show that $\la\mapsto\hf^\la(v)$ is continuous on $(-\infty,\sup\f]$. By concavity, it is continuous on $(-\infty,\sup\f)$. Since $\la\mapsto\hf^\la$ is decreasing and vanishes at $v_\triv$, $\p:=\inf_{\la<\sup\f}\hf^\la$ is $\theta$-psh (see Theorem~\ref{thm:psh12}), and also clearly homogeneous. Further, $\hf^{\max}:=\hf^{\sup\f}\le\p$, and it remains to show that equality holds. For each $\la<\sup\f$, we have $\p\le\hf^\la\le\f-\la$, and hence $\p\le\f-\sup\f$. By the extremal property, we infer $\p\le\hf^{\max}$, which concludes the proof.
\end{proof}

\begin{thm}\label{thm:homogpsh} For any $\theta\in\Num(X)$ and $\f\in\PSH(\theta)$, we have
\begin{equation}\label{equ:homdec}
\f=\sup_{\la<\sup\f}\{\hf^\la+\la\}
\end{equation}
pointwise on $X^\an$. Conversely, if $(\p_\la)_{\la<\sup\f}$ is a concave family of functions in $\PSH_\hom(\theta)$ such that $\f=\sup_{\la<\sup\f}\{\p_\la+\la\}$ holds on $X^\div$, then $\p_\la=\hf^\la$ for all $\la<\sup\f$. 
\end{thm}

\begin{lem}\label{lem:actionconvex} For each $v\in X^\an$ and $\f\in\PSH(\theta)$, $t\mapsto(t\cdot\f)(v)$ is convex on $\R_{>0}$, and is decreasing if $\f\le 0$. 
\end{lem}
\begin{proof} First assume $\f\in\cH^\gf_\R(L)$ with $L\in\Pic(X)_\Q$, and write $\f=m^{-1}\max_i\{\log|s_i|+\la_i\}$ for a finite set $(s_i)$ of nonzero sections of $mL$ and $\la_i\in\R$. Then $t\cdot\f=m^{-1}\max_i\{\log|s_i|+t\la_i\}$, which is a convex function of $t$, decreasing when $\sup\f=\max\la_i\le0$.  The general case follows easily, see Definition~\ref{defi:psh}. 
\end{proof}

\begin{proof}[Proof of Theorem~\ref{thm:homogpsh}] For any $v\in X^\an$, $t\mapsto(t\cdot\f)(v)$ in convex on $\R_{>0}$, with (convex) Legendre transform 
$$
\la\mapsto\sup_{t>0}\{t\la-(t\cdot\f)(v)\}=-\hf^\la(v).
$$
By Legendre duality, we thus have 
$$
(t\cdot\f)(v)=\sup_{\la\in\R}\{\hf^\la(v)+t\la\}
$$
for all $t>0$. For $t=1$, this is precisely~\eqref{equ:homdec}. Conversely, assume we are given a concave family $(\p_\la)_{\la<\sup\f}$ in $\PSH_\hom(\theta)$ such that $\f(v)=\sup_{\la<\sup\f}\{\p_\la(v)+\la\}$ for all $v\in X^\div$. Then $(t\cdot\f)(v)=\sup_{\la<\sup\f}\{\p_\la(v)+t\la\}$ for all $t>0$, and hence 
$$
\p_\la(v)=\inf_{t>0}\{(t\cdot\f)(v)-t\la\}=\hf^\la(v)
$$
for all $\la<\sup\f$, again by Legendre duality. This shows $\p_\la=\hf^\la$ on $X^\div$, and hence on $X^\an$, by Corollary~\ref{cor:Hausdorff}. 
\end{proof}

\begin{rmk} If $\om$ is ample and $\f\in\cE^1(\om)$, then one can view $(t\cdot\f)_{t>0}$ as a geodesic ray in the space $\cE^1(\om)$ (see~\cite{nakstab2,Reb20}), and the above result is then in line with the Legendre transform approach to geodesic rays pioneered in~\cite{RWN}. 
\end{rmk}

%

In analogy with Lemma~\ref{lem:envlsc}, we also prove: 

\begin{lem}\label{lem:hatphi} Assume $\f\in\CPSH(L)$ with $L\in\Pic(X)_\Q$ ample, and pick $\la<\sup\f$. Then
$$
\hf^\la=\sup\left\{\p\in\cH_\hom(L)\mid\p\le\f-\la\right\}
$$
pointwise on $X^\div$. 
\end{lem}

\begin{proof} By Theorem~\ref{thm:homogpsh}~(ii), the right-hand side $\tau$ satisfies $\tau\le\hf^\la$. Pick $\e\in(0,\sup\f-\la)$. By Theorem~\ref{thm:PSH}~(vi), we can find $\p\in\cH(L)$ such that 
$$
\f-(\la+\e)\le\p\le\f-\la.
$$
Since $0\le\sup\f-(\la+\e)\le\sup\p$, Theorem~\ref{thm:homogFS} yields $\widehat\p\in\cH_\hom(L)$. Now $\widehat\p\le\p\le\f-\la$, and hence $\widehat\p\le\tau$. On the other hand, $\hf^{\la+\e}=\widehat{\f-(\la+\e)}\le \widehat\p$, and hence 
$$
\hf^{\la+\e}\le\tau\le\hf^\la.
$$
By Theorem~\ref{thm:homogpsh}~(iii), we have $\lim_{\e\to 0}\hf^{\la+\e}=\hf^\la$ pointwise on $X^\div$, and we are done.
\end{proof}

\begin{exam} Assume $\f=\f_\fa$ for a flag ideal $\fa$, \ie $\f=\max_\la\{\log|\fa_\la|+\la\}$ for a decreasing sequence $(\fa_\la)_{\la\in\Z}$ of ideals on $X$. For any $\la\in\Z$ with $\la\le\sup\f$, Theorem~\ref{thm:homogFS} yields $\hf^\la\in\PL^+_\hom(X)$, \ie $\hf^\la=m^{-1}\log|\fb|$ for an ideal $\fb\subset\cO_X$ and $m\in\Z_{>0}$, and it is natural to wonder whether in fact $\hf^\la=\log|\fa_\la|$. We do have $\hf^\la\ge\log|\fa_\la|$, with equality if $\la\le\inf\f$ or $\la=\sup\f$, by Lemma~\ref{lem:legendre}~(i), but equality however fails in general. Indeed, the concavity of $\la\mapsto\hf^\la$ would otherwise imply that $\fa_{\la-1}\cdot\fa_{\la+1}$ is contained in the integral closure of $\fa_\la^2$, which need not be satisfied. 
\end{exam}

Pick $\theta\in\Num(X)$ and $\f\in\PSH(\theta)$. Besides $\hf^0=\hf$, the case of $\hf^\la$ with $\la=\sup\f$ also plays a special role, and we set
$$
\hf^{\max}:=\hf^{\sup\f}\in\PSH_\hom(\theta).
$$
Thus $\hf^{\max}\le\f-\sup\f$, and $\hf^{\max}$ is the largest homogeneous function with this property. In particular, 
$$
\hf^{\max}=0\Longleftrightarrow\f\text{ constant}. 
$$
Note also that $\hf^{\max}$ is invariant under addition of a constant to $\f$. As the next result shows, $\hf^{\max}$ can be understood as the `G\^ateaux differential' of $\f$ at $v_\triv$. 

\begin{lem}\label{lem:hfmax} For any $\theta\in\Num(X)$, $\f\in\PSH(\theta)$ and $v\in X^\an$, we have 
$$
\hf^{\max}(v)=\lim_{t\to 0_+}\frac{\f(tv)-\f(v_\triv)}{t}. 
$$
\end{lem} 
Recall from Proposition~\ref{prop:pshconvex}  that $t\mapsto\f(tv)$ is convex in $\R_{>0}$ for any $v\in X^\an$, and that $\sup\f=\f(v_\triv)$. 

\begin{proof} Replacing $\f$ with $\f-\sup\f$, we may assume $\sup\f=\f(v_\triv)=0$. Then $t\mapsto t\cdot\f$ is decreasing, by Lemma~\ref{lem:actionconvex}, and hence
$$
\hf^{\max}(v)=\hf(v)=\lim_{t\to+\infty}t\f(t^{-1} v)=\lim_{s\to 0_+} s^{-1}\f(s v), 
$$
which proves the result. 
\end{proof}

\begin{exam} Pick $L\in\Pic(X)_\Q$, $\f\in\cH^\val_\R(L)$, and write $\f=m^{-1}\max_i\{\log|s_i|+\la_i\}$ as in~\eqref{equ:FS}. Then Lemma~\ref{lem:legendre}~(i) yields $\hf^{\max}=m^{-1}\max_i\log|s_i|$. 
\end{exam}

%

Thanks to Theorem~\ref{thm:homogpsh}, we can now extend Theorem~\ref{thm:pshdiv} to arbitrary psh functions, yielding the following analogue of Siu's decomposition theorem. 

\begin{thm}\label{thm:Siu} Assume that $X$ is normal. Pick $\theta\in\Num(X)$, an effective $\Q$-Cartier divisor $E$, and assume $\theta\in\Num(X)$ has the envelope property. For any $\f\in\PSH(\theta+E)$, we then have: 
$$
\f\le\log|s_E|+O(1)\Longleftrightarrow\f-\log|s_E|\in\PSH(\theta).
$$ 
\end{thm}
Recall that we expect that any $\theta\in\Num(X)$ has the envelope property, see Conjecture~\ref{conj:contenv}. Also note the shift by $E$ compared to the notation of Theorem~\ref{thm:pshdiv}
\begin{proof} For each $\la\le\sup\f$, $\hf^\la\in\PSH_\hom(\theta+E)$ satisfies $\hf^\la\le\f-\la\le\log|s_E|+O(1)$, and hence $\hf^\la\le\log|s_E|$, by homogeneity. By Theorem~\ref{thm:pshdiv}, we thus have $\hf^\la=\p_\la+\log|s_E|$ for a unique $\p_\la\in\PSH_\hom(\theta)$, and hence $\f=\tau+\log|s_E|$ pointwise on $X^\an$ with $\tau:=\sup_{\la\le\sup\f}(\p_\la+\la)$. The envelope property guarantees that $\f':=\tau^\star$ lies in $\PSH(\theta)$, and it satisfies $\f'=\f-\log|s_E|$ by Theorem~\ref{thm:divnegl}. 
\end{proof}

%
%
\subsection{Homogeneous PL functions and Cartier $b$-divisors}\label{sec:homognefbdiv}
  
In this section, $X$ is assumed to be \textbf{normal}. A \emph{model} $Y$ of $X$ is a normal projective variety together with a birational map $\pi\colon Y\to X$. Recall that (see for instance~\cite[\S 1]{BdFF}): 
\begin{itemize}
\item a \emph{(rational) $b$-divisor over $X$} is a collection $B=(B_Y)$ of $\Q$-Weil divisors $B_Y\in\Zun(Y)_\Q$ on all models of $X$, compatible under push-forward as cycles, \ie an element of the projective limit
$$
\Zunb(X):=\varprojlim_Y \Zun(Y)_\Q
$$
\item a $b$-divisor $B\in\Zunb(X)$ is ($\Q$-)\emph{Cartier} if there exists a model $Y$, called a \emph{determination} of $B$, such that $B_Y\in\Car(Y)\subset\Zun(Y)$ and $B_{Y'}$ is the pullback of $B_Y$ for all higher birational models $Y'$.
\end{itemize}
The $\Q$-linear subspace of Cartier $b$-divisors $$\Carb(X)\subset\Zunb(X)$$ can thus be identified with the direct limit $\varinjlim_Y \Car(Y)_\Q$.

A Cartier b-divisor is \emph{relatively semiample} if some (equivalently any) determination $B_Y\in\Zun(Y)$ is relatively semiample for the morphism $Y\to X$. We will write $$\Carbplus(X)\subset\Carb(X)$$ for the set of $B\in\Carb(X)$ that are relatively semiample and antieffective, that is, $B\le 0$. Note that $\Carb(X)$ is a birational invariant of $X$, but $\Carbplus(X)$ is not.

\medskip 

Any $b$-divisor $B\in\Zunb(X)$ determines a homogeneous function $\p_B\colon X^\div\to\Q$, such that $\p_B(\ord_E)=\ord_E(B_Y)$ for any model $Y$ and prime divisor $E\subset Y$. The map $B\mapsto\p_B$ sets up a 1--1 correspondence between $\Zunb(X)$ and the space of homogeneous functions $\p\colon X^\div\to\Q$ such that 
$\p(\ord_E)$ is nonzero for only finitely many prime divisors $E\subset X$. Under this correspondence, a net $B_i$ converges to $B$ in the inverse limit topology on $\Zunb(X)$ iff $\p_{B_i}\to\p_B$ pointwise on $X^\div$: we then simply say that $B_i\to B$ pointwise.

If $B\in\Carb(X)$, then $\p_B$ admits a (unique) continuous extension $\p_B\colon X^\val\to\R$, defined by $\p_B(v)=v(B_Y)$ for any determination $Y$ of $B$ and $v\in X^\val\simeq Y^\val$. This yields an injection $\Carb(X)\hookrightarrow\Cz(X^\val)$, and the next result specifies its image.

\begin{thm}\label{thm:bdivPL} The map $B\mapsto\p_B$ induces isomorphisms $$\Carb(X)\simto\PL_\hom(X)\quad\text{and}\quad\Carbplus(X)\simto\PL_\hom^+(X).$$
\end{thm}
\begin{proof}
  We start by proving the second isomorphism.
  First consider $B\in\Carbplus(X)$, and pick a determination $\pi\colon Y\to X$ of $B$. Thus $B_Y$ is $\pi$-semiample and antieffective, so for $m$ sufficiently divisible, the generically trivial ideal $\cO_Y(mB_Y)\subset\cO_Y$ is $\pi$-globally generated, \ie $\cO_Y(mB_Y)=\fb\cdot\cO_Y$ where $\fb:=\pi_\star\cO_Y(mB_Y)$ is a generically trivial ideal. This yields $\p_B=m^{-1}\log|\fb|\in\PL_\hom^+(X)$.

  Conversely, if $\p\in\PL_\hom^+(X)$, then $\p=m^{-1}\log|\fb|$, where $m\ge1$ and $\fb\subset\cO_X$ is a generically trivial ideal. Let $\pi\colon Y\to X$ denotes the normalized blowup of $\fb$, with exceptional divisor $E$. Then $\p=\p_B$ with $B\in\Carb(X)$ determined on $Y$ by $B_Y=-m^{-1}E\le0$.

  Next we prove the first isomorphism above. If $\p\in\PL_\hom(X)$, then we can write $\p=\p_1-\p_2$, with $\p_i\in\PL_\hom^+(X)$. By what precedes, $\p_i=\p_{B_i}$ for $B_i\in\Carbplus(X)$, so $\p=\p_B$, where $B=B_1-B_2\in\Carb(X)$.

  Conversely, any $B\in\Carb(X)$ can be written $B=B_1-B_2$, where $B_i\in\Carb(X)$ are relatively semiample. It would take some work to also arrange $B_i\le0$, so instead we argue as follows. Pick a common determination  $\pi\colon Y\to X$ of the $B_i$. Thus $B_{i,Y}$ is $\pi$-semiample, so for $m$ sufficiently divisible, the generically trivial fractional ideal $\cO_Y(mB_{i,Y})\subset\cO_Y$ is $\pi$-globally generated, \ie $\cO_Y(mB_{i,Y})=\fb_i\cdot\cO_Y$ where $\fb_i:=\pi_\star\cO_Y(mB_{i,Y})$ is a generically trivial fractional ideal. It follows that $\p_B=m^{-1}(\log|\fb_1|-\log|\fb_2|)\in\PL_\hom(X)$, see Example~\ref{exam:homfrac}, which completes the proof.
\end{proof}

\begin{rmk}\label{rmk:hfmax} By Lemma~\ref{lem:hfmax}, $\f\mapsto\hf^{\max}$ can be extended to a linear map $\PL(X)\to\PL_\hom(X)$. Pick $\f\in\PL(X)$, and write $\f=\f_D$ with $D\in\VCar(\cX)_\Q$ for an integrally closed test configuration $\cX$ that dominates $\cX_\triv$ (see Theorem~\ref{thm:VCarPL}). The strict transform $Y$ of $(\cX_\triv)_0=X\times\{0\}$ is an irreducible component of $\cX_0$ that induces the trivial valuation $v_\triv$ on $X$. Thus $\ord_Y(D)=\sup\f$, and one checks that $\hf^{\max}=\p_B$ with $B\in\Carb(X)$ determined on $Y$ by the $\Q$-Cartier divisor $B_Y:=\left(D-\ord_Y(D)\cX_0\right)|_Y$. 
\end{rmk}

%
%
%
%
%

Now consider $L\in\Pic(X)_\Q$. We will describe the image of $\cH_\hom(L)\subset\PL_\hom(X)$ in $\Carb(X)$ under the isomorphism in Theorem~\ref{thm:bdivPL}. To this end, we say that a Cartier b-divisor $B\in\Carb(X)$ is \emph{semiample} if $B_Y$ is.

\begin{lem}\label{lem:FShomsemi} For any $L\in\Pic(X)_\Q$ and $B\in\Carb(X)$, we have 
$$
\p_B\in\cH_\hom(L)\Longleftrightarrow B\le 0\text{ and }L+B\text{ semiample},
$$
where the last condition means that $\pi^\star L+B_Y$ is semiample for some (or, equivalently, any) determination $\pi\colon Y\to X$ of $B$. 
\end{lem} 

\begin{proof} The direct implication follows from Lemma~\ref{lem:homFSideal}. Conversely, assume $B\le 0$ and that there exists a determination $\pi\colon Y\to X$ such that $\pi^\star L+B_Y$ is semiample. Pick $m$ such that $m(\pi^\star L+B_Y)$ is a globally generated line bundle. Since $X$ is normal and $B_Y\le 0$, we have $\pi_\star\cO_Y(mB_Y)\subset\pi_\star\cO_Y=\cO_X$, by Zariski's main theorem. Further, the ideal $\fb_m\subset\cO_X$ locally generated by the image of 
$$
\Hnot\left(Y,m(\pi^\star L+B_Y)\right)\simeq\Hnot\left(X,mL\otimes\pi_\star\cO_Y(-mE)\right)\hookrightarrow\Hnot(X,mL)
$$
satisfies $\fb_m\cdot\cO_Y=\cO_Y(mB_Y)$, and hence $\p_B=m^{-1}\log|\fb_m|$, which lies in $\cH_\hom(L)$, by Lemma~\ref{lem:homFSideal}.
\end{proof}

Similarly, if $\theta\in\Num(X)$ we say that $\theta+B$ is nef if $\pi^\star \theta+B_Y\in\Nef(Y)$ for some (or, equivalently, any) determination $\pi\colon Y\to X$ of $B$. The proof of the next result, which describes the image of $\PL_\hom(X)\cap\PSH(\theta)$ in $\Carb(X)$, is more involved. 

\begin{thm}\label{thm:hPLpsh} Assume that all classes in $\Amp(X)$ have the envelope property, see~\S\ref{sec:envprop}. For any $\theta\in\Num(X)$ and $B\in\Carb(X)$, we then have
$$
\p_B\in\PSH_\hom(\theta)\Longleftrightarrow B\le 0\text{ and }\theta+B\text{ nef}.
$$
\end{thm}

\begin{rmk} The condition $\p_B\in\PSH_\hom(\theta)$ means here that $\p_B|_{X^\div}=\p|_{X^\div}$ for some $\p\in\PSH_\hom(\theta)$, which is necessarily unique by Corollary~\ref{cor:Hausdorff}. If $\charac k=0$, Remark~\ref{rmk:Hausdorff} applies, and shows that $\p_B$ and $\p$ in fact coincide on $X^\val$.
\end{rmk}

The proof relies on the following `homogeneous version' of Corollary~\ref{cor:goodman}. 

\begin{lem}\label{lem:goodmanhom} Pick $\theta\in\Num(X)$. Let $(B_i)$ be a decreasing net in $\Carb(X)$ that converges pointwise to $B\in\Carb(X)$. Assume also that we are given $\theta_i\in\Num(X)$ with $B_i\le 0$, $\theta_i+B_i$ is nef, and $\theta_i\to\theta\in\Num(X)$. Then $\theta+B$ is nef. 
\end{lem}
\begin{proof} Let $\pi\colon Y\to X$ be a determination of $B$, and $C\subset Y$ an irreducible curve. Following~\cite[Proposition 8]{Good}, denote by $\mu\colon Z\to Y$ the normalized blowup of $C$, with exceptional divisor $F$, and pick an ample line bundle $A$ on $Z$. Then $\mu_\star(F\cdot A^{n-2})=aC$ with $a\in\Q_{>0}$ as numerical classes on $Y$, and it will thus be enough to show $(\mu^\star\pi^\star\theta+B_Z)\cdot (F\cdot A^{n-2})\ge 0$. Since $\theta_i+B_i$ is nef, the projection formula yields 
$(\mu^\star\pi^\star\theta_i+B_{i,Z})\cdot (F\cdot A^{n-2})\ge 0$, and we will be done if we show $(B_{i,Z}\cdot F\cdot A^{n-2})\to (B_Z\cdot F\cdot A^{n-2})$. Denote by $(E_\b)$ the finite set of prime components of $B_Z$. Since $B\le B_i\le 0$, $B_{i,Z}$ is also supported in the $E_\b$'s. Thus 
$$
(B_{i,Z}\cdot F\cdot A^{n-2})=\sum_\b \ord_{E_\b}(B_i) (E_\b\cdot F\cdot A^{n-2})\to
(B_Z\cdot F\cdot A^{n-2})=\sum_\b\ord_{E_\b}(B) (E_\b\cdot F\cdot A^{n-2}),
$$
and we are done.
\end{proof}

\begin{proof}[Proof of Theorem~\ref{thm:hPLpsh}] Assume $\p_B\in\PSH_\hom(\theta)$. Then $\p_B\le 0$, and hence $B\le 0$. By Corollary~\ref{cor:pshhom}, we can write $\p_B$ as the pointwise limit of a decreasing net $\p_i\in\cH_\hom(L_i)$ with $L_i\in\Pic(X)_\Q$ and $c_1(L_i)\to\theta$. Denote by $B_i$ the $\Q$-Cartier $b$-divisor associated to $\p_i$. By Lemma~\ref{lem:FShomsemi}, $L_i+B_i$ is semiample. Thus $c_1(L_i)+B_i$ is nef, and Lemma~\ref{lem:goodmanhom} shows that $\theta+B$ is nef. 

Assume, conversely, $B\le 0$ and $\theta+B$ nef. We claim that it is enough to prove $\p_B\in\PSH(\theta)$ when $\theta$ is further ample. To see this, pick an effective Cartier divisor $E$ on $X$ such that $\om:=\theta+E$ is ample. Denote by $\bar E\in\Carb(X)$ the Cartier $b$-divisor determined on $X$ by $E$, and set $B':=B-\bar E$. Any determination $\pi\colon Y\to X$ of $B$ is also a determination of $B'$, and $\pi^\star\om+B'=\pi^\star\theta+B$ in $\Num(Y)$, which shows that $\om+B'$ is nef. On the other hand, for any $v\in X^\val$ we have $\p_{\bar E}(v)=v(E)=-\log|s_E|(v)$, and hence $\p_{B'}=\p_B+\log|s_E|$. Assuming the result in the ample case, we get $\p_{B}+\log|s_E|=\p_{B'}\in\PSH_\hom(\om)=\PSH_\hom(\theta+E)$. By Theorem~\ref{thm:pshdiv}, this implies $\p_B\in\PSH_\hom(\theta)$, which proves the claim. 

From now on we assume $\theta=\om\in\Amp(X)$. Since $\pi$ is isomorphic to the blowup of an ideal of $X$, we can choose a $\pi$-ample divisor $H\le 0$ on $Y$. Pick $L\in\Pic(X)_\Q$ such that $\om':=c_1(L)-\om$ is ample. Since $H$ is $\pi$-ample, $\pi^\star\om'+\e H$ is ample for all $\e\in\Q_{>0}$ small enough. Now $\pi^\star \om+B_Y$ is nef, and hence 
$$
(\pi^\star\om+B_Y)+(\pi^\star\om'+\e H)=\pi^\star L+(B_Y+\e H)
$$
is ample. By Lemma~\ref{lem:FShomsemi}, it follows that $\p_\e:=\p_B+\e\p_H\in\cH_\hom(L)\subset\PSH_\hom(L)$. Since $c_1(L)=\om+\om'$ is ample, it has the envelope property. As $\p_B\colon X^\an\to\R\cup\{-\infty\}$ is continuous and $\p_\e\nearrow\p_B$ on $X^\div$ as $\e\searrow 0$, Lemma~\ref{lem:envppty} thus yields $\p_B\in\PSH_\hom(L)$. As this holds for all $L\in\Pic(X)_\Q$ as above, Theorem~\ref{thm:psh12} yields, as desired, $\p_B\in\PSH_\hom(\om)$.  

\end{proof}

%
%
%
%
%
%

%
 \subsection{Nef $b$-divisors and homogeneous psh functions}
We finally consider general nef $b$-divisors, and describe these in terms of homogeneous psh functions. In this section, we assume that $X$ is \textbf{smooth} and $k$ has \textbf{characteristic zero}. 

If $Y$, $Y'$ are smooth models of $X$ and $Y'$ dominates $Y$, then the corresponding birational morphism $\mu\colon Y'\to Y$ induces a linear map $\mu_\star\colon\Num(Y')\to\Num(Y)$. The space of \emph{$b$-divisor classes} is defined as
$$
\Numb(X):=\varprojlim_Y\Num(Y)
$$
with $Y$ running over all smooth models of $X$, endowed with the projective limit topology (see~\cite{diskant,BdFF,DF}).  Each $B\in\Zunb(X)_\R$ determines a class $[B]\in\Numb(X)$. 
\begin{defi} We say that a (real) $b$-divisor $B\in\Zunb(X)_\R$ is \emph{exceptional} if $B_X=0$. 
\end{defi}
The map $B\mapsto\p_B$ sets up a linear isomorphism between the space of exceptional $b$-divisors and the space $\cE$ of all homogeneous functions $\p\colon X^\div\to\R$ such that $\p(\ord_E)=0$ for all prime divisors $E\subset X$. We equip the $\R$-vector space $\cE$ with the topology of pointwise convergence. 
\medskip

By the negativity lemma, the map $B\mapsto [B]$ is injective on exceptional $b$-divisors. For any $\a\in\Numb(X)$, we can thus find a unique exceptional $b$-divisor $B_\a\in\Zunb(X)_\R$ such that $[(B_\a)_Y]=\a_Y-\pi^\star\a_X$
for all smooth models $\pi\colon Y\to X$ (see~\cite[Lemma 1.11]{BdFF}). We denote by $\p_\a:=\p_{B_\a}\in\cE$ the corresponding homogeneous function. As a consequence of~\cite[Lemma 1.12]{BdFF}, we can now state: 

\begin{lem}\label{lem:classcont} The map $\a\mapsto(\a_X,\p_\a)$ defines a topological vector space isomorphism 
\begin{equation}\label{equ:bdiviso}
\Numb(X)\simto\Num(X)\times\cE. 
\end{equation}
\end{lem}

As recalled above, one says that a Cartier $b$-divisor $B\in\Carb(X)$ is nef if $B_Y$ is nef for some (or, equivalently, any) determination $Y$. The set
$$
\Nefb(X)\subset\Numb(X)
$$
of \emph{nef $b$-divisor classes} is defined as the closure of the set of nef Cartier $b$-divisor classes. Thus a class $\a\in\Numb(X)$ is nef iff there exists a net $(B_i)$ of nef Cartier $b$-divisors such that $[B_{i,Y}]\to\a_Y$ for all smooth models $Y\to X$. Since $\Nef(Y)$ is a closed convex cone for each model $Y$, it is not hard to see that $\Nefb(X)$ is a closed convex cone. 

By~\cite[Lemma 2.12]{BdFF}, a class $\a\in\Numb(X)$ is nef iff, for each smooth model $Y$, $\a_Y\in\Num(Y)$ is nef in codimension one.
A typical example is provided by the `positive part' in the divisorial Zariski decomposition of a Cartier $b$-divisor~\cite{Bou,Naka,diskant}. 

As we next show, nef $b$-divisor classes admit a precise description as homogeneous psh functions. 

\begin{thm}\label{thm:bnefpsh} The isomorphism~\eqref{equ:bdiviso} maps $\Nefb(X)$ onto the set of pairs $(\theta,\p)$ with $\theta\in\Num(X)$ and $\p\in\PSH_\hom(\theta)$ such that $\p(\ord_E)=0$ for all prime divisors $E\subset X$.  
\end{thm}
As a consequence of Theorem~\ref{thm:bnefpsh}, we recover the monotone approximation result of~\cite{DF}. 

\begin{cor}\label{cor:DF} For any $\a\in\Nefb(X)$, there exists a decreasing net $(B_i)$ of nef Cartier $b$-divisors such that $[B_i]\to\a$ in $\Numb(X)$. 
\end{cor}
By Remark~\ref{rmk:homocount}, the result is actually valid with a sequence instead of a net, as in~\cite[Theorem A]{DF}. 

\begin{rmk} Both~\cite[Theorem A]{DF} and Corollary~\ref{cor:DF} ultimately rely on the multiplier ideals technique that goes back to~\cite{ELS}, and was already used in a similar manner in~\cite[Theorem 8.5]{siminag}. However, the proof of~\cite{DF} is much more direct, and the main interest of the present discussion lies rather in Theorem~\ref{thm:bnefpsh}. 
\end{rmk}

\begin{proof}[Proof of Theorem~\ref{thm:bnefpsh}] Assume first that $\a=[B]$ for a nef Cartier $b$-divisor $B$, and pick a determination $\pi\colon Y\to X$ of $B$. Then $\p_\a\in\PL_\hom(X)$ is the function associated to the $\Q$-Cartier $b$-divisor determined on $Y$ by the $\pi$-exceptional divisor $E_Y:=B_Y-\pi^\star B_X$, and 
$$
\pi^\star\theta+[E_Y]=[B_Y]
$$
is nef. By Theorem~\ref{thm:hPLpsh}, we thus have $\p_\a\in\PSH_\hom(\a_X)$. 

Consider next any $\a\in\Nefb(X)$. By assumption, there exists a net $(B_i)$ of nef $b$-divisors such that $\a_i:=[B_i]\to\a$ in $\Numb(X)$, \ie $\theta_i:=(\a_i)_X\to\theta$ and $\p_{\a_i}\to\p_\a$ pointwise on $X^\div$ (see Lemma~\ref{lem:classcont}). For any ample class $\om$, we have $\theta+\om-\theta_i\in\Amp(X)$ for $i$ large enough, and hence $\p_{\a_i}\in\PSH_\hom(\theta_i)\subset\PSH_\hom(\theta+\om)$. By Corollary~\ref{cor:homcompact}, the latter space is compact, and we may thus assume, after passing to a subnet, that $\p_{\a_i}\to\p\in\PSH_\hom(\theta+\om)$. Since $\p_{\a_i}\to\p_\a$ on $X^\div$, it follows that $\p_\a\in\PSH_\hom(\theta+\om)$, and hence $\p_\a\in\PSH_\hom(\theta)$, since this holds for all $\om\in\Amp(X)$ (see Theorem~\ref{thm:psh12}). 

Conversely, pick $\p\in\PSH_\hom(\theta)$ such that $\p(\ord_E)=0$ for all prime divisors $E\subset X$. Its restriction to $X^\div$ determines an element $B\in\Zunb(X)_\R$ such that $\p_B=\p$. We define $\a\in\Numb(X)$ by setting $\a_Y:=\pi^\star\theta+[B_Y]$ for all smooth models $\pi\colon Y\to X$, and we claim that $\a$ is nef. 

By Corollary~\ref{cor:pshhom}, we can write $\p$ as the pointwise limit of a decreasing net $(\p_i)$ such that $\p_i\in\cH_\hom(L_i)$ with $L_i\in\Pic(X)_\Q$ and $\lim_i c_1(L_i)=\theta$. For each $i$, we have $\p\le\p_i\le 0$, and hence $\p_i(\ord_E)=0$ for all prime divisors $E\subset X$, which means that the Cartier $b$-divisor $B_i\in\Carb(X)$ such that $\p_{B_i}=\p_i$ satisfies $(B_i)_X=0$. Choose a $\Q$-Cartier divisor $D_i$ on $X$ representing the linear equivalence class $L_i\in\Pic(X)_\Q$, and set $B'_i:=\bar D_i+B_i$. By Lemma~\ref{lem:FShomsemi}, $B'_i$ is semiample, and hence nef. Since $\p_i\to\p$ and $[B_i]_X=c_1(L_i)\to\theta=\a_X$ in $\Num(X)$, Lemma~\ref{lem:classcont} yields $\a_i\to\a$ in $\Numb(X)$, which proves, as desired, that $\a$ is nef. 
\end{proof}

\begin{proof}[Proof of Corollary~\ref{cor:DF}] Pick $\a\in\Nefb(X)$. Write $\a_X=[D]$ for an $\R$-divisor $D$ on $X$, and pick an effective, ample $\R$-divisor $H$ on $X$ such that $L:=D+H$ is an ample $\Q$-divisor. By Theorem~\ref{thm:bnefpsh}, we have $\p_\a-\p_H\in\PSH(L)$. By Corollary~\ref{cor:pshhom}, $\p_\a-\p_H$ is thus the pointwise limit on $X^\div$ of a decreasing net $(\p_i)$ in $\cH_\hom(L)$. Write $\p_i=\p_{C_i}$ with $C_i\in\Carb(X)$. By Lemma~\ref{lem:FShomsemi}, $B_i:=C_i+L$ is a decreasing net of nef Cartier $b$-divisors, and we have by construction $[B_i]\to\a$. 
\end{proof}

%
%
%
%
 \section{Functions of finite energy and mixed Monge--Amp\`ere measures} 
In this section, $X$ denotes a projective variety of dimension $n$,  with irreducible components $X_\a$.   We extend the energy pairing to arbitrary psh functions, and use this to define functions of finite energy, and the mixed Monge--Amp\`ere operator thereon.

%
%
\subsection{Extending the energy pairing}\label{sec:extint}
In~\S\ref{sec:enpairing}, the energy pairing 
$$
(\theta_0,\f_0)\inter(\theta_n,\f_n)\in\R
$$
was defined for pairs $(\theta_i,\f_i)\in\Num(X)\times\PL_\R$. When the $\f_i$ are $\theta_i$-psh, this is an increasing function of the $\f_i$, by Lemma~\ref{lem:enmono}.   On the other hand, any $\om$-psh function with $\om\in\Amp(X)$ can be written as a decreasing limit of functions in $\cH^\dom(\om)\subset\PL\cap\PSH(\om)$ (see Theorem~\ref{thm:pshample}). This allows to extend the energy pairing by monotonicity, as follows.  

\begin{thm}\label{thm:extint} For each $(n+1)$-tuple $\om_0,\dots,\om_n\in\Amp(X)$, the energy pairing 
$$
(\f_0,\dots,\f_n)\mapsto(\om_0,\f_0)\inter(\om_n,\f_n)
$$
admits a unique extension to a map $\prod_{i=0}^n\PSH(\om_i)\to\R\cup\{-\infty\}$ that is
\begin{itemize}
\item upper semicontinuous;
\item increasing in each variable.
\end{itemize}
Furthermore,  this map is continuous along decreasing nets, and satisfies
\begin{equation}\label{equ:extint}
 (\om_0,\f_0)\inter(\om_n,\f_n)=\inf_{\p_i\in\cH^\dom(\om_i),\,\p_i\ge\f_i} (\om_0,\p_0)\inter(\om_n,\p_n)
\end{equation}
for all $(\f_i)\in\prod_i\PSH(\om_i)$. 
\end{thm}
\begin{proof} Any map that is both usc and increasing is automatically continuous along decreasing nets. Since all functions in $\PSH(\om_i)$ are limits of decreasing nets in $\cH^\dom(\om_i)\subset\PL(X)\cap\PSH(\om_i)$, uniqueness is clear. To prove existence, it suffices to show that~\eqref{equ:extint} has the required properties. 

Monotonicity is obvious. To prove upper semicontinuity, suppose that $\f_i\in\PSH(\om_i)$, $0\le i\le n$, and $t\in\R$ satisfy 
$$
(\om_0,\f_0)\inter(\om_n,\f_n)<t,
$$ 
By definition, we can choose $\p_i\in\cH^\dom(\om_i)$ and $0<\e\ll 1$ such that $\f_i\le\p_i$ and 
$$
 (\om_0,\p_0)\inter(\om_n,\p_n)< t-\e.
$$
Set 
$$
C:=\sum_i (\om_0\inter \widehat{\om_i}\inter \om_n)>0.
$$
By Proposition~\ref{prop:suppsh}, $\f\mapsto\sup(\f-\p_i)$ is continuous on $\PSH(\om_i)$, and we can thus find an open neighborhood $U_i$ of $\f_i$ in $\PSH(\om_i)$ such that $\f'_i\le\p_i+\e C^{-1}$ for all $\f'_i\in U_i$. By monotonicity of the energy pairing and Proposition~\ref{prop:intpairing}, this gives 
$$
 (\om_0,\f'_0)\inter(\om_n,\f'_n)\le (\om_0,\p_0+\e C^{-1})\inter(\om_n,\p_n+\e C^{-1})
$$
$$
= (\om_0,\p_0)\inter(\om_n,\p_n)+\e <t
$$
for all $\f'_i\in U_i$, which proves upper semicontinuity. 
\end{proof}

\begin{prop}\label{prop:intpsh} The energy pairing $ (\om_0,\f_0)\inter(\om_n,\f_n)$ is symmetric and $\R_{>0}$-linear with respect to each variable $(\om_i,\f_i)$ with $\om_i\in\Amp(X)$ and $\f_i\in\PSH(\om_i)$. Furthermore,  
 
\begin{equation}\label{equ:addcst}
(\om_0,\f_0+c_0)\inter(\om_n,\f_n+c_n)=(\om_0,\f_0)\inter(\om_n,\f_n)+\sum_{i=0}^n c_i(\om_0\inter\widehat{\om_i}\inter\om_n)_X
\end{equation}
for all $c_i\in\R$,  
\begin{equation}\label{equ:inthom}
 (\om_0,t\cdot\f_0)\inter(\om_n,t\cdot\f_n)=t (\om_0,\f_0)\inter(\om_n,\f_n)
\end{equation}
for all $t\in\R_{>0}$, 
\begin{equation}\label{equ:intsumcomp}
(\om_0,\f_0)\inter(\om_n,\f_n)=\sum_{\dim X_\a=n}(\om_0,\f_0)|_{X_\a}\inter(\om_n,\f_n)|_{X_\a},
\end{equation}
and
\begin{equation}\label{equ:intnorm}
(\nu^\star\om_0,\nu^\star\f_0)\inter(\nu^\star\om_n,\nu^\star\f_n)=(\om_0,\f_0)\inter(\om_n,\f_n)
\end{equation}
with $\nu\colon X^\nu\to X$ the normalization morphism.
\end{prop}

\begin{proof} Using approximation by decreasing nets in $\cH^\dom(\om_i)$, everything is clear from Proposition~\ref{prop:intpairing}, except that this only provides a proof of~\eqref{equ:inthom} for $t\in\Q_{>0}$. To get the general case, we may replace $\f_i$ with $\f_i-\sup\f_i$ and assume $\f_i\le 0$ for all $i$. Write a given $t\in\R_{>0}$ as the limit of an increasing sequence $t_m\in\Q_{>0}$. For each $i$, $(t_m\cdot\f_i)_m$ is a decreasing sequence in $\PSH(\om_i)$ that converges to $t\cdot\f_i$, and hence 
$$
\lim_{m\to\infty} (\om_0,t_m\cdot\f_0)\inter(\om_n,t_m\cdot\f_n)=(\om_0,t\cdot\f_0)\inter(\om_n,t\cdot\f_n).
$$
The result follows. 
\end{proof}

We also record the following useful monotonicity property with respect to $\om_i$. 

\begin{prop}\label{prop:enmono}
  For $i=0,\dots,n$, pick $\om_i,\om'_i\in\Amp(X)$ such that $\om'_i\ge\om_i$, and $\f_i\in\PSH(\om_i)\subset\PSH(\om'_i)$. Then
\begin{itemize}
\item[(i)] if $\f_i\le 0$ for all $i$, then $(\om'_0,\f_0)\inter(\om'_n,\f_n)\le(\om_0,\f_0)\inter(\om_n,\f_n)\le 0$; 
\item[(ii)] if $\f_i\ge 0$ for all $i$, then $(\om'_0,\f_0)\inter(\om'_n,\f_n)\ge(\om_0,\f_0)\inter(\om_n,\f_n)\ge 0$. 
\end{itemize}
\end{prop} 
\begin{proof} By monotone approximation, we may assume $\f_i\in\PL\cap\PSH(\om_i)$. By assumption, $\theta_i:=\om'_i-\om_i$ is nef. The result thus follows by expanding out
$$
(\om'_0,\f_0)\inter(\om'_n,\f_n)=((\om_0,\f_0)+(\theta_0,0))\inter((\om_n,\f_n)+(\theta_n,0))
$$
and applying Corollary~\ref{cor:enmono}.
\end{proof}

  For any $\theta\in\Num(X)$, recall that $\CPSH(\theta)=\Cz(X)\cap\PSH(\theta)\subset\cE^\infty(\theta)$ respectively denote the sets of continuous and bounded $\theta$-psh functions.   We conclude this section with a general continuity result for the energy pairing involving such functions. It is an analogue of the Chern--Levine--Nirenberg inequality in the complex case. 

\begin{thm}\label{thm:intcont} For all $\om_0,\dots,\om_n\in\Amp(X)$, the energy pairing
$$
(\f_0,\dots,\f_n)\mapsto (\om_0,\f_0)\inter(\om_n,\f_n)
$$
is finite-valued and continuous on $\PSH(\om_0)\times\prod_{i=1}^n\CPSH(\om_i)$ with its natural topology  (pointwise convergence on $X^\div$ for the first factor, and uniform convergence for the other ones). 
\end{thm}
\begin{lem}\label{lem:intcont} Pick   $(\f_0,\dots,\f_n)\in\PSH(\om_0)\times\prod_{i=1}^n\cE^\infty(\om_i)$.   Then: 
\begin{itemize}
\item[(i)] $(\om_0,\f_0)\inter(\om_n,\f_n)\in\R$ is finite; 
\item[(ii)] for any   $(\f'_1,\dots,\f'_n)\in\prod_{i=1}^n\cE^\infty(\om_i)$   we have 
\begin{equation}\label{equ:intlip}
\left| (\om_0,\f_0)\cdot(\om_1,\f_1)\inter(\om_n,\f_n)- (\om_0,\f_0)\cdot(\om_1,\f'_1)\inter(\om_n,\f'_n)\right|\le C\sum_{i=1}^n\sup|\f_i-\f'_i|
\end{equation}
with $C:=\max_{1\le i\le n}(\om_0\inter\widehat{\om_i}\inter\om_n)$. 
\end{itemize}
\end{lem}
 
We emphasize that the estimate is uniform with respect to $\f_0$.  
\begin{proof} For any $t\in\R$ such that $\f_i\ge t$ for $i\ge 1$,   \eqref{equ:addcst} yields  
\begin{equation}\label{equ:CLN}
(\om_0,\f_0)\inter(\om_n,\f_n)\ge(\om_0,\f_0)\cdot(\om_1,0)\inter(\om_n,0)+ntC. 
\end{equation}
To prove~(i) we may therefore assume $\f_i=0$ for $1\le i\le n$.
For $\f_0\in\PL_\R$,~\eqref{equ:intsumcomp} and~\eqref{equ:intnumb} yield
$$
(\om_0,\f_0)\cdot(\om_1,0)\inter(\om_n,0)=\sum_{\dim X_\a=n}\f_0(v_{\triv,\a})(\om_1\inter\om_n)|_{X_\a}
$$
By monotone approximation, this remains true for all $\f_0\in\PSH(\om_0)$, proving (i).

For (ii), we may assume that $\f_i=\f'_i$ for $i>1$, by multilinearity and symmetry. Then $\f_1\le\f'_1+\sup|\f_1-\f'_1|$ implies
$$
(\om_0,\f_0)\inter(\om_n,\f_n)\le(\om_0,\f_0)\cdot(\om_1,\f'_1)\cdot(\om_2,\f_2)\inter(\om_n,\f_n)+C\sup|\f_1-\f'_1|,
$$
  using again~\eqref{equ:addcst}.   The result follows.
\end{proof}

\begin{proof}[Proof of Theorem~\ref{thm:intcont}] Assume first that all $\f_i$ are PL, and pick an integrally closed  test configuration $\cX$ dominating $\cX_\triv$ and $B_i\in\VCar(\cX)_\Q$ such that $\f_i=\f_{B_i}$ for $i=1,\dots,n$. By~\eqref{equ:intnumb}, we have 
\begin{multline}\label{equ:MAdirac}
(\om_0,\f_0)\inter(\om_n,\f_n)=(\om_0,0)\cdot(\om_1,\f_1)\inter(\om_n,\f_n)\\
+\sum_E b_E\,\f_0(v_E)(\om_{1,\cX}+B_1 )|_E\inter(\om_{n,\cX}+B_n )|_E
\end{multline}
with $E$ ranging over all irreducible components of $\cX_0$. By continuity along decreasing nets, this remains true for any $\f_0\in\PSH(\om_0)$, and shows that $\f_0\mapsto(\om_0,\f_0)\cdot(\om_1,\f_1)\inter(\om_n,\f_n)$ is continuous on $\PSH(\om_0)$ when $\f_i$ is PL for $i\ge 1$. 

Assume next $\f_i\in\CPSH(\om_i)$ for $i=1,\dots,n$. For each $i\ge 1$, we can choose a sequence $(\f_{ij})_j$ in $\cH^\dom(\om_i)\subset\PL\cap\PSH(\om_i)$ converging uniformly to $\f_i$. By~\eqref{equ:intlip}, the sequence of continuous functions on $\PSH(\om_0)$
$$
\f_0\mapsto (\om_0,\f_0)\cdot(\om_1,\f_{1j})\inter(\om_n,\f_{nj})
$$
converges uniformly to
$$
\f_0\mapsto (\om_0,\f_0)\cdot(\om_1,\f_1)\inter(\om_n,\f_n),
$$
which is therefore continuous as well. Finally consider an arbitrary convergent net 
$$
(\f_{0j},\f_{1j},\dots,\f_{nj})\to(\f_0,\f_1,\dots,\f_n)
$$
in $\PSH(\om_0)\times\prod_{i=1}^n\CPSH(\om_i)$, so that $\f_{0j}\to\f_0$ in $\PSH(\om_0)$ (\ie pointwise on $X^\div$) and $\f_{ij}\to\f_i$ uniformly for $i\ge 1$. Write
\begin{multline*}
(\om_0,\f_{0j})\inter(\om_n,\f_{nj})- (\om_0,\f_0)\inter(\om_n,\f_n)\\
=\left[(\om_0,\f_{0j})\cdot(\om_1,\f_{1j})\inter(\om_n,\f_{nj})-(\om_0,\f_{0j})\cdot(\om_1,\f_1)\inter(\om_n,\f_n)\right]\\
+\left[(\om_0,\f_{0j})\cdot(\om_1,\f_1)\inter(\om_n,\f_n)-(\om_0,\f_0)\cdot(\om_1,\f_1)\inter(\om_n,\f_n)\right]. 
\end{multline*}
By~\eqref{equ:intlip}, the first term tends to $0$ as $j\to\infty$, while the previous step of the proof shows that the second term tends to $0$ as well, so we are done. 
\end{proof}

%
%
%
%
\subsection{Functions of finite energy}\label{sec:conv}
\begin{defi}\label{defi:E} For each $\om\in\Amp(X)$, we define the \emph{Monge--Amp\`ere energy functional} 
$$
\en_\om\colon\PSH(\om)\to\R\cup\{-\infty\}
$$
by setting, for $\f\in\PSH(\om)$,
\begin{equation}\label{equ:defiE}
\en_\om(\f):=\frac{(\om,\f)^{n+1}}{(n+1)(\om^n)}.
\end{equation}
We say that $\f$ has \emph{finite energy} if $\en_\om(\f)>-\infty$.
\end{defi} 
We denote by 
$$
\cE^1(\om)\subset\PSH(\om)
$$ 
the set of $\om$-psh functions of finite energy.  

\begin{prop}\label{prop:E} For each $\om\in\Amp(X)$, the functional $\en_\om\colon\PSH(\om)\to\R\cup\{-\infty\}$ satisfies:
\begin{itemize}
\item[(i)] $\en_\om$ is increasing, concave, usc, and continuous along decreasing nets; 
\item[(ii)] $\en_\om(\f+c)=\en_\om(\f)+c$ for $\f\in\PSH(\om)$ and $c\in\R$;
\item[(iii)] for each $\f\in\PSH(\om)$ and $t\in\R_{>0}$ we have $\en_\om(t\cdot\f)=t\en_\om(\f)$ and $\en_{t\om}(t\f)=t\en_{\om}(\f)$;
\item[(iv)] for each $\f\in\PSH(\om)$ we have
\begin{equation}\label{equ:ensumcomp}
\en_\om(\f)=\sum_\a c_\a\en_{\om|_{X_\a}}(\f|_{X_\a^\an})=\en_{\nu^\star\om}(\nu^\star\f)
\end{equation}
with $c_\a:=(\om^n)_{X_\a}/(\om^n)_X$ and $\nu\colon X^\nu\to X$ the normalization morphism. 
\end{itemize}
\end{prop}
Note that $c_\a>0$ iff $X_\a$ is top-dimensional, and $\sum_\a c_\a=1$. 
\begin{proof} Concavity follows from Theorem~\ref{thm:enconc}, by monotone approximation. The rest of (i)--(ii) and the first half of (iii) are consequences of Proposition~\ref{prop:intpsh}. Pick $\f\in\PSH(\om)$. For any $t\in\R_{>0}$, $t\f\in\PSH(t\om)$ satisfies 
$$
\en_{t\om}(t\f)=\frac{(t\om,t\f)^{n+1}}{(n+1)((t\om)^n)}=t\frac{(\om,\f)^{n+1}}{(n+1)(\om^n)}=t\en_\om(\f),
$$
which concludes the proof of (iii). Finally, (iv) is a consequence of~\eqref{equ:intsumcomp}. 
\end{proof}

\begin{thm}\label{thm:E1} Given ample classes $\om,\om',\om_0,\dots,\om_n\in\Amp(X)$, we have:
  \begin{itemize}
  \item[(i)] $\cE^1(\om)$ is a convex subset of $\PSH(\om)$, and contains all bounded $\om$-psh functions;
   \item[(ii)] $\cE^1(\om)$ is stable under addition of a constant and the scaling action of $\R_{>0}$;
   \item[(iii)] if $\f\in\cE^1(\om)$, then any $\p\in\PSH(\om)$ such that $\p\ge\f$ is also in 
    $\cE^1(\om)$;
     \item[(iv)] if $\om\le\om'$ then $\cE^1(\om)\subset\cE^1(\om')$; 
     \item[(v)] $\cE^1(t\om)=t\,\cE^1(\om)$ for all $t\in\R_{>0}$, and 
     $$
     \cE^1(\om)+\cE^1(\om')\subset\cE^1(\om+\om');
     $$
   \item[(vi)] if $\f_i\in\cE^1(\om_i)$, $i=0,\dots,n$, then $(\om_0,\f_0)\inter(\om_n,\f_n)$ is finite;
   \item[(vii)] if $\f\in\PSH(\om)$ then 
   \begin{align*}
  &  \f\in\cE^1(\om)\Longleftrightarrow\nu^\star\f\in\cE^1(\nu^\star\om) \\
 &  \Longleftrightarrow\f|_{X_\a^\an}\in\cE^1(\om|_{X_\a})\ \text{for all top-dimensional components $X_\a$}.
   \end{align*} 
    \end{itemize}
\end{thm}

\begin{lem}\label{lem:ennef} Pick $\om,\om'\in\Amp(X)$ and $t\ge 1$ such that $\om\le\om'\le t\om$. For all nonpositive $\f\in\PSH(\om)\subset\PSH(\om')$ we have
$$
0\ge (\om,\f)^{n+1}\ge(\om',\f)^{n+1}\ge t^n(\om,\f)^{n+1}.
$$ 
\end{lem}

\begin{proof} By Proposition~\ref{prop:enmono} we have 
$$
0\ge(\om,\f)^{n+1}\ge(\om',\f)^{n+1}\ge (t\om,\f)^{n+1}=t^{n+1}(\om,t^{-1}\f)^{n+1}. 
$$
Since $\f\in\PSH(\om)$ and $t^{-1}\in[0,1]$, concavity of the energy yields
$(\om,t^{-1}\f)^{n+1}\ge t^{-1}(\om,\f)^{n+1}$, and the result follows. 
\end{proof}

\begin{lem}\label{lem:ensum}
Pick $r\ge 1$, $\om_0,\dots,\om_r\in\Amp(X)$, $0\ge\f_i\in\PSH(\om_i)$ for $i=0,\dots,r$. Assume also given $t\ge 1$ such that $\om_i\le t\om_j$ for all $i,j$. Then 
\begin{equation}\label{equ:ensum}
\left(\sum_i\om_i,\sum_i\f_i\right)^{n+1}\ge C_{r,n} t^{rn}\sum_i(\om_i,\f_i)^{n+1}
\end{equation}
with $C_{r,n}:=\left(2^r r!\right)^n$. 
\end{lem}

\begin{proof} Assume first $r=1$. Set $\tilde\om:=\frac{t}{1+t}(\om_0+\om_1)$, and observe that $\om_0\le\tilde\om\le t\om_0$ and $\om_1\le\tilde\om\le t\om_1$. Thus
\begin{multline*}
\left(\om_0+\om_1,\f_0+\f_1\right)^{n+1}=2^{n+1}\left(\tfrac 12(\om_0+\om_1),\tfrac 12(\f_0+\f_1)\right)^{n+1}\\
  \ge 2^{n+1}\left(\tilde\om,\tfrac 12(\f_0+\f_1)\right)^{n+1}
\ge 2^n\left(\left(\tilde\om,\f_0\right)^{n+1}+\left(\tilde\om,\f_1\right)^{n+1}\right)\\
\ge (2t)^n\left((\om_0,\f_0)^{n+1}+(\om_1,\f_1)^{n+1}\right), 
\end{multline*}
where the first inequality holds by Proposition~\ref{prop:enmono}, the second one by concavity of $\p\mapsto(\tilde\om,\p)^{n+1}$, and the third one by Lemma~\ref{lem:ennef}. 

Assume now $r\ge 2$, and set $\om'_0:=\sum_{i>0}\om_i$, $\f'_0:=\sum_{i>0}\f_i$. Since $t^{-1}\om_0\le\om'_0\le rt\om_0$, the first part of the proof yields 
$$
\left(\sum_i\om_i,\sum_i\f_i\right)^{n+1}=(\om_0+\om'_0,\f_0+\f'_0)^{n+1}\ge (2rt)^n\left((\om_0,\f_0)^{n+1}+(\om'_0,\f'_0)^{n+1}\right)
$$
By induction, we have on the other hand
$$
(\om'_0,\f'_0)^{n+1}\ge C_{r-1,n} t^{(r-1)n} \sum_{i>0}(\om_i,\f_i)^{n+1},
$$
The result follows, since $(\om_0,\f_0)^{n+1}\le 0$ and 
$$
C_{r,n} t^{rn}=(2rt)^n C_{r-1,n}\ge(2rt)^n. 
$$
\end{proof}

\begin{cor}\label{cor:enmixte} With the notation of Lemma~\ref{lem:ensum} we have 
\begin{equation}\label{equ:bounden}
 (\om_0,\f_0)\inter(\om_n,\f_n)\gtrsim t^{n^2} \min_{0\le i\le n}(\om_i,\f_i)^{n+1}
\end{equation}
for all nonpositive $\f_i\in\PSH(\om_i)$. 
\end{cor}
\begin{proof} Expanding out $(\om_0+\dots+\om_n,\f_0+\dots+\f_n)^{n+1}$ yields
$$
 (\om_0+\dots+\om_n,\f_0+\dots+\f_n)^{n+1}\le (n+1)! (\om_0,\f_0)\inter(\om_n,\f_n),
$$
and we conclude by Lemma~\ref{lem:ensum} with $r=n$. 
\end{proof}

\begin{proof}[Proof of Theorem~\ref{thm:E1}] Properties~(i)--(iii) follow from Proposition~\ref{prop:E}, while (iv),~(v) and~(vi), respectively, follow from Lemmas~\ref{lem:ennef},~\ref{lem:ensum} and Corollary~\ref{cor:enmixte}. Finally, (vii) follows from~\eqref{equ:ensumcomp}. 

\end{proof}

It will be convenient to extend the energy pairing to a multilinear map.
 Recall that any $\om$-psh function is finite-valued on $X^\lin$, see Corollary~\ref{cor:pshonXlin}. We now introduce 
\begin{defi} We define a \emph{function of finite energy} as a function $\f\colon X^\lin\to\R$ of the form $\f=\f^+-\f^-$ with $\f^\pm\in\cE^1(\om)$ for some $\om\in\Amp(X)$. 
\end{defi}

By Theorem~\ref{thm:E1}, $\bigcup_{\om\in\Amp(X)}\cE^1(\om)$ forms a convex cone in the $\R$-linear space of all functions $\f\colon X^\lin\to\R$, and the set 
$$
\vec\cE^1=\vec\cE^1(X)
$$
of functions of finite energy thus forms an $\R$-vector space, which contains $\PL(X)_\R$, by \eqref{equ:HdomPL}. Note further that $\vec\cE^1$ is generated by $\cE^1(\om)$ for any given $\om\in\Amp(X)$, by Theorem~\ref{thm:E1}~(iv), (v). 

\begin{lem}\label{lem:pairinglin} For each $\om\in\Amp(X)$, we have $\cE^1(\om)=\vec\cE^1\cap\PSH(\om)$.
\end{lem}
In other words, $\cE^1(\om)$ consists precisely of functions that are $\om$-psh and of finite energy, so that the chosen terminology is consistent. 
\begin{proof} The inclusion $\cE^1(\om)\subset\vec\cE^1\cap\PSH(\om)$ is clear. Conversely, pick $\f\in\vec\cE^1\cap\PSH(\om)$, and write $\f=\f^+-\f^-$ with $\f^\pm\in\cE^1(\om')$ for some $\om'\in\Amp(X)$. After adding constants, we may assume $\f,\f^\pm\le 0$. Then $\left((\om,\f)+(\om',\f^-)\right)^{n+1}\le (\om,\f)^{n+1}$, by expanding out the left-hand side. By Theorem~\ref{thm:E1}~(iv), we infer
$$
-\infty<(\om+\om',\f^+)^{n+1}=\left((\om,\f)+(\om',\f^-)\right)^{n+1}\le (\om,\f)^{n+1}. 
$$
Thus $\f\in\cE^1(\om)$, which proves the result.
\end{proof}

\begin{thm}\label{thm:pairinglin} There exists a unique multilinear symmetric pairing 
$$
\left(\Num(X)\times\vec\cE^1\right)^{n+1}\ni ((\theta_0,\f_0),\dots,(\theta_n,\f_n))\mapsto(\theta_0,\f_0)\inter(\theta_n,\f_n)\in\R
$$
that is compatible with the one defined in Theorem~\ref{thm:extint} for tuples $(\theta_i,\f_i)\in\Num(X)\times\vec\cE^1$ such that $\theta_i\in\Amp(X)$ and $\f_i\in\cE^1(\theta_i)=\vec\cE^1\cap\PSH(\theta_i)$.   Furthermore:
\begin{itemize}
\item[(i)] this pairing is compatible with the one defined in \S\ref{sec:enpairing} for tuples in $\Num(X)\times\PL(X)_\R$;
\item[(ii)] \eqref{equ:inthom}, \eqref{equ:intsumcomp} and~\eqref{equ:intnorm} remain valid on $\Num(X)\times\vec\cE^1$; 
\item[(iii)] for any $\om\in\Amp(X)$ and $\theta_0,\dots,\theta_n\in\Num(X)$, the map 
$$
(\f_0,\dots,\f_n)\mapsto(\theta_0,\f_0)\inter(\theta_n,\f_n)
$$ 
is continuous along decreasing nets in $\cE^1(\om)$.
\end{itemize}
\end{thm}

\begin{proof} Consider the vector space $V:=\Num(X)\times\vec\cE^1$. By Theorem~\ref{thm:E1}~(v), 
$$
C:=\bigcup_{\om\in\Amp(X)}\{\om\}\times\cE^1(\om).
$$
is a convex cone in $V$, and the energy pairing is an $\R_{>0}$-multilinear pairing on $C$, by Proposition~\ref{prop:intpsh}. Since $\Amp(X)$ spans $\Num(X)$, it is straightforward to see, using Theorem~\ref{thm:E1}~(iv), that $C$ spans $V$. It is now a simple general fact that the pairing on $C$ uniquely extends to an $\R$-multilinear pairing on $V$. 

The first compatibility assertion follows from Lemma~\ref{lem:pairinglin}. The pairing on $C$ restricts to the one from~\S\ref{sec:enpairing} on the subcone 
$$
C':=\bigcup_{\om\in\Amp(X)}\{\om\}\times(\PL(X)_\R\cap\PSH(\om)),
$$
which spans $V':=\Num(X)\times\PL(X)_\R$; the extended pairing on $V$ thus coincides with the given one on $V'$,   which proves (i). Next, (ii) is immediate by multilinearity. Finally (iii) holds when $\theta_i\ge\om$, by Theorem~\ref{thm:extint}; the general case follows, again by multilinearity.  
\end{proof}

By multilinearity, \eqref{equ:inthom}, \eqref{equ:intsumcomp} and~\eqref{equ:intnorm} remain valid on $\Num(X)\times\vec\cE^1$. We also note: 

\begin{lem}\label{lem:posnef} If $\theta_0,\dots,\theta_n\in\Nef(X)$, then $(\theta_0,\f_0)\inter (\theta_n,\f_n)$ is an increasing function of $\f_i\in\PSH(\theta_i)\cap\vec\cE^1$. 
\end{lem} 
\begin{proof} It suffices to write $(\theta_0,\f_0)\inter(\theta_n,\f_n)=\lim_{\e\to 0}(\theta_0+\e\theta,\f_0)\inter(\theta_n+\e\theta,\f_n)$ with $\theta\in\Amp(X)$. 
\end{proof}

\begin{rmk}\label{rmk:E1nef} Mimicking the complex analytic case, one can define $\cE^1(\theta)$ for any $\theta\in\Nef(X)$ as the set of $\f\in\PSH(\theta)$ such that $\inf_j(\theta,\f_j)^{n+1}>-\infty$ with $\f_j:=\max\{\f,-j\}\in\cE^\infty(\theta)$. The proof of Lemma~\ref{lem:pairinglin} still yields $\vec\cE^1\cap\PSH(\theta)\subset\cE^1(\theta)$, but the inclusion is strict in general when $\theta$ is not ample (compare~\cite{DiN16}). 
\end{rmk}

%
%
%
%
\subsection{Mixed Monge--Amp\`ere measures}\label{sec:mixedMA}
Recall that a Radon measure $\mu$ on the compact space $X^\an$ is a regular (positive) Borel measure on $X^\an$. By the Riesz representation theorem, Radon measures are in 1--1 correspondence with positive linear forms on $\Cz(X)$, see for instance~\cite[\S7.1--2]{Folland}. Any usc function $f:X^\an\to\R\cup\{-\infty\}$ satisfies 
$$
f=\inf\left\{g\in\Cz(X),\,g\ge f\right\}
$$
pointwise, and 
\begin{equation}\label{equ:intusc}
\int f\,\mu=\inf\left\{\int g\,\mu\mid g\in\Cz(X),\,g\ge f\right\}\in\R\cup\{-\infty\}. 
\end{equation}

We will sometimes need to rely on the following monotone convergence theorem for (possibly uncountable) nets of usc functions, a simple consequence of~\eqref{equ:intusc} and Dini's lemma (see for instance~\cite[Proposition~7.12]{Folland}). 

\begin{lem}\label{lem:mononets}
  If $\mu$ is a Radon measure on $X^\an$ 
  and $(f_j)_j$ a decreasing net of usc functions on $X^\an$,
  converging pointwise to a (usc) function $f$, then 
  $\lim_j\int f_j\mu=\int f\mu$.
\end{lem}
More generally, continuous linear forms in $\Cz(X)^\vee$ correspond to \emph{signed Radon measures} on $X^\an$. Any such measure can be written as a difference of Radon measures, and Lemma~\ref{lem:mononets} thus applies as well when $\mu$ is a signed measure. 
 
With these preliminaries in hand, we now generalize the construction of mixed Monge--Amp\`ere measures, introduced in \S\ref{sec:enpairing} in the PL case. 

\begin{thm}\label{thm:mixedMA} For each $n$-tuple $(\theta_i,\f_i)\in\Num(X)\times\vec\cE^1$, $i=1,\dots,n$, there exists a unique signed measure 
\begin{equation}\label{equ:mixedMA1}
(\theta_1+\ddc \f_1)\winter (\theta_n+\ddc \f_n)=\bigwedge_i(\theta_i+\ddc\f_i)\in\Cz(X)^\vee
\end{equation}
such that 
\begin{equation}\label{equ:mixedMA2}
\int_{X^\an}\f\bigwedge_i(\theta_i+\ddc\f_i)=(0,\f)\cdot(\theta_1,\f_1)\inter(\theta_n,\f_n)
\end{equation}
for all $\f\in\PL(X)_\R$. If further $\theta_i\in\Nef(X)$ and $\f_i\in\PSH(\theta_i)$ for $i=1,\dots,n$, then~\eqref{equ:mixedMA1} is a positive measure. 
\end{thm}

While the notation mimics the one for mixed Monge--Amp\`ere measures in the complex analytic case, we will not define the individual factors $\theta_i+\ddc\f_i$. At least if the $\f_i$ are continuous, this could, however, be done using the approach in~\cite{CLD} (cf.~Remark~\ref{rmk:MACLD} below).

\begin{proof} By density of $\PL(X)$ in $\Cz(X)$ (see Theorem~\ref{thm:PLdense}), the prescription~\eqref{equ:mixedMA2} uniquely determines the measure~\eqref{equ:mixedMA1}. To show existence, we may assume, by multilinearity, that $\theta_i\in\Nef(X)$ and $\f_i\in\vec\cE^1\cap\PSH(\theta_i)$ for $i=1,\dots,n$. Since any $\f\in\PL(X)_\R$ can be written as a difference of functions in $\PL_\R\cap\PSH(\om)$ for some $\om\in\Amp(X)$ (see~\eqref{equ:HdomPL}), Lemma~\ref{lem:posnef} shows that $\f\mapsto (0,\f)\cdot(\theta_1,\f_1)\inter(\theta_n,\f_n)$ is a positive linear form on $\PL(X)_\R$. By density of $\PL(X)_\R$ in $\Cz(X)$, it thus uniquely extends to a positive linear form on $\Cz(X)$, and we are done. 
\end{proof}

\begin{prop}\label{prop:mixedMA} Mixed Monge--Amp\`ere measures satisfy the following properties: 
\begin{itemize}
\item[(i)] the signed measure $\bigwedge_i(\theta_i+\ddc\f_i)$ is a symmetric and multilinear function of the $n$-tuple $(\theta_i,\f_i)\in\Num(X)\times\vec\cE^1$, of total mass
$$
\int_{X^\an}\bigwedge_i(\theta_i+\ddc\f_i)=(\theta_1\inter\theta_n)_X;
$$
\item[(ii)] assume $(\theta_i,\f_i)\in\Num(X)\times\PL_\R$, and pick an integrally closed test configuration $\cX$ dominating $\cX_\triv$ and $D_i\in\VCar(\cX)_\R$ such that $\f_i=\f_{D_i}$; then
$$
\bigwedge_i(\theta_i+\ddc\f_i)=\sum_E c_E\,\d_{v_E},
$$
where $E$ ranges over the irreducible components of $\cX_0$ and
$$
c_E:=\ord_E(\cX_0)(\theta_{1,\cX}+D_1 )|_E\inter(\theta_{n,\cX}+D_n )|_E; 
$$
\item[(iii)] for all $(\theta_i,\f_i)\in\Num(X)\times\vec\cE^1$ and $t\in\R_{>0}$ we have
\begin{equation}\label{equ:mixedMAscaling}
\bigwedge_i\left(\theta_i+\ddc(t\cdot\f_i)\right)=t_\star\bigwedge_i(\theta_i+\ddc\f_i),
\end{equation}
\begin{equation}\label{equ:mixedMAsum}
\bigwedge_i\left(\theta_i+\ddc\f_i\right)=\sum_{\dim X_\a=n}\bigwedge_i\left(\theta_i|_{X_\a}+\ddc\f_i|_{X_\a^\an}\right),
\end{equation}
\begin{equation}\label{equ:mixedMAnorm}
\bigwedge_i\left(\theta_i+\ddc\f_i\right)=\nu_\star\bigwedge_i\left(\nu^\star\theta_i+\ddc\nu^\star\f_i\right),
\end{equation}
with $\nu\colon X^\nu\to X$ the normalization; 
\item[(iv)] for any $\om\in\Amp(X)$, the map $\cE^1(\om)^n\ni(\f_1,\dots,\f_n)\mapsto\bigwedge_i(\theta_i+\ddc\f_i)$ is continuous along decreasing nets; 
\item[(v)] for all $(\theta_i,\f_i)\in\Num(X)\times\vec\cE^1$ and $\f\in\vec\cE^1$, $\f$ is integrable against $\bigwedge_i(\theta_i+\ddc\f_i)$, and~\eqref{equ:mixedMA2} holds. 
\end{itemize}
\end{prop}

When $\f_i=0$ for some $i$, we drop the term $\ddc\f_i$ from the notation. As a special case of Proposition~\ref{prop:mixedMA}~(ii), we then have
\begin{equation}\label{equ:mixedMAzero}
  \theta_1\winter\theta_n
  =\sum_{\dim X_\a=n}\left(\theta_1\inter\theta_n\right)_{X_\a}\d_{v_{\triv,\a}}. 
\end{equation}

\begin{proof}[Proof of Proposition~\ref{prop:mixedMA}] Points (i) and (iii) follow from Proposition~\ref{prop:intpsh}, while (ii) is a reformulation of~\eqref{equ:intnumb}. By density of $\PL(X)$ in $\Cz(X)$, (iv) is equivalent to the fact that, for each $\f\in\PL(X)$, 
$$
(\f_1,\dots,\f_n)\mapsto\int\f\,\bigwedge_i\left(\theta_i+\ddc\f_i\right)=(0,\f)\cdot(\theta_1,\f_1)\inter(\theta_n,\f_n)
$$
is continuous along decreasing nets in $\cE^1(\om)^{n+1}$. Now $\f$ can be written as a difference of functions in $\PL\cap\PSH(\om')$ for some $\om'\in\Amp(X)$, and the desired continuity is thus a consequence of   Theorem~\ref{thm:pairinglin}. To prove (v), we may assume $\theta_i\in\Amp(X)$ and $\f_i\in\cE^1(\theta_i)$, by multilinearity.  Any $\f\in\cE^1(\om)$ with $\om\in\Amp(X)$ can be written as the pointwise limit of a decreasing net in $\PL\cap\PSH(\om)$, and~\eqref{equ:mixedMA2} thus holds for $\f$, by monotone convergence and the continuity of the energy pairing along decreasing limits.  
\end{proof}

\begin{exam}\label{exam:MAcurve} Assume $X$ is a smooth curve and pick $\om\in\Amp(X)$. By Example~\ref{exam:curvepsh}, any $\f\in\PSH(\om)$ determines a positive Radon measure 
$(\deg\om)\d_{v_{\triv}}+\D\f$, with 
$$
\D\f=\sum_{p\in X(k)}\left[\frac{d}{dt}\bigg|_{0+}\f(t\ord_p)\d_{v_\triv}-\frac{d}{dt}\bigg|_{+\infty}\f(t\ord_p)\d_{v_{p,\triv}}+\frac{d^2}{dt^2}\f(t\ord_p)\right]
$$
the tree Laplacian (see~\cite[\S 7]{valtree}). As we saw in Example~\ref{exam:curvePLMA}, we have $\om+\ddc\f=(\deg\om)\d_{v_{\triv}}+\D\f$ for any $\f\in\PL_\R\cap\PSH(\om)$, and this remains true for any $\f\in\cE^1(\om)$, by monotone approximation. For such functions, $\om+\ddc\f$ further puts no mass on the endpoints $v_{p,\triv}$, which means that the convex function $t\mapsto\f(t\ord_p)$ has sublinear growth (see respectively Corollary~\ref{cor:mixedMArestr} and Theorem~\ref{thm:phiTbelow}). 
\end{exam}

\begin{prop}\label{prop:MAmixedconv} Pick $\om_i\in\Amp(X)$ and $\f_i\in\cE^1(\om_i)$, $i=1,\dots,n$, and choose also $\om_0\in\Amp(X)$ and $\f_0\in\PSH(\om_0)$. Then:  
\begin{itemize}
\item[(i)] $\f_0$ is integrable with respect to $\mu:=\bigwedge_i(\om_i+\ddc\f_i)$ iff $(\om_0,\f_0)\inter(\om_n,\f_n)>-\infty$, and this holds as soon as $\f_i$ is bounded for $i=1,\dots,n$;  
\item[(ii)] assume $\f_0\in L^1(\mu)$. Pick for $i=0,\dots,n$ a decreasing net $(\f_{ij})_j$ in $\PSH(\om_i)$ that converges pointwise to $\f_i$, and set $\mu_j:=\bigwedge_i(\om_i+\ddc\f_{ij})$. Then $\f_{0j}\in L^1(\mu_j)$ for all $j$, and 
$\f_{0j}\mu_j\to\f_0\mu$ weakly in $\Cz(X)^\vee$.
\end{itemize}
\end{prop}

\begin{lem}\label{lem:monoweak} Let $K$ be a compact topological space with a weakly convergent net of positive Radon measures 
$\mu_j\to\mu$. Assume also given an increasing net of lsc functions $f_j\colon K\to\R\cup\{+\infty\}$ such that $f_j\in L^1(\mu_j)$ for all $j$, $f:=\lim_j f_j\in L^1(\mu)$ and $\int f_j\,\mu_j\to\int f\,\mu$. Then $f_j\mu_j\to f\mu$ weakly in $\Cz(K)^\vee$. 
\end{lem}
\begin{proof} Being lsc, each $f_j$ is bounded below, and the increasing net $(f_j)$ is thus ultimately bounded below. Since $\mu_j\to\mu$, we may therefore assume, after adding a constant, that $f_j\ge 0$ for all $j$. Since $\int f_j\,\mu_j$ converges, 
the positive measures $\sigma_j:=f_j\mu_j$ stay in a fixed weakly compact subset of $\Cz(X)^\vee$, and it suffices to show that any limit point $\sigma_\infty$ of $\sigma_j$ must coincide with $\sigma:=f\mu$. By assumption, $\int\sigma_j\to\int\sigma$, hence $\int\sigma_\infty=\int\sigma$, and it will thus be enough to show that $\sigma_\infty\ge\sigma$. To this end, pick $0\le g\in\Cz(K)$. For all $j\ge k$ we have $f_j\ge f_k$, and hence
$$
\int g\,\sigma_j=\int g f_j\,\mu_j\ge \int g f_k\,\mu_j. 
$$
By lower semicontinuity of $g f_k$ and the weak convergence $\mu_j\to\mu$, we infer $\int g\,\sigma_\infty\ge \int g f_k\,\mu$. Using monotone convergence, this now yields in turn,
$$
\int g\,\sigma_\infty\ge \int g f\,\mu=\int g\,\sigma, 
$$
which concludes the proof.
\end{proof}

\begin{proof}[Proof of Proposition~\ref{prop:MAmixedconv}] Set $c:=(\om_0,0)\cdot(\om_1,\f_1)\inter(\om_n,\f_n)\in\R$, and pick a decreasing net $(\f_{0j})$ in $\PL\cap\PSH(\om_0)$ such that $\f_{0j}\to\f_0$ pointwise. By construction of $\mu$, we have for each $j$ 
$$
\int\f_{0j}\,\mu=(\om_0,\f_{0j})\cdot(\om_1,\f_1)\inter(\om_n,\f_n)-c. 
$$
By Lemma~\ref{lem:mononets}, $\int\f_{0j}\,\mu\to\int\f_0\,\mu$. By Theorem~\ref{thm:extint}, we have on the other hand $(\om_0,\f_{0j})\cdot(\om_1,\f_1)\inter(\om_n,\f_n)\to(\om_0,\f_{0})\cdot(\om_1,\f_1)\inter(\om_n,\f_n)$. Thus 
$$
\int\f_0\,\mu=(\om_0,\f_0)\cdot(\om_1,\f_1)\inter(\om_n,\f_n)-c, 
$$
and (i) follows. 

We turn to (ii). By Proposition~\ref{prop:mixedMA}~(iv), we have $\mu_j\to\mu$. Further, 
$$
\int\f_{0j}\,\mu_j=(\om_0,0)\cdot(\om_1,\f_{1j})\inter(\om_n,\f_{nj})-(\om_0,\f_{0j})\inter(\om_n,\f_{nj})
$$
is finite for each $j$, and converges to 
$$
\int\f_0\,\mu=(\theta_0,0)\cdot(\om_1,\f_{1})\inter(\om_n,\f_{n})-(\om_0,\f_0)\inter(\om_n,\f_{n}). 
$$
The result is now a consequence of Lemma~\ref{lem:monoweak}. 
\end{proof}

\begin{rmk}\label{rmk:MACLD} For continuous $\om$-psh functions, mixed Monge--Amp\`ere measures can also be defined using the general theory developed by Chambert-Loir and Ducros~\cite{CLD}. By base change invariance of their theory, it follows from Proposition~\ref{prop:mixedMA}~(ii)  and~\cite[Theorem 8.18]{BE} that the present approach is compatible with~\cite{CLD} for continuous $\om$-psh functions. 
\end{rmk}

%
%
%
%
\subsection{The Monge--Amp\`ere operator and energy functionals}\label{sec:MAen}
In this section we fix $\om\in\Amp(X)$. We denote by $V:=(\om^n)$ its volume, and write $\cE^1:=\cE^1(\om)$ and $\en:=\en_\om$. 

\begin{defi} The \emph{Monge--Amp\`ere operator} takes $\f\in\cE^1$ to the Radon probability measure 
$$
\MA(\f)=\MA_\om(\f):=V^{-1}(\om+\ddc\f)^n
$$
on $\Xan$.
\end{defi}

\begin{defi} For any two $\f,\p\in\cE^1$ we set
$$
\jj_\p(\f)=\jj_{\om,\p}(\f)=\en(\p)-\en(\f)+\int(\f-\p)\MA(\p)
$$
and 
\begin{align*}
  \ii(\f,\p)=\ii_\om(\f,\p)
  &=\int (\f-\p)\left(\MA(\p)-\MA(\f)\right)\\
  &=\jj_\p(\f)+\jj_\f(\p).
\end{align*}
\end{defi}
Note that $\jj_\p(\f)$ and $\ii(\f,\p)$ are both invariant under translation of $\f,\p$ by constants. When $\p=0$ we simply write
$$
\ii(\f):=\ii(\f,0)=\int\f\MA(0)-\int\f\MA(\f)
$$
and
$$
\jj(\f):=\jj_0(\f)=\int\f\MA(0)-\en(\f).  
$$

\begin{exam}\label{exam:MAzero}
  By~\eqref{equ:mixedMAzero} with $\theta_1=\dots=\theta_n=\om$ we have 
  \begin{equation}\label{equ:MAzero}
    \MA(0)=\sum_{\dim X_\a=n}c_\a\d_{v_{\triv,\a}},
  \end{equation}
  with $c_\a=(\om^n)_{X_\a}/(\om^n)$. As a consequence, $\int\f\MA(0)=\sum_\a c_\a\f(v_{\triv,\a})=\sum_\a c_\a\sup(\f|_{X_\a})$, see~\eqref{equ:maxprin}. 
  In particular, when $X$ is irreducible, $\MA(0)=\d_{v_\triv}$ and $\int\f\MA(0)=\sup\f$.
\end{exam}

\medskip
Thanks to Theorem~\ref{thm:E1}, we can make sense of the definitions from \S\ref{sec:estimates} for functions in $\cE^1$. In particular, we set, for any two $\f,\p\in\cE^1$, 
$$
d_\om(\f,\p):=\max_{0\le j\le n-1}\|\f-\p\|^2_{(\om,\f)^j\cdot(\om,\p)^{n-1-j}}
$$
with 
\begin{multline}\label{equ:squarenorm}
\|\f-\p\|^2_{(\om,\f)^j\cdot(\om,\p)^{n-1-j}}=-(0,\f-\p)^2\cdot(\om,\f)^j\cdot(\om,\p)^{n-j}\\
=-\int(\f-\p)\ddc(\f-\p)\wedge(\om+\ddc\f)^j\wedge(\om+\ddc\p)^{n-j}. 
\end{multline}

\begin{prop}\label{prop:EIJ} For all $\f,\p\in\cE^1$ we have: 
\begin{equation}\label{equ:Ebis}
\en(\f)-\en(\p)=\frac{1}{n+1}\sum_{j=0}^n V^{-1}\int(\f-\p)\,(\om+\ddc\f)^j\wedge(\om+\ddc\p)^{n-j}; 
\end{equation}
\begin{equation}\label{equ:Ebound}
\int(\f-\p)\MA(\f)\le\en(\f)-\en(\p)\le\int(\f-\p)\MA(\p); 
\end{equation}
\begin{equation}\label{equ:derE}
\frac{d}{dt}\bigg|_{t=0}\en((1-t)\f+t\p)=\int(\p-\f)\MA(\f); 
\end{equation}
\begin{equation}\label{equ:I}
\ii(\f,\p)=V^{-1}\sum_{j=0}^{n-1}\|\f-\p\|^2_{(\om,\f)^j,(\om,\p)^{n-1-j}};
\end{equation}
\begin{equation}\label{equ:Isumcomp}
\ii(\f,\p)=\sum_\a c_\a \ii(\f|_{\Xan_\a},\p|_{\Xan_\a}).
\end{equation}
\begin{equation}\label{equ:J}
\jj_\p(\f)=V^{-1}\sum_{j=0}^{n-1}\frac{j+1}{n+1}\|\f-\p\|^2_{(\om,\f)^j,(\om,\p)^{n-1-j}}; 
\end{equation}
\begin{equation}\label{equ:IJ} 
\tfrac{1}{n+1}\ii(\f,\p)\le\jj_\p(\f)\le\tfrac{n}{n+1}\ii(\f,\p);
\end{equation}
\begin{equation}\label{equ:IJd}
V^{-1} d_\om(\f,\p)\approx\ii(\f,\p)\approx\jj_\p(\f)\approx\jj_\f(\p);
\end{equation}
\begin{equation}\label{equ:Itriangle}
\ii(\f_1,\f_2)\lesssim\ii(\f_1,\f_3)+\ii(\f_3,\f_2).
\end{equation}
\end{prop}

\begin{proof} As in Lemma~\ref{lem:en},~\eqref{equ:Ebis},~\eqref{equ:derE},~\eqref{equ:I} and~\eqref{equ:J} follow from straightforward computations based on the multilinearity and symmetry of the energy pairing on $\cE^1$, while~\eqref{equ:Isumcomp} follows from~\eqref{equ:intsumcomp}. Equation~\eqref{equ:Ebound} follows from~\eqref{equ:enconc} by monotone approximation. Equations~\eqref{equ:I} and~\eqref{equ:J} imply~\eqref{equ:IJ} and~\eqref{equ:IJd}, and~\eqref{equ:Itriangle} is now a consequence of Theorem~\ref{thm:triangle} and monotone approximation. 
\end{proof}

Taking into account~\eqref{equ:derE}, we can write~\eqref{equ:IJ} as
\begin{equation}\label{equ:strictconc}
\en(\p)+\en'(\p)(\f-\p)-\en(\f)\ge\tfrac 1{n+1} \ii(\f,\p),
\end{equation}
which can be seen as a strict concavity property of $\en$ with respect to $\ii$.   Indeed, it yields the following uniform concavity estimate: 
\begin{thm}\label{thm:strictconc} For any two $\f,\p\in\cE^1$ and $t\in[0,1]$ we have 
$$
\en((1-t)\f+t\p)-\left((1-t)\en(\f)+t\en(\p)\right)\gtrsim t(1-t)\ii(\f,\p).
$$
\end{thm}
The proof relies on the following elementary estimate, which is certainly well-known. 

\begin{lem}\label{lem:elemconv} For all $a,b\in\R_{\ge 0}$ and $t\in [0,1]$ we have
$(1-t)a+tb\ge t(1-t)(a+b)$. 
\end{lem}
\begin{proof} By homogeneity we may assume $a+b=1$. Write $t=\frac 12+x$ and $a=\tfrac 12+y$ with $x,y\in[-\tfrac 12,\tfrac 12]$, where we may assume $x\ge 0$, by symmetry. Then 
$$
(1-t)a+tb=(\frac 12-x)(\frac 12+y)+(\frac 12+x)(\frac 12-y)=\frac 12-2xy,\quad t(1-t)(a+b)=(\frac 12+x)(\frac 12-x)=\frac 14-x^2,
$$
and hence
$$
(1-t)a+tb-t(1-t)(a+b)=\frac 14-2xy+x^2\ge\frac 14-\frac 12 x+x^2=(\frac 12-x)^2\ge 0.
$$
\end{proof}

\begin{proof}[Proof of Theorem~\ref{thm:strictconc}] Set $\f_t:=(1-t)\f+t\p$. By~\eqref{equ:strictconc}, we have 
$$
\en(\f_t)+\int(\f-\f_t)\MA(\f_t)-\en(\f)\gtrsim\ii(\f,\f_t)
$$
and
$$
\en(\f_t)+\int(\p-\f_t)\MA(\f_t)-\en(\p)\gtrsim\ii(\f_t,\p),
$$
Since $(1-t)(\f-\f_t)+t(\p-\f_t)=0$, this yields 
$$
\en(\f_t)-\left((1-t)\en(\f)+t\en(\p)\right)\gtrsim (1-t)\ii(\f,\f_t)+t\ii(\p,\f_t)
$$
$$
\ge t(1-t)\left(\ii(\f,\f_t)+\ii(\f_t,\p)\right)\gtrsim t(1-t)\ii(\f,\p), 
$$
by Lemma~\ref{lem:elemconv} and~\eqref{equ:Itriangle}. 
\end{proof}

For later use, we record the following crucial consequence of Corollary~\ref{cor:CS} (compare~\cite[Lemma 3.13]{BBGZ}). 
 
\begin{lem}\label{lem:BBGZ} For all $\f,\f',\p,\p'\in\cE^1$ we have 
$$
\left|\int(\f-\f')\left(\MA(\p)-\MA(\p')\right)\right|\lesssim \ii(\f,\f')^{\a_n}\ii(\p,\p')^{\frac 12}\max\{\jj(\f),\jj(\f'),\jj(\p),\jj(\p')\}^{\frac 12-\a_n}
$$
with $\a_n:=2^{-n}$. 
\end{lem}
 
\begin{proof}
  After regularization, we may assume $\f,\f',\p,\p'\in\PL_\R\cap\PSH(\om)$.
  Note that
\begin{multline*}
V\int(\f-\f')\left(\MA(\p)-\MA(\p')\right)=(0,\f-\f')\cdot(\om,\p)^n-(0,\f-\f')\cdot(\om,\p')^n\\
=\sum_{j=1}^{n-1}(0,\f-\f')\cdot(0,\p-\p')\cdot(\om,\p)^j\cdot(\om,\p')^{n-1-j}.
\end{multline*}
By the Cauchy--Schwarz inequality~\eqref{equ:CS} we infer
\begin{multline*}
\left|(0,\f-\f')\cdot(\om,\p)^n-(\f-\f')\cdot(\om,\p')^n\right|\\
\lesssim\max_j\left(\|\f-\f'\|_{(\om,\p)^j\cdot(\om,\p')^{n-1-j}}\|\p-\p'\|_{(\om,\p)^j\cdot(\om,\p')^{n-1-j}}\right)\\
\le\left(\max_j\|\f-\f'\|_{(\om,\p)^j\cdot(\om,\p')^{n-1-j}}\right) d_\om(\p,\p')^{\frac 12}
\end{multline*}
and Corollary~\ref{cor:CS} together with~\eqref{equ:IJd} yield the desired estimate. 
\end{proof}

\begin{cor}\label{cor:jjhold} For all $\f,\f',\p\in\cE^1$ we have 
$$
|\jj_\p(\f)-\jj_\p(\f')|\lesssim \ii(\f,\f')^{\a_n}\max\{\jj(\f),\jj(\f'),\jj(\p)\}^{1-\a_n}
$$
with $\a_n=2^{-n}$. 
\end{cor} 
\begin{proof} Set $M:=\max\{\jj(\f),\jj(\f'),\jj(\p)\}$. We have 
$$
\jj_\p(\f)-\jj_\p(\f')=\en(\f')-\en(\f)+\int(\f-\f')\MA(\p)=\jj_{\f'}(\f)+\int(\f-\f')(\MA(\p)-\MA(\f')),
$$
and~\eqref{equ:IJ} and Lemma~\ref{lem:BBGZ} thus yield
$$
|\jj_\p(\f)-\jj_\p(\f')|\lesssim\ii(\f,\f')+\ii(\f,\f')^{\a_n}\ii(\p,\f')^{\frac 12} M^{\frac 12-\a_n}. 
$$
By the quasi-triangle inequality, we have 
$$
\ii(\f,\f')=\ii(\f,\f')^{\a_n}\ii(\f,\f')^{1-\a_n}\lesssim \ii(\f,\f')^{\a_n} M^{1-\a_n}
$$
and $\ii(\p,\f')\lesssim M$, and we obtain, as desired, 
$|\jj_\p(\f)-\jj_\p(\f')|\lesssim \ii(\f,\f')^{\a_n} M^{1-\a_n}$. 
\end{proof}

By Proposition~\ref{prop:MAmixedconv}~(i), Monge--Amp\`ere measures of bounded $\om$-psh functions integrate all $\om$-psh functions. In that setting,~\eqref{equ:intlip} yields the following simpler variant of Lemma~\ref{lem:BBGZ}, which can be viewed as a version of the classical Chern--Levine--Nirenberg inequality. 

\begin{lem}\label{lem:CLN} If $\p,\p'\in\PSH(\om)$ are bounded, then 
$$
\left|\int\f\left(\MA(\p)-\MA(\p')\right)\right|\le n\sup|\p-\p'|
$$
for all $\f\in\PSH(\om)$. 
\end{lem}

We conclude this section with two useful additional estimates. 
\begin{lem}\label{lem:quadJ} For each $\f,\p\in\cE^1$ and $t\in[0,1]$ we have 
$$
\ii(t\f+(1-t)\p,\p)\le(1-(1-t)^n)\ii(\f,\p)\le n t^2\ii(\f,\p)
$$
and
$$
\jj_\p(t\f+(1-t)\p)\le t^{1+\frac1n}\jj_\p(\f).
$$
\end{lem}
\begin{proof} Adding a constant to $\f$, we may assume $\int(\f-\p)\MA(\p)=0$. 
  Set $\f_t:=t\f+(1-t)\p$.
  Then $\f_t-\p=t(\f-\p)$, so $\int(\f_t-\p)\MA(\p)=0$, and
  \begin{multline*}
    \ii(\f_t,\p)
    =-t\int(\f-\p)\MA(\f_t)\\
    =-tV^{-1}\sum_{j=0}^n\binom{n}{j}t^j(1-t)^{n-j}\int(\f-\p)(\om+\ddc\f)^j\wedge(\om+\ddc\p)^{n-j}.
  \end{multline*}
  Here the integral vanishes for $j=0$ and is bounded below by $\int(\f-\p)\MA(\f)$
  for $j>0$, so 
  \begin{equation*}
    \ii(\f_t,\p)
    \le -t(1-(1-t)^n)\int(\f-\p)\MA(\f)
    \le nt^2\ii(\f,\p),
  \end{equation*}
  by the concavity of $t\mapsto(1-(1-t)^n)$.
  
  To prove the last inequality, note that 
  \begin{equation*}
    \jj_\p(\f_t)=\en(\p)-\en(\f_t)+t\int(\f-\p)\MA(\p).
  \end{equation*}
  Differentiating this with respect to $t$, and using~\eqref{equ:IJ}, gives
  \begin{equation*}
    \frac{d}{dt}\jj_\psi(\f_t)
    =-\int(\f-\p)\MA(\f_t)+\int(\f-\p)\MA(\p)
    =t^{-1}\ii(\f_t,\p)
    \ge\tfrac{n+1}{n}t^{-1}\jj_\p(\f_t)
  \end{equation*}
  for $0<t\le 1$, from which the desired estimate follows easily.
\end{proof}
%
%
\subsection{H\"older continuity of the energy pairing}
In this section, $X$ is \textbf{irreducible}.
Using the estimates of \S\ref{sec:estimates}, we are going to establish the following general H\"older continuity property of the energy pairing. 
For $\om\in\Amp(X)$ we write $$\cE^1_{\sup}(\om):=\{\f\in\cE^1(\om)\mid\sup\f=0\}.$$

\begin{thm}\label{thm:energyholder} Pick $\om_0,\dots,\om_n\in\Amp(X)$, and assume we are given $t\ge 1$ such that $\om_i\le t\om_j$ for all $i,j$. For all tuples of pairs $\f_i,\f'_i\in\cE^1_{\sup}(\om_i)$, we then have 
$$
\left|(\om_0,\f_0)\inter(\om_n,\f_n)-(\om_0,\f'_0)\inter(\om_n,\f'_n)\right|
$$
$$
\lesssim t^{n^2}\max_i d_{\om_i}(\f_i,\f'_i)^{\a_n}\max\{\max_i d_{\om_i}(\f_i),\max_i d_{\om_i}(\f_i')\}^{1-\a_n}.
$$
with $\a_n\in(0,1]$ only depending on $n$. 
\end{thm}
Recall that we have set $d_\om(\f)=d_\om(\f,0)$. To each $\om\in\Amp(X)$ we associate a norm on $\Num(X)$ by setting 
$$
\|\theta\|_\om:=\inf\left\{C\ge 0\mid -C\theta\le\om\le C\theta\right\}. 
$$

\begin{cor}\label{cor:energybounded} Pick $\om\in\Amp(X)$, $\f_0,\dots,\f_n\in\cE^1(\om)$, $\theta_0,\dots,\theta_n\in\Num(X)$. Then 
$$
\left|(\theta_0,\f_0)\inter(\theta_n,\f_n)\right|\lesssim\max_i\{1,\|\theta_i\|_\om,d_\om(\f_i)\}^{n+1}+\left(\max_i|\sup\f_i|\right)\left(\max_i\|\theta_i\|_\om^n\right)(\om^n). 
$$
\end{cor}

\begin{lem}\label{lem:sumstrong} Pick $r\ge 1$, $\om_0,\dots,\om_r\in\Amp(X)$, and $t\ge 1$ such that $\om_i\le t\om_j$ for all $i,j$. Pick $\f_i,\f'_i\in\cE^1_{\sup}(\om)$, $i=0,\dots,r$, and set 
$$
\om:=\sum_i\om_i,\quad\f:=\sum_i\f_i,\quad\f':=\sum_i \f'_i.
$$
Then 
$$
d_\om(\f,\f')\le C_{r,n} t^{rn}\max_i d_{\om_i}(\f_i,\f'_i)^{\b_n}\max\{\max_i d_{\om_i}(\f_i),\max_i d_{\om_i}(\f'_i)\}^{1-\b_n}
$$
with $\b_n\in(0,1]$ and $C_{r,n}\in(0,+\infty)$ only depending on $n$ and $r,n$, respectively. 
\end{lem}

\begin{proof} Set $M:=\max\{\max_i d_{\om_i}(\f_i),\max_i d_{\om_i}(\f'_i)\}$. By definition we have 
$$
d_{\om}(\f,\f')=\max_{a+b=n-1}\|\sum_{i=1}^r (\f_i-\f'_i)\|^2_{(\om,\f)^a,(\om,\f')^b}. 
$$
The triangle inequality thus yields 
$$
d_{\om}(\f,\f')\le r^2\max_{a+b=n-1}\max_i\|\f_i-\f'_i\|^2_{(\om,\f)^a,(\om,\f')^b}. 
$$
and hence
$$
d_{\om}(\f,\f')\lesssim r^2\max_i d_\om(\f_i,\f'_i)^{2^{1-n}}\max\{d_\om(\f),d_\om(\f'),M)\}^{1-2^{1-n}},
$$
by Corollary~\ref{cor:CS}. On the one hand, we have for each $i$ 
$$
d_\om(\f_i,\f'_i)=\max_{a+b=n-1}\|\f_i-\f'_i\|^2_{(\om,\f_i)^a\cdot (\om,\f'_i)^b}
$$
$$
\lesssim\max_{a+b+c=n-1}\|\f_i-\f'_i\|^2_{(\om_i,\f_i)^a\cdot(\om_i,\f'_i)^b(\om_i',0)^c}
$$
with $\om'_i:=\om-\om_i=\sum_{j\ne i}\om_j$. Since $\om'_i\le rt\om_i$, we infer 
$$
d_\om(\f_i,\f'_i)\lesssim (rt)^n\max_{a+b+c=n-1}\|\f_i-\f'_i\|^2_{(\om_i,\f_i)^a\cdot(\om_i,\f'_i)^b(\om_i,0)^c}
$$
$$
\lesssim (rt)^n d_{\om_i}(\f_i,\f'_i)^{2^{1-n}}M^{1-2^{1-n}},
$$
using once more Corollary~\ref{cor:CS}. On the other hand, \eqref{equ:IJd} implies $d_\om(\f)\approx-(\om,\f)^{n+1}$, $d_\om(\f')\approx-(\om,\f')^{n+1}$, and Lemma~\ref{lem:ensum} thus yields 
$$
\max\{d_\om(\f),d_\om(\f'),M\}\le  C_{r,n} t^{rn} M. 
$$
All in all we infer
$$
d_\om(\f,\f')\le C'_{r,n} t^{(n-1) 2^{1-n}+rn(1-2^{1-n})} \max_i d_{\om_i}(\f_i,\f'_i)^{2^{2-2n}} M^{\left(1-2^{1-n}\right)\left(1+2^{1-n}\right)}
$$
which yields the desired estimate with $\b_n:=2^{2-2n}$. 
\end{proof}

\begin{proof}[Proof of Theorem~\ref{thm:energyholder}] Set $M:=\max\{\max_i d_{\om_i}(\f_i),\max_i d_{\om_i}(\f'_i)\}$. Since the energy pairing is symmetric and multilinear, the general polarization formula yields
$$
(\om_0,\f_0)\inter(\om_n,\f_n)=\sum_{I\subset\{0,\dots,n\}}(-1)^{n+1-|I|}(\om_I,\f_I)^{n+1},
$$
with $\om_I=\sum_{i\in I}\om_i$, $\f_I=\sum_{i\in I}\f_i$. Thus 
$$
\left|(\om_0,\f_0)\inter(\om_n,\f_n)-(\om_0,\f'_0)\inter(\om_n,\f'_n)\right|
$$
$$
\lesssim\max_I \left|(\om_I,\f_I)^{n+1}-(\om_I,\f'_I)^{n+1}\right|.
$$
By Corollary~\ref{cor:jjhold}, we have for each $I$
$$
\left|(\om_I,\f_I)^{n+1}-(\om_I,\f'_I)^{n+1}\right|\lesssim d_{\om_I}(\f_I,\f')^{2^{-n}}\max\{d_{\om_I}(\f_I),d_{\om_I}(\f'_I)\}^{1-2^{-n}},
$$
since $\sup\f_I=\sup\f'_I=0$. Now Lemma~\ref{lem:sumstrong} implies 
$$
d_{\om_I}(\f_I,\f'_I)\lesssim s^{n^2}\max_i d_{\om_i}(\f_i,\f'_i)^{\b_n} M^{1-\b_n}
$$
and 
$$
\max\{d_{\om_I}(\f_I),d_{\om_I}(\f'_I)\}\lesssim t^n M.
$$
Combining all this we get, as desired, 
$$
\left|(\om_0,\f_0)\inter(\om_n,\f_n)-(\om_0,\f'_0)\inter(\om_n,\f'_n)\right|\lesssim s^{n^2} \max_i d_{\om_i}(\f_i,\f'_i)^{\a_n} M^{1-\a_n},
$$
with $\a_n:=\b_n 2^{-n}$. 
\end{proof}

\begin{proof}[Proof of Corollary~\ref{cor:energybounded}] Set first $c_i:=\sup\f_i$, $\f'_i:=\f_i-c_i$. Then 
\begin{align*}
(\theta_0,\f_0)\inter(\theta_n,\f_n) & =(\theta_0,\f'_0+c_0)\inter(\theta_n,\f'_n+c_n)\\
& =(\theta_0,\f'_0)\inter(\theta_n,\f'_n)+\sum_{i=0}^n c_i(\theta_0\inter\widehat{\theta_i}\inter\theta_n). 
\end{align*}
Now 
$$
\left|\sum_{i=0}^n c_i(\theta_0\inter\widehat{\theta_i}\inter\theta_n)\right|\lesssim \left(\max_i|\sup\f_i|\right)\left(\max_i\|\theta_i\|_\om^n\right)(\om^n). 
$$
and we may thus assume wlog that $\sup\f_i=0$. Set $C:=\max_i\{1,\|\theta_i\|_\om,d_\om(\f_i)\}$. Since $C\ge 1$, $C^{-1}\f_i\in\cE^1(\om)$ satisfies $d_\om(C^{-1}\f_i)\lesssim 1$, by the quasi-convexity estimate~\eqref{equ:dconv}. After replacing $\theta_i$ and $\f_i$ with $C^{-1}\theta_i$ and $C^{-1}\f_i$, we may thus assume as well that $C\lesssim 1$, and we then need to prove $\left|(\theta_0,\f_0)\inter(\theta_n,\f_n)\right|\lesssim 1$. 

Set $\om_i:=\theta_i+(C+1)\om$, so that $\om\le\om_i\le(2C+1)\om$ and $\theta_i=\om+\om_i-(C+2)\om$. Expanding out 
$$
(\theta_0,\f_0)\inter(\theta_n,\f_n)
$$
$$
=\left[(\om,\f_0)+(\om_0,0)-(C+2)(\om,0)\right]\inter\left[(\om,\f_n)+(\om_n,0)-(C+2)(\om,0)\right]
$$
now yields the desired estimate, thanks to Theorem~\ref{thm:energyholder}. 
\end{proof}

As a further consequence of Theorem~\ref{thm:energyholder}, we show:

\begin{lem}\label{lem:MAweakcont} Pick $(\theta_i,\f_i)\in\Num(X)\times\vec\cE^1$, $i=1,\dots,n$, and set $\mu:=\bigwedge_i (\theta_i+\ddc\f_i)$. Suppose also that we are given $\om\in\Amp(X)$, and a convergent net $\p_j\to\p$ in $\PSH(\om)$ with $\jj_\om(\p_j)$ uniformly bounded. Then $\lim_j\int\p_j\,\mu=\int\p\,\mu$. 
\end{lem}
\begin{proof} By multilinearity, we may assume wlog $\theta_i\in\Amp(X)$ and $\f_i\in\cE^1(\theta_i)$. By Theorem~\ref{thm:pshample}, each $\f_i$ can be written as the limit of a decreasing net $(\f_{il})_l$ in $\cH^\dom(\theta_i)$. For each $l$, the measure $\mu_l:=\bigwedge_i (\theta_i+\ddc\f_{il})$ has finite support in $X^\div$ (see Proposition~\ref{prop:mixedMA}~(ii)), and hence $\lim_j\int\p_j\,\mu_l=\int\p\,\mu$, since $\p_j\to\p$ pointwise on $X^\div$ by assumption. To conclude the proof, it will thus be enough to show that $\lim_l \int\tau\,\mu_l=\int\tau\,\mu$ uniformly for $\tau\in\cE^1(\om)$ with $\jj_\om(\tau)\le C$. Now 
$$
\int\tau\,\mu_l=(\om,\tau)\cdot(\theta_1,\f_{1l})\inter(\theta_n,\f_{nl})-(\om,0)\cdot(\theta_1,\f_{1l})\inter(\theta_n,\f_{nl}), 
$$
$$
\int\tau\,\mu=(\om,\tau)\cdot(\theta_1,\f_1)\inter(\theta_n,\f_n)-(\om,0)\cdot(\theta_1,\f_1)\inter(\theta_n,\f_n). 
$$
Since $(\f_{il})_l$ is decreasing, we have $\lim_l d_{\theta_i}(\f_{il},\f_i)=0$, by~\eqref{equ:squarenorm} and Theorem~\ref{thm:pairinglin}. The desired uniform convergence is now a consequence of Theorem~\ref{thm:energyholder}. 
\end{proof}

This yields in turn the following monotone convergence theorem. 

\begin{thm}\label{thm:monostrong} Pick $\om\in\Amp(X)$. For $i=0,\dots,n$, assume we are given $\theta_i\in\Num(X)$ and an increasing net $(\f_{ij})_j$ in $\cE^1(\om)$ that converges pointwise on $X^\div$ to $\f_i\in\cE^1(\om)$. Then 
$$
\lim_j (\theta_1,\f_{1j})\inter(\theta_n,\f_{nj})=(\theta_1,\f_1)\inter(\theta_n,\f_n).
$$
\end{thm}
\begin{proof} By multilinearity, we can assume wlog $\theta_i\ge\om$ for all $i$, and hence $\cE^1(\om)\subset\cE^1(\theta_i)$. We proceed by induction on $p=0,\dots,n$ such that $\f_{ij}=\f_i$ is a constant net for $i\ge p$. The case $p=0$ is trivial, so assume $p\ge1$. By monotonicity of the energy pairing on $\prod_i\PSH(\theta_i)$ (see Theorem~\ref{thm:extint}), $j\mapsto(\theta_0,\f_{0j})\inter(\theta_n,\f_{nj})$ is increasing, and 
$$
\lim_j(\theta_0,\f_{0j})\inter(\theta_n,\f_{nj})\le(\theta_0,\f_0)\inter(\theta_n,\f_n).
$$
Conversely, pick $j\ge l$. Then $\f_{0j}\ge\f_{0l}$, and hence 
$$
 (\theta_0,\f_{0j})\inter(\theta_n,\f_{nj})\ge (\theta_0,\f_{0l})\cdot(\theta_1,\f_{1j})\inter(\theta_n,\f_{nj}). 
$$
Using the inductive assumption, we infer
$$
\lim_j (\theta_0,\f_{0j})\inter(\theta_n,\f_{nj})\ge (\theta_0,\f_{0l})\cdot(\theta_1,\f_1)\inter(\theta_n,\f_n)
$$
$$
=(\theta_0,0)\cdot(\theta_1,\f_1)\inter(\theta,\f_n)+\int\f_{0l}\,(\theta_1+\ddc\f_1)\winter(\theta_n+\ddc\f_n).
$$
Since $\jj_{\theta_0}(\f_{0l})=\sup\f_{0l}-\en_{\theta_0}(\f_{0l})$ is ultimately bounded, Lemma~\ref{lem:MAweakcont} next shows that 
$$
\lim_l\int\f_{0l}\,(\theta_1+\ddc\f_1)\winter(\theta_n+\ddc\f_n)=\int\f_0\,(\theta_1+\ddc\f_1)\winter(\theta_n+\ddc\f_n).
$$
Thus
$$
\lim_j (\theta_0,\f_{0j})\inter(\theta_n,\f_{nj})\ge (\theta_0,0)\cdot(\theta_1,\f_1)\inter(\theta_n,\f_n)
$$
$$+\int\f_0\,(\theta_1+\ddc\f_1)\winter(\theta_n+\ddc\f_n)=(\theta_0,\f_0)\cdot(\theta_1,\f_1)\inter(\theta_n,\f_n),
$$
and we are done. 
\end{proof}

\begin{rmk}\label{rmk:weakcontpairing} Lemma~\ref{lem:MAweakcont} and Theorem~\ref{thm:monostrong} remain valid when $X$ is reducible, by~\eqref{equ:mixedMAsum} and~\eqref{equ:intsumcomp}, respectively.  
\end{rmk}

%
%
%
%
\subsection{Locality and the comparison principle}
The next result and its consequences play a crucial role in analyzing deep properties of $\om$-psh functions and the Monge-Amp\`ere operator.
\begin{thm}\label{thm:locality}
  If $\om\in\Amp(X)$, then
  \begin{equation}\label{equ:locality}
    \one_{\{\f>\f'\}}\MA(\max\{\f,\f'\})=\one_{\{\f>\f'\}}\MA(\f)
  \end{equation}
  for all $\f,\f'\in\cE^1(\om)$.
\end{thm}
A first consequence is the fact that the mixed Monge--Amp\`ere operator is local in nature, something that
is not an immediate consequence of our definition in~\S\ref{sec:mixedMA}.
\begin{cor}\label{cor:locality}
  Let $G\subset\Xan$ be open set. If $\om_i\in\Amp(X)$ and $\f_i,\p_i\in\cE^1(\om_i)$, $1\le i\le n$, are such that $\f_i=\p_i$ on $G$, then
  $$
  (\om_1+\ddc\f_1)\wedge\dots\wedge(\om_n+\ddc\f_n)
  =
  (\om_1+\ddc\p_1)\wedge\dots\wedge(\om_n+\ddc\p_n)
  \ \text{on $G$}.
  $$
  In particular, if $\om\in\Amp(X)$ and $\f,\f'\in\cE^1(\om)$ are such that $\f=\f'$ on $G$, then $\MA(\f)=\MA(\f')$ on $G$.
\end{cor}
\begin{proof}
  By multilinearity, it suffices to prove the final statement.
  As in~\cite[Corollary 5.2]{nama}, given $\e>0$ we apply Theorem~\ref{thm:locality} to 
  $\f+\e$ and $\f'$. This gives $\MA(\max\{\f+\e,\f'\})=\MA(\f)$ 
  on $G\subseteq\{\f+\e>\f'\}$. Letting $\e\to 0$ gives $\MA(\max\{\f,\f'\})=\MA(\f)$ on $G$,
  since the Monge--Amp\`ere operator is continuous under decreasing limits. Exchanging the roles of $\f$ and $\f'$ completes the proof.
\end{proof}
\begin{rmk}
  As mentioned in Remark~\ref{rmk:MACLD}, when all the functions involved are continuous, the mixed Monge--Amp\`ere measure coincides with the one in~\cite{CLD}. Since the latter is local in nature, Theorem~\ref{thm:locality} and Corollary~\ref{cor:locality} follow in this case. However, for applications it is important to consider functions $\f,\f'$ that are not continuous, and in this case the Borel set $\{\f>\f'\}$ can be quite complicated. By definition, it is an open subset in the \emph{plurifine topology}, see~\cite{BT87}. 
\end{rmk}

\begin{lem}\label{lem:locality} Let $\f,\f'\in\PL(X)$, and pick an integrally closed test configuration $\cX$ such that $\f$, $\f'$ and $\f'':=\max\{\f,\f'\}$ are associated to elements $D$, $D'$, $D''$ of $\VCar(\cX)_\Q$. Pick an irreducible component $E$ of $\cX_0$, and assume that $\f(v_E)>\f'(v_E)$. Then $\f(v_F)\ge\f'(v_F)$ for any irreducible component $F$ of $\cX_0$ that intersects $E$. 
\end{lem}
\begin{proof} First, $\f_{D''-D}=\f_{D''}-\f_D\ge 0$ implies $D''-D\ge 0$, by Lemma~\ref{lem:PLvanishing}, and similarly $D''-D'\ge 0$. Furthermore,  $0<\f''(v_E)-\f'(v_E)=\sigma(v_E)(D''-D')$, and $E$ is thus in the support of $D''-D'$. Now pick any $k^\times$-invariant divisorial valuation $w$ on $\cX$ with center $\xi\in E\cap F$, normalized by $w(\cX_0)=1$; denote by $v\in X^\div$ its restriction to $k(X)\hookrightarrow k(\cX)$, so that $w=\sigma(v)$ (see~\S\ref{sec:Gauss}). Then $w(D''-D')=\f''(v)-\f'(v)>0$, and hence $\f''(v)=\f(v)$. This means that $D''=D$ at $\xi$; hence $0=\sigma(v_F)(D''-D)=\f''(v_F)-\f(v_F)$, which proves the result.  
\end{proof}

\begin{proof}[Proof of Theorem~\ref{thm:locality}]
  First assume that $\f,\p\in\PL\cap\PSH(\om)$. As in Lemma~\ref{lem:locality}, pick an integrally closed test configuration $\cX$ dominating the trivial one, such that $\f=\f_D$, $\f'=\f_{D'}$ and $\max\{\f,\f'\}=\f_{D''}$ with $D,D',D''\in\VCar(\cX)_\Q$. As in Proposition~\ref{prop:mixedMA}~(ii)  we have
\begin{equation*}
\MA(\f_D)
  =V^{-1}\sum_E b_E\left((\om_\cX+D )|_E\right)^n\d_{v_E}, 
\end{equation*}
and a similar formula holds for $\MA(\f_{D''})$. We thus need to show that $\left((\om_\cX+D )|_E\right)^n=\left((\om_\cX+D''])|_E\right)^n$ for any $E$ with $\f_D(v_E)>\f_{D'}(v_E)$. Now Lemma~\ref{lem:locality} yields $\f_D(v_F)=\max\{\f(v_F),\f'(v_F)\}=\f_{D''}(v_F)$ for all components $F$ of $\cX_0$ that intersect $E$. This implies that $D$ and $D''$ coincide in a neighborhood of $E$. Hence $D |_E=D''|_E$, which implies the result.

  Now consider the general case. Set 
  $$
  f:=\max\{\f-\f',0\}=\max\{\f,\f'\}-\f'.
  $$ 
  Then $\{\f>\f'\}=\{f>0\}$, and it thus suffices to prove that
  \begin{equation}\label{e602}
    f\MA(\max\{\f,\f'\})=f\MA(\f). 
  \end{equation}
  Indeed, multiplying by $f^{-1}$ on $\{f>0\}$ will then yield the result. 
  Let $(\f_j)_j$ and $(\f'_j)_j$ be decreasing nets in $\PL\cap\PSH(\om)$ converging
  to $\f$ and $\f'$, respectively, and set 
  $$
  f_j:=\max\{\f_j-\f'_j,0\}=\max\{\f_j,\f'_j\}-\f'_j.
  $$
  By the first part of the proof we have $f_j\MA(\max\{\f_j,\f'_j\})=f_j\MA(\f_j)$, and~\eqref{e602} is now a consequence of Proposition~\ref{prop:mixedMA} (iii), since $\max\{\f_j,\f'_j\}\in\cE^1(\om)$ decreases to $\max\{\f,\f'\}$.
  \end{proof}
Using Theorem~\ref{thm:locality} we obtain the very useful \emph{comparison principle}.
\begin{thm}\label{thm:compprinc}
  If $\om\in\Amp(X)$ and $\f,\p\in\cE^1(\om)$, then \[\int\limits_{\{\f<\p\}}\MA(\f)\ge\int\limits_{\{\f<\p\}}\MA(\p).\]
\end{thm}
\begin{proof}
  We follow~\cite[Theorem 1.5]{GZ2}.  For any $\e>0$ we have
  \begin{multline*}
    1
    =\int\MA(\max\{\f,\p-\e\})
    \ge\int\limits_{\{\f<\p-\e\}}\MA(\max\{\f,\p-\e\})
    +\int\limits_{\{\f>\p-\e\}}\MA(\max\{\f,\p-\e\})\\
    =\int\limits_{\{\f<\p-\e\}}\MA(\p)
    +\int\limits_{\{\f>\p-\e\}}\MA(\f)
    =1+\int\limits_{\{\f<\p-\e\}}\MA(\p)
    -\int\limits_{\{\f\le \p-\e\}}\MA(\f),
  \end{multline*}
  where the second equality follows from Theorem~\ref{thm:locality}
  and from $\MA(\f-\e)=\MA(\f)$.
  We complete the proof by letting $\e\to0$.
\end{proof}
Another simple consequence of Theorem~\ref{thm:locality} is the following
formula for the $\ii$-functional.
\begin{prop}\label{prop:Imaxgeod}
  If $\om\in\Amp(X)$, and $\f,\p\in\cE^1(\om)$, then $\max\{\f,\p\}\in\cE^1(\om)$, and
  \begin{equation}\label{equ:iipyth}
    \ii(\f,\p)=\ii\left(\f,\max\{\f,\p\}\right)+\ii\left(\max\{\f,\p\},\p\right).
  \end{equation}
\end{prop}
\begin{proof}
  The first assertion follows from Theorem~\ref{thm:E1}~(iii). For the second, note that
  \begin{align*}
    \ii(\f,\p)
    &=\int(\f-\p)(\MA(\f)-\MA(\p))\\
    &=\int\limits_{\{\f>\p\}}(\f-\p)(\MA(\f)-\MA(\p))
    +\int\limits_{\{\p>\f\}}(\f-\p)(\MA(\f)-\MA(\p))\\
    &=\int(\max\{\f,\p\}-\p)(\MA(\max\{\f,\p\})-\MA(\p))\\
    &+\int(\max\{\f,\p\}-\f)(\MA(\max\{\f,\p\})-\MA(\f))\\
    &=\ii\left(\f,\max\{\f,\p\}\right)+\ii\left(\max\{\f,\p\},\p\right),
  \end{align*}
where the third inequality follows from Theorem~\ref{thm:locality}.
\end{proof}

%
%
\subsection{Bedford--Taylor capacity of sublevel sets}\label{sec:BT1}
The following notion goes back to Bedford and Taylor~\cite{BT82} in the
complex-analytic case. A thorough study in our context will be conducted in Section~\ref{sec:neglpp}. Unless otherwise specified, $\om\in\Amp(X)$ is a fixed class. 

\begin{defi} The \emph{Bedford--Taylor capacity} of a Borel set $E\subset X^\an$ is defined by 
  \begin{equation}\label{equ:BTcap}
    \Capa(E)=\Capa_\om(E):=\sup\left\{\int_E\MA(\p)\mid \p\in\PSH, -1\le \p\le 0\right\}.
  \end{equation}
\end{defi}
Note that $0\le\Capa(E)\le 1$ for all Borel sets $E$, and that $\Capa(X^\an)=1$. 

\begin{lem}\label{lem:capsublevel1}
Every $\f\in\PSH_{\sup}(\om)$ satisfies 
\begin{itemize}
\item[(i)] $\Capa(\f\le-t)\le n t^{-1}$; 
\item[(ii)] $\jj(\f)\le 1+\int_1^\infty t^n\Capa(\f\le-t)\,dt$. 
\end{itemize}
Thus 
$$
\int_1^\infty t^n\Capa(\f\le-t)dt<\infty\Longrightarrow\f\in\cE^1(\om).
$$
\end{lem}
\begin{proof}
Pick any $\p\in\PSH(\om)$ with $-1\le\p\le 0$. Then
$$
    \int_{\{\f\le-t\}}\MA(\p)
    \le t^{-1}\int(-\f)\MA(\p)\le t^{-1} n
  $$
  where the last inequality follows from Lemma~\ref{lem:CLN}. 
  Taking the supremum over $\p$ gives~(i).
  
  To prove~(ii), set $\f_s:=\max\{\f,-s\}$ and $\mu_s:=\MA(\f_s)$ for $s\ge1$. By~\eqref{equ:Ebound} we have 
$$
    \jj(\f_s)
    \le\int(-\f_s)\mu_s
    =\int_0^s\mu_s\{\f_s\le-t\}\,dt
  $$
  $$  
    =\int_0^s\mu_s\{\f\le-t\}\,dt\\
    =\int_0^s\mu_t\{\f\le -t\}
    \le 1+\int_1^s\mu_t\{\f\le-t\}\,dt,
$$
  where the third equality holds since $\mu_s$ and $\mu_t$
  are probability measures that agree on the set $\{\f>-t\}$,
  by the locality principle, see Theorem~\ref{thm:locality}.
  For $t\ge 1$, we have $t^{-1}\f_t\in\PSH(\om)$ and $-1\le t^{-1}\f_t\le 0$,
  so $\Capa\ge\MA(t^{-1}\f_t)\ge t^{-n}\MA(\f_t)=t^{-n}\mu_t$. Thus
  $$
  \jj(\f_s)\le 1+\int_1^\infty\Capa(\f\le -t)dt,
  $$
 and~(ii) follows since $\jj(\f_s)\to\jj(\f)$ as $s\to\infty$.  
\end{proof}

\begin{cor}\label{cor:ppandE1}
For any $\f\in\PSH(\om)$ with $\f\le -1$ and $\a\in (0,\tfrac 1n)$ we have $-(-\f)^\a\in\cE^1(\om)$. In particular, every pluripolar set $E\subset X^\an$ is contained in $\{\p=-\infty\}$ for some $\p\in\cE^1(\om)$. 
\end{cor}
\begin{proof} The function $\chi(t)=-(-t)^\a$ is convex with $0\le\chi'(t)\le1$ on $(-\infty,-1]$, and $\p:=\chi(\f)$ is thus $\om$-psh, by Corollary~\ref{cor:pshcvx}. Further, Lemma~\ref{lem:capsublevel1}~(i) yields $\Capa\{\p\le-t\}=O(t^{-\a^{-1}})$, and hence $\p\in\cE^1(\om)$, by Lemma~\ref{lem:capsublevel1}~(ii). 
\end{proof}

\begin{cor}\label{cor:mixedMArestr} Pick $\om_i\in\Amp(X)$ and $\f_i\in\cE^1(\om_i)$, $i=1,\dots,n$, with mixed Monge--Amp\`ere measure $\mu:=(\om_1+\ddc\f_1)\winter(\om_n+\ddc\f_n)$. Then: 
\begin{itemize}
\item[(i)] $\mu$ puts no mass on pluripolar sets; 
\item[(ii)] for each irreducible component $X_\a$ of $X$ we have 
$$
\one_{X_\a}\mu=\left\{
\begin{array}{ll}
\left(\om_1|_{X_\a}+\ddc\f_1|_{X_\a^\an}\right)\winter\left(\om_n|_{X_\a}+\ddc\f_n|_{X_\a^\an}\right) & \text{ if }\dim X_\a=n\\
0 & \text{ otherwise}
\end{array}
\right. 
$$
\end{itemize}
\end{cor}
\begin{proof} By Proposition~\ref{prop:mixedMA}~(v), $\mu$ integrates all functions $\p\in\cE^1(\om)$, and (i) thus follows from Corollary~\ref{cor:ppandE1}. 

By Lemma~\ref{lem:Zarpp}, $X_\a^\an\cap X_{\a'}^\an$ is pluripolar for all $\a\ne\a'$, being a nowhere dense Zariski closed subset. By (i), we thus have $\mu(X_\a^\an\cap X_{\a'}^\an)=0$, and (ii) is now a consequence of Proposition~\ref{prop:mixedMA}~(ii). 
\end{proof}

%
%
%
%
  
 \section{Differentiability of the extended energy functional}\label{sec:endiff}
 As before, $X$ is a projective variety of dimension $n$, with irreducible components $X_\a$. Fix an ample class $\om\in\Amp(X)$. Above we considered the Monge--Amp\`ere energy functional $\en\colon\PSH(\om)\to\R\cup\{-\infty\}$, defined the set $\cE^1=\{\en>-\infty\}$, and associated a Monge-Amp\`ere measure  $\MA(\f)$ to each $\f\in\cE^1$. Formally, $\MA$ is the differential of $\en$, but $\cE^1$ is not a real vector space. In this section, we extend the Monge--Amp\`ere energy to various spaces of functions on $X^\an$, and prove quite precise differentiability results for these extensions.
  
%
%
\subsection{Extending the Monge--Amp\`ere energy}\label{sec:encont}
Unless stated otherwise, $\om\in\Amp(X)$ is an ample class, and we set $\PSH:=\PSH(\om)$, $\cE^1:=\cE^1(\om)$ etc.
The increasing functional 
$$
\en\colon\PSH\to\R\cup\{-\infty\}
$$ 
admits a natural extension to arbitrary functions $\f\colon X^\an\to\R\cup\{\pm\infty\}$ by setting 
\begin{equation}\label{equ:enmono}
\en(\f)=\en_\om(\f):=\sup\left\{\en(\p)\mid\p\in\PSH,\,\p\le\f\right\},
\end{equation}
 with the convention $\sup\emptyset=-\infty$.  
For $\p\in\PSH(\om)$, $\p\le\f\Leftrightarrow\p\le\env(\f)$, and hence 
\begin{equation}\label{equ:enenv}
\en(\f)=\en\left(\env(\f)\right). 
\end{equation}

We will mainly be concerned with the case when $\f$ is bounded, or even continuous. 
\begin{prop}\label{prop:tencont}
  The extended functional $\en\colon\Cz(X)\to\R$ satisfies the following properties: 
\begin{itemize}
\item[(i)] it is increasing, concave, and $1$-Lipschitz continuous; 
\item[(ii)] for each $\f\in\Cz(X)$, $c\in\R$ and $t\in\R_{>0}$ we have 
$$
\en(\f+c)=\en(\f)+c,\quad\en(t\cdot\f)=t\en(\f),\quad{and}\quad\en_{t\om}(t\f)=t\en_\om(\f); 
$$ 
\item[(iii)] for any $\f\in\Cz(X)$, $\om\mapsto\en_\om(\f)$ is continuous on $\Amp(X)$;  
\item[(iv)] for each $\f\in\Cz(X)$ we have 
\begin{equation}\label{equ:extensumcomp}
\en(\f)=\sum_\a c_\a\en(\f|_{X_\a})
\end{equation}
with $c_\a:=(\om^n)_{X_\a}/(\om^n)_X$.  
\end{itemize} 
\end{prop}

\begin{proof} Properties (i) and (ii) are straightforward consequences of Proposition~\ref{prop:E}.

To prove (iii), we may replace $\f$ with $\f-\inf\f$ and assume $\f\ge 0$, by (ii). Then
\begin{equation}\label{equ:tensuppos}
\en(\f)=\sup\left\{\en(\p)\mid\p\in\PSH,\,0\le\p\le\f\right\}.
\end{equation}
Now pick a convergent sequence $\om_i\to\om$ in $\Amp(X)$. We can then find a sequence $(t_i)$ in $\R_{>0}$ such that $t_i^{-1}\om\le\om_i\le t_i\om$ and $t_i\to 1$. By Proposition~\ref{prop:enmono}, any 
$$
\p\in\PSH(\om_i)\subset\PSH(t_i\om)
$$
such that $0\le\p\le\f$ satisfies $(\om_i,\p)^{n+1}\le \left(t_i\om,\p\right)^{n+1}$, and~\eqref{equ:tensuppos} thus yields
$$
\en_{\om_i}(\f)\le c_i\en_{t_i\om}(\f)=c_it_i\en_{\om}(t_i^{-1}\f)
$$
with $c_i:=t_i^n(\om^n)/(\om_i^n)\to 1$, where the right-hand equality follows from (ii). Similarly, 
$$
t_i^{-1}\en_\om(t_i\f)=\en_{t_i^{-1}\om}(\f)\le c'_i\en_{\om_i}(\f)
$$
with $c'_i:=t_i^n(\om_i^n)/(\om^n)\to 1$. Since $t_i^{-1}\f\to\f$ and $t_i\f\to\f$ uniformly, we conclude, as desired, $\en_{\om_i}(\f)\to\en_\om(\f)$ by continuity of $\en_\om$, proving (iii).

  The proof of (iv) is slightly more involved. Introduce $Y:=\coprod_\a X_\a$, with its canonical birational morphism $\pi\colon Y\to X$. The data of a family of functions $\p_\a\in\PSH(\om|_{X_\a})$ with $\p_\a\le\f|_{X_\a^\an}$ is equivalent to that of $\p\in\PSH(\pi^\star\om)$ with $\p\le\pi^\star\f$, and~\eqref{equ:ensumcomp} implies that the right-hand side of~\eqref{equ:extensumcomp} coincides with $\en_{\pi^\star\om}(\pi^\star\f)$. Since $\pi^\star\PSH(\om)\subset\PSH(\pi^\star\om)$, \eqref{equ:ensumcomp} further yields $\en_\om(\f)\le\en_{\pi^\star\om}(\pi^\star\f)$. 

Conversely, pick $\p\in\PSH(\pi^\star\om)$ such that $\p\le\pi^\star\f$. After replacing $\p$ with $\max\{\p,c\}$ for a constant $c\le\f$, we may assume that $\p$ is bounded.  Choose $\om_X\in\Amp(X)$ and $\f_X\in\PSH(\om_X)$ as in Theorem~\ref{thm:descentpsh}. Replacing $\om_X$ and $\f_X$ with small enough multiples, we may assume $\om_X\le\om$, and hence $\f_X\in\PSH(\om)$. Since $\f$ and $\p$ are both bounded, we may further arrange, after adding a constant to $\f_X$, that $\f_X\le\f$ and $\pi^\star\f_X\le\p$. For each $m\in\Z_{>0}$, we get $\f_m\in\PSH(\om)$ such that 
$$
\pi^\star\f_m=(1-\tfrac 1m)\p+\tfrac 1m\pi^\star\f_X.
$$
Since $\f_X\le\f$, we have $\pi^\star\f_m\le\pi^\star\f$, and hence $\f_m\le\f$, see Theorem~\ref{thm:suppsh}. On the other hand, since $\pi^\star\f_X\le\p$, Theorem~\ref{thm:suppsh} shows that $(\f_m)$ is an increasing sequence in $\PSH(\om)$, such that $\pi^\star\f_m\to\p$ pointwise on $Y^\div$. Now set $c:=\min\{\inf\f,\inf\p\}$ and $\f'_m:=\max\{\f_m,c\}$. Then $(\f'_m)$ is also an increasing sequence, this time in $\cE^\infty(\om)\subset\cE^1(\om)$, with $\pi^\star\f'_m\to\p$ pointwise on $Y^\div$. Further, $\f'_m\le\f$, and hence 
$$
\sum_\a c_\a\en_{\om|_{X_\a}}(\f'_m|_{X_\a^\an})=\en_\om(\f'_m)\le\en_\om(\f),
$$
by~\eqref{equ:ensumcomp}. Since $\f'_m|_{X_\a^\an}\to\p|_{X_\a^\an}$ pointwise on $X_\a^\div$, Theorem~\ref{thm:monostrong} yields  $\en_{\om|_{X_\a}}(\f'_m|_{X_\a^\an})\to\en_{\om|_{X_\a}}(\p|_{X_\a^\an})$, and hence 
$$
\en_{\pi^\star\om}(\p)=\sum_\a c_\a\en_{\om|_{X_\a}}(\p|_{X_\a^\an})\le\en_\om(\f).
$$
Taking the supremum over all $\p$, we conclude, as desired, $\en_{\pi^\star\om}(\pi^\star\f)\le\en_\om(\f)$. 
\end{proof}

%
 
\subsection{Further extensions}\label{sec:ensc}
Using the extension of the Monge--Amp\`ere energy $\en$ to continuous functions, we now go further and extend it to (upper or lower) semicontinuous functions.   
\begin{defi} For any $\f\colon X^\an\to\R\cup\{\pm\infty\}$ we set
\begin{equation}\label{equ:enminus}
\enh(\f):=\sup\left\{\en(\p)\mid\p\in\Cz(X),\,\p\le\f\right\}
\end{equation}
and
\begin{equation}\label{equ:enplus}
\enb(\f):=\inf\left\{\en(\p)\mid\p\in\Cz(X),\,\p\ge\f\right\}.
\end{equation}
\end{defi}
 
While $\enh$ and $\enb$ are defined for arbitrary functions, they are mainly of interest when restricted to lsc and usc functions, respectively.
 
\begin{prop}\label{prop:ensc} The functionals $\enh$ and $\enb$ satisfy the following  properties: 
\begin{itemize}
\item[(i)] they are increasing, concave, and satisfy the algebraic properties of Proposition~\ref{prop:tencont}~(ii); 
\item[(ii)] for any $\f\colon X^\an\to\R\cup\{\pm\infty\}$, 
$$
\enh(\f)\le\en(\f)\le\enb(\f),
$$
and $\enh(\f)>-\infty$ (resp.~$\enb(\f)<+\infty$) iff $\f$ is bounded below (resp.~bounded above); 
\item[(iii)] we have 
$$
\enh(\f)=\en(\qq(\f))=\sup\left\{\en(\p)\mid\rho\in\CPSH, \p\le\f\right\};
$$
we further have
$$
\enh(\f)=\en(\f_\star)=\en(\env(\f_\star))
$$ 
if $\f$ is bounded below; 
\item[(iv)] if $\f$ is bounded above, then 
\begin{equation*}
\enb(\f)=\enb(\f^\star)
\end{equation*}
we further have 
\begin{equation}\label{equ:enusc}
\enb(\f)=\en\left(\f^\star\right)=\en\left(\env(\f^\star)\right)
\end{equation}
if $\f\in\Cz(X)$, or $\f\in\PSH(\om)$, and for any $\f$ bounded above if the envelope property holds for $\om$; 
\item[(v)] $\enh$ (resp.~$\enb$) is continuous along increasing (resp.~decreasing) nets of bounded-below, lsc (resp.~bounded above, usc) functions. 
\end{itemize}

\end{prop}

\begin{proof} (i) and (ii) are obvious. We prove (iii). For any $\p\in\Cz(X)$, we have $\env(\p)=\qq(\p)$ by Lemma~\ref{lem:envlsc}, and~\eqref{equ:enenv} yields
$$
\en(\p)=\en(\env(\p))=\en(\qq(\p))=\sup\left\{\en(\rho)\mid\rho\in\CPSH,\,\rho\le\p\right\}. 
$$
This implies in turn
$$
\enh(\f)=\sup\left\{\en(\p)\mid\p\in\Cz(X),\,\p\le\f\right\}=\sup\left\{\en(\rho)\mid\rho\in\CPSH,\,\rho\le\f\right\}=\en(\qq(\f)),
$$
by monotonicity of $\en$. If $\f$ is bounded below, then $\qq(\f)=\env(\f_\star)$ by Lemma~\ref{lem:envlsc}, and hence $\en(\qq(\f))=\en(\env(\f_\star))=\en(\f_\star)$, which proves (iii). 

Next we show that $\enh$ is continuous along increasing nets of bounded-below, lsc functions. Let thus $(\f_i)$ be an increasing net of such functions, converging to the lsc function $\f=\sup_i\f_i$. On the one hand, $\f\ge\f_i$ implies
$$
\enh(\f)\ge S:=\sup_i\enh(\f_i)=\lim_i\enh(\f_i).
$$
On the other hand, for each $\e>0$ we can find $\p\in\Cz(X)$ such that $\p\le\f$ and $\en(\p)\ge\enh(\f)-\e$. Setting, for each $i$, $V_i:=\{\p<\f_i+\e\}$ defines an increasing family of open sets with $\bigcup_i V_i=X^\an$. By compactness, we get $V_i=X^\an$ for all $i$ large enough, \ie $\p\le\f_i+\e$ on $X^\an$, and hence $\en(\p)\le\enh(\f_i)+\e\le S+\e$.  We thus have $\enh(\f)\le S+2\e$ for all $\e>0$, and we get, as desired, $\enh(\f)=S$. 

The proof that $\enb$ is continuous along decreasing nets of bounded-above usc functions is entirely similar; hence (v).


If $\f\in\Cz(X)$ then~\eqref{equ:enusc} follows from~\eqref{equ:enenv}. If $\f\in\PSH(\om)$, write it as the limit of a decreasing net in $\CPSH(\om)$. Then $\en(\f_i)\to\enb(\f)$, by (v), and $\en(\f_i)\to\en(\f)$, by continuity of $\en$ along decreasing nets in $\PSH(\om)$. This proves~\eqref{equ:enusc} in that case. 

Now consider the general case, assuming the envelope property. Write $\f^\star$ as the decreasing limit of a net $(\f_i)$ in $\Cz(X)$. For each $i$ we have $\enb(\f_i)=\en(\env(\f_i))$. On the one hand, $\enb(\f_i)\to\enb(\f)$, by (v). On the other hand, Corollary~\ref{cor:envusc} implies that $\env(\f_i)$ decreases to $\env(\f^\star)$ in $\PSH(\om)$, and hence $\en(\env(\f_i))\to\en(\env(\f^\star))=\en(\f^\star)$, by continuity of $\en$ along decreasing nets in $\PSH(\om)$ and~\eqref{equ:enenv}. This concludes the proof of (iv). 
\end{proof}

\begin{cor}\label{cor:extensumcomp} If $\f\colon X^\an\to\R\cup\{-\infty\}$ is usc, then 
$$
\enb(\f)=\sum_\a c_\a\enb(\f|_{X_\a^\an}). 
$$
If $\f\colon X^\an\to\R\cup\{+\infty\}$ is lsc, then 
$$
\enh(\f)=\sum_\a c_\a\enh(\f|_{X_\a^\an}). 
$$
\end{cor}
\begin{proof} Assume $\f\colon X^\an\to\R\cup\{-\infty\}$ is usc, and pick a decreasing net $(\f_i)$ in $\Cz(X)$ such that $\f_i\searrow\f$ pointwise. For each $i$, \eqref{equ:extensumcomp} yields
$\en(\f_i)=\sum_\a c_\a \en(\f_i|_{X_\a^\an})$. By Proposition~\ref{prop:ensc}~(v), we have 
$$
\en(\f_i)=\enb(\f_i)\searrow\enb(\f)\quad\text{and}\quad \en(\f_i|_{X_\a^\an})=\enb(\f_i|_{X_\a^\an})\searrow\enb(\f|_{X_\a^\an}),
$$
and the first point follows. In the lsc case, the proof is similar, using an increasing net instead.
\end{proof}

%
\subsection{Uniform differentiability of the energy}\label{sec:diff}
Building on~\cite{BGM} we prove the following uniform differentiability result, which is crucial for what follows.  

\begin{thm}\label{thm:unifdiff} Pick $f\in\PL(X)_\R$. For any $\f\in\cE^1$, we then have 
$$
\enb(\f+\e f)=\en(\f)+\e\int f\MA(\f)+O(\e^2)
$$
as $\e\to 0$, where the implicit constant in the $O$ is uniform with respect to $\f$. 
\end{thm}
If $\f\in\CPSH(\om)$ (and for any $\f\in\cE^1$, if the envelope property holds), the left-hand side satisfies 
$$
\enb(\f+\e f)=\en(\f+\e f)=\en(\env(\f+\e f)),
$$
by~\eqref{equ:enusc}. 

\begin{cor}\label{cor:diff} For any $\f\in\cE^1$ and $f\in\Cz(X)$ we have 
$$
\frac{d}{dt}\bigg|_{t=0}\enb(\f+t f)=\int f\MA(\f).
$$
\end{cor}
Again, one can replace $\enb$ with $\en=\en\circ\env$ if $\f\in\CPSH(\om)$, or if the envelope property holds. 
 
For any $\om\in\Amp(X)$ we set $V_\om:=(\om^n)$. 

\begin{lem}\label{lem:BGM1} Pick $\om,\om'\in\Amp(X)$, $\p_1,\p_2\in\PL_\R\cap\PSH(\om')$, and set $f:=\p_1-\p_2$. For any $\f\in\CPSH(\om)$ we then have 
$$
\left|\en_\om(\f+f)-\en_\om(\f)-V_\om^{-1}\int f\,\left(\om+\om'+\ddc(\f+\p_1)\right)^n\right|\le\left(V_\om^{-1}V_{\om+\om'}-1\right)\sup|f|.
$$
\end{lem}

\begin{proof}
  We may assume $\om=c_1(L)$, $\om'=c_1(L')$ for ample line bundles $L,L'$ on $X$, Indeed, by homogeneity, we then obtain the case when $L,L'$ are ample $\Q$-line bundles, and the general case by a simple perturbation argument based on Proposition~\ref{prop:tencont}~(iii).

Consider first any continuous function $\rho\in\Cz(X)$. For each $m\in\N$, $\rho$ defines a sup-norm $\|\cdot\|_{m\rho}$ on $\Hnot(X,mL)$, defined by $$
\|s\|_{m\rho}=\sup_{X^\an}|s| e^{-m\rho}.
$$
This norm induces a norm $\det\|\cdot\|_{m\rho}$ on the determinant line $\det\Hnot(X,mL)$. Comparing with the trivial norm, we can and will think of $\det\|\cdot\|_{m\rho}$ as a number. It now follows from~\cite{CMac} (see~\cite[Theorem 9.5]{BE}) that the \emph{volume}
$$
\vol_L(\rho):=-\lim_{m\to\infty}\frac{1}{m h^0(mL)}\log\det\|\cdot\|_{m\rho}
$$
exists in $\R$.
We claim that $\vol_L(\rho)=\en_L(\rho)$.   Using~\cite[(9.7)]{BE} and~\eqref{equ:extensumcomp}, it is enough to prove this when $X$ is irreducible. When $\rho\in\CPSH(L)$, the equality holds by~\cite[Theorem~A]{BE} (or~\cite[Lemma 4.5]{nakstabold}). In the general case, this yields   
$$
\en_L(\rho)=\sup\left\{\en_L(\p)\mid\p\in\CPSH(L),\,\p\le\rho\right\}\le\vol_L(\rho),
$$ 
by monotonicity of $\vol_L$. On the other hand, Theorem~5.1.1 of~\cite{Reb21} (which is valid for arbitrary non-Archimedean fields; see also~\cite[Theorem 4.13]{nakstabold} for the trivially valued case) shows that
$$
\rho_m:=m^{-1}\log\sup_{s\in\Hnot(X,mL)\setminus\{0\}}\frac{|s|}{\|s\|_{m\rho}}
$$
satisfies $\rho_m\in\PL_\R\cap\PSH(\om)$, $\rho_m\le\rho$ and $\en_L(\rho_m)\to\vol_L(\rho)$; the claim follows. 

In particular, $\en_L(\f+f)=\vol_L(\f+f)$ and $\en_L(\f)=\vol_L(\f)$. The desired estimate is now a consequence of Lemma~\ref{lem:BGM2} below, itself a reformulation of~\cite[Lemma 3.2]{BGM}.
\end{proof}
\begin{lem}\label{lem:BGM2}
  Let $L$ and $L'$ be ample line bundles on $X$, and set $\om=c_1(L)$, $\om'=c_1(L')$.
  Consider functions $\f\in\CPSH(\om)$, and $\p_1, \p_2\in\CPSH(\om')$. Set $f:=\p_1-\p_2\in\Cz$ and $C:=V_{\om+\om'}-V_\om>0$. Then
  \begin{equation*}
    C\inf_{X^\an}f
    \le\int_{X^\an}f(\om+\om'+\ddc(\f+\p_1))^n
    -V_\om(\vol_L(\f+f)-\vol_L(\f))
    \le C\sup_{X^\an}f.
  \end{equation*}
\end{lem}
\begin{proof}
  This is a reformulation in our context of~\cite[Lemma 3.2]{BGM} in the special case of a trivially valued field $k$. The fact that $k$ is trivially valued means that we can view $\f$, $\p_1,\p_2$ as continuous psh metrics on the Berkovich analytifications of the line bundles $L$ and $L'$, respectively. The difference $\vol_L(\f+f)-\vol_L(\f)$ can be viewed as the relative volume of the metrics $\f+f$ and $\f$ on $L$. Note that our current normalization of the volume differs from the one in \emph{loc.cit.} by a factor $V_\om=(L^n)$. 
\end{proof}
\begin{proof}[Proof of Theorem~\ref{thm:unifdiff}] We argue along the lines of the proof of~\cite[Theorem 3.1]{BGM}. Write $\f\in\cE^1$ as the limit of a decreasing net $(\f_i)$ in $\CPSH(\om)$. Then $\en(\f_i+\e f)\to\enb(\f+\e f)$ (by Proposition~\ref{prop:ensc}~(v)), $\en(\f_i)\to\en(\f)$ (by continuity of $\en$ along decreasing nets in $\PSH(\om)$), and $\MA(\f_i)\to\MA(\f)$ weakly (by Proposition~\ref{prop:mixedMA}). 

It is therefore enough to prove the result for $\f\in\CPSH(\om)$. Since $f\in\PL(X)_\R$, we can choose $\om'\in\Amp(X)$ and $\p_1,\p_2\in\PL_\R\cap\PSH(\om')$ such that $f=\p_1-\p_2$. By Lemma~\ref{lem:BGM1} we have, for any $\e>0$,
\begin{multline*}
\left|\en_\om(\f+\e f)-\en_\om(\f)-\e V_\om^{-1}\int f\,\left(\om+\ddc\f+\e(\om'+\ddc\p_1)\right)^n\right|\\
\le\e\left(V_\om^{-1}V_{\om+\e\om'}-1\right)\sup|f|.
\end{multline*}
Now $V_\om^{-1}V_{\om+\e\om'}-1=O(\e)$, while 
$$
V_\om^{-1} \left(\om+\ddc\f+\e(\om'+\ddc\p_1)\right)^n=\MA_\om(\f)+V_\om^{-1}\sum_{j=1}^n{n\choose j}\e^j (\om+\ddc\f)^{n-j}\wedge (\om'+\ddc\p_1)^j,
$$
where $(\om+\ddc\f)^{n-j}\wedge (\om'+\ddc\p_1)^j$ is a Radon measure of mass $(\om^{n-j}\cdot \om'^j)$. The result follows. 
\end{proof}

\begin{proof}[Proof of Corollary~\ref{cor:diff}] For $f\in\PL(X)_\R$ the result follows directly from Theorem~\ref{thm:unifdiff}. For an arbitrary $f\in\Cz(X)$ we argue as in the proof of~\cite[Theorem 3.1]{BGM}. By Theorem~\ref{thm:PLdense}, we can pick a sequence $(f_m)$ in $\PL(X)_\R$ such that $\e_m:=\sup|f-f_m|\to 0$. For each $t>0$, we have $f_m-t\e_m\le f\le f_m+t\e_m$; Proposition~\ref{prop:ensc}~(i) thus yields 
$$
\enb(\f+t f_m)-t\e_m\le\enb(\f+t f)\le\enb(\f+t f_m)+t\e_m,
$$
and hence 
$$
\int f_m\MA(\f)-\e_m\le\liminf_{t\to 0_+} t^{-1}\left(\enb(\f+t f)-\enb(\f)\right)
$$
$$
\le\limsup_{t\to 0+} t^{-1}\left(\enb(\f+t f)-\enb(\f)\right)\le\int f_m\MA(\f)-\e_m, 
$$
by Theorem~\ref{thm:unifdiff} applied to $f_m$. Letting $m\to\infty$, we get
$$
\lim_{t\to 0+} t^{-1}\left(\enb(\f+t f)-\enb(\f)\right)=\int f\MA(\f),
$$
and replacing $f$ with $-f$ concludes the proof.
\end{proof}

%
%
%
%
 
\section{Measures of finite energy}\label{sec:M1}
 Denote as above by $X$ a projective variety of dimension $n$, and fix an ample class $\om\in\Amp(X)$.
 We define the Monge--Amp\`ere energy of a probability measure on $X^\an$, and begin a study of the space of measures of finite energy.
 
%
%
\subsection{The energy of a measure}
Denote by $\cM=\cM(X)$ the space of Radon probability measures on $X^\an$. It is a compact convex subset of the dual $\Cz(X)^\vee$ for the weak topology.
 Fix a class $\om\in\Amp(X)$. 

\begin{defi} The \emph{energy} of a Radon probability measure $\mu\in\cM$ is defined by 
\begin{equation}\label{equ:enmeas}
  \en^\vee(\mu)
  :=\en^\vee_\om(\mu)
  :=\sup_{\f\in\cE^1}\left(\en(\f)-\int\f\,\mu\right)\in[0,+\infty]. 
\end{equation}
We say that $\mu$ has \emph{finite energy} if $\en^\vee(\mu)<+\infty$, and denote by $\cM^1\subset\cM$ the set of such measures. 
\end{defi}
  Here we write $\cE^1=\cE^1(\om)$ for simplicity.  
As we shall see, measures in $\cM^1=\cM^1(\om)$ have a prescribed mass on each $X_\a^\an$, determined by $\om$ (see Corollary~\ref{cor:massMA}). However, when $X$ is irreducible, $\cM^1$ turns out to be independent of $\om\in\Amp(X)$ cf.~Theorem~\ref{thm:M1ind} below. 


By definition, measures of finite energy integrate all functions in $\cE^1$. Combined with Corollary~\ref{cor:ppandE1}, this implies: 

\begin{lem}\label{lem:M1npp} Measures in $\cM^1$ put no mass on pluripolar sets. 
\end{lem}
 
We will later study the dependence $\Amp(X)\ni\om\to\en^\vee_\om(\mu)$ for a fixed measure $\mu\in\cM$. Already now, we note that Proposition~\ref{prop:E}~(iii) implies
\begin{lem}\label{lem:enveehom}
  For $\om\in\Amp(X)$, $t\in\R_{>0}$ and $\mu\in\cM$, we have $\en^\vee_{t\om}(\mu)=t\en^\vee_\om(\mu)$.
\end{lem}

\begin{prop}\label{prop:Evee} The energy functional $\en^\vee\colon\cM\to[0,+\infty]$ satisfies:
\begin{itemize}
\item[(i)] for each $\mu\in\cM$ we have 
\begin{equation}\label{equ:enFS}
\en^\vee(\mu)=\sup_{\f\in\cH^\dom(\om)}\left(\en(\f)-\int\f\,\mu\right);
\end{equation}
\item[(ii)] $\en^\vee$ is convex, lsc, and homogeneous with respect to the scaling action of $\R_{>0}$, \ie  
$$
\en^\vee(t_\star\mu)=t\en^\vee(\mu)
$$ 
for all $\mu\in\cM$ and $t\in\R_{>0}$;
\end{itemize}
\end{prop}

\begin{proof} Denote by $S\in[0,+\infty]$ the right-hand side of~\eqref{equ:enFS}. Trivially, $\en^\vee(\mu)\ge S$. Conversely pick $\f\in\cE^1$, and choose a decreasing net $(\f_j)$ in $\cH^\dom(\om)$ such that $\f_j\to\f$. Then $\en(\f)\le\en(\f_j)$, and hence
$$
\en(\f)-\int\f\,\mu\le\left(\en(\f_j)-\int\f_j\,\mu\right)+\int(\f_j-\f)\mu\le S+\int(\f_j-\f)\mu. 
$$
By Lemma~\ref{lem:mononets}, $\int\f_j\,\mu\to\int\f\,\mu$; hence
$$
\en^\vee(\mu)=\sup_{\f\in\cE^1}\left(\en^\vee(\f)-\int\f\,\mu\right)\le S,
$$
which proves (i). 

Convexity and lower semicontinuity of $\en^\vee$ follow directly from (i), since $\mu\mapsto\en(\f)-\int\f\,\mu$ is affine and continuous for every $\f\in\PL\cap\PSH(\om)$. Now pick $\mu\in\cM$ and $t\in\R_{>0}$. Note that $\f\mapsto t\cdot\f=t\f(t^{-1}\cdot)$ is a bijection of $\cE^1$. By Proposition~\ref{prop:E} we infer
   \begin{multline*}
    \en^\vee(t_\star\mu)
    =\sup_{\f\in\cE^1}\left(\en(t\cdot\f)-\int(t\cdot\f) (t_\star\mu)\right)\\
    =\sup_{\f\in\cE^1}\left(t\en(\f)-\int t\f\mu\right)
    =t\sup_{\f\in\cE^1}\left(\en(\f)-\int\f\mu\right)
    =t\en^\vee(\mu),
  \end{multline*}
which concludes the proof of (ii). 
\end{proof}

\begin{defi} For each $C>0$ we set 
$$
\cM^1_C:=\left\{\mu\in\cM\mid\en^\vee(\mu)\le C\right\}.
$$
\end{defi}
By (ii) of Proposition~\ref{prop:Evee}, $\cM^1_C$ is a (weakly) compact, convex subset of $\cM$.

The next result provides an important source of examples of measures of finite energy. In fact, we will see in~\S\ref{sec:CY} that when the envelope property holds for $\om$, every measure $\mu\in\cM^1$ is of the form $\MA(\f)$ for $\f\in\cE^1$.  
\begin{prop}\label{prop:envar} For any $\f\in\cE^1$, the measure $\MA(\f)$ has finite energy, and $\f$ achieves the supremum in~\eqref{equ:enmeas}, \ie
\begin{equation}\label{equ:enma}
\en^\vee(\MA(\f))=\en(\f)-\int\f\MA(\f)=\ii(\f)-\jj(\f). 
\end{equation}
Furthermore,  
\begin{equation}\label{equ:enveeJ}
n^{-1}\jj(\f)\le\en^\vee\left(\MA(\f)\right)\le n\jj(\f),
\end{equation}
and 
\begin{equation}\label{equ:Eleg}
\en(\f)=\inf_{\mu\in\cM^1}\left(\en^\vee(\mu)+\int\f\,\mu\right),
\end{equation}
where the infimum is achieved for $\mu=\MA(\f)$. 
\end{prop}

\begin{proof}  By~\eqref{equ:Ebound}, we have, for each $\p\in\cE^1$,
$$
\en(\f)+\int(\p-\f)\MA(\f)\ge\en(\p),
$$
and hence 
$$
\en^\vee\left(\MA(\f)\right)=\sup_{\p\in\cE^1}\left(\en(\p)-\int\p\MA(\f)\right)=\en(\f)-\int\f\MA(\f).
$$
For any $\mu\in\cM$, we have $\en(\f)\le\en^\vee(\mu)+\int\f\,\mu$, with equality for $\mu=\MA(\f)$, proving~\eqref{equ:Eleg}. Finally~\eqref{equ:enveeJ} follows from~\eqref{equ:IJ}.
\end{proof}

\begin{exam} If $\mu=\MA(0)$, then $\en^\vee(\mu)=\ii(0)-\jj(0)=0$. Conversely, we will show in Corollary~\ref{cor:MA} that if $\en^\vee(\mu)=0$, then $\mu=\MA(0)$.
\end{exam}

%
%
 
\subsection{Legendre duality}\label{sec:moreleg}
The functional $\en^\vee\colon\cM\to\R\cup\{+\infty\}$ is defined as the Legendre transform of the Monge-Amp\`ere energy functional $\en\colon\cE^1\to\R$.
Here were prove a couple of additional duality formulations, involving the extensions of the Monge--Amp\`ere energy considered in~\S\ref{sec:endiff}
 
\begin{prop}\label{prop:encontleg}
  For any $\mu\in\cM$ we have
  \begin{equation}\label{equ:enveeleg}
    \en^\vee(\mu)=\sup_{\f\in\Cz(X)}\left(\en(\f)-\int\f\,\mu\right)
  \end{equation}
  and for any $\f\in\Cz(X)$ we have
  \begin{equation}\label{equ:encontleg}
    \en(\f)=\inf_{\mu\in\cM}\left(\en^\vee(\mu)+\int\f\,\mu\right)
  \end{equation}
\end{prop}

\begin{proof} Note that~\eqref{equ:enFS} implies
$$
\en^\vee(\mu)=\sup_{\f\in\PL\cap\PSH}\left(\en(\f)-\int\f\,\mu\right)\le\sup_{\f\in\Cz(X)}\left(\en(\f)-\int\f\,\mu\right).
$$
Now pick $\f\in\Cz(X)$ and $\p\in\cE^1$ such that $\p\le\f$. Then 
$$
\en(\p)-\int\f\,\mu\le\en(\p)-\int\p\,\mu\le\en^\vee(\mu),
$$
and taking the supremum over $\p$ yields
$$
\sup_{\f\in\Cz(X)}\left(\en(\f)-\int\f\,\mu\right)\le\en^\vee(\mu).
$$
 
Thus~\eqref{equ:enveeleg} holds. Now use this equation to define $\en^\vee(\mu)$ for any signed measure $\mu\in\Cz(X)^\vee$. If $\mu(1)\ne 1$, then $\en(c)=c$ for $c\in\R$ implies that 
$$
\en^\vee(\mu)\ge\sup_{c\in\R}\left(c(1-\mu(1))\right)=+\infty. 
$$
Similarly, if $\mu\in\Cz(X)^\vee$ and $\f\in\Cz(X)$ satisfy $\f\ge 0$ and $\int\f\,\mu<0$, then $\en(t\f)\ge 0$ for $t\ge 0$ yields
$$
\en^\vee(\mu)\ge\sup_{t\ge 0}\left(-t\int\f\,\mu\right)=+\infty.
$$ 
Thus $\en^\vee\equiv+\infty$ on $\Cz(X)^\vee\setminus\cM^1$. As $\en$ is concave,
Legendre duality now yields~\eqref{equ:encontleg}.
 
\end{proof}
The next duality statement will be used in the proof of Theorem~\ref{thm:preCY}.
\begin{prop}\label{prop:enuscleg}
  For any usc function $\f\colon X^\an\to\R\cup\{-\infty\}$, we have 
  \begin{equation}\label{equ:enuscleg}
    \enb(\f)=\inf_{\mu\in\cM}\left(\en^\vee(\mu)+\int\f\,\mu\right).
  \end{equation}
\end{prop}
\begin{proof}
Write $\f$ as the limit of a decreasing net $(\f_i)$ in $\Cz(X)$. On the one hand, for each $\mu\in\cM$ we have 
$$
\enb(\f)\le\en(\f_i)\le\en^\vee(\mu)+\int\f_i\,\mu,
$$
by~\eqref{equ:encontleg}, and hence $\enb(\f)\le\en^\vee(\mu)+\int\f\,\mu$, by Lemma~\ref{lem:mononets}. On the other hand, 
$$
\inf_\mu\left(\en^\vee(\mu)+\int\f\,\mu\right)\le\inf_\mu\left(\en^\vee(\mu)+\int\f_i\,\mu\right)=\en(\f_i),
$$
where the right-hand side converges to $\enb(\f)$, by Proposition~\ref{prop:ensc}~(v). This proves~\eqref{equ:enuscleg}. 
\end{proof}
 
%
%
\subsection{Maximizing nets}
With a view towards solving the Monge--Amp\`ere equation $\MA(\f)=\mu$, we introduce the following notion.

\begin{defi} Given a measure $\mu$ of finite energy, we shall say that a net $(\f_i)$ in $\cE^1$ is \emph{maximizing for $\mu$} if it computes $\en^\vee(\mu)=\sup_{\f\in\cE^1}\left(\en(\f)-\int\f\,\mu\right)$
  in the sense that
$$
\en(\f_i)-\int\f_i\,\mu\to\en^\vee(\mu).
$$ 
\end{defi} 

As a key consequence of the uniform differentiability property of Theorem~\ref{thm:unifdiff}, we prove: 

\begin{thm}\label{thm:preCY} Let $\mu$ be a measure of finite energy. For any maximizing net $(\f_i)$ for $\mu$ we have $\MA(\f_i)\to\mu$ weakly in $\cM$. 
\end{thm}
\begin{proof} It is enough to show $\int f\MA(\f_i)\to\int f\,\mu$ for any $f\in\PL(X)$, by Theorem~\ref{thm:PLdense}. Set $\d_i:=\en^\vee(\mu)-\en(\f_i)+\int\f_i\,\mu$, so that $\d_i\to 0$. For any $\e>0$,~\eqref{equ:enuscleg} shows that
$$
\enb(\f_i+\e f)-\int(\f_i+\e f)\,\mu\le\en^\vee(\mu)=\en(\f_i)-\int\f_i\,\mu+\d_i.
$$ 
By Theorem~\ref{thm:unifdiff}, we infer
$$
\e\int f\MA(\f_i)\le\e\int f\,\mu+\d_i+C\e^2
$$
for a constant $C>0$ independent of $i$ and $\e$. Dividing by $\e$, letting $i\to\infty$, and then $\e\to 0$, we get $\limsup_i\int f\MA(\f_i)\le\int f\,\mu$, and replacing $f$ with $-f$ yields the result. 
\end{proof}

This yields a variational characterization of solutions to Monge--Amp\`ere equations:

\begin{cor}\label{cor:MA} Pick $\f\in\cE^1$ and $\mu\in\cM^1$. Then $\MA(\f)=\mu$ iff $\f$ computes 
$\en^\vee(\mu)=\sup_{\p\in\cE^1}\left(\en(\p)-\int\p\,\mu\right)$. In particular, $\en^\vee(\mu)=0$ iff $\mu=\MA(0)$. 
\end{cor}
\begin{proof} If $\MA(\f)=\mu$ then $\en(\f)-\int\f\,\mu=\en^\vee(\mu)$ by~\eqref{equ:enma}. If the converse holds then the constant net $\f_i=\f$ is maximizing, and hence $\MA(\f)=\mu$, by Theorem~\ref{thm:preCY}. 
\end{proof}

As another consequence, we prove that measures of finite energy have a prescribed mass on the irreducible components of $X^\an$. 
\begin{cor}\label{cor:massMA} Pick $\mu\in\cM^1$. Then $\mu(X_\a)=(\om^n)_{X_\a}/(\om^n)_X$ for all $\a$, and $\mu(Z)=0$ for all Zariski closed subsets $Z\subset X$ with $\dim Z<n$.  
\end{cor}
\begin{proof} Set $c_\a:=(\om^n)_{X_\a}/(\om^n)_X$ (which vanishes if $\dim X_\a<n$). Pick a maximizing net $(\f_i)$ for $\mu$. By Theorem~\ref{thm:preCY}, $\mu_i:=\MA(\f_i)$ converges weakly to $\mu$. By Corollary~\ref{cor:mixedMArestr}~(ii), we have $\mu_i(X_\a)=c_\a$ for all $i$ and $\a$, and hence 
$$
\mu(X_\a)\ge\limsup_i \mu_i(X_\a)=c_\a.
$$
As in the proof of Corollary~\ref{cor:mixedMArestr}, when $\a\ne\a'$ $\mu$ further puts no mass on the pluripolar set $X_\a^\an\cap X_{\a'}^\an$, by Lemma~\ref{lem:M1npp}. Thus
$$
1=\mu(X)=\sum_\a\mu(X_\a)\ge\sum_\a c_\a=1, 
$$
which forces $\mu(X_\a)=c_\a$ for all $\a$. Finally let $Z\subset X$ be a Zariski closed subset with $\dim Z<n$, and pick an irreducible component $Y$ of $Z$. If $Y$ is not a component of $X$ then $Y^\an$ is pluripolar, and hence $\mu(Y)=0$, by Lemma~\ref{lem:M1npp}. If $Y=X_\a$ for some $\a$ then $\mu(Y)=c_\a$, which vanishes as well since $\dim X_\a<n$. Since this holds for all components of $Z$, we get, as desired, $\mu(Z)=0$. 
\end{proof}

\begin{defi}\label{defi:jmu} For any $\mu\in\cM^1$ we define a functional $\jj_\mu\colon\cE^1\to[0,+\infty)$ by 
$$
\jj_\mu(\f):=\en^\vee(\mu)-\en(\f)+\int\f\,\mu.
$$
\end{defi}
The notation is justified by the fact that for $\f,\p\in\cE^1$ we have
\begin{equation}\label{equ:jmuMA}
\jj_{\MA(\p)}(\f)=\en(\p)-\int\p\MA(\p)-\en(\f)+\int\f\MA(\p)=\jj_\p(\f).
\end{equation}
Note also that 
\begin{equation}\label{equ:jmuen}
\jj_\mu(0)=\en^\vee(\mu),
\end{equation}
and that a net $(\f_i)$ in $\cE^1$ is maximizing for $\mu$ iff $\jj_\mu(\f_i)\to 0$. 

\begin{lem}\label{lem:Iminmax} For all $\f,\f'\in\cE^1$ we have 
$$
\ii(\f,\f')\approx\inf_{\mu\in\cM^1}(\jj_\mu(\f)+\jj_\mu(\f')). 
$$
In particular, we have for all $\f\in\cE^1$ and $\mu\in\cM^1$
\begin{equation}\label{equ:jjmu}
\jj(\f)\lesssim\jj_\mu(\f)+\en^\vee(\mu).
\end{equation} 
\end{lem}
\begin{proof} On the one hand, $\inf_\mu(\jj_\mu(\f)+\jj_\mu(\f'))\le\jj_\f(\f')\approx\ii(\f,\f')$, thanks to~\eqref{equ:IJ}. Conversely, set $\tau:=\tfrac 12(\f+\f')$. By strict concavity of $\en$ (Theorem~\ref{thm:strictconc}), 
$$
\en^\vee(\mu)\ge\en(\tau)-\int\tau\,\mu\ge\frac 12(\en(\f)+\en(\f'))-\frac 12\int(\f+\f')\,\mu+C_n^{-1}\ii(\f,\f'), 
$$
and hence 
$$
\inf_\mu(\jj_\mu(\f)+\jj_\mu(\f'))\gtrsim\ii(\f,\f').
$$
In particular, $\jj(\f)\approx\ii(\f,0)\lesssim \jj_\mu(\f)+\jj_\mu(0)=\jj_\mu(\f)+\en^\vee(\mu)$, which proves the final point. 
\end{proof}

\begin{cor}\label{cor:enmaxim} Pick $\mu\in\cM^1$, and choose a maximizing net $(\f_i)$ for $\mu$. Then $\jj(\f_i)\lesssim\en^\vee(\mu)+o(1)$ for all $i$ large enough.
\end{cor}


\begin{lem}\label{lem:BBGZ2} For each $C>0$ we have 
$$
\left|\int\f\left(\mu-\MA(\p)\right)\right|\lesssim C\jj_\mu(\p)^{\frac 12}
$$
for all $\f,\p\in\{\jj\le C\}\subset\cE^1$ and $\mu\in\cM^1_C$. 
\end{lem}

\begin{proof} By monotone approximation, we may assume $\f\in\PL_\R$. Pick a maximizing sequence $(\p_i)$ for $\mu$. By Corollary~\ref{cor:enmaxim}, $\jj(\p_i)\lesssim C+o(1)$ for all $i$ large enough. By Lemma~\ref{lem:BBGZ}, we infer
$$
\left|\int\f(\MA(\p_i)-\MA(\p))\right|\lesssim (C+o(1))\ii(\p_i,\p)^{\frac 12}\lesssim C\max\{\jj_\mu(\p_i),\jj_\mu(\p)\}^{\frac 12},
$$
where the right-hand inequality follows from Lemma~\ref{lem:Iminmax}. By Theorem~\ref{thm:preCY}, we have $\MA(\f_i)\to\mu$ weakly; hence $\int\f\MA(\p_i)\to\int\f\,\mu$, while $\jj_\mu(\p_i)\to 0$, so we are done. 
\end{proof}

\begin{cor}\label{cor:maxen} Pick $\mu\in\cM^1$. If $(\f_i)$ is maximizing for $\mu$ then $\lim_i\int\f_i(\MA(\f_i)-\mu)=0$ and $\lim_i\en^\vee(\MA(\f_i))=\en^\vee(\mu)$. 
\end{cor}
\begin{proof} By Corollary~\ref{cor:enmaxim}, $\jj(\f_i)$ is bounded for $i$ large enough. By Lemma~\ref{lem:BBGZ2}, it follows that 
$$
c_i:=\int\f_i\mu-\int\f_i\MA(\f_i)=O(\jj_\mu(\f_i)^{\frac 12})
$$
tends to $0$, and hence 
$$
\en^\vee(\MA(\f_i))=\en(\f_i)-\int\f_i\MA(\f_i)
=\en(\f_i)-\int\f_i\,\mu+c_i\to\en^\vee(\mu),
$$
completing the proof.
\end{proof}

As another useful consequence, we get: 

\begin{prop}\label{prop:weakcont} Let $\f_i\to\f$ be a convergent net in $\PSH(\om)$, and assume that $\jj(\f_i)$ is eventually bounded. For any $\mu\in\cM^1$ we then have $\f_i\to\f$ in $L^1(\mu)$. 
\end{prop}

\begin{proof}  We first claim that $\p\mapsto\int\p\,\mu$ is continuous on $\{\jj\le C\}\subset\PSH(\om)$ for any $C>0$. To see this, pick a maximizing net $(\rho_j)_j$ for $\mu$ in $\PL\cap\PSH(\om)$. For each $j$, $\f\mapsto\int\f\MA(\rho_j)$ is continuous on $\PSH(\om)$, by Proposition~\ref{prop:mixedMA}~(ii). By Lemma~\ref{lem:BBGZ2}, $\int\p\MA(\rho_j)\to\int\p\,\mu$ uniformly for $\p\in\{\jj\le C\}\subset\PSH(\om)$, and the claim follows since continuity is preserved under uniform convergence. 

  The functional $\jj$ is lsc, so $\jj(\f)\le\liminf\jj(\f_i)<\infty$. Set $\f'_i:=\max\{\f_i,\f\}$, so that $|\f_i-\f|=2(\f'_i-\f)+(\f_i-\f)$. Then $\f'_i\to\f$ pointwise on $X^\div$, \ie $\f'_i\to\f$ in $\PSH(\om)$. By~\eqref{equ:iipyth} we further have $\ii(\f'_i,\f)\le\ii(\f_i,\f)$, and $\jj(\f'_i)\approx\ii(\f'_i,0)$ is thus eventually bounded, by the quasi-triangle inequality for $\ii$. The first part of the proof thus yields $\int\f_i\,\mu\to\int\f\,\mu$, $\int\f'_i\,\mu\to\int\f\,\mu$, and hence $\int|\f_i-\f|\mu\to 0$, which concludes the proof. 
\end{proof}

%
%
\subsection{Dependence on the ample class}\label{sec:energyandclass}
In this section, $X$ is assumed to be \textbf{irreducible}. As in~\cite[Proposition~3.4]{BBGZ} we then have the following useful characterization of measures of finite energy. 

\begin{thm}\label{thm:charM1}
 A Radon probability measure $\mu\in\cM$ has finite energy iff $\cE^1\subset L^1(\mu)$.  
 \end{thm}
 
\begin{proof} Assume that $\mu$ is finite on $\cE^1$. We first claim that for each $C>0$ there exists $C'>0$ such that $\int\f\,\mu\ge -C'$ for all $\f\in\cE^1_{\sup}$ such that $-\en(\f)=\jj(\f)\le C$. Arguing by contradiction, we may assume 
  there exist $C>0$ and a sequence $(\f_j)_1^\infty$ in $\cE^1_{\sup}$ 
  such that $\jj(\f)\le C$ and $\int\f_j\mu\le-2^j$ for all $j$. 
  Set $\p_m:=\sum_{j=1}^m2^{-j}\f_j$ for $m\ge1$.
  Then $\p_m$ is a decreasing sequence in $\PSH(\om)$ converging to 
  $\p:=\sum_{j=1}^\infty2^{-j}\f_j$, which is thus either $\om$-psh or identically $-\infty$ (since $X$ is irreducible).
  By concavity of $\en$ we have $\en(\p_m)\ge-C$ for all $m$; 
  hence $\p$ is $\om$-psh, with $\en(\p)=\lim_{m\to\infty}\en(\p_m)\ge-C$. 
  On the other hand, monotone convergence gives 
  \begin{equation*}
    \int\p\,\mu
    =\lim_{m\to\infty}\int\p_m\,\mu
    =\lim_{m\to\infty}\sum_{j=1}^m2^{-j}\int\f_j\,\mu
    =-\infty,
  \end{equation*}
a contradiction. We next claim that there exist $A,B>0$ such that 
$$
\left|\int\f\,\mu\right|\le A\jj(\f)^{1/2}+B
$$
for all $\f\in\cE^1_{\sup}$, which will imply, as desired, that
$$
\en^\vee(\mu)=\sup_{\f\in\cE^1_{\sup}}\left(-\int\f\,\mu-\jj(\f)\right)
$$
is finite. If $\jj(\f)\le 1$ then the estimate follows from the first part of the proof. We may thus assume $t:=\jj(\f)^{-1/2}\le 1$. By Lemma~\ref{lem:quadJ}, $\jj(t\f)\lesssim  t^2\jj(\f)=1$, and the first part of the proof yields a uniform constant $A>0$ such that $t\left|\int\f\,\mu\right|\le A$, which gives the desired estimate. 
\end{proof}

As in~\cite[Proposition 4.1]{DiN16}, we infer: 

\begin{cor}\label{cor:M1ind} The set $\cM^1(\om)$ of measures of finite energy is independent of the choice of $\om\in\Amp(X)$, and it is dense in the space $\cM$ of all Radon probability measures on $\Xan$.
\end{cor}
%
%
\begin{proof} For any $\om,\om'\in\Amp(X)$, we can find $s\gg 1$ such that $s^{-1}\om\le\om'\le s\om$. By Theorem~\ref{thm:E1} we then have $s^{-1}\cE^1(\om)\subset\cE^1(\om')\subset s\,\cE^1(\om)$, and $\mu$ is thus finite on $\cE^1(\om)$ iff it is finite on $\cE^1(\om')$. We conclude by Theorem~\ref{thm:charM1} that $\cM^1(\om)$ is independent of $\om$.

For the second part,  we note that finite atomic measures are dense in $\cM$. As $\Xdiv$ is dense in $\Xan$, finite atomic measures with support in $\Xdiv$ are also dense in $\cM$, and Theorem~\ref{thm:charM1} shows that any such measure has finite energy. Indeed, any function in $\PSH(\om)$ is finite on $\Xdiv$.
\end{proof}

\begin{rmk} In view of Corollary~\ref{cor:massMA}, the above results fail when $X$ has more than one top-dimensional component. 
\end{rmk}

By Proposition~\ref{prop:mixedMA}~(v) and Theorem~\ref{thm:charM1}, we also have: 

\begin{cor}\label{cor:M1mixed} For $i=1,\dots,n$, pick $\om_i\in\Amp(X)$ and $\f_i\in\cE^1(\om_i)$. Then the Radon probability measure 
$$
(\om_1\inter\om_n)^{-1}(\om_1+\ddc\f_1)\wedge\dots\wedge(\om_n+\ddc\f_n)
$$
has finite energy. 
\end{cor}

 For further reference, we establish a more precise version of Corollary~\ref{cor:M1ind}. 
 
\begin{thm}\label{thm:M1ind} Pick $\om,\om'\in\Amp(X)$ and $s\ge 1$ such that $s^{-1}\om\le\om'\le s\om$. Then 
$$
s^{-C_n}\en^\vee_\om\le\en^\vee_{\om'}\le s^{C_n}\en^\vee_{\om}
$$ 
on $\cM$, with $C_n:=1+2n^2$.  
\end{thm}
The main ingredient in the proof is the following estimate.
 
\begin{lem}\label{lem:enandtheta}
  Suppose $\om,\om'\in\Amp(X)$ and $\om\le\om'\le s\om$, where $s\ge1$. For any nonpositive $\f\in\PSH(\om)\subset\PSH(\om')$, we then have 
 \begin{equation}\label{equ:enandtheta}
  0\ge s^{-n}\en_\om(\f)\ge\en_{\om'}(\f)\ge s^n\en_{\om}(\f), 
  \end{equation}
whereas
  \begin{equation}\label{equ:enveeandtheta}
    \en^\vee_{\om'}(\mu)\ge\left((n+1)-n s^n\right)\en^\vee_{\om}(\mu)
  \end{equation}
for all $\mu\in\cM$. 
\end{lem}
 
Of course,~\eqref{equ:enveeandtheta} is useful only when $s\le (1+\tfrac 1n)^{1/n}$. In order to upgrade this to the global estimate in Theorem~\ref{thm:M1ind}, we use the \emph{Thompson metric} $\d_T$ of the open convex cone $\Amp(X)\subset\Num(X)$. As in~\cite{Tho}, this is defined by
$$
\d_T(\om,\om')=\sup\{\d\in\R\mid e^{-\d}\om\le\om'\le e^\d\om\}.
$$
It is not hard to see that $(\Amp(X),\d_T)$ is a complete metric space in which line segments are (constant speed) geodesics, see~\cite[Lemma~2.6.2]{LN}.
Theorem~\ref{thm:M1ind} is equivalent to the $C_n$-Lipschitz continuity of $\om\mapsto\log\en^\vee_\om(\mu)$ for any $\mu\ne\d_{v_\triv}$. 

\begin{lem}\label{lem:lip} Let $(Z,d)$ be a geodesic metric space. Pick a function $f\colon Z\to\R$, and suppose we are given $\e>0$ and $\rho:[0,\e]\to[0,+\infty)$ of class $C^1$ with $\rho(0)=0$,  such that for all $x,y\in Z$ we have 
\begin{equation}\label{equ:smalllip}
d(x,y)\le\e\Longrightarrow |f(x)-f(y)| \le \rho\left(d(x,y)\right)
\end{equation}
Then $f$ is Lipschitz continuous, with Lipschitz constant $|\rho'(0)|$. 
\end{lem}
\begin{proof}[Proof of Theorem~\ref{thm:M1ind}] Set $\d:=\d_T(\om,\om')$. Then $e^{-\d}\om\le\om'\le e^{\d}\om=e^{2\d}(e^{-\d}\om)$, so Lemma~\ref{lem:enveehom} and~\eqref{equ:enveeandtheta} yield
$\en^\vee_{\om'}(\mu)\ge e^{-\d}((n+1)-n e^{2n\d})\en^\vee_{\om}(\mu)$. As a result, the function $f\colon\Amp(X)\to\R$ defined by $f(\om):=\log\en^\vee_\om(\mu)$ satisfies the assumptions of Lemma~\ref{lem:lip} with $\e:=\frac{1}{2n}\log(1+\frac{1}{n})$ and $\rho\colon[0,\e]\to[0,+\infty)$ defined by 
$$
\rho(t)=t-\log\left((n+1)-n e^{2nt}\right).
$$
Now $\rho'(0)=C_n:=1+2n^2$. Lemma~\ref{lem:lip} shows that $f$ is $C_n$-Lipschitz continuous, which is equivalent to the desired estimate. 
\end{proof}
 
\begin{proof}[Proof of Lemma~\ref{lem:enandtheta}]
  By Lemma~\ref{lem:ennef}, we have 
$$
0\ge(\om,\f)^{n+1}\ge(\om',\f)^{n+1}\ge s^n (\om,\f)^{n+1},
$$
and~\eqref{equ:enandtheta} follows since $(\om^n)\le(\om'^n)\le s^n(\om^n)$. 

Now pick $\mu\in\cM$. If $\mu\notin\cM^1$, then $\en^\vee_{\om'}(\mu)=+\infty$, and~\eqref{equ:enveeandtheta} is trivial. Now assume $\mu\in\cM^1$, and pick a maximizing sequence $(\f_i)$ in $\cE^1(\om)$ for $\mu$, normalized by $\sup\f_i=0$. If we set  $\mu_i:=\MA_\om(\f_i)$, then $\en_\om(\f_i)-\int\f_i\,\mu$ by definition, and 
$\int\f_i\,\mu=\int\f_i\,\mu_i+o(1)$, by Corollary~\ref{cor:maxen}. Thus
  \begin{align*}
    \en^\vee_{\om'}(\mu)
    &\ge\en_{\om'}(\f_i)-\int\f_i\,\mu\\
    &\ge s^n\en_{\om}(\f_i)-\int\f_i\,\mu\\
    &= s^n(\en_\om(\f_i)-\int\f_i\,\mu)+(s^n-1)\int\f_i\,\mu\\
    &=s^n\en^\vee_\om(\mu)+(s^n-1)\int\f_i\,\mu_i+o(1),
   \end{align*}
  where the first inequality is definitional, and the second follows from~\eqref{equ:enandtheta}. Finally,~\eqref{equ:enma},~\eqref{equ:enveeJ} and~\ref{equ:IJ} yield 
  $$
  -\int\f_i\,\mu_i=\ii_\om(\f_i)\le (n+1)\en^\vee_\om(\mu_i)=(n+1)\en^\vee_\om(\mu)+o(1),
  $$
  and~\eqref{equ:enveeandtheta} follows. 
\end{proof}
 
\begin{proof}[Proof of Lemma~\ref{lem:lip}]
  Pick $C>|\rho'(0)|$, and choose $\eta\in(0,\e]$ such that $|\rho'(\d)|\le C$ for $\d\in[0,\eta]$, and hence $0\le \rho(\d)\le C\d$ for $\d\in[0,\eta]$. For all $x,y\in Z$ we thus have  
$$
d(x,y)\le\eta\Longrightarrow |f(x)-f(y)|\le C d(x,y). 
$$
Now pick $x,y\in Z$ at arbitrary distance, and choose a geodesic $\g\colon[0,1]\to Z$ connecting $x$ to $y$. By compactness of $[0,1]$, we can find a chain $0=t_0<t_1<\dots<t_N=1$ in $[0,1]$ such that the $x_i:=\g(t_i)$ satisfy $d(x_i,x_{i+1})\le\eta$ for $i<N$. Since $\g$ is a geodesic, we have $d(x,y)=\sum_{i<N} d(x_i,x_{i+1})$. Now~\eqref{equ:smalllip} yields $|f(x_i)-f(x_{i+1})|\le C d(x_i,x_{i+1})$ for all $i<N$, and hence 
$$
|f(x)-f(y)|\le\sum_{i<N}|f(x_i)-f(x_{i+1})|\le C d(x,y).  
$$
This holds for any $C>|\rho'(0)|$, and the result follows. 
\end{proof}


 
%
%
%
%
  
\section{The strong topology on $\cM^1$}\label{sec:M1strtop}
As in \S\ref{sec:M1}, we denote by $X$ a projective variety of dimension $n$, and fix $\om\in\Amp(X)$. The set $\cM^1\subset\cM$ of measures of finite energy comes equipped with the weak topology of Radon probability measures. Here we introduce and study a stronger topology on $\cM^1$.
   
  %
%
%
%
\subsection{A quasimetric on $\cM^1$} 
Recall from Definition~\ref{defi:jmu} that 
$$
\jj_\mu(\f)=\en^\vee(\mu)-\en(\f)+\int\f\,\mu\ge 0
$$
for all $\mu\in\cM^1$ and $\f\in\cE^1$. Dualizing Lemma~\ref{lem:Iminmax} we introduce: 

\begin{defi} For any two $\mu,\mu'\in\cM^1$ we set
$$
\ii^\vee(\mu,\mu'):=\inf_{\f\in\cE^1}(\jj_\mu(\f)+\jj_{\mu'}(\f)).
$$
\end{defi}


\begin{thm}\label{thm:trianglevee} The functional $\ii^\vee$ is a quasi-metric on $\cM^1$. Furthermore, 
\begin{equation}\label{equ:iveeen}
\ii^\vee(\mu,\MA(0))\approx\en^\vee(\mu)
\end{equation}
and
\begin{equation}\label{equ:ivjmu}
\ii^\vee\left(\mu,\MA(\f)\right)\approx\jj_\mu(\f)
\end{equation}
for all $\mu\in\cM^1$ and $\f\in\cE^1$. 
\end{thm}
Theorem~\ref{thm:trianglevee} will be proved below, together with the following crucial estimates. 

\begin{thm}\label{thm:BBGZ3} For all $\f,\f'\in\cE^1(\om)$ and $\mu,\mu'\in\cM^1$ we have 
$$
\left|\int(\f-\f')\left(\mu-\mu'\right)\right|\lesssim \ii(\f,\f')^{\a_n}\ii^\vee(\mu,\mu')^{\frac 12}\max\{\jj(\f),\jj(\f'),\en^\vee(\mu),\en^\vee(\mu')\}^{\frac 12-\a_n}
$$
with $\a_n:=2^{-n}$. In particular, 
$$
0\le\int\f\MA(0)-\int\f\,\mu\lesssim\max\{1,\jj(\f)\}^{\frac 12}\max\{1,\en^\vee(\mu)\}^{1-\a_n}.
$$
\end{thm}
This last estimate should be compared with the trivial bound
$$
0\le\int\f\MA(0)-\int\f\,\mu\le\jj(\f)+\en^\vee(\mu),
$$
which holds by definition of $\en^\vee$. 

As an important consequence of Theorem~\ref{thm:BBGZ3}, we have:

\begin{cor}\label{cor:MAinj}
Assume that $X$ is equidimensional. If $\f,\p\in\cE^1$, then the following properties are equivalent:
  \begin{itemize}
  \item[(i)]
    $\MA(\f)=\MA(\p)$;
  \item[(ii)]
    $\ii(\f,\p)=0$;
 \item[(iii)]
   $\f-\p$ is locally constant.  
    \end{itemize}
\end{cor}
 
\begin{proof} The implications~(iii)$\Rightarrow$(i)$\Rightarrow$(ii) are clear. If $\ii(\f,\p)=0$, then $\ii\left(\f|_{X_\a^\an},\p|_{X_\a^\an}\right)=0$ for each irreducible component $X_\a$ of $X$, by~\eqref{equ:Isumcomp}.
In order to prove that $\f-\p$ is locally constant, we may thus assume that $X$ is irreducible. After adding constant, assume $\sup\f=\sup\p$. Then $\f(v)=\p(v)$ for every $v\in X^\div$, as a consequence of Theorem~\ref{thm:BBGZ3} with $\mu=\delta_v$ and $\mu'=\delta_{v_\triv}$, since $\d_v\in\cM^1$ by Proposition~\ref{prop:enT}. By Theorem~\ref{thm:suppsh}, we infer $\f=\p$ on $X^\an$, which proves (ii)$\Rightarrow$(iii). 
\end{proof}

Corollary~\ref{cor:MAinj} in turn implies the following useful results. 

\begin{cor}\label{cor:sameenergy}
Assume that $X$ is equidimensional. If $\f,\p\in\cE^1$ satisfy $\f\ge\p$ and $\en(\f)=\en(\p)$, then $\f=\p$.
\end{cor}  
\begin{proof} By~\eqref{equ:intsumcomp}, we may assume that $X$ is connected. From~\eqref{equ:Ebis} and the assumptions, it follows that
 $$
\int(\f-\p)\MA(\f)=\int(\f-\p)\MA(\p)=0.
$$
Thus $\ii(\f,\p)=\int(\f-\p)(\MA(\p)-\MA(\f))=0$, and hence $\p=\f+c$ for a constant $c$, by Corollary~\ref{cor:MAinj}. But then $\en(\f)=\en(\p)=0$ gives $c=0$ and we are done. 
\end{proof}

\begin{cor}[Domination principle]\label{cor:domprinc} Assume that $X$ is equidimensional. Let $\f\in\PSH(\om)$, $\p\in\cE^1$, and assume $\f\le \p$ a.e.\ with respect to $\MA(\p)$. Then $\f\le \p$ on $\Xan$.
\end{cor}
\begin{proof} Arguing on each connected component, we may assume that $X$ is connected. After replacing $\f$ with $\max\{\f,\p\}$, we may assume $\f\ge\p$, and hence $\f\in\cE^1$. The assumption then becomes $\f=\p$ a.e for $\MA(\p)$, and we need to show $\f=\p$. Note that 
$$
\ii(\f,\p)=\int(\f-\p)(\MA(\p)-\MA(\f))=\int(\p-\f)\MA(\f)\le 0.
$$
By Corollary~\ref{cor:MAinj}, we infer $\f=\p+c$ with $c\in\R$, and $\f=\p$ a.e.~for $\MA(\p)$ implies $c=0$. 
\end{proof}
%
%
%
%
\subsection{Proof of Theorems~\ref{thm:trianglevee} and~\ref{thm:BBGZ3}}

\begin{lem}\label{lem:trianglevee} If $\p,\p'\in\cE^1$ and $\mu=\MA(\p)$, $\mu'=\MA(\p')$, then 
$$
\ii(\p,\p')\approx\ii^\vee(\mu,\mu').
$$
\end{lem}
\begin{proof} By~\eqref{equ:jmuMA} we have 
$$
\ii^\vee(\mu,\mu')=\inf_{\f\in\cE^1}(\jj_\p(\f)+\jj_{\p'}(\f)). 
$$
By~\eqref{equ:IJ}, this implies
$$
\ii^\vee(\mu,\mu')\approx \inf_{\f\in\cE^1}(\ii(\p,\f)+\ii(\p',\f)),
$$
and the result follows thanks to the quasi-triangle inequality for $\ii$, see~\eqref{equ:Itriangle}. 
\end{proof}

\begin{lem}\label{lem:Iveeunif} For each $C>0$ there exists $C'\lesssim C$ such that 
$$
\ii^\vee(\mu,\mu')=\inf\left\{\jj_\mu(\f)+\jj_{\mu'}(\f)\mid\f\in\cE^1,\,\jj(\f)\le C'\right\}
$$
for all $\mu,\mu'\in\cM^1_C$. 
\end{lem}
\begin{proof} By definition of $\ii^\vee$ we have   
$$
\ii^\vee(\mu,\mu')\le\jj_\mu(0)+\jj_{\mu'}(0)=\en^\vee(\mu)+\en^\vee(\mu')\le 2C,
$$
thanks to~\eqref{equ:jmuen}. Thus 
$$
\ii^\vee(\mu,\mu')=\inf\left\{\jj_\mu(\f)+\jj_{\mu'}(\f)\mid\f\in\cE^1_{\sup},\,\jj_\mu(\f)\le 3C\right\}.
$$
Now $\jj_\mu(\f)\le 3C$ implies $\jj(\f)\lesssim\jj_\mu(\f)+\en^\vee(\mu)\le 4C$, by~\eqref{equ:jjmu}, and we are done. 
\end{proof}

\begin{lem}
  Given $\mu,\mu'\in\cM^1$ with maximizing nets $(\p_i)$, $(\p'_i)$, respectively, we have
\begin{equation}\label{equ:convIvee}
  \ii^\vee\left(\MA(\p_i),\MA(\p'_i)\right)\to\ii^\vee(\mu,\mu')
\end{equation}
and 
\begin{equation}\label{equ:convjmu}
\jj_{\p_i}(\f)\to\jj_\mu(\f)
\end{equation}
for all $\f\in\cE^1$. 
\end{lem}
\begin{proof}
By Corollary~\ref{cor:maxen}, the measures $\MA(\p_i),\MA(\p'_i)$ have uniformly bounded energy for $i$ large enough. By Lemma~\ref{lem:Iveeunif}, we can thus find $C>0$ such that
$$
\ii^\vee(\mu,\mu')=\inf\left\{\jj_\mu(\f)+\jj_{\mu'}(\f)\mid\f\in\cE^1,\,\jj(\f)\le C\right\}
$$ 
and
$$
\ii^\vee\left(\MA(\p_i),\MA(\p'_i)\right)=\inf\left\{\jj_{\p_i}(\f)+\jj_{\p'_i}(\f)\mid \f\in\cE^1,\,\jj(\f)\le C\right\}
$$
for all $i$ large enough. By Corollary~\ref{cor:maxen} and Lemma~\ref{lem:BBGZ2}, 
$$
\jj_{\p_i}(\f)=\en^\vee(\MA(\p_i))-\en(\f)+\int\f\MA(\p_i)
$$ 
converges to 
$$
\jj_\mu(\f)=\en^\vee(\mu)-\en(\f)+\int\f\,\mu,
$$ 
uniformly with respect to $\f\in\{\jj\le C\}\subset\cE^1$; hence~\eqref{equ:convjmu}. The same holds for $\jj_{\p'_i}(\f)\to\jj_{\mu'}(\f)$, and~\eqref{equ:convIvee} follows. 
\end{proof}

\begin{proof}[Proof of Theorem~\ref{thm:BBGZ3}] Pick maximizing sequences $(\p_i)$, $(\p'_i)$ for $\mu,\mu'$. By~\eqref{equ:convjmu} we have 
$$
\jj(\p_i)\approx\en^\vee(\mu)+o(1),\quad\jj(\p'_i)\approx\en^\vee(\mu')+o(1),
$$
while Lemma~\ref{lem:trianglevee} and~\eqref{equ:convIvee} give
$$
\ii(\p_i,\p'_i)\approx\ii^\vee(\mu,\mu')+o(1).
$$
Thanks to Lemma~\ref{lem:BBGZ}, we infer
$$
\left|\int(\f-\f')\left(\MA(\p_i)-\MA(\p_i')\right)\right|
$$
$$
\lesssim \ii(\f,\f')^{\a_n}\left(\ii^\vee(\mu,\mu')+o(1)\right)^{\frac 12}\max\{\jj(\f),\jj(\f'),\en^\vee(\mu)+o(1),\en^\vee(\mu')+o(1)\}^{\frac 12-\a_n}.
$$
By Lemma~\ref{lem:BBGZ2} we have 
$$
\int(\f-\f')\MA(\p_i)\to\int(\f-\f')\mu,\quad\int(\f-\f')\MA(\p'_i)\to\int(\f-\f')\mu',
$$
and the result follows. 
\end{proof}

\begin{proof}[Proof of Theorem~\ref{thm:trianglevee}]
Let $\mu_1,\mu_2,\mu_3\in\cM^1$, and for $i=1,2,3$ pick a maximizing sequence $(\f_{ij})_j$ for $\mu_i$. By Theorem~\ref{thm:triangle}, we have for each $i$
$$
\ii(\f_{1j},\f_{2j})\lesssim \max\{\ii(\f_{1j},\f_{3j}),\ii(\f_{3j},\f_{2j})\}. 
$$
As in the proof of Theorem~\ref{thm:BBGZ3}, this yields
$$
\ii^\vee(\mu_1,\mu_2)\lesssim\max\{\ii^\vee(\mu_1,\mu_3),\ii^\vee(\mu_3,\mu_2)\},
$$
thanks to Lemma~\ref{lem:trianglevee} and~\eqref{equ:convIvee}. 

Suppose now that $\ii^\vee(\mu,\mu')=0$. By Theorem~\ref{thm:BBGZ3} and \eqref{equ:HdomPL}, it follows that $\int\f\,\mu=\int\f\,\mu'$ for all $\f\in\PL(X)$, and hence $\mu=\mu'$, by density of $\PL(X)$ in $\Cz(X)$ (Theorem~\ref{thm:PLdense}). 

To establish~\eqref{equ:iveeen} and~\eqref{equ:ivjmu}, choose again a maximizing sequence $(\p_i)$ for $\mu$, and set $\mu_i:=\MA(\p_i)$. For each $i$, Lemma~\ref{lem:trianglevee} and~\eqref{equ:enveeJ} yield
$$
\ii^\vee\left(\mu_i,\MA(0)\right)\approx\ii(\p_i)\approx\en^\vee(\mu_i)
$$
and 
$$
\ii^\vee\left(\mu_i,\MA(\f)\right)\approx\ii(\p_i,\f)\approx\jj_{\p_i}(\f). 
$$
By~\eqref{equ:convIvee}, Corollary~\ref{cor:maxen} and Lemma~\ref{lem:BBGZ2}, we have 
$$
\ii^\vee\left(\mu_i,\MA(0)\right)\to\ii^\vee(\mu,\MA(0)),\quad \en^\vee(\mu_i)\to\en^\vee(\mu),\quad \ii^\vee\left(\mu_i,\MA(\f)\right)\to\ii^\vee\left(\mu,\MA(\f)\right),
$$
and 
$$
\jj_{\p_i}(\f)=\en^\vee(\mu_i)-\en(\f)+\int\f\,\mu_i\to\en^\vee(\mu)-\en(\f)+\int\f\,\mu=\jj_\mu(\f),
$$
which proves the result. 
\end{proof}

%
\subsection{Strict convexity of the dual energy}

The dual energy functional $\en^\vee\colon\cM\to[0,+\infty]$ is convex (see~Proposition~\ref{prop:Evee}). As we next show, its restriction to $\cM^1$ is even uniformly convex with respect to the quasi-metric $\ii^\vee$ (compare Theorem~\ref{thm:strictconc} for the energy $\en$). 

\begin{prop}\label{prop:enconvex} For all $\mu,\mu'\in\cM^1$ and $t\in[0,1]$ we have 
$$
\en^\vee((1-t)\mu+t\mu')\le (1-t)\en^\vee(\mu)+t\en^\vee(\mu')-t(1-t)\ii^\vee(\mu,\mu'). 
$$
\end{prop}

\begin{proof} Pick a maximizing sequence $(\f_i)$ in $\cE^1$ for $\mu_t:=(1-t)\mu+t\mu'$. Then 
$$
\jj_{\mu}(\f_i)=\en^\vee(\mu)-\en(\f_i)+\int\f_i\,\mu,\quad\jj_{\mu'}(\f_i)=\en^\vee(\mu')-\en(\f_i)+\int\f_i\,\mu',
$$
and hence 
$$
(1-t)\jj_\mu(\f_i)+t\jj_{\mu'}(\f_i)=(1-t)\en^\vee(\mu)+t\en^\vee(\mu')-\en(\f_i)+\int\f_i\,\mu_t.
$$
By Lemma~\ref{lem:elemconv}, we have
\begin{align*}
(1-t)\jj_\mu(\f_i)+t\jj_{\mu'}(\f_i) & \ge  t(1-t)\left(\jj_\mu(\f_i)+\jj_{\mu'}(\f_i)\right) \\
& \ge  t(1-t)\ii^\vee(\mu,\mu'),
\end{align*}
while $\en(\f_i)-\int\f_i\,\mu_t\to\en^\vee(\mu_t)$; the result follows. 
\end{proof}

%
\subsection{The strong topology of $\cM^1$}

Following~\cite[Definition 2.5]{BBEGZ} we introduce: 

\begin{defi} The \emph{strong topology} on $\cM^1$ is defined as the coarsest refinement of the weak topology for which $\en^\vee\colon\cM^1\to\R$ becomes continuous. 
\end{defi}
Thus a net $(\mu_i)$ in $\cM^1$ converges strongly to $\mu\in\cM^1$ iff $\mu_i\to\mu$ weakly and $\en^\vee(\mu_i)\to\en^\vee(\mu)$. 

When $X$ is irreducible, Corollary~\ref{cor:M1ind} shows that the set $\cM^1\subset\cM$ is independent of the choice of $\om\in\Amp(X)$. As we shall see below, this is then also true of the strong topology of $\cM^1$, cf.~Proposition~\ref{prop:strongind}.

\begin{thm}\label{thm:strongquasi} For a net $(\mu_i)$ and $\mu$ in $\cM^1$, the following are equivalent: 
\begin{itemize}
\item[(i)] $\mu_i\to\mu$ strongly in $\cM^1$; 
\item[(ii)] $\ii^\vee(\mu_i,\mu)\to 0$; 
\item[(iii)]   for each $C>0$ we have $\int\f\,\mu_i\to\int\f\,\mu$ uniformly for $\f\in\cE^1$ with $\jj(\f)\le C$.  . 
\end{itemize}
\end{thm}

By (ii), the strong topology is associated to a canonical (metrizable) uniform structure defined by the quasi-metric $\ii^\vee$.

\begin{cor}\label{cor:maxstrong} Pick $\mu\in\cM^1$, and a net $(\f_i)$ in $\cE^1$. Then $(\f_i)$ is maximizing for $\mu$ iff $\MA(\f_i)\to\mu$ strongly in $\cM^1$. 
\end{cor}
\begin{proof} By definition, $(\f_i)$ is maximizing for $\mu$ iff $\jj_\mu(\f_i)=\en^\vee(\mu)-\en(\f_i)+\int\f_i\,\mu$ tends to $0$. Now~\eqref{equ:ivjmu} yields $\jj_\mu(\f_i)\approx\ii^\vee(\mu,\MA(\f_i))$, and the result is thus a consequence of Theorem~\ref{thm:strongquasi}. 
\end{proof}

\begin{proof}[Proof of Theorem~\ref{thm:strongquasi}] Assume (i), \ie $\mu_j\to\mu$ weakly and $\en^\vee(\mu_j)\to\en^\vee(\mu)$, and pick $\e>0$. By~\eqref{equ:enFS}, we can choose $\f\in\PL_\R\cap\PSH(\om)$ such that 
$$
\jj_\mu(\f)=\en^\vee(\mu)-\en(\f)+\int\f\,\mu\le\e.
$$ 
For all $i$ large enough we have $\int\f\,\mu_i\le\int\f\,\mu+\e$ and $\en^\vee(\mu_i)\le\en^\vee(\mu)+\e$, and hence
$$
\jj_{\mu_i}(\f)=\en^\vee(\mu_i)-\en(\f)+\int\f\,\mu_i\le 3\e. 
$$
Thus $\ii^\vee(\mu_i,\mu)\le\max\{\jj_\mu(\f),\jj_{\mu_i}(\f)\}\le 3\e$ for all $i$ large enough, and we have proved (i)$\Rightarrow$(ii). 

Next assume (ii). Theorem~\ref{thm:BBGZ3} shows that for each $C>0$ we have $\int\f\,\mu_i\to\int\f\,\mu$ uniformly for $\f\in\{\jj\le C\}\subset\cE^1$, hence (ii)$\Rightarrow$(iii). 

Finally assume (iii). First, $\int\f\,\mu_i\to\int\f\,\mu$ for all $\f\in\cE^1$ implies $\int\f\,\mu_i\to\int\f\,\mu$ for all $\f\in\PL(X)$, by \eqref{equ:HdomPL}, and hence $\mu_i\to\mu$ weakly, by density of $\PL(X)$ in $\Cz(X)$.   Next, we claim that $\en^\vee(\mu_i)$ is eventually bounded. By Corollary~\ref{cor:maxen}, for each $i$ we can choose $\f_i\in\cE^1$, normalized by $\int\f_i\MA(0)=0$, such that 
$$
\left|\en^\vee(\MA(\f_i))-\en^\vee(\mu_i)\right|\le 1,\quad\left|\int\f_i(\mu_i-\MA(\f_i))\right|\le 1,
$$
and hence
\begin{equation}\label{equ:enmui}
\en^\vee(\mu_i)\le\en^\vee(\MA(\f_i))+1\lesssim\ii(\f_i)+1=\int(-\f_i)\MA(\f_i)+1\le\int(-\f_i)\,\mu_i+2. 
\end{equation}
Pick $A>0$, to be determined in a moment, and set 
$$
t_i:=\min\{1,A/\en^\vee(\mu_i)\}\in[0,1],\quad s_i:=t_i\en^\vee(\mu_i)=\min\{\en^\vee(\mu_i),A\}.
$$
By concavity of $\en$, $t_i\f_i\in\cE^1$ satisfies 
$$
\jj(t_i\f_i)=t_i\sup\f_i-\en(t_i\f_i)\le t_i\jj(\f_i)\approx t_i\en^\vee(\MA(\f_i))\le s_i+1\le A+1. 
$$
The condition in (iii) thus yields $\int(-t_i\f_i)\mu_i\le\int(-t_i\f_i)\mu+1$ for $i$ large enough (depending on $A$). By Theorem~\ref{thm:BBGZ3}, we have on the other hand 
$$
\int(-t_i\f_i)\mu\le C\max\{1,\jj(t_i\f_i)\}^{1/2}\lesssim C (s_i+1)^{1/2}
$$ 
for a constant $C=C(\mu)>0$ only depending on $\mu$. Combining these estimates with~\eqref{equ:enmui}, we get 
$s_i=t_i\en^\vee(\mu_i)\le C' (s_i^{1/2}+1)$ for a constant $C'=C'(\mu)>0$, and hence $s_i\le C''=C''(\mu)$. Choosing $A>C''$, this yields, as desired, $\en^\vee(\mu_i)\le C''$ for $i$ large enough. By Corollary~\ref{cor:enmaxim}, we can now find a uniform constant $B>0$ such that 
$$
\en^\vee(\mu_i)=\sup_{\jj(\f)\le B}\left(\en(\f)-\int\f\,\mu_i\right)
$$
for all $i$ large enough, and the condition in (iii) now implies $\en^\vee(\mu_i)\to\en^\vee(\mu)$, which proves (iii)$\Rightarrow$(i).  
\end{proof}

We finally show that the uniform structure of $\cM^1$ is complete. 
\begin{thm}\label{thm:M1complete} Let $(\mu_i)$ be a Cauchy net in $\cM^1$, \ie $\lim_{i,j}\ii^\vee(\mu_i,\mu_j)=0$. Then $(\mu_i)$ converges in the strong topology.
\end{thm}
 
\begin{rmk}\label{rmk:d1barcomplete}
  In~\cite{nakstab1} we exhibit a natural complete metric $\dd_1$ on $\cM^1$ that defines the strong topology. (The case $\om=c_1(L)$ with $L\in\Pic(X)_\Q$ is handled in~\cite{nakstab1} and the general case in~\cite{nakstab2}.)
\end{rmk}
 
\begin{proof} By the quasi-triangle inequality, it is enough to show that some subnet of $(\mu_i)$ converges strongly in $\cM^1$. After passing to a subnet, we can thus assume wlog that $(\mu_i)$ converges weakly to a Radon probability measure $\mu\in\cM$, by weak compactness of $\cM$. By Theorem~\ref{thm:trianglevee}, $\en^\vee(\mu_i)$ is eventually bounded, and hence $\en^\vee(\mu)<+\infty$, \ie $\mu\in\cM^1$, by lower semicontinuity of $\en^\vee$ in the weak topology of measures. We thus have $\mu\in\cM^1$, and it remains to prove that for each $C>0$ we have $\int\f\,\mu_i\to\int\f\,\mu$ uniformly for $\f\in\{\jj\le C\}\subset\cE^1$, by Theorem~\ref{thm:strongquasi}. To see this, pick $\e>0$. For all $i,j$ large enough, Theorem~\ref{thm:BBGZ3} yields $\left|\int\f\,\mu_i-\int\f\,\mu_j\right|\le\e$ for all $\f\in\{\jj\le C\}\cap\PL_\R$. Since $\mu_j\to\mu$ weakly, it follows that for all $i$ large enough we have $\left|\int\f\,\mu_i-\int\f\,\mu\right|\le\e$ for all $\f\in\{\jj\le C\}\cap\PL_\R$, and hence also for $\f\in\{\jj\le C\}\subset\cE^1$, by monotone convergence. 
\end{proof}

%
%
%
%
%

%
%
\subsection{Dependence on the ample class} 
In this section, $X$ is assumed to be \textbf{irreducible}. Recall from Corollary~\ref{cor:M1ind} that the set $\cM^1=\cM^1(\om)$ of measures of finite energy is independent of $\om\in\Amp(X)$. 

\begin{prop}\label{prop:strongind}
The strong topology on $\cM^1=\cM^1(\om)$ does not depend on the choice of $\om\in\Amp(X)$.
\end{prop}

\begin{proof} Pick $\mu\in\cM^1$, and a net $(\mu_i)$  in $\cM^1$. By Theorem~\ref{thm:M1ind}, the condition that $\en_\om^\vee(\mu_i)$ is eventually bounded is independent of $\om\in\Amp(X)$. On the other hand, pick $\om,\om'\in\Amp(X)$ and choose $s\ge 1$ such that $s^{-1}\om\le\om'\le s\om$. By Lemma~\ref{lem:enandtheta}, for each $C>0$ there exists $C'>0$ such that $\{\jj_\om\le C\}\subset s\{\jj_{\om'}\le C'\}$. As a result, the condition that $\int\f\,\mu_i\to\int\f\,\mu$ uniformly for $\f\in\{\jj_\om\le C\}$ for all $C>0$ is also independent of $\om$, and Theorem~\ref{thm:strongquasi} yields, as desired, that the strong convergence of $(\mu_i)$ to $\mu$ is independent of the choice of $\om$. 
\end{proof}

\begin{rmk}
  In view of Theorem~\ref{thm:M1ind}, one may ask whether for any two $\om,\om'\in\Amp(X)$ we have $C^{-1}\ii^\vee_\om\le\ii^\vee_{\om'}\le C\ii^\vee_\om$ for some $C=C(\om,\om')\ge 1$, in which case the uniform structure on $\cM^1=\cM^1(\om)$ would also independent of $\om$. 
Theorem~\ref{thm:M1ind} at least shows that the class of \emph{bounded subsets} of $\cM^1$, \ie subsets of $\cM^1_C(\om)=\{\mu\in\cM^1\mid\en_\om^\vee(\mu)\le C\}$ for some $C>0$, is independent of the choice of $\om$. 
\end{rmk}

\begin{lem}\label{lem:Iveecont}
For all $\mu,\mu'\in\cM^1$, the map $\om\mapsto\ii^\vee_\om(\mu,\mu')$ is continuous on $\Amp(X)$, uniformly for $\mu,\mu'$ in any bounded subset of $\cM^1$. 
\end{lem}
\begin{proof} By Lemma~\ref{lem:enandtheta}, there exists a function $s:[1,1+\e]\to [1,2]$ with $\lim_{t\to 1} s(t)=1$, only depending on $n$, such that for all $\om,\om'\in\Amp(X)$ with 
\begin{equation}\label{equ:boundtheta}
t^{-1}\om\le\om'\le t\om,\,\,\,t\in[1,1+\e]
\end{equation}
we have
$$
s(t)^{-1}\en^\vee_{\om}\le\en^\vee_{\om'}\le s(t)\en^\vee_\om,\quad\en_{\om'}(t^{-1}\f)\ge s(t)\en_\om(\f)
$$ 
for all $\f\in\cE^1_{\sup}(\om)\subset t\cE^1_{\sup}(\om')$. 

Pick $\om,\om'\in\Amp(X)$ satisfying~\eqref{equ:boundtheta}. Let $C>0$ and $\mu,\mu'\in\cM^1_C(\om)$. By Lemma~\ref{lem:Iveeunif}, we can find $C'\lesssim C$ such that 
\begin{equation}\label{equ:iveetheta}
\ii^\vee_{\om}(\mu,\mu')=\inf\left\{\en^\vee_\om(\mu)+\en^\vee_\om(\mu')
    -2\en_\om(\f)
    +\int\f(\mu+\mu')\ \big|\ \f\in\cE^1_{C'}(\om)\right\}
 \end{equation}
For each $\f\in\cE^1_{C'}(\om)$, we have 
$$
\int\f(\mu+\mu')\ge2\en_{\om}(\f)-\en^\vee_\om(\mu)-\en^\vee_\om(\mu')\ge-2(C+C').
$$ 
Thus
$$
    \ii^\vee_{\om'}(\mu,\mu')
    \le\en^\vee_{\om'}(\mu)+\en^\vee_{\om'}(\mu')
    -2\en_{\om'}(t^{-1}\f)
    +\int t^{-1} \f(\mu+\mu')
    $$
    $$
    \le s(t)\left(\en^\vee_\om(\mu)+\en^\vee_\om(\mu')
    -2\en_\om(\f)
    +\int\f(\mu+\mu')\right)
    +2(s(t)-t^{-1})(C+C'),
    $$
    and~\eqref{equ:iveetheta} yields
$$
 \ii^\vee_{\om'}(\mu,\mu')\le s(t)\ii^\vee_\om(\mu,\mu')+2(s(t)-t^{-1})(C+C')\le\ii^\vee_\om(\mu,\mu')+\e(t)C
 $$
 with $\lim_{t\to 1}\e(t)=0$, since $\ii^\vee(\mu,\mu')\lesssim C$ by Theorem~\ref{thm:trianglevee}. Since $\mu,\mu'\in\cM^1_C(\om)\subset\cM^1_{2C}(\om')$, the desired result follows by symmetry.
 \end{proof}

 %
%
 \section{Valuations of linear growth}
 Assume that $X$ is of dimension $n$ and \textbf{irreducible}, and fix $\om\in\Amp(X)$ ample. The purpose of this section is to show that the set $X^\lin$ of valuations of linear growth is endowed with a natural metric induced by $\om$-psh functions, with respect to which it sits as a bi-Lipschitz subspace of $(\cM^1,\ii^\vee)$.
%
%
\subsection{The energy of a Dirac mass}
Recall from \S\ref{sec:pluripolar} that a point $v\in X^\an$ is nonpluripolar iff 
\begin{equation}\label{equ:tee}
\te(v)=\sup_{\f\in\PSH}\left(\sup\f-\f(v)\right)=\sup_{\f\in\PSH}\left(\f(v_\triv)-\f(v)\right)
\end{equation}
is finite, and that the set of such points coincides with the set $X^\lin$ of valuations of linear growth. 

\begin{prop}\label{prop:enT} For any $v\in X^\an$, the Dirac mass $\d_v\in\cM$ has finite energy iff $v\in X^\lin$. Moreover, 
\begin{equation}\label{equ:TvsE}
\frac 1{n+1}\te(v)\le\en^\vee(\d_v)\le\frac{n}{n+1}\te(v).
\end{equation}
\end{prop}

\begin{proof} Let us first show that
\begin{equation}\label{equ:TvsE2}
\frac{1}{n+1}\te(v)\le\en^\vee(\d_v)\le\te(v),
\end{equation}
which will already imply $\d_v\in\cM^1\Longleftrightarrow v\in X^\lin$. The right-hand inequality is trivial, since 
$$
\en^\vee(\d_v)
=\sup_{\f\in\cE^1}\left(\en(\f)-\f(v)\right)
\le\sup_{\f\in\cE^1}\left(\sup\f-\f(v)\right)
=\sup_{\f\in\PSH(\om)}\left(\sup\f-\f(v)\right)=\te(v). 
$$
Here the second equality follows since every function in $\PSH(\om)$ is a decreasing limit of functions in $\cH^\dom(\om)\subset\cE^1$.

The left-hand inequality in~\eqref{equ:TvsE2} is equivalent to $\f(v)\ge -(n+1)\en^\vee(\d_v)$ for all $\f\in\PSH_{\sup}$. For each $m\in\N$ set $\p_m:=\max\{\f,\f(v),-m\}$. Then $\p_m\in\PSH_{\sup}$ and 
$$
\p_m\ge\p_m(v)=\max\{\f(v),-m\}.
$$
Since $X$ is irreducible, $\int\p_m\MA(0)=\p_m(v_\triv)=\sup\p_m=0$, and~\eqref{equ:Ebis} thus yields
$$
-\frac{1}{n+1}\max\{\f(v),-m\}=\frac{1}{n+1}\int(\p_m-\p_m(v))\MA(0)
$$
$$
\le\en(\p_m-\p_m(v))=\en(\p_m)-\p_m(v)\le\en^\vee(\d_v). 
$$
As $m\to\infty$, this yields, as desired, $-\f(v)\le (n+1)\en^\vee(\d_v)$. 

It remains to establish the stronger right-hand inequality in~\eqref{equ:TvsE} for $v\in X^\lin$. Pick a maximizing sequence $(\f_i)$ for $\d_v$ in $\cE^1$. By~\eqref{equ:enma} and~\eqref{equ:IJ} we have
$$
\en^\vee(\MA(\f_i))
=\ii(\f_i)-\jj(\f_i)
\le\frac{n}{n+1}\ii(\f_i),
$$
with 
$$
\ii(\f_i)=\sup\f_i-\int\f_i\MA(\f_i)\le\te(v)+\int\f_i\left(\d_v-\MA(\f_i)\right).
$$
By Corollary~\ref{cor:maxen}, $\en^\vee(\MA(\f_i))\to\en^\vee(\d_v)$, and $\int\f_i(\d_v-\MA(\f_i))\to 0$. The result follows. 
\end{proof}

\begin{rmk} When $X$ has at least two top-dimensional components, point masses never have finite energy, by Corollary~\ref{cor:massMA}.
\end{rmk}

As we next show, the image of the embedding $X^\lin\hto\cM^1$ can further be characterized as the set of extremal points.  

\begin{prop}\label{prop:extpoint} A measure $\mu\in\cM^1$ is an extremal point of the convex set $\cM^1$ iff $\mu=\d_v$ with $v\in X^\lin$. 
\end{prop}
\begin{proof} As is well-known, the extremal points of the space $\cM$ of all Radon probability measures are the Dirac masses $\d_v$, $v\in X^\an$; if $v\in X^\lin$, then $\d_v$ is \emph{a fortiori} an extremal point of $\cM^1$. Conversely, assume $\mu\in\cM^1$ is not a Dirac mass. Then $\mu$ is not an extremal point of $\cM$, and hence can be written as $\mu=(1-t)\mu_0+t\mu_1$ with $\mu_0\ne\mu_1\in\cM$ and $t\in (0,1)$. For any $\f\in\cE^1_{\sup}$ we then have $(1-t)\int\f\,\mu_0\ge\int\f\,\mu>-\infty$, $t\int\f\,\mu_1\ge\int\f\,\mu>-\infty$. By Theorem~\ref{thm:charM1}, this implies $\mu_0,\mu_1\in\cM^1$, which contradicts the extremality of $\mu$ in $\cM^1$. 
\end{proof}

%
%
\subsection{Weak convergence of psh functions}
Every $\f\in\PSH(\om)$ is finite-valued on $X^\lin$, see Corollary~\ref{cor:pshonXlin}, but the topology of $\PSH(\om)$ is defined as the topology of pointwise convergence on the strict subset $X^\div\subset X^\lin$. However we show: 

\begin{thm}\label{thm:weaklin} The topology of $\PSH(\om)$ coincides with the topology of pointwise convergence on $X^\lin$. 
\end{thm}
See also Corollary~\ref{cor:charweaktop} for other characterizations of the topology of $\PSH(\om)$, assuming the envelope property.

\begin{proof} Given a convergent net $\f_i\to\f$ in $\PSH(\om)$ and $v\in X^\lin$, we need to show $\f_i(v)\to\f(v)$. Since $\sup\f_i=\f_i(v_\triv)$ converges to $\sup\f=\f(v_\triv)$ and $\sup\f_i\le\f_i(v)+\te(v)$, the net $(\f_i(v))$ is eventually bounded. Replacing $\f_i$ and $\f$ with $\max\{\f_i,-t\}$ and $\max\{\f,-t\}$ for $t\gg 1$, we may thus assume wlog that $(\f_i)$ is uniformly bounded, and hence that $\jj(\f_i)$ is bounded. The result is now a consequence of Proposition~\ref{prop:weakcont}, since $\d_v$ has finite energy by Proposition~\ref{prop:enT}. 
\end{proof}

\begin{rmk}\label{rmk:weaklin} Theorem~\ref{thm:weaklin} remains true for a general projective variety $X$, simply by restricting to each of its irreducible components. 
\end{rmk}

%
%
\subsection{The $\dd_\infty$-metric}\label{sec:Izumi}
For any two $v,w\in X^\lin$, we set 
\begin{equation}\label{equ:Izumi}
  \dd_\infty(v,w):=\sup_{\f\in\PSH}\left|\f(v)-\f(w)\right|.
\end{equation}

\begin{prop}\label{prop:Izumi} The following properties hold. 

\begin{itemize}
\item[(i)] $\dd_\infty$ is a metric on $X^\lin$, and is the smallest one with respect to which $\f\colon X^\lin\to\R$ is $1$-Lipschitz for each $\f\in\PSH(\om)$; 
\item[(ii)] $\te(v)=\dd_\infty(v,v_\triv)$ for all $v\in X^\lin$; 
\item[(iii)] $\dd_\infty(tv,tw)=t \dd_\infty(v,w)$ for all $v,w\in X^\lin$ and $t\in\R_{>0}$; 
\item[(iv)] $\dd_\infty(v,w)=\sup_{\f\in\PL\cap\PSH}|\f(v)-\f(w)|$; 
\item[(v)] $\dd_\infty$ is lsc on $X^\lin\times X^\lin$. 
\end{itemize}
\end{prop}
\begin{proof} The function $\dd_\infty$ is finite-valued, since $\dd_\infty(v,w)\le \te(v)+\te(w)$ by~\eqref{equ:tee}. It is further obviously symmetric, and satisfies the triangle inequality. If $\dd_\infty(v,w)=0$, then $\f(v)=\f(v)$ for all $\f\in\PL_\R\cap\PSH(\om)$. By \eqref{equ:HdomPL}, this implies $\f(v)=\f(w)$ for all $\f\in\PL_\R$, and hence $v=w$, since $\PL(X)$ separates points by Lemma~\ref{lem:FSdense}. The second half of (i) is tautological. (ii) and (iii) follow directly from the definition. An easy approximation argument yields (iv), which in turn implies (v).  \end{proof}

We refer to the metric space topology defined by $\dd_\infty$ as the \emph{strong topology} of $X^\lin$, while the \emph{weak topology} of $X^\lin$ means the topology inherited from $X^\an$. By (iv), the weak topology is coarser than the strong one. 

If $\om\le\om'$, then $\dd_{\infty,\om}\le \dd_{\infty,\om'}$. The metrics $\dd_{\infty,\om}$ with $\om\in\Amp(X)$ are thus all Lipschitz equivalent, and the strong topology of $X^\lin$ is independent of the choice of $\om$. 

\medskip

In the case of classes of $\Q$-line bundles, the $\dd_\infty$ metric admits the following alternative description.
\begin{lem}\label{lem:Izumilim}
  If $L$ is an ample $\Q$-line bundle, then $\dd_\infty=\dd_{\infty,c_1(L)}$ satisfies 
\begin{align*}
  \dd_{\infty}(v,w)
  &=\sup\left\{m^{-1}|v(s)-v(w)|\mid m\in\Z_{>0},\,s\in\Hnot(X,mL)\setminus\{0\}\right\}\\
  &=\lim_{m\to\infty}\sup\left\{m^{-1}\left|v(s)-w(s)\right|\mid s\in\Hnot(X,mL)\setminus\{0\}\right\}
\end{align*}
for any $v,w\in X^\lin$.
\end{lem}
\begin{proof} For each $s\in\Hnot(X,mL)$, $\f:=m^{-1}\log|s|$ is $\om$-psh, and $|v(s)-w(s)|=|\f(v)-\f(w)|$, which proves 
$$
\dd_\infty(v,w)\ge S:=\sup\left\{m^{-1}|v(s)-v(w)|\mid m\in\Z_{>0},\,s\in\Hnot(X,mL)\setminus\{0\}\right\}.
$$
Next pick $\f\in\cH(L)$, \ie $\f=m^{-1}\max_i\{\log|s_i|+\la_i\}$ for a basepoint free, finite set of sections $(s_i)$ of $\Hnot(X,mL)$ and $\la_i\in\Q$. For each $i$ we have $m^{-1}\log|s_i|(v)\le m^{-1}\log|s_i|(w)+S$, thus $\f(v)\le \f(w)+S$. Assume $\f\in\PSH(L)$. By Theorem~\ref{thm:pshample}, $\f$ is the pointwise limit of a decreasing net $(\f_i)$ in $\cH(L)$, and hence $\f(v)\le\f(w)+S$, which proves $\dd_\infty(v,w)\le S$, by symmetry. 

Finally, the second equality follows from 
$$
|v(s)-w(s)|=\max\{v(s)-w(s),w(s)-v(s)\}
$$ 
and the superadditivity of 
$$
m\mapsto\sup_{s\in\Hnot(X,mL)\setminus\{0\}}(v(s)-w(s)),
$$
thanks to Fekete's lemma. 
\end{proof}

\begin{exam}\label{exam:Izumicurve} If $X$ is a smooth curve, then the parametrizations described in \S\ref{sec:curve} endow $X^\lin=X^\val=X^\an\setminus X(k)=\bigcup_{p\in X}\iota_p([0,+\infty))$ with a metric which equals $\dd_{\infty,\om}$ up to a factor $\deg\om$. In this case, the metric space $(\Xlin,\dd_{\infty,\om})$ is an $\R$-tree.
\end{exam}

\begin{thm}\label{thm:complete} The metric space $(X^\lin,\dd_\infty)$ is complete. 
\end{thm}
\begin{proof} Let $(v_i)$ be a Cauchy net for $(X^\lin,\dd_\infty)$. Upon passing to a subnet, we may assume that $v_i$ admits a limit $v\in X^\an$ in the topology of $X^\an$. Pick $\e>0$, and choose $i_0$ such that $\dd_\infty(v_i,v_j)\le\e$ for all $i,j\ge i_0$. We claim that $v\in X^\lin$ and $\dd_\infty(v_i,v)\le\e$ for all $i\ge i_0$, which will prove, as desired, that the Cauchy net $(v_i)$ admits a limit in $(X^\lin,\dd_\infty)$. Indeed, for all $\f\in\PSH(\om)$ and $i,j\ge i_0$, we have $|\f(v_i)-\f(v_j)|\le\e$. Letting $j\to\infty$, this first shows that $\f(v)<+\infty$ for all $\f\in\PSH(\om)$, and hence $v\in X^\lin$. Furthermore,  $|\f(v_i)-\f(v)|\le\e$, and taking the supremum over $\f$ yields the claim. 
\end{proof}

As noted above, the Lipschitz equivalence class of $(X^\lin,\dd_{\infty})$ is independent of $\om\in\Amp(X)$. It is also a birational invariant of $X$: 

\begin{prop}\label{prop:dinftybir} Every birational morphism $\pi\colon Y\to X$ induces a bi-Lipschitz isomorphism $(Y^\lin,\dd_{\infty})\simeq(X^\lin,\dd_\infty)$. 
\end{prop}
This refines~\cite[Lemma 2.8]{BKMS}.

\begin{proof} To see that $Y^\lin\to X^\lin$ is Lipschitz, pick $\om_X\in\Amp(X)$ and $\om_Y\in\Amp(Y)$ such that $\pi^\star\om_X\le\om_Y$. Then $\pi^\star\PSH(\om_X)\subset\PSH(\om_Y)$, which implies $\dd_{\infty,\om_X}(v,w)\le  \dd_{\infty,\om_Y}(v,w)$ for all $v,w\in Y^\lin$.

  In order to prove that the inverse $X^\lin\to Y^\lin$ is Lipschitz, we use Corollary~\ref{cor:descentpsh} to find $\om_X\in\Amp(Y)$, $\om_Y\in\Amp(Y)$ and $\f_X\in\PSH(\om_X)$ such that 
$$
\PSH(\om_Y)+\pi^\star\f_X\subset\pi^\star\PSH(\om_X). 
$$
For all $v,w\in Y^\lin$, this implies $\dd_{\infty,\om_Y}(v,w)\le 2 \dd_{\infty,\om_X}(v,w)$, and the result follows. 
\end{proof}

%

\subsection{Growth and weak continuity of psh functions}

By definition of $\te$, any $\om$-psh function $\f$ satisfies a linear growth estimate 
$$
|\f(v)|\le \te(v)+O(1)=\dd_\infty(v,v_\triv)+O(1)
$$ 
on $X^\lin$, and this cannot be improved in general. For functions in $\cE^1$, the growth is sublinear: 

\begin{thm}\label{thm:phiTbelow}
Each $\f\in\cE^1$ satisfies 
$$
|\f|\le A\te^{1-\a_n}+B
$$
on $X^\lin$, with $\a_n:=2^{-n}$ and $A,B>0$ only depending on $\jj(\f)$ and $\sup\f$, respectively.  
\end{thm}
\begin{proof} For each $v\in X^\lin$, Proposition~\ref{prop:enT} shows that $\en^\vee(\d_v)\approx \te(v)$, and Theorem~\ref{thm:BBGZ3} thus gives
$$
\sup\f-\f(v)\lesssim\max\{1,\jj(\f)\}^{\frac 12}\max\left\{1,\te(v)\right\}^{1-\a_n}. 
$$
The result follows. 
\end{proof}

As already noticed, the restriction of any $\f\in\PSH(\om)$ to $X^\lin$ is strongly continuous, and even $1$-Lipschitz with respect to $\dd_\infty$. With respect to the weak topology we have: 

\begin{thm}\label{thm:pshweakcont} The restriction of any $\f\in\PSH(\om)$ to a bounded subset of $(X^\lin,\dd_\infty)$ is weakly continuous.
\end{thm}

\begin{proof} We may assume $\f\le 0$. Then $\exp(\f)\in\cE^\infty\subset\cE^1$, see Corollary~\ref{cor:exppsh}, and it is thus enough to prove the result for $\f\in\cE^1$. Pick a decreasing net $(\f_i)$ in $\PL_\R\cap\PSH(\om)$ converging pointwise to $\f$. Then 
$$
0\le \ii(\f_i,\f)=\int(\f_i-\f)(\MA(\f)-\MA(\f_i))\le\int(\f_i-\f)\MA(\f)\to 0, 
$$
by monotone convergence. On the other hand, $\dd_\infty(v,v_\triv)=\te(v)\le C$ implies $\en^\vee(\d_v)\lesssim C$, by Proposition~\ref{prop:enT}, and Theorem~\ref{thm:BBGZ3} applied to $\mu=\d_v$, $\mu'=\d_{v_\triv}$ implies that $\f_i(v)\to\f(v)$ uniformly for $v$ in the ball $\{T\le C\}$. Since each $\f_i$ is weakly continuous on $\PSH(\om)$, we conclude, as desired, that $\f$ is weakly continuous on $\{T\le C\}$. 
\end{proof}

%

\subsection{Bi-Lipschitz embedding into $\cM^1$}

By Proposition~\ref{prop:enT}, $v\mapsto\d_v$ defines an injection $X^\lin\hookrightarrow\cM^1$ such that 
$\te(v)\approx\en^\vee(\d_v)$, \ie $\dd_\infty(v,v_\triv)\approx\ii^\vee(\d_v,\d_{v_\triv})$, by Theorem~\ref{thm:trianglevee}. More generally, we prove: 

\begin{thm}\label{thm:Izumi} For all $v,w\in X^\lin$ we have $\dd_\infty(v,w)\approx\ii^\vee(\d_v,\d_w)$. 
In particular, $X^\lin\hookrightarrow\cM^1$ is a topological embedding with respect to the strong topologies. 
\end{thm}
 
\begin{rmk}\label{rmk:Xlinmetrics}
  In~\cite{nakstab1,nakstab2} we explore additional natural metrics on $X^\lin$; these are all equivalent to $\dd_\infty$.
\end{rmk}
 
By Theorem~\ref{thm:strongquasi}, we infer: 

\begin{cor}\label{cor:Izumi} The strong topology of $X^\lin$ is the coarsest refinement of the weak topology in which $v\mapsto\en^\vee(\d_v)$ becomes continuous. 
\end{cor} 
In other words, a net $(v_i)$ of $X^\lin$ converges strongly to $v\in X^\lin$ iff $v_i\to v$ weakly and $\en^\vee(\d_{v_i})\to\en^\vee(\d_v)$.

\begin{qst} Do we have a similar characterization with $\te(v)$ in place of $\en^\vee(\d_v)$?
\end{qst}

\begin{lem}\label{lem:Izumi} Assume $\om=c_1(L)$ with $L$ an ample line bundle, and pick $a\gg 1$ such that $mL$ is globally generated for all $m\ge a$. Pick $v\in X^\lin$, and set for each $m\ge a$
$$
\f_{v,m}:=m^{-1}\max_{s\in\Hnot(X,mL)\setminus\{0\}}\left(\log|s|+v(s)\right)\in\cH_\R(L). 
$$
Then $(\f_{v,m})_{m\ge a}$ is a maximizing sequence for $\d_v$,   and hence $\MA(\f_{v,m})\to\d_v$ strongly in $\cM^1$.  
\end{lem}
\begin{proof} We first claim that $\lim_{m\to\infty}\en(\f_{v,m})=\sup_{m\ge a}\en(\f_{v,m})$. To see this, observe that
$$
(m+m')\f_{v,m+m'}\ge m\f_{v,m}+m'\f_{v,m'}
$$ 
for all $m,m'\ge a$. Set $S:=\sup_m\en(\f_{v,m})$, pick $\e>0$, and choose $m_0\gg 1$ such that $\en(\f_{v,m_0})\ge S-\e$. For each $m\ge m_0+a$ write $m-a=qm_0+r$ with $q\in\N$ and $r\in\{0,\dots,m_0-1\}$. The above superadditivity property yields
$$
m\f_{v,m}\ge qm_0\f_{v,m_0}+(r+a)\f_{v,r+a}. 
$$
Thus $\f_{v,m}\ge\f_{v,m_0}-O(m^{-1})$, and hence 
$$
\en(\f_{v,m})\ge S-\e-O(m^{-1})\ge S-2\e
$$
for all $m\gg 1$, proving the claim. 

Since $\log|s|(v)=-v(s)$ for $s\in\Hnot(X,mL)\setminus\{0\}$, we have $\f_{v,m}(v)=0$. Thus
$$
\lim_{m\to\infty}\en(\f_{v,m})=\sup_{m\ge a}\en(\f_{v,m})\le\en^\vee(\d_v),
$$
and we need to show that this is an equality. Pick $\f\in\cH(L)$, so that 
$$
\f=m^{-1}\max_i\{\log|s_i|+\la_i\}
$$
for some $m\in\Z_{>0}$, a basepoint free, finite set $(s_i)$ of $\Hnot(X,mL)$ and $\la_i\in\Q$. We may assume wlog $m\ge a$, since we have for $r\in\Z_{>0}$
$$
\f=(rm)^{-1}\max_i\{\log|s_i^r|+r\la_i\}. 
$$
For each $i$ we have 
$$
-v(s_i)+\la_i=\log|s_i|(v)+\la_i\le m\f(v). 
$$
This yields $\f\le\f_{v,m}+\f(v)$, and hence $\en(\f)-\f(v)\le\en(\f_{v,m})$. We infer
$$
\en^\vee(\d_v)=\sup_{\f\in\cH(L)}\left(\en(\f)-\f(v)\right)\le\sup_{m\ge a}\en(\f_m),
$$
which proves, as desired, that $(\f_{v,m})$ is a maximizing sequence for $\d_v$.   The final assertion now follows from Corollary~\ref{cor:maxstrong}.  
\end{proof}

\begin{proof}[Proof of Theorem~\ref{thm:Izumi}]
By Lemma~\ref{lem:Iveecont}, we may assume that $\om=c_1(L)$ with $L\in\Pic(X)_\Q$ ample. By homogeneity, we can even assume that $L$ is a globally generated line bundle. For each $v\in X^\lin$ consider the sequence $(\f_{v,m})$ from Lemma~\ref{lem:Izumi}. For $w\in X^\lin$ we have
$$
\f_{v,m}(w)=m^{-1}\max_{s\in\Hnot(X,mL)\setminus\{0\}}\left(v(s)-w(s)\right),
$$
and Lemma~\ref{lem:Izumilim} thus shows that 
\begin{equation}\label{equ:dlinphiv}
\dd_\infty(v,w)=\lim_{m\to\infty}\max\{\f_{v,m}(w),\f_{w,m}(v)\}.
\end{equation}
On the other hand, Lemma~\ref{lem:trianglevee} and~\eqref{equ:convIvee} yield
\begin{equation}\label{equ:iveephiv}
\ii(\f_{v,m},\f_{w,m})\approx\ii^\vee\left(\MA(\f_{v,m}),\MA(\f_{w,m})\right)\to\ii^\vee(\d_v,\d_w).
\end{equation}
Write
$$
\ii(\f_{v,m},\f_{w,m})=\int(\f_{v,m}-\f_{w,m})\left(\MA(\f_{w,m})-\MA(\f_{v,m})\right).
$$
By Lemma~\ref{lem:Izumi}, $(\f_{v,m})$ and $(\f_{w,m})$ are maximizing sequences for $\d_v$ and $\d_w$. By Corollary~\ref{cor:maxen}, 
$$
\int\f_{v,m}\MA(\f_{v,m})=\int\f_{v,m}(\MA(\f_{v,m})-\d_v)\to 0
$$
and
$$
\int\f_{w,m}\MA(\f_{v,m})-\f_{w,m}(v)\to 0.
$$
Exchanging the roles of $v$ and $w$, we infer
$$
\ii(\f_{v,m},\f_{w,m})=\f_{v,m}(w)+\f_{w,m}(v)+o(1)\approx\max\{\f_{v,m}(w),\f_{w,m}(v)\}+o(1). 
$$
Combining this with~\eqref{equ:dlinphiv} and~\eqref{equ:iveephiv}, we conclude, as desired,
$\dd_\infty(v,w)\approx\ii^\vee(\d_v,\d_w)$. 
\end{proof}

As an application of the above results, we prove: 
\begin{thm}\label{thm:char0dense} Assume $\charac k=0$. Then $X^\div$ is dense in $X^\lin$ with respect to the strong topology.
\end{thm}

\begin{lem}\label{lem:char0dense} Let $v_i\to v$ be a weakly convergent net in $X^\lin$. If $v_i\le v$, then $v_i\to v$ strongly. 
\end{lem}
\begin{proof} Every $\om$-psh function is usc and decreasing on $X^\an$, and 
$$
w\mapsto\en^\vee(\d_w)=\sup_{\f\in\cE^1}(\en(\f)-\f(w))
$$
is thus increasing and lsc on $X^\an$. This yields $\en^\vee(\d_{v_i})\le\en^\vee(\d_v)$ and 
$$
\en^\vee(\d_v)\le\liminf_i\en^\vee(\d_{v_i}),
$$
and hence $\en^\vee(\d_{v_i})\to\en^\vee(\d_v)$. By Corollary~\ref{cor:Izumi}, this proves $v_i\to v$ strongly.
\end{proof}

\begin{proof}[Proof of Theorem~\ref{thm:char0dense}] By Corollary~\ref{cor:apprbelow}, any $v\in X^\an$ is the limit of a net $(v_i)_i$ in $X^\div$ such that $v_i\le v$. When $v\in X^\lin$ this implies $v_i\to v$ strongly, by Lemma~\ref{lem:char0dense}, and thus proves that $X^\div$ is strongly dense in $X^\lin$. 
\end{proof}

%
%
%
%
 \section{The strong topology on $\cE^1$ and the Calabi--Yau theorem}
In this section, we assume until further notice that $X$ is \textbf{irreducible}, and fix an ample class $\om\in\Amp(X)$. 

Having analyzed the strong topology on the space $\cM^1=\cM^1(\om)$ of Radon probability measures of finite energy, we now perform the corresponding analysis on the space $\cE^1=\cE^1(\om)$ of $\om$-psh functions of finite energy.

%
%
\subsection{The strong topology of $\cE^1$} 
Following~\cite{BBGZ,BBEGZ} we introduce: 

\begin{defi} The \emph{weak topology} of $\cE^1$ is the topology inherited from $\PSH(\om)$. The \emph{strong topology} is the coarsest refinement of the weak topology 
for which $\en\colon\cE^1\to\R$ becomes continuous.
\end{defi}

Thus a net $(\f_j)$ in $\cE^1$ converges weakly to $\f\in\cE^1$ iff 
$\f_j\to\f$ pointwise on $X^\div$ (or, equivalently, on $X^\lin$, cf.~Remark~\ref{rmk:weaklin}); it converges strongly iff we further have $\en(\f_j)\to\en(\f)$. By Proposition~\ref{prop:E}, $\en$ is continuous along decreasing nets, and hence: 

\begin{exam}\label{exam:decrstrong} For a decreasing net $(\f_i)$ in $\cE^1$, weak and strong convergence coincide. 
\end{exam}

\begin{defi} We define the quasi-metric $\tii$ on $\cE^1$ by setting 
\begin{equation}\label{equ:quasiE1}
\tii(\f,\f')=\tii_\om(\f,\f'):=\ii(\f,\f')+|\sup\f-\sup\f'|
\end{equation}
for $\f,\f'\in\cE^1$.
\end{defi}
By~\eqref{equ:Itriangle} and Corollary~\ref{cor:MAinj}, $\tii$ is indeed a quasi-metric. It is further immediate to check that it satisfies the analogue of~\eqref{equ:iipyth}, \ie 
\begin{equation}\label{equ:tiipyth}
\tii(\f,\p)=\tii\left(\f,\max\{\f,\p\}\right)+\tii\left(\max\{\f,\p\},\p\right)
\end{equation}
for all $\f,\p\in\cE^1$. As we now show, the quasi-metric $\tii$ defines the strong topology of $\cE^1$. 

\begin{thm}\label{thm:strongunif} A net $(\f_j)$ in $\cE^1$ converges strongly to $\f\in\cE^1$ iff $\tii(\f_j,\f)\to0$. 
\end{thm}

\begin{proof} First assume that $\f_j\to\f$ strongly in $\cE^1$. Then
$\sup\f_j\to\sup\f$, $\en(\f_j)\to \en(\f)$, and $\jj(\f_j)=\sup\f_j-\en(\f_j)$ is thus eventually bounded. Since $\f_j\to\f$ weakly,
Proposition~\ref{prop:weakcont} yields $\int(\f_j-\f)\MA(\f)\to 0$.
Thus 
 \begin{equation}\label{e509}\
 \ii(\f_j,\f)\approx \jj_\f(\f_j)=\en(\f)-\en(\f_j)+\int(\f_j-\f)\MA(\f)\to 0, 
\end{equation}
in view of~\eqref{equ:IJd}, and hence $\tii(\f_j,\f)\to 0$. 

Assume, conversely, $\tii(\f_j,\f)\to 0$, \ie $\sup\f_j\to\sup\f$ and $\ii(\f_j,\f)\to 0$. Since $\d_v$ has finite energy for $v\in X^\div$ (Proposition~\ref{prop:enT}), Theorem~\ref{thm:BBGZ3} yields $\f_j(v)\to\f(v)$, \ie $\f_j\to\f$ weakly. By the quasi-triangle inequality, $\jj(\f_j)\approx\ii(\f_j,0)$ is further eventually bounded, and Proposition~\ref{prop:weakcont} therefore implies $\int\f_j\MA(\f)\to\int\f\MA(\f)$, which in turn shows that $\en(\f_j)\to\en(\f)$, by~\eqref{e509}. Thus $\f_j\to\f$ strongly, which concludes the proof. 
\end{proof}

By Example~\ref{exam:decrstrong}, any weakly convergent decreasing net in $\cE^1$ is strongly convergent.   The next result, which is a direct consequence of Theorem~\ref{thm:monostrong}, shows that this also holds for increasing nets:  

\begin{thm}\label{thm:incrstrong} Let $(\f_j)$ be an increasing net in $\cE^1$, and assume that $\f_j\to\f$ weakly in $\cE^1$. Then $\f_j\to\f$ strongly as well. 
\end{thm}

We end this section with a useful quantitative version of Corollary~\ref{cor:sameenergy}.
\begin{prop}\label{prop:convdom} Let $(\f_i)$ be a net in $\cE^1$ and $\f\in\cE^1$ such that $\f_i\le\f$ for all $i$. Then $\f_i\to\f$ strongly iff $\en(\f_i)\to\en(\f)$. 
\end{prop}
\begin{proof} The `only if' part follows directly from the definition of strong convergence. Conversely, assume $\en(\f_i)\to\en(\f)$. Since $\f_i\le\f$, 
$$
\ii(\f_i,\f)\approx\jj_\f(\f_i)=\en(\f)-\en(\f_i)+\int(\f_i-\f)\MA(\f)\le \en(\f)-\en(\f_i)
$$
tends to $0$, and hence $\f'_i:=\f_i-\sup\f_i$ converges strongly to $\f':=\f-\sup\f$, by Theorem~\ref{thm:strongunif}. In particular, $\en(\f'_i)=\en(\f_i)-\sup\f_i$ converges to $\en(\f')=\en(\f)-\sup\f$, hence $\sup\f_i\to\sup\f$, and we conclude, as desired, $\f_i\to\f$ strongly. 
\end{proof}

%
%
%
%
\subsection{Strong continuity and surjectivity of the Monge--Amp\`ere operator}\label{sec:CY}
We define the strong topology on $\cE^1/\R\simeq\cE^1_{\sup}$ as the one induced by the strong topology on $\cE^1$. It is defined by the quasi-metric $(\f,\f')\mapsto\ii(\f,\f')$, by Theorem~\ref{thm:strongunif}. 

As a direct consequence of Lemma~\ref{lem:trianglevee}, we have: 

\begin{prop}\label{prop:MAemb} The Monge--Amp\`ere operator $\MA\colon\cE^1\to\cM^1$ induces a bi-Lipschitz embedding
$$
(\cE^1/\R,\ii)\hookrightarrow (\cM^1,\ii^\vee).
$$
\end{prop}

Recall from Theorem~\ref{thm:envppty} that the envelope property for $\om\in\Amp(X)$ (which holds if $X$ is smooth and $k$ has characteristic zero, cf.~Theorem~\ref{thm:contenvlisse}) is equivalent to the compactness of $\PSH(\om)/\R$. In a similar vein we have: 

\begin{thm}\label{thm:contenvE1} For any $\om\in\Amp(X)$, the following are equivalent:
\begin{itemize}
\item[(i)] the envelope property holds for $\om$; 
\item[(ii)] the Monge--Amp\`ere operator $\MA_\om\colon\cE^1(\om)\to\cM^1$ is onto; 
\item[(iii)] the Monge--Amp\`ere operator induces a bi-Lipschitz isomorphism 
$$
(\cE^1(\om)/\R,\ii_\om)\simeq(\cM^1,\ii_\om^\vee);
$$
\item[(iv)] $(\cE^1(\om)/\R,\ii_\om)$ is complete; 
\item[(v)] $(\cE^1(\om),\tii_\om)$ is complete.
\end{itemize}
\end{thm}
In particular, $\cE^1(\om)$ can only be complete if $X$ is unibranch, cf.~Theorem~\ref{thm:unibranch}. 

\begin{proof} Assume (i). By Corollary~\ref{cor:MA}, (ii) holds iff the supremum defining
$$
\en^\vee(\mu)=\sup_{\f\in\cE^1}\left(\en(\f)-\int\f\,\mu\right)
$$
is achieved for any $\mu\in\cM^1$. Thus choose a maximizing sequence $(\f_i)$ for $\mu$ in $\cE^1_{\sup}$. By Corollary~\ref{cor:contenvcomp}, $\PSH_{\sup}$ is weakly compact, and we may thus assume, after passing to a subnet,  that $(\f_i)$ converges weakly to $\f\in\PSH_{\sup}$. Since $\jj_{\mu}(\f_i)\to 0$,~\eqref{equ:jjmu} shows that $\jj(\f_i)$ is eventually bounded, and hence that $\f\in\cE^1(\om)$, since $\jj$ is weakly lsc. By Proposition~\ref{prop:weakcont}, the weak convergence $\f_i\to\f$ implies $\int\f_i\,\mu\to\int\f\,\mu$, and hence 
$$
\en(\f)-\int\f\,\mu\ge\limsup_i\left(\en(\f_i)-\int\f_i\,\mu\right)=\en^\vee(\mu),
$$
since $\en$ is weakly usc. This proves (i)$\Rightarrow$(ii). Proposition~\ref{prop:MAemb} gives (ii)$\Leftrightarrow$(iii), and (iii)$\Rightarrow$(iv) follows from the completeness of $(\cM^1,\ii^\vee_\om)$ (Theorem~\ref{thm:M1complete}). 

Assume (iv), and pick a Cauchy sequence $(\f_i)$ in $(\cE^1(\om),\tii_\om)$. Then $(\f_i-\sup\f_i)$ is Cauchy in $\cE^1_{\sup}(\om)\simeq\cE^1(\om)/\R$, and hence strongly convergent in $\cE^1(\om)$. Now $(\sup\f_i)$ is Cauchy as well, and hence convergent in $\R$, so we conclude that $(\f_i)$ is strongly convergent, proving (iv)$\Rightarrow$(v). 

Finally we prove (v)$\Rightarrow$(i). Assume that $\cE^1(\om)$ is complete. By Theorem~\ref{thm:envppty}, (i) is equivalent to the fact that any bounded, increasing net $(\f_i)$ in $\PSH(\om)$ converges weakly in $\PSH(\om)$. Since the increasing net $(\sup\f_i)$ converges in $\R$, it will be enough to show that $(\f_i)$ is Cauchy in $\cE^1(\om)/\R$, \ie that
$$
\ii_\om(\f_i,\f_j)\approx\jj_{\om,\f_i}(\f_j)=\en_\om(\f_i)-\en_\om(\f_j)+\int(\f_j-\f_i)\MA_\om(\f_i)
$$
tends to $0$ as $i,j\to\infty$, where by symmetry of $\ii_\om$ we may assume $i\le j$. The increasing net $(\en_\om(\f_i))$ converges in $\R$. Given $\e>0$, choose $i_0$ such that $0\le\en_\om(\f_j)-\en_\om(\f_i)\le\e$ for $j\ge i\ge i_0$.
As $\f_i\le\f_j$,~\eqref{equ:Ebis} yields
$$
0\le\int(\f_j-\f_i)\MA_\om(\f_i)\le(n+1)\left(\en_\om(\f_j)-\en_\om(\f_i)\right)\le(n+1)\e.
$$
Thus 
$$
\en_\om(\f_i)-\en_\om(\f_j)+\int(\f_j-\f_i)\MA_\om(\f_i)\le n\e,
$$
which proves, as desired, that $(\f_i)$ is Cauchy in $\cE^1(\om)/\R$.  
\end{proof}
 
As a direct consequence of Theorem~\ref{thm:contenvE1} and Theorem~\ref{thm:contenvlisse}, we have the following result, \cf Theorem~A.
\begin{thm}\label{thm:CY}
  Suppose that $X$ is smooth, that $\charac k=0$ or $\dim X\le 2$, and fix $\om\in\Amp(X)$. Then the Monge--Amp\`ere equation $\MA(\f)=\mu$ admits a solution $\f\in\cE^1(\om)$, unique up to translation,  for every $\mu\in\cM^1$.
\end{thm}

\begin{exam}\label{exam:phiv} Assume the envelope property for $\om\in\Amp(X)$. For each $v\in X^\lin$, there exists a unique $\f_v\in\cE^1(\om)$ such that $\MA(\f_v)=\d_v$ and $\f_v(v)=0$, which we call the \emph{Green's function} of $v$. See Proposition~\ref{prop:phiv} for more on these functions. 
\end{exam}

\begin{exam}\label{exam:nodal} Let $X$ be an irreducible curve with a nodal singularity at $p\in X(k)$, and denote by $p_1,p_2\in X^\nu(k)$ the preimages of $p$ (see Figure~\ref{fig:curve}). Then $v:=\ord_{p_1}\in X^{\nu,\,\div}\simeq X^\div$, and $\d_v\in\cM^1(X)$ is not in the image of $\MA\colon\cE^1\to\cM^1$. Using Example~\ref{exam:MAcurve} and $\nu_\star\MA(\nu^\star\f)=\MA(\f)$ (see~\eqref{equ:mixedMAnorm}), one sees indeed that a function $\f\in\PSH(\om)$ such that $\MA(\f)=\d_v$ must satisfy $\nu^\star\f(v_\triv)=\nu^\star\f(v_{p_2,\triv})=\f(v_{p,\triv})=\nu^\star\f(v_{p_1,\triv})$, and hence must be constant on the ray through $v$, a contradiction. 
\end{exam}

%
%
\subsection{Continuity of solutions to Monge--Amp\`ere equations}

\begin{thm}\label{thm:contMA}
Let $\f\in\cE^1$, and suppose that $\mu:=\MA(\f)$ is supported in a bounded subset of $(X^\lin,\dd_\infty)$, \ie $\supp\mu\subset\{\te\le C\}$ for some $C>0$. Then $\f\in\CPSH(\om)$. 
\end{thm}

\begin{exam} The condition on $\mu$ is satisfied if $\charac k=0$, $X$ is smooth, and $\mu$ is supported on the dual complex of an snc test configuration for $X$ (see Theorem~\ref{thm:pshdual} below). 
\end{exam}

\begin{proof}[Proof of Theorem~\ref{thm:contMA}] Write $\f$ as the limit of a decreasing net $(\f_i)$ in $\CPSH(\om)$. By Theorem~\ref{thm:pshweakcont}, $\f$ is continuous on the (weakly) compact set $K:=\supp\mu\subset X^\an$, and hence $\f_i\to\f$ uniformly on $K$, by Dini's lemma. For any $\e>0$ we therefore have $\f_i\le\f+\e$ on $K$ for $i$ large enough. Thus $\f_i\le\f+\e$ a.e.~for $\MA(\f+\e)=\mu$, and hence $\f_i\le\f+\e$ on $X^\an$, by the domination principle (Corollary~\ref{cor:domprinc}). This shows $\f_i\to\f$ uniformly on $X^\an$, and hence $\f\in\Cz(X)$. \end{proof}

\begin{prop}\label{prop:phiv} Assume the envelope property for $\om$. Pick $v\in X^\lin$, and recall from Example~\ref{exam:phiv} that the Green's function $\f_v\in\cE^1$ is the unique function such that $\MA(\f_v)=\d_v$ and $\f_v(v)=0$. Then:
\begin{itemize}
\item[(i)] $\f_v\in\CPSH(X)$;
\item[(ii)] $\f_v=\sup\{\f\in\PSH(\om)\mid\f(v)\le 0\}$; 
\item[(iii)] for all $v,w\in X^\lin$ we have $\dd_\infty(v,w)=\max\{\f_v(w),\f_w(v)\}=\sup_{X^\an}|\f_v-\f_w|$;
\item[(iv)] $\sup\f_v=\te(v)$. 
\end{itemize}
If $\om=c_1(L)$ for an ample bundle $L$, then $\f_v$ is the uniform limit as $m\to\infty$ of the functions
$$
\f_{v,m}:=m^{-1}\max_{s\in\Hnot(X,mL)\setminus\{0\}}\left(\log|s|+v(s)\right)\in\cH_\R(L)
$$
from Lemma~\ref{lem:Izumi}. 
\end{prop}
Thus (iii) yields an isometric embedding $(X^\lin,\dd_\infty)\hookrightarrow\CPSH(\om)$ with respect to the supnorm. 

\begin{proof} Assertions (i) and (ii) are direct consequences of Theorem~\ref{thm:contMA} and the domination principle (Corollary~\ref{cor:domprinc}), respectively. To establish (iii), note first that $\f_v(w)=\f_v(w)-\f_v(v)\le \dd_\infty(v,w)$, and hence $\dd_\infty(v,w)\ge\max\{\f_v(w),\f_w(v)\}$, by symmetry. Conversely pick $\f\in\PSH(\om)$. By (ii), $\f-\f(v)\le\f_v$. Thus $\f(w)-\f(v)\le\f_v(w)$, and hence 
$$
\dd_\infty(v,w)=\sup_{\f\in\PSH}|\f(v)-\f(w)|\le\max\{\f_v(w),\f_w(v)\},
$$
which proves (iii), of which (iv) is the special case $w=v_\triv$, see Proposition~\ref{prop:Izumi}~(ii).

To establish the last statement in the proposition, we first note, as in the proof of Lemma~\ref{lem:Izumi}, that the sequence $m\mapsto m\f_{v,m}$ is superadditive, and hence 
$$
\lim_{m\to\infty}\f_{v,m}=\sup_{m\to\infty}\f_{v,m}
$$
pointwise on $X^\an$, by Fekete's lemma. On the one hand, $\f_{v,m}(v)=0$ yields $\sup_{m\to\infty}\f_{v,m}\le\f_v$, by (ii). On the other hand, the proof of Lemma~\ref{lem:Izumi} shows that any $\f\in\cH(L)$ with $\f(v)=0$ satisfies $\f\le\f_{v,m}$ for some $m\gg 1$, and hence 
$$
\f_v=\sup\{\f\in\cH(L)\mid \f(v)=0\}\le\sup_m\f_{v,m},
$$
where the left-hand equality is an easy consequence of (ii). We conclude $\f_{v,m}\to\f_v$ pointwise on $X^\an$, and a simple variant of Dini's lemma, based on the superadditivity of $m\mapsto m\f_{v,m}$, shows that the convergence is uniform. 
\end{proof}

%
%
%
%
%
We conclude this section with the following consequence of Proposition~\ref{prop:phiv}: 

\begin{cor}\label{cor:charweaktop} Assume that $\om\in\Amp(X)$ has the envelope property.  
 Let $\f\in\PSH(\om)$, and let $(\f_j)$ be a net in $\PSH(\om)$. Then the following are equivalent:
  \begin{itemize}
  \item[(i)]
    $\f_j\to\f$ in $\PSH(\om)$; 
  \item[(ii)]
    $\lim_j\int\f_j\MA(\p)=\int\f\MA(\p)$ for every $\p\in\PL\cap\PSH(\om)$; 
  \item[(iii)] 
  $\lim_j\int\f_j\MA(\p)=\int\f\MA(\p)$ for every $\p\in\CPSH(\om)$.
   \end{itemize}
\end{cor}

\begin{proof} By Proposition~\ref{prop:mixedMA}~(ii) , for each $\p\in\PL\cap\PSH(\om)$ the measure $\MA(\p)$ is supported in a finite subset of $X^\div$, and the pointwise convergence $\f_i\to\f$ on $X^\div$ thus implies $\lim_i\int\f_i\MA(\p)\to\int\f\MA(\p)$, \ie (i)$\Rightarrow$(ii). Assume (ii), pick $\p\in\CPSH(\om)$, and choose a sequence $(\p_m)$ in $\PL\cap\PSH(\om)$. By Lemma~\ref{lem:CLN}, we have 
$$
\left|\int\f\MA(\p_m)-\int\f\MA(\p)\right|\le n\sup|\p-\p_m|
$$
and
$$
\left|\int\f_i\MA(\p_m)-\int\f_i\MA(\p)\right|\le n\sup|\p-\p_m|
$$
for all $i$, and (ii)$\Rightarrow$(iii) follows easily. Finally, (iii) implies in particular 
$$
\f_i(v)=\int\f_i\MA(\f_v)\to\f(v)=\int\f\MA(\f_v)
$$ 
for each $v\in X^\div$ , by  Proposition~\ref{prop:phiv}. Hence (iii)$\Rightarrow$(i). 
 \end{proof}
 
 %
%
%
%
\subsection{Countable regularization and Choquet's lemma}\label{sec:countreg}
In this final section, $X$ is again allowed to be \textbf{reducible}. Using the energy functionals, we establish the following \emph{countable convergence} result.

\begin{thm}\label{thm:moncount} Pick $\theta\in\Num(X)$. Let $(\f_i)_{i\in I}$ be a monotone (\ie increasing or decreasing) net in $\PSH(\theta)$, and assume that it converges to $\f\in\PSH(\theta)$. Then there exists an increasing map $\N\to I$ $m\mapsto i(m)$ such that $\lim_{m\to\infty}\f_{i(m)}=\f$ in $\PSH(\theta)$. 
\end{thm}

\begin{rmk}\label{rmk:moncount} Note that $(\f_{i(m)})_{m\in\N}$ is not claimed to be a subnet of $(\f_i)_{i\in I}$, as the directed set $I$ might not admit any countable cofinal subset. On the other hand, by monotonicity of $(\f_i)$, the sequence $(\f_{i'(m)})$ also converges to $\f$ for any increasing map $m\mapsto i'(m)$ such that $i'(m)\ge i(m)$. 
\end{rmk}

\begin{proof} Assume first that $X$ is irreducible. Pick $\theta'\in\Amp(X)$ such that $\theta'-\theta$ is nef. Then $\PSH(\theta)\subset\PSH(\theta')$, and it is enough to show the result when $\om:=\theta$ is ample, which we now assume. 

Suppose first that $\f$ and all the $\f_i$ lie in $\cE^1(\om)$. By Example~\ref{exam:decrstrong} (resp.~Theorem~\ref{thm:incrstrong}), we have $\f_i\to\f$ strongly, \ie $\lim_{i\in I}\tii(\f_i,\f)= 0$. We can thus inductively construct an increasing map $m\mapsto i(m)$ such that $\tii(\f_{i(m)},\f)\le 1/m$. Then $\lim_m\f_{i(m)}=\f$ strongly in $\cE^1(\om)$, and hence in $\PSH(\om)$ (\ie pointwise on $X^\div$). 

Consider now the general case. Observe first that we may assume wlog $\f\le 0$ and $\f_i\le 0$ for all $i$. If $(\f_i)$ is increasing, it suffices to replace $\f_i$ with $\f_i-\sup\f$. If $(\f_i)$ is decreasing, replace $I$ with $\{i\in I\mid i\ge i_0\}$ for some $i_0\in I$, and $\f_i$ with $\f_i-\sup\f_{i_0}$. 

By Corollary~\ref{cor:exppsh}, we then have 
  $$
  \tilde\f:=\exp(\f),\,\tilde\f_i:=\exp(\f_i)\in\cE^\infty(\om)\subset\cE^1(\om). 
  $$
 By what precedes, we can find an increasing map $m\mapsto i(m)$ such that $\tilde\f_{i(m)}\to\tilde\f$ pointwise on $X^\div$. Since the $\f_i$ and $\f$ are finite-valued on $X^\div$, this implies $\f_{i(m)}\to\f$ on $X^\div$, and we are done.
 
Assume finally that $X$ is reducible. For each $\a$, the previous step yields an increasing $m\mapsto i_\a(m)$ such that $\f_{i_\a(m)}\to\f$ on $X_\a^\an$. By Remark~\ref{rmk:moncount}, it remains to pick any increasing map $m\mapsto i(m)$ such that $i(m)\ge\max_\a i_\a(m)$ for all $m$; such a map is easily constructed by induction on $m$. 
\end{proof}

Combining Theorem~\ref{thm:moncount} with Theorem~\ref{thm:pshample}, we infer: 

\begin{cor}\label{cor:countreg} Pick $\theta\in\Num(X)$ and $\f\in\PSH(\theta)$. Then:
\begin{itemize}
\item[(i)] $\f$ can be written as the limit of a decreasing sequence $\f_m\in\cH^\gf_\Q(L_m)$ with $L_m\in\Pic(X)_\Q$ such that $c_1(L_m)-\theta$ is ample and $\lim_m c_1(L_m)=\theta$; 
\item[(ii)] if $\theta\in\Nef(X)$ then (i) holds with $\f_m\in\cH(L_m)$ and $L_m$ ample; 
\item[(iii)] if $\theta\in\Amp(X)$, then $\f$ can be written as the limit of a decreasing sequence in $\cH^\dom(\theta)$ (see Definition~\ref{defi:Hdom}). 
\end{itemize}
\end{cor}

\begin{rmk} A version of (iii) was proved in~\cite[Proposition 4.7]{nama}, using the Bedford--Taylor capacity instead of the energy functionals.
\end{rmk}

This implies in turn that psh functions are non-constant only on relatively small subsets of $X^\an$, when the ground field $k$ is uncountable. 

\begin{thm}\label{thm:almostconstant} For any $\theta\in\Num(X)$ and $\f\in\PSH(\theta)$, there exists a countable intersection of non-empty Zariski open subsets $U\subset X$ such that $\f\equiv\sup\f$ on $c^{-1}(U)$.  
\end{thm}
Recall that $c\colon X^\an\to X$ denotes the center map, which is anticontinuous. 

\begin{proof} By Corollary~\ref{cor:countreg}, we can write $\f$ as the limit of a decreasing sequence of generically finite Fubini--Study functions $\f_m$. By Lemma~\ref{lem:almostconstant}, for each $m$ there exists a non-empty Zariski open subset $U_m\subset X$ such that $\f_m\equiv\sup\f_m$ on $c^{-1}(U_m)$. Now $\sup\f_m\to\sup\f$ (see Example~\ref{exam:supcont}), and the result follows with $U:=\bigcap_m U_m$. 
\end{proof}

\begin{cor}\label{cor:ppnotdense} Any pluripolar set $E\subset X^\an$ is contained in a countable union of open sets of the form $c^{-1}(Z)$, with $Z\subset X$ a strict Zariski closed subset. If $k$ is uncountable, then $E$ cannot be dense. 
\end{cor}
\begin{proof} There exists $\theta\in\Num(X)$ and $\f\in\PSH(\theta)$ such that $E\subset\{\f=-\infty\}\subset\{\f<\sup\f\}$, and the first point thus follows from Theorem~\ref{thm:almostconstant}. When $k$ is uncountable, $X$ cannot be written as a countable union of strict Zariski closed subsets, and we can thus find $p\in X(k)$ such that $E$ is disjoint from $c^{-1}(\{p\})\subset X^\an$. Now the latter is open and non-empty (as it contains, for instance, the trivial semivaluation $v_{p,\triv}$). This prevents $E$ from being dense. 
\end{proof}

As a further consequence of Theorem~\ref{thm:moncount}, we also get a version of \emph{Choquet's lemma}, which will be put to use in Section~\ref{sec:neglpp}. 

\begin{cor}\label{cor:Choquet}
  Assume that $\theta\in\Num(X)$ have the envelope property. Let $(\f_\a)_{\a\in A}$ be a family in $\PSH(\theta)$ that is uniformly bounded from above. 
  Then there exists an at most countable subset $B\subset A$ such that 
  \begin{equation*}
    \supstar_{\a\in A}\f_\a = \supstar_{\a\in B}\f_\a.
  \end{equation*}
\end{cor}
\begin{proof} By assumption, $\f:=\supstar_{\a\in A}\f_\a$ is $\theta$-psh. Denote by $(A_i)_{i\in I}$ the net of all finite subsets of $A$, and for each $i\in I$ set $\p_i:=\max_{\a\in A_i}\f_\a$. Then $(\p_i)_{i\in I}$ is an increasing net in $\PSH(\om)$, which converges to $\f$. By Theorem~\ref{thm:moncount}, we can find an   increasing map $m\mapsto i(m)$ such that $\p_{i(m)}\to\f$. Then $B:=\bigcup_m A_{i(m)}$ is at most countable , and we have $\sup_{\a\in B}\f_\a=\sup_{\a\in A}\f_\a$ on $X^\div$, and hence $\supstar_{\a\in B}\f_\a=\supstar_{\a\in A}\f_\a$ on $X^\an$ by Corollary~\ref{cor:Hausdorff}, since both functions are $\theta$-psh. 
\end{proof}
%
%
%
%
\section{Capacities, pluripolar sets, and negligible sets}\label{sec:neglpp}
Assuming the envelope property, we make a more detailed study of pluripotential theory, adapting to our setting classical arguments from the complex analytic case. In particular, we prove the fundamental result that negligible sets are the same as pluripolar sets.

In what follows, $X$ is of dimension $n$ and \textbf{irreducible}. We fix an ample class $\om\in\Amp(X)$, and we assume that \textbf{the envelope property} holds for $\om$, see~\S\ref{sec:envprop}. As before, we write $V=V_\om=(\om^n)$ and $\MA(\f):=V^{-1}(\om+\ddc\f)^n$ for $\f\in\cE^1=\cE^1(\om)$.
%
%
\subsection{Envelopes of lsc functions}
Consider an arbitrary function $\f:X^\an\to\R\cup\{+\infty\}$, and assume that $\f$ is bounded below. By Lemma~\ref{lem:contenv}, the usc regularization 
$$
\env(\f)^\star=\supstar\left\{\p\in\PSH\mid\p\le\f\right\} 
$$
of the $\om$-psh envelope $\env(\f)$ is either $\om$-psh and bounded, or satisfies $\env(\f)^\star\equiv+\infty$.

\begin{thm}\label{thm:ortho} Pick any lsc function $\f:X^\an\to\R\cup\{+\infty\}$, and assume $\env(\f)^\star\not\equiv+\infty$. Then 
$$
\env(\f)^\star=\env(\f)=\f\text{ a.e.~for }\MA(\env(\f)^\star).
$$
\end{thm}
In particular, the \emph{orthogonality property}
\begin{equation}\label{equ:ortho}
\int_{X^\an}\left(\f-\env(\f)^\star\right)\MA(\env(\f)^\star)=0
\end{equation}
holds (compare~\cite{nama,BGM}).

\begin{lem}\label{lem:neglnegl} If $N\subset X^\an$ is a negligible Borel set (see Definition~\ref{defi:negl}), then $\mu(N)=0$ for all measures of finite energy $\mu\in\cM^1$. 
\end{lem}
This strengthens in particular Lemma~\ref{lem:M1npp}, since pluripolar sets are negligible (Proposition~\ref{prop:plurineg}), 

\begin{proof} By Choquet's lemma (see Corollary~\ref{cor:Choquet}), there exists an increasing, uniformly bounded
  sequence $(\f_m)_1^\infty$ in $\PSH(\om)$ such that $N\subset\{\f<\f^\star\}$, where
  $\f:=\lim_m\f_m$ pointwise. Recall that $\f^\star\in\PSH(\om)$, thanks to the standing assumption that the envelope property holds.  On $X^\div$ we have $\f^\star=\f$, by Theorem~\ref{thm:divnegl}. Thus $\f_m\to\f^\star$ weakly in $\cE^1$ (and in fact strongly as well, by Theorem~\ref{thm:incrstrong}), and hence 
 $\int\f_m\,\mu\to\int\f^\star\,\mu$, by Proposition~\ref{prop:weakcont}. On the other hand, the monotone convergence theorem shows that $\int\f_m\,\mu\to\int\f\,\mu$. Thus $\int(\f^\star -\f)\mu=0$, which implies $\mu(N)=0$.
 \end{proof}

\begin{proof}[Proof of Theorem~\ref{thm:ortho}] Set $\p:=\env(\f)^\star$ and $\mu:=\MA(\p)$. Assume first $\f\in\Cz(X)$, so that $\p=\env(\f)$ is continuous and $\om$-psh. Arguing as in the proof of~\cite[Lemma 3.5]{BGM}, set $f:=\f-\p\ge 0$, and observe that for $t\in[0,1]$, 
$$
\f_t:=\p+tf=(1-t)\env(\f)+t\f
$$ 
satisfies $\env(\f)\le\f_t\le\f$, and hence $\env(\f_t)=\env(\f)$. Thus
$\en(\p+t f)=\en(\env(\f))$ for $t\in[0,1]$, and 
Corollary~\ref{cor:diff} yields 
$$
\int(\f-\p)\MA(\p)=\frac{d}{dt}\bigg|_{t=0}\en(\p+t f)=0,
$$
and hence $\f=\p$ $\mu$-a.e.

Now consider the general case. Since the Borel set $\{\p>\env(\f)\}$ is negligible, Lemma~\ref{lem:neglnegl} shows that $\p=\env(\f)$ $\mu$-a.e. In particular, $\p\le\f$ $\mu$-a.e., and it thus remains to show that $\mu$ puts no mass on the open set $U:=\{\p<\f\}$. 

Since $\f$ is lsc and $X^\an$ is compact, there exists an increasing net $(\f_j)$ in $\Cz(X)$ that converges pointwise to $\f$. By Lemma~\ref{lem:envlsc}, the increasing net $\p_j:=\env(\f_j)$ converges pointwise to $\env(\f)$. By Theorem~\ref{thm:divnegl}, we thus have $\p_j\to\p$ on $X^\div$, and Theorem~\ref{thm:incrstrong} implies that $\mu_j:=\MA(\p_j)$ converges weakly to $\mu$. 
  By the first part of the proof, $\mu_j$ puts no mass on the open set $U_j:=\{\p_j<\f_j\}$, and we want to take the limit and deduce that $\mu(U)=0$; this will
  conclude the proof.
  To do so, suppose $\mu(U)>0$, and pick a point $v\in U\cap\supp\mu$. 
  Thus $\p(v)<\f(v)$, so there exists $j_0$ such that $\p(v)<\f_j(v)$
  for $j\ge j_0$. Set $V:=\{\p<\f_{j_0}\}$. Then $V$ is open, $v\in V$ and $V\subset U_j$
  for $j\ge j_0$. Since $v\in\supp\mu$, we have $\mu(V)>0$.
  As $\mu_j\to\mu$ weakly, we must have $\mu_j(V)>0$ for $j\gg0$, which contradicts $\mu_j(U_j)=0$. 
  \end{proof}
%
%
\subsection{The Bedford--Taylor capacity, reprise}\label{sec:BT2}
Recall from~\S\ref{sec:BT1} that we have set, for each Borel set $E\subset X^\an$, 
$$
    \Capa(E)=\Capa_\om(E):=\sup\left\{\int_E\MA(\p)\mid \p\in\PSH, -1\le \p\le 0\right\}.
$$
Note that $0\le\Capa(E)\le 1$ for all Borel sets $E\subset X^\an$, and that $\Capa(X)=1$.

\begin{prop}\label{prop:BT}
  The function $\Capa$ is a \emph{precapacity}, \ie it satisfies: 
\begin{equation}\label{equ:capamono} 
 E\subset E'\Longrightarrow\Capa(E)\le\Capa(E')
 \end{equation}
 for all Borel subsets $E,E'$, and 
\begin{equation}\label{equ:capacont}
\Capa\left(\bigcup_m E_m\right)=\sup_m\Capa(E_m)
\end{equation}
for any increasing sequence $(E_m)$ of Borel subsets. 

It is further \emph{subadditive}, \ie 
\begin{equation}\label{equ:capasub}
\Capa\left(\bigcup E_m\right)\le\sum_j\Capa(E_m)
\end{equation}
for any sequence $(E_m)$ of Borel subsets, and \emph{inner regular}, \ie
\begin{equation}\label{equ:capainner}
 \Capa(E)=\sup\left\{\Capa(K)\mid K\subset E\ \text{compact}\right\}
\end{equation}
for any Borel subset $E$. 
  \end{prop}
\begin{proof}
  The first three properties follow immediately from 
  the definition of the capacity as the supremum of 
  the Borel measures $\MA(\p)$, and the last one likewise follows
  from the fact that these Borel measures are Radon
  measures, and hence inner regular on all Borel sets.
\end{proof}

We define the \emph{outer capacity} of any subset $E\subset\Xan$ as
\begin{equation*}
  \Capa^\star (E):=\inf\left\{\Capa(U)\mid U\supset E\ \text{open}\right\}.
\end{equation*}
Trivially, $\Capa^\star(U)=\Capa(U)$ for any open $U\subset X^\an$, and $\Capa^\star$ is outer regular, \ie
\begin{equation}\label{equ:capaouter}
\Capa^\star(E)=\inf_{U\ \text{open $\supset E$}}\Capa^\star(U)
\end{equation}
for all $E$. Proposition~\ref{prop:BT} immediately yields: 

\begin{prop}\label{prop:outcapprops}
The function $\Capa^\star$ is monotone and subadditive, \ie the analogues of~\eqref{equ:capamono} and~\eqref{equ:capasub} hold for arbitrary subsets of $X^\an$. 
\end{prop}
As we shall see, $\Capa^\star$ also satisfies the continuity condition~\eqref{equ:capacont},  but this is much more involved, see Theorem~\ref{thm:outcapincr} below.

\begin{exam}\label{exam:capaopen} Every nonempty open subset $U\subset X^\an$ has $\Capa(U)>0$. Indeed, $U$ must contain a divisorial valuation $v$, by density of $X^\div$. For $0<\e\ll 1$, $\p:=\e\f_v-1$ is a candidate in~\eqref{equ:BTcap}, and hence 
\begin{equation}\label{equ:capaopen}
\Capa(U)\ge\int_U\MA(\p)\ge\e^n\int_U\MA(\f_v)\ge\e^n. 
\end{equation}
\end{exam}

By Proposition~\ref{prop:envar} and Lemma~\ref{lem:neglnegl}, we have, on the other hand: 
\begin{exam}\label{exam:negcapzero} Every negligible Borel set $N\subset X^\an$ satisfies $\Capa(N)=0$. 
\end{exam}

As we shall see below, this holds in fact for $\Capa^\star$ as well, see Lemma~\ref{lem:neglpp}. 

To conclude this section, we consider the dependence of $\Capa=\Capa_\om$ on $\om\in\Amp(X)$. 

\begin{thm}\label{thm:BTclass}
  For any ample classes $\om,\om'\in\Amp(X)$ there exists a constant $C=C(\om,\om')\ge 1$ such that $C^{-1}\Capa_\om\le\Capa_{\om'}\le C\Capa_\om$.
\end{thm}
This is a direct consequence of the following more precise result. 

\begin{lem}\label{lem:BTclass}
  Let $\om,\om'\in\Amp(X)$. Then:
  \begin{itemize}
  \item[(i)]
    if $\om\le\om'$, then $V_\om\Capa_\om\le V_{\om'}\Capa_{\om'}$;
  \item[(ii)]
    if $t\ge1$, then $t^{-n}\Capa_\om\le\Capa_{t\om}\le\Capa_\om$;
  \end{itemize}
\end{lem}
\begin{proof}
  If $\om\le\om'$, then $\PSH(\om)\subset\PSH(\om')$ and
  \begin{equation*}
    V_\om\MA_\om(\p)
    =(\om+\ddc \p)^n
    \le(\om'+\ddc \p)^n
    =V_{\om'}\MA_{\om'}(\p)
  \end{equation*}
  for any $\p\in\PSH(\om)$ with $-1\le \p\le 0$.
  This proves~(i).

  The first inequality in~(ii) follows from~(i). For the second, suppose $\f\in\PSH(\om)$ and $-1\le \p\le 0$. Then $t^{-1}\p\in\PSH(\om)$, $-1\le t^{-1}\f\le 0$,  and the Radon measures $\MA_\om(t^{-1}\p)$ and $\MA_{t\om}(\p)$ coincide. Thus $\Capa_{t\om}\le\Capa_\om$.
\end{proof}

%
%
\subsection{Quasicontinuity}
While $\om$-psh functions are not continuous in general, they satisfy the following quasicontinuity property, reminiscent of Lusin's theorem in measure theory. 

\begin{thm}\label{thm:qcont}
For each $\f\in\PSH(\om)$ and $\e>0$, there exists a compact $K\subset X^\an$ such that $\Capa(K^c)\le\e$ and $\f|_K\in\Cz(K)$.  
\end{thm}
Here $K^c:=X^\an\setminus K$.

\begin{cor}\label{cor:qcont} Given a sequence $(\f_m)$ in $\PSH(\om)$ and $\e>0$ there exists a compact $K\subset X^\an$ such that $\Capa(K^c)\le \e$ and $\f_m|_K\in\Cz(K)$ for all $m$.
\end{cor}

\begin{proof}[Proof of Theorem~\ref{thm:qcont}] After adding a constant, we may assume $\sup\f=0$, \ie $\f\in\PSH_{\sup}$.
  First suppose that $\f$ is bounded. Let $(\f_j)_j$ be a decreasing net in $\PL(X)\cap\PSH(\om)$ converging to $\f$, and pick $\p\in\PSH(\om)$ with $-1\le \p\le 0$. By Theorem~\ref{thm:BBGZ3} $\int(\f_j-\f)\MA(\p)\to 0$ uniformly with respect to $\p$. Using Chebyshev's inequality and the definition of the Bedford--Taylor capacity we can thus find, for every integer $m\ge 1$,
  an index $j_m\in I$ such that $j_{m+1}>j_m$ and
  the compact set 
  $$
  K_m:=\{\f_{j_m}-\f\le m^{-1}\}
  $$
  satisfies $\Capa(K_m^c)\le 2^{-m}\e$. If we set $K:=\bigcap_mK_m$, then $\Capa(K^c)\le\sum_m\Capa(K_m^c)\le\e$ by countable subadditivity of $\Capa$, and $\f_{j_m}\to\f$ uniformly on $K$, so that $\f|_K$ is continuous. 

  Consider now a possibly unbounded $\f\in\PSH_{\sup}$ and $\e>0$.
  By Lemma~\ref{lem:capsublevel1}, we can find $t>0$ such that $K':=\{\f\ge -t\}$ satisfies $\Capa(K'^c)\le \e$.
  Applying what precedes to the bounded function $\max\{\f,-t\}$ gives a compact $K\subset K'$ such that $\Capa(K^c)\le 2\e$ and $\f|_K=\max\{\f,-t\}|_K$ is continuous, and we are done. 
\end{proof}

\begin{proof}[Proof of Corollary~\ref{cor:qcont}] Pick $\e>0$. For each $m$ we find a compact $K_m\subset X^\an$ such that $\Capa(K^c_m)\le \e 2^{-m}$ and $\f_m|_{K_m}\in\Cz(X)$. Then $K:=\bigcap_m K_m$ is a compact such that $\Capa(K^c)\le\sum_m\Capa(K^c_m)\le\e$, and $\f_m|_K$ is continuous for all $m$. 
\end{proof}

\begin{rmk} Corollary~\ref{cor:qcont} implies Theorem~\ref{thm:pshweakcont}.
  Indeed, arguing as in Example~\ref{exam:capaopen} (or using the Alexander--Taylor inequality~\eqref{equ:AT} below) shows that given $C>0$ there exists $\e>0$ such that any $v\in X^\lin$ with $\te(v)\le C$ satisfies $\Capa(\{v\})\ge 2\e$. For any $\f\in\PSH(\om)$, the compact $K$ provided by Theorem~\ref{thm:qcont} must then contain $\{T\le C\}$, and $\f$ must thus be weakly continuous thereon. 
\end{rmk}

%
%
\subsection{Extremal functions}
We now look for functions achieving the supremum in~\eqref{equ:BTcap}.
Let $E\subset\Xan$ be any subset.
The \emph{extremal function} of $E$ is 
\begin{equation}\label{equ:extfcn}
  \f_E:=\sup\{\f\in\PSH\mid \f\le 0, \f|_E\le -1\}=\env(-\one_E).
\end{equation}
Note that $-1\le\f_E\le 0$, and $\f_E\equiv-1$ on $E$, since the function $\f\equiv-1$ is
a competitor in~\eqref{equ:extfcn}.
Its usc regularization $\f_E^\star$ lies in $\PSH(\om)$, by the envelope property,
and also satisfies $-1\le\f_E^\star\le 0$, but it may happen that $\f_E^\star\not\equiv-1$ on $E$.

\begin{thm}\label{thm:capformula}
  If $K\subset X$ is compact, then
  \begin{equation}\label{equ:capformula}
    \Capa(K)
    =\int\limits_K\MA(\f_K^\star )
    =\int\limits_K(-\f_K^\star )\MA(\f_K^\star )
    =\int_{X^\an}(-\f_K^\star )\MA(\f_K^\star ).
  \end{equation}
\end{thm}
We refer to Theorem~\ref{thm:outcapformula} below for a version of this result for arbitrary subsets of $X^\an$.

\begin{proof} Since $\f_K$ is the psh envelope of the lsc function $-\one_K$, Theorem~\ref{thm:ortho} implies that $\f_K^\star=-\one_K$ a.e.~for $\MA(\f_K^\star)$, and hence  
\begin{equation*}
    \int_{X^\an}(-\f_K^\star )\MA(\f_K^\star )
    =\int_K(-\f_K^\star )\MA(\f_K^\star )
    =\int_K\MA(\f_K^\star ).
  \end{equation*}
 
  Since $-1\le \f_K^\star \le 0$, we have $\Capa(K)\ge\int_K\MA(\f_K^\star )$.
  To prove the reverse inequality, pick any $\p\in\PSH(\om)$ with $-1\le \p\le 0$, and $t\in(0,1)$. By Lemma~\ref{lem:neglnegl}, $t\p> -1=\f_K^\star$ on $K$ a.e.~for $\MA(t\p)$, and hence
$$
t^n\int_K\MA(\p)\le\int\limits_K\MA(t\p)\le\int\limits_{\{\f_K^\star <t\p\}}\MA(t\p)
    \le\int\limits_{\{\f_K^\star <t\p\}}\MA(\f_K^\star )
    \le\int\limits_K\MA(\f_K^\star ), 
$$
where the third inequality follows from the comparison principle,
  and the last inequality from the fact that $t\p<0=\f_K^\star$ outside $K$ a.e.~for $\MA(t\p)$. Letting $t\to 1$ yields the reverse inequality $\Capa(K)\le\int_K\MA(\f_K^\star )$.
\end{proof}

\begin{cor}\label{cor:capaouter}
The Bedford--Taylor capacity is outer regular on compact sets, \ie $\Capa^\star(K)=\Capa(K)$ for all compact $K\subset X^\an$. 
\end{cor}
\begin{proof} On the one hand, $\Capa^\star(K)$ is the decreasing limit of $\Capa(L)$ for $L$ ranging over the directed sets of all compact neighborhoods of $K$. On the other hand, $(-\one_L)$ forms an increasing net of lsc functions converging pointwise to $-\one_K$. By Lemma~\ref{lem:envlsc}, $\f_L^\star=\envstar(-\one_L)$ therefore converges strongly to $\f_K^\star=\envstar(-\one_K)$ in $\cE^1$, and hence 
$$
\Capa(L)=\int(-\f_{L}^\star)\MA(\f_{L}^\star)\to\Capa(K)=\int(-\f_K^\star)\MA(\f_K^\star),
$$
by Theorem~\ref{thm:BBGZ3}. Thus $\Capa^\star(K)=\Capa(K)$. 
\end{proof}

%
%
\subsection{Negligible sets are pluripolar}
Recall from~\S\ref{sec:pluripolar} that 
$$
\te(E)=\sup\left\{\sup\f-\sup_E\f\mid\f\in\PSH\right\}\in[0,+\infty]
$$
for any subset $E\subset X^\an$. The next result is a direct analogue of~\cite[Proposition 7.1]{GZ1}, itself an adaptation of the Alexander--Taylor inequality~\cite{AT84}. 

\begin{thm}\label{thm:AT}
  For any subset $E\subset\Xan$ we have 
  \begin{equation}\label{equ:AT}
   \min\{1,\te(E)^{-n}\}\le\Capa^\star (E)\le n\te(E)^{-1}. 
    \end{equation}
\end{thm}

Before entering the proof, we attach to each subset $E$ another extremal function $V_E\colon X^\an\to[0,+\infty]$, defined as 
\begin{equation*}
  V_E:=\sup\{\f\in\PSH\mid\f\le0\ \text{on $E$}\}=\env(f_E)
\end{equation*}
with $f_E\colon X\to\R\cup\{+\infty\}$ such that $f_E\equiv 0$ on $E$ and $f_E\equiv+\infty$ on $X\setminus E$.
Note that
\begin{equation}\label{equ:VT}
\sup V_E=\sup V_E^\star=\te(E).
\end{equation}

\begin{lem}\label{lem:AT}
  Let $E\subset\Xan$ be any subset. 
  \begin{itemize}
  \item[(i)]
    If $E$ is pluripolar then $V_E^\star \equiv+\infty$;
  \item[(ii)]
    if $E$ is not pluripolar, then $V_E^\star \in\PSH(\theta)$; 
  \item[(iii)]
    if $E$ is open, then $V_E^\star =V_E$;
  \item[(iv)]
    if $E$ is compact and non-pluripolar, then $\MA(V_E^\star )$ is supported in $E$.
  \end{itemize}
\end{lem}

\begin{proof}
  Suppose that $E$ is pluripolar, and pick $\p\in\PSH(\om)$ such that
  $\p|_E\equiv-\infty$. For any $t\ge 0$, we have $\p_t:=\p+t\in\PSH(\om)$ and $\p_t|_E\le0$, so
  $V_E\ge\p_t$. Thus $V_E^\star=V_E=+\infty$ on $X^\div$, and hence $V_E^\star\equiv+\infty$, by density of $X^\div$. This proves (i), while (ii) holds by the envelope property. If $E$ is open, then $f_E$ is usc, and hence $V_E^\star =V_E$, implying~(iii). 
  If $E$ is compact and nonpluripolar, then $f_E$ is lsc, so Theorem~\ref{thm:ortho}
  shows that $\MA(V_E^\star )$ is supported on the set $\{V_E^\star\ge f_E\}\subset E$, which proves (iv) 
\end{proof}

%

\begin{lem}\label{lem:ATreg}
For each subset $E\subset X^\an$ we have 
\begin{equation}\label{equ:teouter}
    \te(E)=\sup\{\te(U)\mid U\supset E\ \text{open}\}, 
\end{equation}
and
\begin{equation}\label{equ:teinner}
  \te(U)=\inf\{\te(K)\mid K\subset U\ \text{compact}\}
\end{equation}
 for every open $U\subset X^\an$. 
  \end{lem}
\begin{proof} Denote the right-hand side of~\eqref{equ:teouter} by $S\in[0,+\infty]$. Clearly, $\te(E)\ge S$. Pick $\f\in\PSH(\om)$ and $\e>0$. Since $\f$ is usc, $U:=\{\f<\sup_E\f+\e\}$ is an open set containing $E$, and hence 
$$
\sup\f\le\sup_U\f+\te(U)\le\sup_E\f+\e+S.
$$
It follows that $\te(E)\le S+\e$ for all $\e>0$, which proves~\eqref{equ:teouter}. 

To prove~\eqref{equ:teinner}, it is enough to show that the decreasing net $(V^\star_K)$ with $K$ in the directed set of compact subsets of $U$ converges to $V^\star_U=V_U$. By Theorem~\ref{thm:PSH}, the decreasing limit $\f:=\lim_K V^\star_K$ is $\om$-psh. For each compact $K\subset U$, we have  $V^\star_K\ge V_U$, and hence $\f\ge V_U\ge 0$. On the other hand, $\f\le V_K^\star\le V_{\mathring{K}}^\star=V_{\mathring{K}}$, and hence $\f=0$ on $\mathring{K}$. As this holds for all compact $K\subset U$, we infer $\f=0$ on $U$, and hence $\f\le V_U$.  
\end{proof}

\begin{proof}[Proof of Theorem~\ref{thm:AT}]
 
 Assume first that $E=K$ is compact. In this case $\Capa^\star (K)=\Capa(K)$, by Corollary~\ref{cor:capaouter}. 
  Set $T:=\te(K)=\sup V_K^\star$, and suppose first $T\le 1$. We claim that $\Capa(K)=1$, which implies~\eqref{equ:AT}. 
  Indeed, we always have $\Capa(K)\le1$, and $V_K^\star -1$ is a candidate in the
  definition of $\Capa(K)$. By Lemma~\ref{lem:AT}, $\MA(V_K^\star -1)=\MA(V_K^\star )$
  is supported on $K$, and hence
    \begin{equation*}
    \Capa(K)
    \ge \int_K\MA(V^\star _K-1)
    =\int\MA(V_K^\star )=1. 
  \end{equation*}

  Now suppose $T\ge1$. Then $T^{-1}V_K^\star -1\in\PSH(\om)$
  is a competitor in the definition of $\Capa(K)$, so
  \begin{equation*}
    \Capa(K)
    \ge\int_K\MA(T^{-1}V_K^\star -1)
    \ge T^{-n}\int_K\MA(V_K^\star )
    =T^{-n}, 
  \end{equation*}
  which proves the left-hand inequality in~\eqref{equ:AT}.
  On the other hand, we have $T^{-1}V_K^\star -1\le \f_K^\star $, so
  by 
  Theorem~\ref{thm:outcapformula} we have
  \begin{equation*}
    \Capa(K)
    =\int(-\f_K^\star )\MA(\f_K^\star )
    \le T^{-1}\int(T-V_K^\star)\MA(\f_K^\star )
    \le n T^{-1}
  \end{equation*}
where the last inequality follows from Lemma~\ref{lem:CLN}, since $V_K^\star -T\in\PSH_{\sup}$, which implies $\int(V_K^\star-T)\MA(0)=0$.
  This proves the right-hand inequality in~\eqref{equ:AT} when $E=K$ is compact. 
  
  Assume next that $E=U$ is open. By inner regularity of $\Capa$, we have
  $$
  \Capa^\star(U)=\Capa(U)=\sup_{K\subset U}\Capa(K)
  $$
  with $K$ ranging over the compact subsets of $U$, while $\te(U)=\inf_{K\subset U}\te(K)$, by~\eqref{equ:teinner}. 
  Thus~\eqref{equ:AT} for compact sets implies the case of open sets. 
  
  Finally for an arbitrary subset $E$ we have $\Capa^\star(E)=\inf_{U\supset E}\Capa(U)$ by definition of the outer capacity, and $\te(E)=\sup_{U\supset E}\te(U)$ by~\eqref{equ:teouter}. Thus~\eqref{equ:AT} for open sets implies the general case. 
\end{proof}

We are finally in a position to establish the converse of Proposition~\ref{prop:plurineg}.

\begin{thm}\label{thm:neglpp} Every negligible subset $E\subset X^\an$ is pluripolar.
\end{thm}
 
Recall that our assumptions here are that $X$ is irreducible and that $\om\in\Amp(X)$ is a class for which the envelope property holds. However, the conclusion of Theorem~\ref{thm:neglpp} holds in other cases too.
\begin{cor}\label{cor:neglpp}
  Let $X$ be any projective variety, and assume that $\charac k=0$ or $\dim X\le 2$. Then any negligible subset of $X$ is pluripolar.
\end{cor}
\begin{proof}
  By Corollaries~\ref{cor:ppcomp} and~\ref{cor:neglcomp} we may assume $X$ is irreducible. Our assumptions imply that $X$ admits a resolution of singularities, so by Lemmas~\ref{lem:ppbir} and~\ref{lem:neglbir} we may assume that $X$ is smooth. In this case, Theorem~\ref{thm:contenvlisse} implies that any ample class has the envelope property, and we conclude using Theorem~\ref{thm:neglpp}.
\end{proof}
To prove Theorem~\ref{thm:neglpp} we need
 
\begin{lem}\label{lem:neglpp}
  If $E\subset\Xan$ is a negligible subset, then $\Capa^\star (E)=0$.
\end{lem}

\begin{proof}
   By Choquet's lemma (see Corollary~\ref{cor:Choquet})
  there exists bounded, countable family $(\f_m)$ in $\PSH(\om)$ such that 
  $E\subset\{\f<\f^\star \}$, where $\f:=\sup_m\f_m$ pointwise.
  Pick any $\e>0$. By Corollary~\ref{cor:qcont}, we can
  find a compact $K\subset\Xan$ such that $\Capa^\star (K^c)=\Capa(K^c)\le\e$
  and such that $\f_m|_K$ is continuous for all $m$. Thus $\f|_K$ is lsc. 
  For each $l\in\Z_{>0}$ define
  \begin{equation*}
    K_l:=K\cap\{\f+l^{-1}\le\f^\star\}
  \end{equation*}
  Since $(\f-\f^\star)|_K$ is lsc, $K_l$ is compact. It is also negligible, so
  $\Capa^\star (K_l)=\Capa(K_l)=0$, by Corollary~\ref{cor:capaouter} and Example~\ref{exam:negcapzero}.  Furthermore, $E\cap K\subset\bigcup_l K_l$, and hence
  \begin{equation*}
    \Capa^\star (E)
    \le\Capa^\star (K^c)+\sum_l\Capa^\star (K_l)\le \e,
    \end{equation*}
    by subadditivity of $\Capa^\star$; the result follows. 
\end{proof}

\begin{proof}[Proof of Theorem~\ref{thm:neglpp}] In view of Theorem~\ref{thm:ppT}, we need to show that every negligible subset $E\subset X^\an$ satisfies $\te(E)=\infty$. This follows from Lemma~\ref{lem:neglpp} and the Alexander--Taylor inequality~\eqref{equ:AT}. 
\end{proof}

%
%
\subsection{More on envelopes} 
As a consequence of Theorem~\ref{thm:neglpp}, we have: 

\begin{thm}\label{thm:pshenv} If $\f\colon X^\an\to\R\cup\{+\infty\}$ is bounded below, then $\env(\f)^\star$ is the largest function $\p\in\PSH(\om)$ such that $\p\le\f$ outside a pluripolar set. 
\end{thm}
\begin{proof} On the one hand, $\env(\f)^\star=\env(\f)\le\f$ outside the set $\{\env(\f)<\env(\f)^\star\}$, which is negligible, and hence pluripolar by Theorem~\ref{thm:neglpp}. Assume conversely that $\p\in\PSH(\om)$ satisfies $\p\le\f$ outside a pluripolar set $E$. Pick $\rho\in\PSH_{\sup}$ with $\rho\equiv-\infty$ on $E$. For each $\e>0$, we have $(1-\e)\p+\e\rho\le (1-\e)\f\le\f+C\e$ on the whole of $X^\an$, with $C:=-\inf\f$, and  hence 
$$
(1-\e)\p+\e\rho-C\e\le\env(\f)\le\env(\f)^\star.
$$
Letting $\e\to 0$ yields $\p\le\env(\f)^\star$ outside the pluripolar set $\{\rho=-\infty\}$. In particular, $\p\le\env(\f)^\star$ on $X^\div$, and hence $\p\le\env(\f)^\star$ on $X^\an$, by Theorem~\ref{thm:suppsh}. 
\end{proof}

As a consequence, we obtain the following partial generalization of Corollary~\ref{cor:envusc}. 

\begin{cor}\label{cor:envdecr} Consider a decreasing sequence of functions $\f_m\colon X^\an\to\R\cup\{+\infty\}$ that is uniformly bounded below, and set $\f:=\lim_m\f_m$. Then $\env(\f_m)^\star\searrow\env(\f)^\star$ in $\PSH(\om)$.
\end{cor}

\begin{proof} By Theorem~\ref{thm:PSH}, $\p_m:=\env(\f_m)^\star$ converges in $\PSH(\om)$ to $\p:=\inf_m\p_m\ge\env(\f)^\star$, and we need to show that equality holds. For each $m$, we have $\p_m\le\f_i$ outside a pluripolar set $E_m$, by Theorem~\ref{thm:pshenv}. By Lemma~\ref{lem:ppcount}, $E:=\bigcup_m E_m$ is pluripolar, and $\p\le\f$ outside $E$. By Theorem~\ref{thm:pshenv}, we infer $\p\le\env(\f)^\star$, and we are done.
\end{proof}

\begin{defi} We say that a bounded $\om$-psh function $\f\in\PSH(\om)$ is \emph{regularizable from below} if there exists an increasing net $(\f_j)$ of functions in $\CPSH(\theta)$ converging to $\f$ in $\PSH(\om)$. 
\end{defi}
Thus $\f_j\nearrow\f$ pointwise on $X^\div$, and $\f=\supstar_j\f_j$. It is equivalent to demand the same property with $\f_j\in\PL\cap\PSH(\om)$. 


\begin{exam}\label{exam:regbelow} Pick a pluripolar subset $E\subset X^\an$ that contains $v_\triv$ in its closure (see Example~\ref{exam:ppdense}). Choose $\rho\in\PSH_{\sup}$ such that $\rho\equiv-\infty$ on $E$, and set $\f:=e^\rho$. Then $\f\in\PSH(\om)$ is bounded, but not regularizable from below. Indeed, any $\p\in\CPSH(\om)$ such that $\p\le\f$ satisfies $\p\le 0$ on $E$, and hence $\sup\p=\p(v_\triv)\le 0<\sup\f=1$. 
\end{exam}


Inspired by the main result in~\cite{Bed80}, we now prove the following characterization of psh functions regularizable from below, reminiscent of that of Riemann integrable functions among Lebesgue integrable functions. 

\begin{thm}\label{thm:regbelow} For a bounded function $\f\in\PSH(\om)$, the following are equivalent:
\begin{itemize}
\item[(i)] $\f$ is regularizable from below; 
\item[(ii)] $\f=\qq(\f)^\star$; 
\item[(iii)] the discontinuity locus of $\f$ is pluripolar. 
\end{itemize}
\end{thm}
Here $\qq(\f)^\star$ denotes the usc regularization of
$$
\qq(\f)=\sup\left\{\p\in\CPSH\mid\p\le\f\right\}=\env(\f_\star),
$$
cf.~Lemma~\ref{lem:envlsc}. 

\begin{proof} That (i)$\Leftrightarrow$(ii) is straightforward. Further, (ii) holds iff $\f\le\qq(\f)^\star=\env(\f_\star)^\star$, which is equivalent to $\f\le\f_\star$ outside a pluripolar set, by Theorem~\ref{thm:pshenv}. This is also equivalent to $\f=\f_\star$ outside a pluripolar set, which is a reformulation of (iii), since $\f$ is usc. Thus (ii)$\Leftrightarrow$(iii). 
\end{proof}

In dimension $n=1$, any $\f\in\PSH(\om)$  is continuous outside $v_\triv$, so since $\{v_\triv\}$ is non-pluripolar, any $\f\in\cE^\infty(\om)$ that is regularizable from below is automatically continuous. In higher dimension, the situation is different.
\begin{exam} Assume $\dim X>1$, let $C\subset X$ be an irreducible curve, and pick an ample $\Q$-line bundle. We can find $m\ge 1$ and sections $s_1,\dots,s_r,t_0,t_1\in\Hnot(X,mL)$ such that $C=\bigcap_j\{s_j=0\}$, $t_0|_C,t_1|_C\not\equiv0$, and $t_1/t_0$ defines a nonconstant rational function on $C$.
  Pick a sequence $(a_l)_1^\infty$ of distinct elements of $k$.
  Then the function
  \begin{equation*}
    \f:=\frac1m\sum_{l=1}^\infty2^{-j}\log\max\{|s_1|,\dots,|s_r|,|t_1-a_jt_0|\}
  \end{equation*}
  is $L$-psh, and it is continuous outside $C^\an$. However, it is not continuous at $v_{C,\triv}$, since $\f(v_{C,\triv})=0$ but $\f=-\infty$ on a Zariski dense subset of $C^\an$. As $\f\le0$, and $C^\an$ is pluripolar, the function $\exp(\f)$ is $L$-psh (see Corollary~\ref{cor:exppsh} ) and regularizable from below, but not continuous.
\end{exam}

%
%
\subsection{More on the outer capacity} 
The next result generalizes Theorem~\ref{thm:capformula}. 

\begin{thm}\label{thm:outcapformula}
  For any subset $E\subset\Xan$ we have
  \begin{equation}\label{equ:capaextMA}
    \Capa^\star (E)=\int(-\f_E^\star )\MA(\f_E^\star ).
  \end{equation}
\end{thm}

\begin{lem}\label{lem:extfcnopen}
 For any open $U\subset\Xan$, $\f_U=\f_U^\star$ is the limit of the decreasing net $(\f_K^\star )_K$, where $K$
  runs over the directed set of compact subsets of $U$.
\end{lem}

\begin{proof}
  We obviously have $\f_K^\star \ge \f_U^\star=\f_U$ for all $K\subset U$. 
  The limit $\p:=\lim_K\f_K^\star$ 
  therefore satisfies $\p\ge \f_U$, and it remains to show $\p(v)\le -1$ for all $v\in U$. As $\Xan$ is compact and Hausdorff,  we can find an open $V\subset\Xan$ such that 
  $v\in V\Subset U$. 
  Since $\f_U\equiv-1$ on $U$, we must have $\f_U^\star \equiv-1$ on $U$, and hence
  \begin{equation*}
    \p(v)
    \le \f_{\overline{U}}^\star (v)
    \le \f_U (v)=-1,
  \end{equation*}
  as claimed.
\end{proof}

\begin{lem}\label{lem:extfcn}
  Let $E\subset\Xan$ be any subset. Then there exists a decreasing sequence $(U_m)_{m=1}^\infty$
  of open neighborhoods of $E$ such that $\Capa(U_m)\to\Capa^\star(E)$ and $\f_{U_m}=\f_{U_m}^\star\nearrow\f_E^\star$ as $m\to\infty$.
\end{lem}

\begin{proof}
  We first claim that there exists a decreasing sequence $(U'_m)_{m=1}^\infty$
  of open neighborhoods of $E$ such that $\f_{U'_m}\nearrow\f_E^\star $ 
  as $m\to\infty$. By Choquet's lemma, we can find an increasing sequence $(\f_m)_{m=1}^\infty$ in 
  $\PSH(\om)$ such that $\f_m=-1$ on $E$ and $\f_m$ converges weakly to $\f_E^\star $.
  Set
  \begin{equation*}
    U'_m:=\{\f_m<-1+\tfrac1m\}.
  \end{equation*}
  Then $\f_m-\frac1m\le \f_{U'_m}\le \f_E$, and hence 
  $\f_m-\frac1m\le \f_{U'_m}^\star \le \f_E^\star $, which proves the claim. 

We can also, evidently, find a
  decreasing sequence $(U''_m)_{m=1}^\infty$
  of open neighborhoods of $E$ such that $\Capa(U''_m)$ decreases to $\Capa^\star (E)$.
  If we set $U_m:=U'_m\cap U''_m$, then $\Capa^\star (E)\le\Capa(U_m)\le\Capa(U''_m)$
  and $\f_{U'_m}^\star \le \f_{U_m}^\star \le \f_E^\star $ for all $m$, and the result follows.
\end{proof}

\begin{proof}[Proof of Theorem~\ref{thm:outcapformula}]
  When $E=K$ is compact, $\Capa^\star(K)=\Capa(K)$ (Corollary~\ref{cor:capaouter}), and the result thus amounts to Theorem~\ref{thm:capformula}. 
 
 Next consider the case when $E=U$ is open. By inner regularity of $\Capa$, $\Capa^\star(U)=\lim_K\Capa^\star(K)$, where $K$ runs over the directed set of compact subsets of $U$. On the other hand, Lemma~\ref{lem:extfcnopen} implies that $\f_K^\star\to\f_U$ strongly in $\cE^1$, and hence $\int(-\f_K^\star )\MA(\f_K^\star )\to\int(\f_U)\MA(\f_U)$, which proves~\eqref{equ:capaextMA} for open sets. A similar reasoning, based on Lemma~\ref{lem:extfcn}, yields the case of a general subset $E$. 
\end{proof}

\begin{thm}\label{thm:outcapincr}
  If $(E_m)$ is an increasing sequence of subsets of $X^\an$ and
  $E:=\bigcup_m E_m$, then:
  \begin{itemize}
  \item[(i)] $\f_{E_m}^\star\searrow\f_E^\star $, and $\Capa^\star (E_m)\nearrow\Capa^\star (E)$; 
  \item[(ii)] $V_{E_m}^{\star}\searrow V_E^\star$, and $\te(E_m)\searrow\te(E)$.
  \end{itemize}
\end{thm}

\begin{proof} The sequence $m\mapsto\one_{E_m}$ is increasing, and converges pointwise to $\one_E$. By Corollary~\ref{cor:envdecr}, the decreasing sequence $\f_{E_m}^\star=\env(-\one_{E_m})^\star$ thus converges to $\f_E^\star=\env(-\one_E)^\star$, and hence $\Capa^\star(E_m)\to\Capa^\star(E)$, by~\eqref{equ:capaextMA} and the continuity of Monge--Amp\`ere integrals along decreasing nets (Theorem~\ref{thm:extint}). This proves (i). 
 
The proof of (ii) is entirely similar, and left to the reader. 
\end{proof}

%
%
%
\appendix
%
%
%
\section{Dual complexes and PL functions}\label{sec:dual} 
%
%
%
%
In this section, we assume $\charac k=0$, and show how the well-known description of Berkovich spaces over discretely valued fields as limits of dual complexes carries over to the trivially valued case. In what follows, $X$ is a projective variety of dimension $n$.

%
%
\subsection{Snc test configurations and dual complexes} 
We use~\cite{jonmus} as a reference for what follows. An \emph{snc pair} $(Y,B)$ over $X$ is defined as a smooth birational model $\pi\colon Y\to X$ together with a reduced snc divisor $B=\sum_{i\in I} B_i$ on $Y$.  The \emph{dual cone complex} $\widehat\D(Y,B)$ is the simplicial cone complex whose faces are in 1--1 correspondence with the \emph{strata} $Z$ of $B$, \ie connected components of a non-empty intersection $B_J:=\bigcap_{i\in J} B_i$ for some $J\subset I$, the cone $\widehat\sigma_Z$ attached to $Z$ being identified with $\R_+^J$. In the simpler case where all $B_J$ are connected, $\widehat\D(Y,B)$ can be identified with the rational fan 
$$
\bigcup_{B_J\ne\emptyset} \R_+^J\subset\R_+^I.
$$ 
There is a canonical $\R_{>0}$-equivariant topological embedding $\widehat\D(Y,B)\hookrightarrow X^\val$, which maps $\a\in\widehat\sigma_Z\simeq\R_+^J$ to the unique valuation $\val_\a$ on $X$ that is monomial with respect to local equations of the $E_i$, centered on $Z$, and satisfies $\val_\a(E_i)=\a_i$. There is also a continuous retraction 
$$
p_{(Y,B)}\colon X^\val\to\widehat\D(Y,B), 
$$
defined as follows. If the center $c(v)\in X$ lies outside $B$, \ie $v(B)=0$, $p_{(Y,B)}$ maps $v$ to the apex of $\widehat\D(Y,B)$. Assume now $c(v)\in B$, let $(B_j)_{j\in J}$ be the set of irreducible components of $B$ containing $c(v)$, and let $Z$ be the connected component of $B_J$ containing $c(v)$. Then $p_{(Y,B)}(v)$ is the element of $\widehat\sigma_Z\simeq\R_+^J$ with coordinates $(v(B_i))_{i\in J}$. 

Snc pairs over $X$ form a directed poset, and the retractions $p_{(Y,B)}$ induce a homeomorphism 
$$
X^\val\simeq\varprojlim\hD(Y,B),
$$
see~\cite[Theorem 4.9]{jonmus}. This is further compatible with the PL structure of $X^\val$, in the sense that a function $\f$ on $X^\val$ is in $\PL_\hom(X)$ iff $\f$ is the pullback of a usual homogeneous PL function on $\hD(Y,B)$ for some $(Y,B)$; indeed, this follows from the description of $\PL_\hom(X)$ in terms of $\Q$-Cartier $b$-divisors. 

\smallskip

Assume now that $X$ is smooth, and consider an \emph{snc test configuration} $\cX$ for $X$, \ie a test configuration such that $\cX$ is nonsingular and $\cX_{0,\redu}$ is snc. By Hironaka's theorem, snc test configurations are cofinal in the directed set of all test configurations for $X$. 

Let $\cX_0=\sum_{i\in I} b_i E_i$ be the irreducible decomposition. Applying the above considerations to the reduced snc divisor $\cX_{0,\redu}$ provides a natural realization of the dual cone complex 
$$
\widehat\D_\cX:=\widehat\D(\cX,\cX_{0,\redu})
$$ 
as a set of monomial valuations $\val_\a\in\cX^\val$, which are further $k^\times$-invariant, by $\Gm$-invariance of $\cX_0$. The condition $\val_\a(\unipar)=1$ cuts out a finite simplicial complex 
$$
\D_\cX\subset\widehat\D_\cX\subset\cX^\val,
$$ 
whose faces 
$$
\sigma_Z\simeq\left\{\a\in\R_+^J\mid \sum_{i\in J} b_i \a_i=1\right\}\subset\widehat\sigma_Z\simeq\R_+^J
$$
are equipped with the integral affine structure inherited from $\Z^J\subset\R^J$ (see \eg~\cite[\S1.3]{konsoib} for details). In particular, the vertices $(e_i)$ of $\D_\cX$ are in 1--1 correspondence with the irreducible components $(E_i)$ of $\cX_0$. For any $\a\in\D_\cX$, the restriction of $\val_\a$ to 
$$
k(X)\hookrightarrow k(\cX)\simeq k(X)(\unipar)
$$
is a valuation $v_\a\in X^\val$ with Gauss extension $\sigma(v_\a)=\val_\a$, and the map $\a\mapsto v_\a$ thus provides a factorization 
$$
\D_\cX\hookrightarrow X^\val\hookrightarrow\cX^\val.
$$
As above, there is also a natural continuous retraction 
$$
p_\cX\colon X^\an\to\D_\cX,
$$
and we can now state the following trivially valued version of the well-known description of the Berkovich analytification as the limit of dual complexes of snc models, something that goes back to the fundamental work of Berkovich~\cite{Berk99}. 

\begin{thm}\label{thm:projlim} For any smooth projective variety $X$ over $k$ of characteristic $0$, the retraction maps $p_\cX\colon X^\an\to\D_\cX$ induce a homeomorphism
$$
p\colon X^\an\simeq\varprojlim_\cX \D_\cX,
$$
where the limit is over the directed set of all snc test configurations $\cX$ for $X$. For each $v\in X^\an$, we further have 
$p_\cX(v)\le v$, and $\lim_\cX p_\cX(v)=v$. 
\end{thm}

For each snc test configuration $\cX$, denote by $\Aff_\Q(\D_\cX)\subset\Cz(\D_\cX)$ the $\Q$-vector space of functions on $\D_\cX$ that are rational affine on each face (with respect to the canonical integral affine structure). 

\begin{lem}\label{lem:projlim} The space $\PL(X)\subset\Cz(X)$ satisfies $\PL(X)=\bigcup_\cX p_\cX^\star\Aff_\Q(\D_\cX)$.
\end{lem}
\begin{proof} Mapping a function $f\in\Aff_\Q(\D_\cX)$ to $D_f:=\sum_i f(e_i) b_i E_i\in\VCar(\cX)_\Q$ defines an isomorphism $\Aff_\Q(\D_\cX)\simeq \VCar(\cX)_\Q$, and it is easy to see from the definition of $p_\cX$ that $f\circ p_\cX=\f_{D_f}$. The rest follows from Theorem~\ref{thm:VCarPL}. 
\end{proof}

\begin{proof}[Proof of Theorem~\ref{thm:projlim}] The map $p$ is continuous, and $X^\div$ maps onto the dense subset $\varprojlim_\cX\D_\cX(\Q)$ of $\varprojlim_\cX \D_\cX$. Since $X^\an$ is compact (Hausdorff), it is thus remains to show that $p$ is injective. This is a simple consequence of Lemma~\ref{lem:projlim}, since $\PL(X)$ separates the points of $X^\an$ (Lemma~\ref{lem:FSdense}). 

Now pick $v\in X^\an$. That $p_\cX(v)\le v$ is immediate from the definition of $p_\cX$. Since $\PL(X)$ is dense in $\Cz(X)$, it remains to see that $\lim_\cX\f(p_\cX(v))=\f(v)$ for all $\f\in\PL(X)$, which again follows from Lemma~\ref{lem:projlim}. 
\end{proof}

\begin{cor}\label{cor:apprbelow} Let $X$ be a projective variety over $k$ with $\charac k=0$. Then any $v\in X^\an$ is 
the limit of a net $(v_i)$ in $X^\div$ such that $v_i\le v$ for all $i$. 
\end{cor}
\begin{proof}
  It suffices to prove that for any $t_0>1$ and any neighborhood $U$ of $v$ in $\Xan$, there exists $v_0\in\Xdiv$ with $v_0\le t_0v$.

First assume $X$ is smooth. 
By Theorem~\ref{thm:projlim}, we can find an snc test configuration $\cX$ such that $w:=p_\cX(v)\in\Delta_\cX\cap U$. Moreover, $w\le v$. If $w\in\Xdiv$, then we can take $v_0=w$, so suppose $w\not\in\Xdiv$, and let $\sigma\subset\Delta_\cX$ be the unique simplex containing $w$ in its interior. Then $\dim\sigma>0$, or else $w\in\Xdiv$.

We claim there exists a continuous function $t\colon\mathring{\sigma}\to[1,\infty)$ such that $t(w)=1$ and $w'\le t(w')w$ for all $w'\in\mathring{\sigma}$. In view of the embedding $\Delta_\cX\hookrightarrow X^\val\hookrightarrow\cX^\val$ above, this follows from the elementary fact that if $\alpha\in\R_{>0}^J$, then there exists a continuous function $t\colon\R_{>0}^J\to\R_{>0}$ such that $t(\a)=1$ and $\a'_j\le t(\a')\a_j$ for all $\a'\in\R_{>0}^J$ and all $j\in J$. This function can be chosen as $t(\a')=\max_j\a'_j/\a_j$, for example.

Now $\Xdiv$ is dense in $\sigma$, so we can pick $w'\in\Xdiv\cap\sigma$ close enough to $w$ so that $t(w')\le t_0$, and then we can pick $v_0=w'$.

  In the general case, let $\mu\colon X'\to X$ be a resolution of singularities, and pick $v'\in X^{\prime\,\mathrm{an}}$ with $\mu^{\mathrm{an}}(v')=v$. By what precedes, there exists $v'_0\in X^{\prime\,\mathrm{div}}\cap(\mu^{\mathrm{an}})^{-1}(U)$ such that $v'_0\le t_0v'$. We can then choose $v_0=\mu^{\mathrm{an}}(v'_0)$.
\end{proof}

%
%
\subsection{Psh functions and dual complexes}

Building on the uniform Izumi-type estimates of~\cite{siminag}, we show: 

\begin{thm}\label{thm:pshdual} Let $\cX$ be an snc test configuration, with dual complex $\D_\cX\hookrightarrow X^\an$. Then: 
\begin{itemize}
\item[(i)] for each $\f\in\PSH(\om)$, $\f|_{\D_\cX}$ is finite-valued, continuous, and convex on each face of $\D_\cX$;
\item[(ii)] the set $\{\f|_{\D_\cX}\mid \f\in\PSH(\om)\}$ is equi-Lipschitz continuous;
\item[(iii)] $\D_\cX$ is a strongly compact subset of $X^\lin$.
\end{itemize}
\end{thm}

\begin{proof} Set $K:=k\lau{\unipar}$. This is a non-Archimedean field with valuation ring $K^\circ=k\cro{\unipar}$. Consider the Berkovich analytification $X_K^\an$. Gauss extension can be viewed as a continuous section $\sigma\colon X^\an\to X_K^\an$ of the natural projection $\pi\colon X_K^\an\to X^\an$.

  The pullback of $\om$ is an ample class $\om_K\in\Num(X_{K^\circ}/K^{\circ})$, which is interpreted as a closed $(1,1)$-form with ample de Rham class in~\cite{siminag}, and it follows immediately from the definitions that $\f\circ\pi$ is a $\om_K$-psh model function on $X_K^\an$ for each $\f\in\PL\cap\PSH(\om)$. By~\cite[Theorem 6.1]{siminag}, $\f|_{\D_\cX}$ is thus continuous, convex on each face of $\D_\cX$, and Lipschitz continuous with Lipschitz constant $C>0$ only depending on $\om$.

Since each $\f\in\PSH(\om)$ is a decreasing limit of functions in $\PL\cap\PSH(\om)$, (i) and (ii) follow. By (i), $\D_\cX$ is contained in $X^\lin$, and (ii) precisely says that the embedding $\D_\cX\hookrightarrow\Xlin$ is Lipschitz continuous with respect to $\dd_\infty$, see Definition~\ref{equ:Izumi}.
By compactness of $\D_\cX$, it follows that the inclusion is a homeomorphism onto its image with respect to the strong topology, which proves (iii). 
\end{proof}

%
%
 \section{The toric case}\label{sec:toric}
The goal of this section is to provide a brief description of various objects considered in this paper in the context of toric varieties~\cite{FultonToric,BPS}.

\begin{itemize}
\item  Consider an algebraic torus $T\simeq\Gm^n$, with associated dual lattices $M:=\Hom(T,\Gm)$ and $N:=\Hom(\Gm,T)$. We have a canonical embedding $M\hookrightarrow k(T)^\times$ onto the set $T$-invariant functions, and a dual canonical embedding $N_\R\hookrightarrow T^\val$ onto the set of $T(k)$-invariant valuations, such that $v(u)=\langle v,u\rangle$ for all $v\in N_\R$ and $u\in M\hookrightarrow k(T)^\times$. There is also a canonical retraction $\rho\colon T^\val\to N_\R$, which maps $v$ to the linear form on $M$ given by $u\mapsto v(u)$.  

\item A proper (normal) toric variety $X$ corresponds to a rational fan decomposition $\Sigma$ of $N_\R$, and an ample class $\om\in\Amp(X)$ to a polyhedron $P\subset M_\R$ (up to translation) with normal fan $\Sigma$. The support function of $P$ is the convex PL function $f_P\colon N_\R\to\R$ defined by 
$$
f_P(v):=\sup_{u\in P}\langle v,u\rangle,
$$
which is linear precisely on the cones of $\Sigma$. 

\item For any $\om$-psh function $\f$, the function $f_\f:N_\R\to\R$ defined by
$$
f_\f:=(\f+f_P)|_{N_\R}
$$
is convex. This sets up a 1--1 correspondence between the set $\PSH_\tor(\om)$ of $\te(k)$-invariant $\om$-psh functions $\f$ on $X^\an$ and the set of all convex functions $f:N_\R\to\R$ such that $f\le f_P+O(1)$, the inverse being given by 
$$
f\mapsto\f:=\rho^\star(f-f_P),
$$
with $\rho\colon X^\val=T^\val\to N_\R$ the retraction. 

\item The latter set of convex functions is in turn in 1--1 correspondence with the set $\cP$ of all lsc convex functions $g:P\to\R\cup\{+\infty\}$, via the Legendre transform
$$
g(u)=f^\vee(u)=\sup_{v\in N_\R}\left(\langle v,u\rangle-f(v)\right),
$$
$$
f(v)=g^\vee(v)=\sup_{u\in P}\left( \langle v,u\rangle-g(u)\right).
$$
\item For any bounded $\f\in\PSH_\tor(\om)$, the real Monge--Amp\`ere measure $\MA_\R(f_\f)$ of the convex function $f_\f=f_P+O(1)$ is a positive measure on $N_\R$ of total mass $V_P:=\vol(P)=(\om^n)/n!$, and the non-Archimedean Monge--Amp\`ere measure of $\f$ satisfies
$$
\MA_\om(\f)=V_P^{-1}\iota_\star\MA_\R(f_\f)
$$
with
$$
\iota:N_\R\hookrightarrow T^\val=X^\val\subset X^\an
$$ 
the inclusion. Furthermore,  $f_\f^\vee$ is bounded on $P$, and 
\begin{equation}\label{equ:Etoric}
\en_\om(\f)=-\fint_P f_\f^\vee.
\end{equation}

\item Equation~\eqref{equ:Etoric} remains valid for all $\f\in\PSH_\tor(\om)$, and shows that 
$$
\f\in\cE^1_\tor(\om)\Longleftrightarrow f_\f^\vee\in L^1(P).
$$
Furthermore, $\f_i\to\f$ strongly in $\cE^1_\tor(\om)$ iff $f_{\f_i}^\vee\to f_\f^\vee$ in $L^1(P)$. 

\item For any $v\in N_\R\subset X^\lin$, the function $\f_v\in\CPSH(\om)$ satisfies 
$$
f_{\f_v}(w)=f_P(w-v)+f_P(v).
$$
In particular, $\te_\om(v)=f_P(v)+f_P(-v)=f_{P+(-P)}(v)$.

\item As in~\cite[Proposition 5.7]{Berm13}, the energy of a probability measure $\mu$ on $N_\R\subset X^\an$ coincides with the optimal cost $\mathrm{C}_P(\mu)$ of transporting the measure $\mu$ to $\la_P$ with respect to the cost function $c_P:N_\R\times P\to\R_+$ given by 
$$
c_P(v,u):=f_P(v)-\langle v,u\rangle.
$$  
Indeed, 
$$
\en^\vee_\om(\mu)=\sup_{\f\in\cE^1_{\tor}(\om)}\left(-\fint_P f_\f^\vee-\int_{N_\R} \f\,\mu\right)=\sup_{\f\in\mathrm{C}^0_b(N_\R)}\left(-\fint_P (f_P+\f)^\vee-\int_{N_\R}\f\,\mu\right),
$$
and Monge--Kantorovich duality yields
$$
\en^\vee_\om(\mu)=\inf_{\La}\int_{N_\R\times P} c_P(v,u)\,\La=\mathrm{C}_P(\mu),
$$
where $\La$ ranges over all probability mesures on $N_\R\times P$ with marginals $\mu$ and $\la_P$. In particular, we have for any $v\in N_\R$ 
$$
\en^\vee_\om(\d_v)=f_P(v)-\langle v,u_P\rangle=\sup_Pv-\fint_P v, 
$$
with $u_P\in P$ the center of mass.

\end{itemize}
%
%

%
%
%
%
%
%
\end{document}